\documentclass[11pt,english]{article}
\usepackage[T1]{fontenc}
\usepackage[latin9]{inputenc}
\usepackage{geometry}
\usepackage{color}
\usepackage{babel}
\usepackage{array}
\usepackage{rotating}
\usepackage{multirow}
\usepackage{amstext}
\usepackage{amssymb}
\usepackage{setspace}
\usepackage[unicode=true]
 {hyperref}
\usepackage{breakurl}

\makeatletter

\providecommand{\tabularnewline}{\\}

\date{}

\makeatother

\begin{document}

\title{High-Dimensional Semiparametric Selection Models: Estimation Theory
with an Application to the Retail Gasoline Market%
\thanks{First, I thank James Powell, Martin Wainwright, Miguel Villas-Boas,
and Demian Pouzo for useful suggestions and comments. I also thank
Ganesh Iyer for the data set and helpful discussions on the empirical
application. I am also grateful to Michael Jansson, Bryan Graham,
Valentin Verdier, Mikkel Soelvsten, Przemyslaw Jeziorski, and other
participants at the Econometrics Seminar of Department of Economics
and the Shansby Marketing Seminar of the Haas School of Business,
UC Berkeley. All errors are my own. This work was supported by Haas
School of Business at UC Berkeley. %
} }

\author{Ying Zhu\\
\textbf{(Job Market Paper)}\\
November 4, 2014}

\maketitle
\textit{Haas School of Business, UC Berkeley. 2220 Piedmont Ave.,
Berkeley, CA 94720. ying\_zhu@haas.berkeley.edu. Tel: 406-465-0498.
Fax: 510-643-4255}
\begin{abstract}
This paper proposes a multi-stage projection-based \textit{\textcolor{black}{Lasso}}
procedure for the semiparametric sample selection model in high-dimensional
settings under a weak nonparametric restriction on the form of the
selection correction. In particular, the number of regressors in the
main equation, $p$, and the number of regressors in the selection
equation, $d$, can grow with and exceed the sample size $n$. The
analysis considers the \textit{exact sparsity} case where the number
of non-zero components in the vectors of coefficients is bounded above
by some integer which is allowed to grow with $n$ but slowly compared
to $n$, and also considers the \textit{approximate sparsity} case,
where the vectors of coefficients can be approximated by exactly sparse
vectors. The main theoretical results of this paper are finite-sample
bounds from which sufficient scaling conditions on the sample size
for estimation consistency and variable-selection consistency (i.e.,
the multi-stage high-dimensional estimation procedure correctly selects
the non-zero coefficients in the main equation with high probability)
are established.\textcolor{black}{{} A technical issue related to a
set of high-level assumptions on the regressors for estimation consistency
and selection consistency arises in the multi-stage estimation procedure
from allowing the number of regressors in the main equation to exceed
}$n$\textcolor{black}{{} and this paper provides analysis to verify
these conditions.} These verifications also provide a finite-sample
guarantee of the population identification condition required by the
semiparametric sample selection models. Statistical efficiency of
the proposed estimators is studied via lower bounds on minimax risks
and the result shows that, for a family of models with exactly sparse
structure on the coefficient vector in the main equation, one of the
proposed estimators attains the smallest estimation error up to the
$(n,\, d,\, p)-$scaling among a class of procedures in worst-case
scenarios. Inference procedures for the coefficients of the main equation,
one based on a pivotal Dantzig selector to construct non-asymptotic
confidence sets and one based on a post-selection strategy (when perfect
or near-perfect selection of the high-dimensional coefficients is
achieved), are discussed. Other theoretical contributions of this
paper include establishing the non-asymptotic counterpart of the familiar
asymptotic {}``oracle'' type of results from previous literature:
the estimator of the coefficients in the main equation behaves as
if the unknown nonparametric component were known, provided the nonparametric
component is sufficiently smooth. Small-sample performance of the
high-dimensional multi-stage estimation procedure is evaluated by
Monte-Carlo simulations and illustrated with an empirical application
to the retail gasoline market in the Greater Saint Louis area. Proofs
are included in the online supplementary material (\href{https://sites.google.com/site/yingzhu1215/home/JobMar_Proofs.pdf}{https://sites.google.com/site/yingzhu1215/home/JobMar\_{}Proofs.pdf}).

\newpage{}
\end{abstract}

\section{Introduction }

The past decade has witnessed research activities in high-dimensional
statistics that considers inference for models in which the dimension
of the parameters of interests is comparable to or even larger than
the sample size. The rapid advance of data collection technology is
a major driving force of the development of high-dimensional statistics:
it allows for not only more observations but also more explanatory
variables to be collected. Recently, high-dimensional estimation techniques
have been studied in several popular econometric models and some first
applications of these techniques are now available in economics. However,
a very important class of models, sample selection models, have not
been considered in high-dimensional settings even though they are
central in many economic and marketing applications. For example,
on the demand side, consumers often face choosing a service or brand
followed by the amount of utilization or the number of quantities
to purchase conditional on the chosen service or brand. On the supply
side, firms first decide on the product positioning and then a pricing
scheme based on the chosen product type. Selection models are also
seen in the auction literature. In estimating the underlying selection
models to study these empirical problems, only a low-dimensional set
of explanatory variables has been considered even though the actual
information available to these empirical problems can be far richer
than what has been used by the researchers. The lack of estimation
methods that deal with these {}``data-rich'' selection problems
may have limited the use of high-dimensional techniques in many economics
and marketing problems. This paper aims to provide estimation tools
together with their theoretical guarantees for this important but
little studied topic.

Observational studies are rarely based on pure random samples. When
a sample, intentionally or unintentionally, is based in part on values
taken by a dependent variable (e.g., Gronau, 1973; Heckman, 1974),
parameter estimates without corrective measures may be inconsistent.
Such samples can be broadly defined as selected samples. Selection
may be due to self-selection, with the outcome of interest determined
in part by individual choice of whether or not to participate in the
activity of interest. It can also result from endogenous stratification,
with those who participate in the activity of interest deliberately
oversampled - an extreme case being sampling only participants. 

In the classical low-dimensional selection models, parameter estimates
obtained from OLS may be inconsistent unless corrective measures are
taken. For the parametric case where the error terms are jointly normally
distributed and homoskedastic, the most well-known estimator is Heckman's
two-step procedure (1974, 1976). For semiparametric estimation of
the parameters of selection models when the joint distribution of
the error terms is of unknown form, many estimators have adopted the
two-step estimation strategy similar to Heckman's, under an additional
{}``single-index'' restriction on the form of the selection equation.
Several estimators of the parameters of the selection equation have
been proposed in the literature on semiparametric estimation; while
some of these methods sidestep estimation of the unknown distribution
function of the errors (e.g., Manski, 1975 and 1985; Han, 1987), others
use nonparametric regression methods to estimate this distribution
function along with the parameters of the underlying regression function
(e.g., Cosslett, 198l; Ichimura, 1987; Klein and Spady, 1987). Similarly,
the methods for estimation of the parameters in the second stage also
involve nonparametric regression methods, which are applied either
to estimation of the selection correction function directly (Lee,
1982; Cosslett, 1991; Gallant and Nychka, 1987; Ichimura and Lee,
1991; Newey, 1991) or to estimation of other regression functions
which depend upon the estimated single index (Powell, 1989; Ahn and
Powell, 1993). 

The object of interest of this paper is a class of \textit{high-dimensional}
selection models under a weak nonparametric restriction on the form
of the selection correction. Statistical estimation and variable selection
in the high-dimensional setting concerns models in which the dimension
of the parameters of interests is comparable or even larger than the
sample size. In the past decade, an increase of research activities
in this field has been facilitated by the advances in data collection
technology. In the literature on high-dimensional sparse linear regression
models, a great deal of attention has been given to the $l_{1}-$penalized
least squares. In particular, the Lasso and the Dantzig selector are
the most studied techniques (see, e.g., Tibshirani, 1996; Candès and
Tao, 2007; Bickel, Ritov, and Tsybakov, 2009; Belloni, Chernozhukov,
and Wang, 2011; Belloni and Chernozhukov, 2011b; Loh and Wainwright,
2012; Negahban, Ravikumar, Wainwright, and Yu, 2012). Variable selection
when the dimension of the problem is larger than the sample size has
also been studied in the likelihood method setting with penalty functions
other than the $l_{1}-$norm (see, e.g., Fan and Li, 2001; Fan and
Lv, 2011). Lecture notes by Koltchinskii (2011), as well as recent
books by Bühlmann and van de Geer (2011) and Wainwright (2015) have
given a more comprehensive introduction to high-dimensional statistics. 

Recently, these $l_{1}-$penalized techniques have been applied in
a number of econometric papers. Caner (2009) studies a Lasso-type
GMM estimator. Rosenbaum and Tsybakov (2010) study the high-dimensional
errors-in-variables problem where the non-random regressors are observed
with additive error and they present an application to hedge fund
portfolio replication. Belloni and Chernozhukov (2011a) study the
$l_{1}-$penalized quantile regression and illustrate its use on an
international economic growth application. Fan, Lv, and Li (2011)
review the literature on sparse high-dimensional econometric models
including the vector autoregressive model for measuring the effects
of monetary policy, panel data model for forecasting home price, and
volatility matrix estimation in finance. Their discussion is not restricted
to $l_{1}-$based regularization methods. Manresa (2014) considers
settings where outcomes depend on an agent's own characteristics and
on the characteristics of other agents in the data and applies a Lasso
type estimator to study individuals generating spillovers and their
strength using panel data on outcomes and characteristics. Bonaldi,
Hortacsu, and Kastl (2014) propose a new measure of systemic risk
based on estimating spillovers between funding costs of individual
banks with a Lasso type procedure applied to the panel of each individual
bank to recover the financial network. Lecture notes by Belloni and
Chernozhukov (2011b) discuss the $l_{1}-$based penalization methods
with various econometric problems including earning regressions and
instrumental selection in Angrist and Krueger data (1991). Belloni,
Chen, Chernozhukov, and Hansen (2012) estimate the optimal instruments
using the Lasso and in an empirical example dealing with the effect
of judicial eminent domain decisions on economic outcomes, they find
the Lasso-based instrumental variable estimator outperforms an intuitive
benchmark. Belloni, Chernozhukov, and Hansen (2014) propose robust
methods for inference on the effect of a treatment variable on a scalar
outcome in the presence of many controls with an application to abortion
and crime. In many economic applications, the number of endogenous
regressors is also large relative to the sample size. The case of
many endogenous regressors and many instrumental variables has been
studied by Gautier and Tsybakov (2011), Zhu (2013), and Fan and Liao
(2014).

While previous literature has extended the estimation theories and
applications of several popular econometric models from the classical
low-dimensional settings to the high-dimensional settings, selection
models have not been considered in high-dimensional settings even
though many economic applications actually fit into this setup. On
the demand side, selection models are used in the context where a
consumer faces choosing a service (such as electricity, cell phone
service, etc.) or brand followed by the amount of utilization or the
number of quantities to purchase conditional on the chosen service
or brand (e.g., Krishnamurth and Raj, 1988; Chintagunta, 1993; Fox,
Kim, and Yang, 2013). On the firms' side, selection models are useful
for situations where a firm first decides on its product positioning
and then a pricing scheme based on the chosen product type. For example,
a grocery store sometimes needs to choose which products to put on
sales or promotions and then the amount of discount on these chosen
products; a gas station first chooses to be either a two-product station
offering both self-service and full-service gasoline or a single-product
station offering only full-service or self-service gasoline, and then
decides on a pricing scheme conditional on the choice of the station
type (Iyer and Seetharaman, 2003). Selection models are also seen
in the auction literature (e.g., Roberts and Sweeting, 2011, 2012);
in particular, by estimating a Heckman selection model with the exclusion
restriction that potential competition affects a bidder's decision
to enter an auction, but has no direct effect on the values of the
bids, Roberts and Sweeting (2011) presents reduced form evidence that
the auction data are best explained by a selection model. 

In estimating the underlying selection models to study these empirical
problems, analysis has been restricted to only a low-dimensional set
of explanatory variables in both the selection equation and the main
equation. However, the actual information available to these empirical
problems can be far richer than what has been used by the researchers.
More importantly, economic theory is not always explicit about the
variables that belong to the true model (e.g., Sala-i-Martin, 1997
concerning development economics). In the selection models used for
consumer demand estimation, the number of explanatory variables formed
by the characteristics (and the transformations of these characteristics)
of a service or brand can be very large. In the grocery store example,
when choosing whether to put a product on sale and the amount of discount,
the store often considers not only the own characteristics of this
product but also characteristics of other products. All these characteristics
can potentially exceed the number of products chosen to be on sale
(namely, the sample size of the observations in the main-equation),
which makes it a high-dimensional selection problem. Similarly, in
the bidder example, when deciding whether to enter an auction, a bidder
considers potential competition from other bidders; when deciding
on the values of its bid upon the entry decision, the bidder may still
consider competition from the set of other {}``enters''. Consequently,
the number of explanatory variables entering the selection equation
and the main equation may grow with the number of bidders. 

In the gasoline example mentioned above, besides the large number
of station characteristics and demographic characteristics which amount
to approximately 400 regressors with only 700 gas stations in the
data studied by Iyer and Seetharaman (2008), more explanatory variables
can be obtained by utilizing the geographic information and spatial
data. In particular, geographic information may be used to analyze
the interaction between different gas stations and identify the competitive
market structure, as will be shown in Section 6 of this paper. Despite
that the explanatory variables in the examples above are high-dimensional,
it is plausible that only a small set of these variables (relative
to the sample size) matter to the underlying response variables but
which variables constitute the relevant regressors are unknown to
the researchers.

The following sets up the models of interests and highlights the major
contributions made by this paper. In particular, we consider estimation
and selection of regression coefficients in the class of selection
models captured by the following system: for all $i=1,\,...,\, n$,
\begin{eqnarray}
y_{1i} & =\mathbb{I} & \left\{ w_{i}^{T}\theta^{*}+\epsilon_{1i}>0\right\} ,\nonumber \\
y_{2i} & = & y_{1i}\left(x_{i}^{T}\beta^{*}+\epsilon_{2i}\right),\label{eq:1}\\
\mathbb{E}\left(\epsilon_{2i}\vert w_{i},\, x_{i},\,\, y_{1i}=1\right) & = & g(w_{i}^{T}\theta^{*}),\nonumber 
\end{eqnarray}
where $x_{i}$ is a $p-$dimensional vector of explanatory variables
and the dimension $p$ of $\beta^{*}$ is large relative to the sample
size $n$ (namely, $p\asymp n$ or even $p\gg n$). Furthermore, $g(\cdot)$
is an unknown function and $w_{i}$ is a $d-$dimensional vector of
explanatory variables with an unknown coefficient vector $\theta^{*}$.
Note here the dimension $d$ of $\theta^{*}$ can also be large relative
the sample size $n$ (namely, $d\asymp n$ or even $d\gg n$). The
third equation in (\ref{eq:1}) is known as the {}``single-index''
restriction used in Powell (1989), Newey (1991), and Ahn and Powell
(1993). Newey (1991) and Powell (1994) discuss sufficient conditions
for this restriction. In particular, it is implied by independence
of the errors $(\epsilon_{1i},\,\epsilon_{2i})$ and the regressors
$(w_{i},\, x_{i})$. Note that the second equation of model (\ref{eq:1})
implies 
\begin{equation}
y_{2i}=x_{i}^{T}\beta^{*}+g(w_{i}^{T}\theta^{*})+\eta_{i}\quad\textrm{whenever}\: y_{1i}=1,\label{eq:2}
\end{equation}
where by construction $\mathbb{E}\left[\eta_{i}\vert w_{i},\, x_{i},\, y_{1i}=1\right]=0$.
Throughout the rest of this paper, when it is clear from the context
that only the selected sample is of our interests, the notation $y_{1i}=1$
will be suppressed. In addition, the values of $n$ will vary according
to whether we are working with the whole sample (the observations
in the selection equation) or the selected sample (the observations
in the main equation). Motivated by the Frisch-Waugh Theorem, applying
a projection idea used in Robinson's semilinear models (1988) yields
the following equivalent model 
\begin{equation}
v_{i0}=v_{i}\beta^{*}+\eta_{i},\label{eq:3}
\end{equation}
where 
\begin{eqnarray*}
v_{i} & = & \left(x_{i1}-\mathbb{E}\left(x_{i1}\vert w_{i}^{T}\theta^{*}\right),\,...\,,x_{ip}-\mathbb{E}\left(x_{ip}\vert w_{i}^{T}\theta^{*}\right)\right),\\
v_{i0} & = & y_{2i}-\mathbb{E}\left(y_{2i}\vert w_{i}^{T}\theta^{*}\right).
\end{eqnarray*}
For convenience, the first equation in (\ref{eq:1}) is referred to
as the selection equation and the second equation in (\ref{eq:1})
as the main equation. 

High dimensionality arises in selection model (\ref{eq:1}) when the
dimension $p$ of $\beta^{*}$ is large relative to the sample size
$n$ (namely, $p\asymp n$ or $p\gg n$) in the main equation. In
addition, a weak nonparametric restriction is imposed on the form
of the selection correction. Specifically, the selection effect is
assumed to depend on the linear combination of some observable selection
variables. The selection model under this nonparametric restriction
on the form of the selection correction when $p\geq n$ has apparently
not been studied in the literature. As in classical low-dimensional
selection models where parameter estimates obtained from OLS may be
inconsistent, direct implementation of the Lasso or the Dantzig selector
fails as sparsity of coefficients in the main equation in (\ref{eq:1})
may not correspond to sparsity of linear projection coefficients and
{}``bias'' from parameter estimates by the direct Lasso procedure
without corrective measures is found to be only exacerbated in the
high-dimensional setting. This evidence is given by the Monte-Carlo
simulation results in Section 5. The selection equation in (\ref{eq:1})
is a linear latent variable model and the selection bias $g(\cdot)$
is assumed to be an unknown function of the single index $w_{i}^{T}\theta^{*}$.
This setup allows us to consider special cases where the dimension
$d$ of $\theta^{*}$ is also large relative to the sample size $n$
(namely, $d\asymp n$ or $d\gg n$) in the selection equation described
by some of the most popular binary response models. It is worth noting
that the general results provided by this paper also hold for the
more general structure where $\mathbb{E}\left(\epsilon_{2i}\vert w_{i},\, y_{1i}=1\right)=g(h(w_{i}^{T},\,\theta^{*}))$
and $h(w_{i}^{T},\,\theta^{*})$ is a scalar unobservable index, under
appropriate identification assumptions. 

The proposed estimation procedure for the high-dimensional linear
coefficients in the main equation in this paper is the penalized version
of a projection-type strategy. In the first-stage, given consistent
estimates $\hat{\theta}$ of $\theta^{*}$ in the selection equation
obtained using one of several methods recently proposed in the high-dimensional
statistics literature, estimates $w_{i}^{T}\hat{\theta}$ of the \textquotedblleft{}single
index\textquotedblright{} variables $w_{i}^{T}\theta^{*}$ are formed.
In the second-stage, nonparametric regression is performed to obtain
estimate $\hat{\mathbb{E}}\left(x_{ij}\vert w_{i}^{T}\hat{\theta}\right)$
of $\mathbb{E}\left(x_{ij}\vert w_{i}^{T}\theta^{*}\right)$ for $j=1,...,p$
and $\hat{\mathbb{E}}\left(y_{2i}\vert w_{i}^{T}\hat{\theta}\right)$
of $\mathbb{E}\left(y_{2i}\vert w_{i}^{T}\theta^{*}\right)$; then
the estimated residuals $\hat{v}_{i}=\left(x_{i1}-\hat{\mathbb{E}}\left(x_{i1}\vert w_{i}^{T}\hat{\theta}\right),\,...\,,x_{ip}-\hat{\mathbb{E}}\left(x_{ip}\vert w_{i}^{T}\hat{\theta}\right)\right)$
of $v_{i}$ and $\hat{v}_{i0}=y_{2i}-\hat{\mathbb{E}}\left(y_{2i}\vert w_{i}^{T}\hat{\theta}\right)$
of $v_{i0}$ are formed. This step is motivated by the estimator of
Robinson (1988) for semilinear models. The second-stage estimation
in this paper involves $p+1$ nonparametric regressions where $p\asymp n$
or $p\gg n$, and in contrast to the classical low-dimensional settings
(e.g., Robinson 1988), a more careful control for the noise from the
$p+1$ nonparametric regressions is required. In particular, the prediction
errors of the nonparametric procedures are shown in this paper to
satisfy 
\[
\mathbb{P}\left\{ \sqrt{\frac{1}{n}\sum_{i=1}^{n}\left[\hat{\mathbb{E}}\left(z_{ij}\vert w_{i}^{T}\hat{\theta}\right)-\mathbb{E}\left(z_{ij}\vert w_{i}^{T}\hat{\theta}\right)\right]^{2}}\geq t\right\} \leq c\exp\left(-nt^{2}\right)
\]
where $z_{ij}=x_{ij}$ for $j=1,...,p$ and $z_{i0}=y_{2i}$, and
as a consequence, 
\begin{eqnarray*}
\mathbb{P}\left\{ \max_{j=0,...,p}\sqrt{\frac{1}{n}\sum_{i=1}^{n}\left[\hat{\mathbb{E}}\left(z_{ij}\vert w_{i}^{T}\hat{\theta}\right)-\mathbb{E}\left(z_{ij}\vert w_{i}^{T}\hat{\theta}\right)\right]^{2}}\geq t\right\}  & \leq & c\exp\left(-nt^{2}+\log p\right)\\
 & = & O\left(\frac{1}{p}\right)
\end{eqnarray*}
where the last equality holds provided $n$ is sufficiently large.
The tail bounds above can be ensured by considering the family of
nonparametric least squares estimators or regularized nonparametric
least squares estimators defined in van de Geer (2000). This family
of estimators include linear regression as the simplest case, sparse
linear regressions, convex regression, Lipschitz and Isotonic regression,
kernel ridge regression based on reproducing kernel Hilbert spaces,
estimators based on series expansion, sieves and spline methods. A
procedure based on Lipschitz regression for the second-stage nonparametric
estimation is illustrated in this paper for a leading case example. 

In the third-stage, regressing $\hat{v}_{i0}$ on $\hat{v}_{i}$ with
$l_{1}-$regularization to estimate the main-equation coefficients
$\beta^{*}$. In particular, for the third-stage estimation, this
paper considers a non-pivotal Lasso procedure whose regularization
parameter depends on the unknown variance of $\eta_{i}$, and a pivotal
Dantzig selector whose regularization parameter does not involve the
unknown variance of $\eta_{i}$. This pivotal Dantzig selector was
originally proposed by Gautier and Tsybakov (2011) in the context
of instrumental variables regression. A by-product of the pivotal
procedure is a set of non-asymptotic confidence intervals (which also
do not involve the unknown variance of $\eta_{i}$). Upon the availability
of estimates of the high-dimensional linear coefficients, two different
estimation strategies for the selection bias function are proposed:
one is a closed form estimator and the other is a nonparametric least
squares estimator. Despite that the nonparametric least squares estimator
of $g(w_{i}^{T}\theta^{*})$ is computationally more involved relative
to the closed-form estimator, its rate of convergence turns out to
be faster. In particular, when $\beta^{*}$ is approximately sparse
with $q_{2}=1$, the closed-form estimator cannot achieve \textit{MSE}-consistency
even if $n\rightarrow\infty$ while the nonparametric least squares
estimator is consistent in \textit{MSE }when $q_{2}=1$. 

While existing semiparametric estimation techniques for the selection
models limit the number of regressors entering the selection equation
and the main equation, the multi-stage estimation procedure with $l_{1}-$regularization
in the first- and third-stage are more flexible and particularly powerful
for applications in which the vector of parameters of interests is
high-dimensional but sparse and there is lack of information about
the relevant explanatory variables. Moreover, this above-mentioned
high-dimensional multi-stage estimation procedure is intuitive and
can be easily implemented using existing software packages. In particular,
it decomposes the joint search of the optimal values for the high-dimensional
linear coefficients and the nonparametric selection bias component
into several sequential searches with each search defined over a much
smaller parameter space. In particular, the second-stage estimation
incurs a computational cost linear in $p$ as it involves solving
$p+1$ independent subproblems and each subproblem can be in general
solved with a polynomial-time algorithm. The computational efficiency
of the first-stage and third-stage estimations is guaranteed by existing
algorithms developed for solving the Lasso or the Dantzig program.
Upon the availability of estimates of the high-dimensional linear
coefficients, the estimator for the selection bias function is simply
a closed form estimator or a nonparametric least squares estimator.
In additional to the computational efficiency, as we will see in Section
4.4 that, under some conditions and when $\beta^{*}$ is exactly sparse,
the proposed procedures for estimating $\beta^{*}$ and $g(\cdot)$
are overall statistically efficient up to the $(n,\, d,\, p)-$factors,
relative to any procedure constructed based on model (\ref{eq:2})
for estimating model (\ref{eq:1}), regardless of its computational
cost.

The main theoretical results of this paper are finite-sample bounds
from which sufficient scaling conditions on the sample size for estimation
consistency in $l_{2}-$norm and variable-selection consistency (i.e.,
the multi-stage high-dimensional estimator correctly selects the non-zero
coefficients in the main equation with high probability) are established.
\textcolor{black}{These results imply that the estimate from performing
the Lasso-type procedures in the third-stage estimation is }$l_{2}-$\textcolor{black}{consistent
as long as }$\beta^{*}$\textcolor{black}{{} is }$l_{q_{2}}-$sparse
with\textit{\textcolor{black}{{} }}$q_{2}\in[0,\,1]$\textcolor{black}{{}
but inconsistent when }$q_{2}>1$.\textcolor{black}{{} A technical issue
related to a set of high-level assumptions on the regressors for estimation
consistency and selection consistency arises in the multi-stage estimation
procedure from allowing the number of regressors in the main equation
to exceed }$n$\textcolor{black}{{} and this paper provides analysis
to verify these conditions.} These verifications also provide a finite-sample
guarantee of the population identification condition required by the
semiparametric selection models. It is worth mentioning that the multi-stage
estimator and the general results in this paper can be applied to
other high-dimensional sparse semiparametric models. Section 4.5 discusses
estimation of a certain type of high-dimensional semilinear models
with the proposed multi-stage strategy when the number of parametric
and (additive) nonparametric components are large relative to the
sample size (the details are left to one of the PhD thesis chapters
by Zhu, 2015).\textcolor{black}{{} }Statistical efficiency of the proposed
estimators is studied via lower bounds on minimax risks and the result
shows that, for a class of models with exactly sparse $\beta^{*}$,
the overall convergence rate of the estimator of the high-dimensional
linear coefficients in the main equation and the nonparametric least
squares estimator of the selection bias function matches the theoretical
lower bound up to the $(n,\, d,\, p)-$factors, and exceeds it at
most by a factor of $k_{2}^{3/2}$. This statistical efficiency result,
however, does not apply to the case where $\beta^{*}$ is approximately
sparse. 

Other theoretical contributions of this paper include establishing
the non-asymptotic counterpart of the familiar asymptotic {}``oracle''
type of results from previous literature: the estimator of the coefficients
in the main equation behaves as if the unknown nonparametric component
were known, provided the nonparametric component is sufficiently smooth.
This new {}``oracle'' result holds for a unified framework of nonparametric
least squares estimators and regularized nonparametric least squares
estimators considered in the second-stage estimation. In general,
for a semiparametric model with two additive components one parametric
and the other nonparametric, when the prediction error or the $\sqrt{MSE}$
(the square root of the mean squared error) of the nonparametric estimation
\textit{per se} is $O_{p}(t_{n})$, this paper shows that the error
arising from not knowing the functional form contributes $O_{p}\left(t_{n}^{2}\right)$
in the $l_{2}-$error of the estimator of $\beta^{*}$. The driver
behind this {}``oracle'' result lies in the projection strategy.
An application of this general result to classical low-dimensional
semilinear models would imply that the nonparametric component needs
to be estimated at a rate \textit{no slower} than $O\left((\frac{1}{n})^{\frac{1}{4}}\right)$
in order for the estimator of the parametric component to achieve
the rate of $O\left(\sqrt{\frac{1}{n}}\right)$. In contrast to the
semilinear models, the low-dimensional selection models require the
rate of the nonparametric component to be at least $O\left((\frac{1}{n})^{\frac{1}{3}}\right)$
because the nonparametric component in the selection model involves
the unknown parameters $\theta^{*}$ that also need to be estimated. 

The high-dimensional multi-stage procedure is illustrated with an
application to the retail gasoline market in the Greater Saint Louis
area. Gasoline stations choose to be one of the two types: a two-product
station offering both self-service and full-service gasoline or a
single-product station offering only full-service or self-service
gasoline. A two-product station, by charging different prices for
full- and self-service gasoline, induces consumers with different
valuations to self-select the product that is consistent with their
preferences. In other words, a two-product station engages in price-discrimination.
A single-product station, on the other hand, is unable to price discriminate.
Similar to Iyer and Seetharaman (2003), this paper models a retailer\textquoteright{}s
incentive to price discriminate by choosing either single-product
or multi-product as a function of market and station characteristics
and then models the retailer\textquoteright{}s pricing decision, conditional
on the choice of the product type. However, Iyer and Seetharaman (2003)
did not account for interactions between the gas stations in their
empirical analysis. This paper uses geographic information and spatial
data to introduce, in the main equation related to the retailers'
pricing decisions, a set of variables that are high-dimensional to
control for interactions between the gas stations and employ a proposed
estimator to identify the competitive market structure. In contrast
to other heuristic ways of defining competitive markets as typically
seen in the retail gasoline industry literature, the proposed method
in this paper is natural and data-driven. The empirical finding highlights
the importance of accounting for potential interactions between stations
and suggests that competition effects from retailers that are not
in the same local market should not be overlooked.

Section 2 presents identification assumptions required for model (\ref{eq:1})
in high-dimensional settings. The estimation procedures are introduced
in Section 3. Theoretical results are established in Section 4. Small-sample
performance of the proposed multi-step high-dimensional estimator
is evaluated with Monte-Carlo simulations in Section 5 and applied
to the retail gasoline market in Section 6. Section 7 concludes this
paper. Proofs of the main results are collected in Appendix I, with
the remaining proofs of technical lemmas contained in Appendix II.
The appendices are included in the online supplementary material (\href{https://sites.google.com/site/yingzhu1215/home/JobMar_Proofs.pdf}{https://sites.google.com/site/yingzhu1215/home/JobMar\_{}Proofs.pdf}).

\section{Identification assumptions}

\textbf{Notation}. The $l_{q}$ norm of a vector $v\in p\times1$
is denoted by $\left|v\right|{}_{q}$, $1\leq q\leq\infty$ where
$\left|v\right|{}_{q}:=\left(\sum_{i=1}^{p}|v_{i}|^{q}\right)^{1/q}$
when $1\leq q<\infty$ and $\left|v\right|{}_{q}:=\max_{i=1,...,p}|v_{i}|$
when $q=\infty$. For a matrix $A\in\mathbb{R}^{p\times p}$, write
$\left|A\right|{}_{\infty}:=\max_{i,j}|a_{ij}|$ to be the elementwise
$l_{\infty}-$norm of $A$. The $l_{2}-$operator norm, or spectral
norm of the matrix $A$ corresponds to its maximum singular value;
i.e., it is defined as $\left\Vert A\right\Vert {}_{2}:=\sup_{v\in S}\left|Av\right|{}_{2}$,
where $S=\{v\in\mathbb{R}^{p}\,\vert\,\left|v\right|{}_{2}=1\}$.
The $l_{\infty}$ matrix norm (maximum absolute row sum) of $A$ is
denoted by $\left\Vert A\right\Vert {}_{\infty}:=\max_{i}\sum_{j}|a_{ij}|$
(note the difference between $\left|A\right|{}_{\infty}$ and $\left\Vert A\right\Vert {}_{\infty}$).
For a square matrix $A$, denote its minimum eigenvalue and maximum
eigenvalue by $\lambda_{\min}(A)$ and $\lambda_{\max}(A)$, respectively.
The $\mathcal{L}_{2}(\mathbb{P})-$error of a vector $\Delta(x)$,
denoted by $\left|\Delta\right|_{\mathcal{L}_{2}(\mathbb{P})}$, is
given by $\left[\mathbb{E}_{\mathbb{X}}(\Delta(x))^{2}\right]^{\frac{1}{2}}$.
Define $\mathbb{P}_{n}:=\frac{1}{n}\sum_{i=1}^{n}\delta_{x_{i}}$
that places a weight $\frac{1}{n}$ on each observation $x_{i}$ for
$i=1,...,n$, and the associated $\mathcal{L}_{2}(\mathbb{P}_{n})-$norm
of the vector $\Delta:=\left\{ \Delta(x_{i})\right\} _{i=1}^{n}$,
denoted by $\left|\Delta\right|_{n}$, is given by $\left[\frac{1}{n}\sum_{i=1}^{n}\left(\Delta(x_{i})\right)^{2}\right]^{\frac{1}{2}}$.
For a vector $v\in\mathbb{R}^{p}$, let $J(v)=\{j\in\{1,...,p\}\,\vert\, v_{j}\neq0\}$
be its support, i.e., the set of indices corresponding to its non-zero
components $v_{j}$. The cardinality of a set $J\subseteq\{1,...,p\}$
is denoted by $|J|$. For functions $f(n)$ and $g(n)$, write $f(n)\succsim g(n)$
to mean that $f(n)\geq cg(n)$ for a universal constant $c\in(0,\,\infty)$
and similarly, $f(n)\precsim g(n)$ to mean that $f(n)\leq c^{'}g(n)$
for a universal constant $c^{'}\in(0,\,\infty)$, and $f(n)\asymp g(n)$
when $f(n)\succsim g(n)$ and $f(n)\precsim g(n)$ hold simultaneously.
Also denote $\max\{a,\, b\}$ by $a\vee b$ and $\min\{a,\, b\}$
by $a\wedge b$.\\
\\
The following assumptions are imposed on model (\ref{eq:1}). \\
\\
\textbf{Assumption 2.1} (Sampling): The data \textbf{$\{y_{1i},\, y_{2i},\, w_{i},\, x_{i}\}$
}are \textit{i.i.d.} with finite second moments.\\
\\
\textbf{Remark}. The identicalness of \textbf{$\{y_{1i},\, y_{2i},\, w_{i},\, x_{i}\}$}
in Assumption 2.1 can be relaxed with the condition that \textbf{$\{\epsilon_{1i},\,\eta_{i}\}$}
are identically distributed but $\{w_{i},\, x_{i}\}$ are not. \\
\\
\textbf{Assumption 2.2} (Sparsity): The coefficient vector $\beta^{*}\in\mathbb{R}^{p}$
belongs to the $l_{q_{2}}-${}``balls'' $\mathcal{B}_{q_{2}}^{p}(R_{q_{2}})$
for a {}``radius'' of $R_{q_{2}}$ and some $q_{2}\in[0,\,1]$,
where the $l_{q}-${}``balls'' of {}``radius'' $R$ for $q\in[0,\,1]$
are defined by 
\begin{eqnarray*}
\mathcal{B}_{q}^{p}(R) & := & \left\{ \beta\in\mathbb{R}^{p}\,\vert\,\left|\beta\right|_{q}^{q}=\sum_{j=1}^{p}\left|\beta_{j}\right|^{q}\leq R\right\} \;\textrm{for }q\in(0,\,1]\\
\mathcal{B}_{0}^{p}(R) & := & \left\{ \beta\in\mathbb{R}^{p}\,\vert\,\left|\beta\right|_{0}=\sum_{j=1}^{p}\mathbb{I}\left[\beta_{j}\neq0\right]\leq R\right\} \;\textrm{for }q=0.
\end{eqnarray*}
\textbf{Remark}.\textbf{ }Assumption 2.2 requires the coefficient
vector to be {}``sparse''. As one might expect, if the high-dimensional
model lacks any additional structure, then there is no hope of obtaining
consistent estimators of $\beta^{*}$ when the ratio $\frac{p}{n}$
stays bounded away from $0$. For this reason, when working in settings
in which $p>n$, it is necessary to impose some type of sparsity assumptions
on the unknown coefficient vector $\beta^{*}\in\mathbb{R}^{p}$. Assumption
2.2 formalizes the sparsity condition by considering the $l_{q}-${}``balls''
$\mathcal{B}_{q}^{p}(R_{q})$ of {}``radius'' $R_{q}$ where $q\in[0,\,1]$.
The exact sparsity on $\beta^{*}$ corresponds to the case of $q=q_{2}=0$
with $R_{q_{2}}=k_{2}$ (in this paper, the subscript {}``$2$''
is generally reserved for the main-equation related parameters and
the subscript {}``$1$'' for the selection-equation related parameters),
which says that $\beta^{*}$ has at most $k_{2}$ non-zero components,
where the sparsity parameter $k_{2}$ is also allowed to increase
to infinity with $n$ but slowly compared to $n$. In the more general
setting $q_{2}\in(0,\,1]$, membership in $\mathcal{B}_{q_{2}}^{p}(R_{q_{2}})$
has various interpretations and one of them involves how quickly the
ordered coefficients decay. When $q_{2}\in[0,\,1),$ the set $\mathcal{B}_{q_{2}}^{p}(R_{q_{2}})$
is non-convex and the $l_{1}-$ball is the \textit{closest convex}
approximation of these non-convex sets. In terms of algorithm design,
the idea of approximating non-convex problems with their closest convex
member (so called {}``convex relaxation'') provides a tremendous
computational advantage. This is one of the reasons for favoring the
$l_{1}-$penalization techniques such as the Lasso (in solving high-dimensional
regression problems with sparsity described by the $l_{q}-$constraint
where $q\in[0,\,1]$) over estimators based on the $l_{q}-$penalty
with $q\in[0,\,1)$ which are computationally more difficult (see
the Bridge estimator in Huang, Horowitz, and Ma, 2008 as an example
of these nonconvex penalization procedures) and estimators based on
$l_{q}-$penalty with $q>1$ (such as the ridge-penalty) which are
not the closest convex approximations. On the other hand, if the coefficient
vector belongs to an Euclidean ball (the $l_{2}-$ball), then it would
make more sense to apply a ridge penalty. The focus of this paper
is on high-dimensional sparse $\beta^{*}$ that belongs to $\mathcal{B}_{q_{2}}^{p}(R_{q_{2}})$
for $q_{2}\in[0,\,1]$.\textbf{}\\
\textbf{}\\
\textbf{Assumption 2.3} (Restricted Identifiability): For a subset
$S\subseteq\{1,\,2,...,p\}$ and all non-zero $\Delta\in\mathbb{C}(S;\, q_{2},\,\varphi)\cap\mathbb{S}_{\delta}$
where 
\[
\mathbb{C}(S;\, q_{2},\,\varphi):=\left\{ \Delta\in\mathbb{R}^{p}\,:\,|\Delta_{S^{c}}|_{1}\leq\varphi|\Delta_{S}|_{1}+(\varphi+1)|\beta_{S^{c}}^{*}|_{1}\right\} \quad\textrm{for some constant \ensuremath{\varphi}}\geq1,
\]
(with $\Delta_{S}$ denoting the vector in $\mathbb{R}^{p}$ that
has the same coordinates as $\Delta$ on $S$ and zero coordinates
on the complement $S^{c}$ of $S$) and 
\[
\mathbb{S}_{\delta}:=\left\{ \Delta\in\mathbb{R}^{p}\,:\,\left|\Delta\right|_{2}\geq\delta\right\} ,
\]
the matrix $\mathbb{E}\left[y_{1i}v_{i}^{T}v_{i}\right]$ satisfies
\[
\frac{\Delta^{T}\mathbb{E}\left[y_{1i}v_{i}^{T}v_{i}\right]\Delta}{\left|\Delta\right|_{2}^{2}}\geq\kappa_{L}>0,
\]
where 
\[
v_{i}=\left(x_{i1}-\mathbb{E}\left(x_{i1}\vert w_{i}^{T}\theta^{*}\right),\,...\,,x_{ip}-\mathbb{E}\left(x_{ip}\vert w_{i}^{T}\theta^{*}\right)\right).
\]
\textbf{Remark}. Assumption 2.3 is the high-dimensional counterpart
of the familiar identification assumption in the low-dimensional selection
model literature (e.g., Powell 1989; Newey, 1991; Ahn and Powell,
1993), which assumes the matrix $\mathbb{E}\left[y_{1i}v_{i}^{T}v_{i}\right]$
is positive definite uniformly over all $\Delta\in\mathbb{R}^{p}\backslash\{\mathbf{0}\}$.
When $v_{i}$ is a zero-mean Gaussian matrix with covariance $\mathbb{E}\left[y_{1i}v_{i}^{T}v_{i}\right]=\sigma^{2}I_{p\times p}$,
the smallest eigenvalue of $\mathbb{E}\left[y_{1i}v_{i}^{T}v_{i}\right]$
is $\sigma^{2}$, so the traditional identification condition in the
low-dimensional case naturally carries to the high-dimensional case.
However, for more general structures on $\mathbb{E}\left[y_{1i}v_{i}^{T}v_{i}\right]$,
while this traditional identification condition is plausible for small
$p$, it may become harder to be satisfied when $p$ is large. Assumption
2.3 relaxes the uniform positive definiteness but only requires it
to hold over a restricted set $\mathbb{C}(S;\, q_{2},\,\varphi)\cap\mathbb{S}_{\delta}$
so that the special case of $x_{i}\subset w_{i}$ is allowed even
in the high-dimensional settings (the choices of $\delta$ and $S$
will be made clear in Section 4 when the theoretical results are presented.).
If $x_{i}\subset w_{i}$, Assumption 2.3 says that for any non-zero
vector $\lambda\in\mathbb{C}(S;\, q_{2},\,\varphi)\cap\mathbb{S}_{\delta}$,
there is no measurable function $f(w_{i}^{T}\theta^{*})$ such that
$x_{i}^{T}\lambda=f(w_{i}^{T}\theta^{*})$ when $y_{1i}=1$. Consequently,
there is at least one component $w_{ij}$ with $\theta_{j}^{*}$ in
the support set of $\theta^{*}$ (namely, the set of non-zero components
in $\theta^{*}$) such that $w_{ij}$ is excluded from $x_{i}$. This
necessary condition is the high-dimensional extension of the familiar
{}``exclusion restriction'' condition in the low-dimensional selection
model literature. 

When $\beta^{*}$ is exactly sparse (namely, $q_{2}=0$), we can take
$\delta=0$ and choose $S=J(\beta^{*})$ (where $J(\beta^{*})$ denotes
the support of $\beta^{*}$), which reduces the set $\mathbb{C}(S;\, q_{2},\,\varphi)\cap\mathbb{S}_{\delta}$
to the following cone: 
\[
\mathbb{C}(J(\beta^{*});\,0,\,\varphi):=\left\{ \Delta\in\mathbb{R}^{p}\,:\,|\Delta_{J(\beta^{*})^{c}}|_{1}\leq\varphi|\Delta_{J(\beta^{*})}|_{1}\right\} .
\]
The sample analog of Assumption 2.3 over the cone $\mathbb{C}(J(\beta^{*});\,0,\,\varphi)$
is the so-called \textit{restricted eigenvalue} condition on the Gram
matrix $\frac{v^{T}v}{n}$, studied in Bickel, et. al. (2009), Meinshausen
and Yu (2009), Raskutti, et al. (2010), Bühlmann and van de Geer (2011),
Loh and Wainwright (2012), Negahban, et. al. (2012), etc. Note that
in the low-dimensional setting where $p<n$, as long as $\textrm{rank}(v)=p$,
we are guaranteed that the Gram matrix $\frac{v^{T}v}{n}$ is positive
definite. In the high-dimensional setting with $p>n$, the matrix
$\frac{v^{T}v}{n}$ is a $p\times p$ matrix with rank at most $n$,
so it is impossible to have the uniform positive definiteness. It
is well-known that the \textit{restricted eigenvalue }assumption,
defined more precisely below, is a sufficient condition for the $l_{2}$-
consistency of the Lasso estimator for the sparse linear regression
models in high-dimensional settings. To motivate the restricted set
$\mathbb{C}(J(\beta^{*});\,0,\,\varphi)$, note that the vectors $\Delta$
in this cone have a substantial part of their {}``mass'' concentrated
on a set of the cardinality of $J(\beta^{*})$. The vectors $\Delta$
of interests often concern the error $\hat{\beta}-\beta^{*}$ where
$\hat{\beta}$ is some estimate of $\beta^{*}$. When the high-dimensional
sparse linear regression models are estimated by the $l_{1}-$penalized
techniques, an appropriate choice of the regularization parameter
would generally ensure the error $\hat{\beta}-\beta^{*}$ to be in
this restricted set. 

The following discussion provides a review of the \textit{restricted
eigenvalue} condition in literature for consistent estimation of both
exactly sparse and approximately sparse regression models in high-dimensional
settings with the Lasso or Dantzig selector. Consider the high-dimensional
sparse linear models 
\begin{equation}
y_{i}=x_{i}^{T}\beta^{*}+\epsilon_{i}=\sum_{j=1}^{p}x_{ij}\beta_{j}^{*}+\epsilon_{i},\quad i=1,\,...,\, n,\label{eq:4}
\end{equation}
where $\mathbb{E}(x_{i}\epsilon_{i})=\mathbf{0}$ for $i=1,...,n$.
Assume $p$, the number of regressors, in the above equation grows
with and exceeds the sample size $n$. Again, the focus here is the
class of sparse models with $\beta^{*}\in\mathcal{B}_{q}^{p}(R_{q})$
for $q\in[0,\,1]$. The Lasso procedure is a combination of the residual
sum of squares and a $l_{1}-$regularization defined by the following
program 
\begin{equation}
\hat{\beta}_{Las}\in\arg\min_{\beta\in\mathbb{R}^{p}}\left\{ \frac{1}{2n}|y-X\beta|_{2}^{2}+\lambda_{n}|\beta|_{1}\right\} ,\label{eq:5}
\end{equation}
where $\lambda_{n}>0$ is some regularization or tuning parameter.
Denote the minimizer to the above program by $\hat{\beta}_{Las}$.
A necessary and sufficient condition of $\hat{\beta}_{Las}$ is that
$0$ belongs to the subdifferential of the convex function $\beta\mapsto\frac{1}{2n}|y-X\beta|_{2}^{2}+\lambda_{n}|\beta|_{1}$.
This implies that the Lasso solution $\hat{\beta}_{Las}$ satisfies
the constraint 
\[
\left|\frac{1}{2n}X^{T}(y-X\hat{\beta}_{Las})\right|_{\infty}\leq\lambda_{n}.
\]
The Dantzig selector of the linear regression function is defined
as a vector having the smallest $l_{1}-$norm among all $\beta$ satisfying
the above constraint, i.e., 
\[
\hat{\beta}_{Dan}\in\arg\min\left\{ |\beta|_{1}:\left|\frac{1}{2n}X^{T}(y-X\beta)\right|_{\infty}\leq\lambda_{n}\right\} .
\]
Under the exact sparsity assumption, Bickel et al., 2009 shows that
the Lasso and the Dantzig selector exhibit similar behavior. 

In the high-dimensional setting, a sufficient condition for the $l_{2}$-
consistency of the Lasso estimator $\hat{\beta}_{Las}$ is the \textit{restricted
eigenvalue} (RE) condition related to the positive definiteness of
the Gram matrix $\frac{X^{T}X}{n}$ over a restricted set (see, e.g.,
Bickel, et. al., 2009; Meinshausen and Yu, 2009; Raskutti, et al.,
2010; Bühlmann and van de Geer, 2011; Loh and Wainwright 2012; Negahban,
et. al., 2012; etc.). Consider the following definition of the RE
condition given by Negahban, et. al. (2012) and Wainwright (2015).\\
\textbf{}\\
\textbf{Definition 2.1} (RE condition). For $q\in[0,\,1]$, the matrix
$X\in\mathbb{R}^{n\times p}$ satisfies the RE condition over a subset
$S\subseteq\{1,\,2,...,p\}$ with parameters $(q,\,\delta,\,\kappa,\,\varphi)$
if 
\begin{equation}
\frac{\frac{1}{n}|X\Delta|_{2}^{2}}{|\Delta|_{2}^{2}}\geq\kappa>0\qquad\textrm{for all nonzero}\;\Delta\in\mathbb{C}(S;\, q,\,\varphi)\cap\mathbb{S}_{\delta},\label{eq:6}
\end{equation}
where $\mathbb{C}(S;\, q,\,\varphi)\cap\mathbb{S}_{\delta}$ is defined
in Assumption 2.3. \\
\\
As discussed previously, when the unknown vector $\beta^{*}\in\mathbb{R}^{p}$
is exactly sparse, the set $\mathbb{C}(S;\, q,\,3)\cap\mathbb{S}_{\delta}$
is reduced to the cone $\mathbb{C}(J(\beta^{*});\,0,\,3)$. When $\beta^{*}$
is approximately sparse (namely, $q\in(0,\,1]$), in sharp contrast
to the case of exact sparsity, the set $\mathbb{C}(S;\, q,\,3)$ is
no longer a cone but rather contains a ball centered at the origin.
As a consequence, it is never possible to ensure that $\frac{\left|X\Delta\right|_{2}^{2}}{n}$
is bounded from below for all vectors $\Delta$ in the set $\mathbb{C}(S;\, q,\,3)$
(see Negahban, et. al., 2012 for a geometric illustration of this
issue). For this reason, in order to obtain a general applicable theory,
it is crucial to further restrict the set $\mathbb{C}(S;\, q,\,3)$
for $q\in(0,\,1]$ by introducing the set 
\[
\mathbb{S}_{\delta}:=\left\{ \Delta\in\mathbb{R}^{p}\,:\,\left|\Delta\right|_{2}\geq\delta\right\} ,
\]
where $\delta>0$ is some parameter depending on the choice of the
regularization parameter $\lambda_{n}$ in the Lasso program (\ref{eq:5}).
Provided the parameter $\delta$ and the set $S$ are suitably chosen,
the intersection $\mathbb{C}(S;\, q,\,3)\cap\mathbb{S}_{\delta}$
excludes many {}``flat'' directions (with eigenvalues of $0$) in
the space for the case of $q\in(0,\,1]$. To the best of my knowledge,
the necessity of this additional set $\mathbb{S}_{\delta}$, essential
for the approximately sparse case of $q\in(0,\,1]$, is first recognized
explicitly in Negahban, et. al. (2012).

Raskutti et al. (2010) shows that the RE condition (\ref{eq:6}) is
satisfied by the design matrix $X\in\mathbb{R}^{n\times p}$ formed
by independently sampling each row $X_{i}\sim N(0,\,\Sigma)$. Rudelson
and Zhou (2011) as well as Loh and Wainwright (2012) extend the verification
of the RE condition from the case of Gaussian designs to the case
of sub-Gaussian designs. The sub-Gaussian assumption says that the
explanatory variables need to be drawn from distributions with well-behaved
tails like Gaussian. In contrast to the Gaussian assumption, sub-Gaussian
variables constitute a more general family of distributions. In this
paper, we make use of the following definition for a sub-Gaussian
matrix. \textbf{}\\
\textbf{}\\
\textbf{Definition 2.2} (Sub-Gaussian variables and matrices). A random
variable $X$ with mean $\mu=\mathbb{E}[X]$ is sub-Gaussian if there
is a positive number $\sigma$ such that 
\[
\mathbb{E}[\exp(t(X-\mu))]\leq\exp(\sigma^{2}t^{2}/2)\qquad\textrm{for all}\, t\in\mathbb{R},
\]
and a random matrix $A\in\mathbb{R}^{n\times p}$ is sub-Gaussian
with parameters $(\Sigma_{A},\,\sigma_{A}^{2})$ if (a) each row $A_{i}\in\mathbb{R}^{p}$
is sampled independently from a zero-mean distribution with covariance
$\Sigma_{A}$, (b) for any unit vector $u\in\mathbb{R}^{p}$, the
random variable $u^{T}A_{i}^{T}$ is sub-Gaussian with parameter at
most $\sigma_{A}^{2}$.\\
\\
When applying the proposed multi-stage procedure in this paper to
estimate the high-dimensional selection models, there is no guarantee
that the random matrix $\frac{\hat{v}^{T}\hat{v}}{n}$ (where $\hat{v}_{i}$
are the estimates of $v_{i}=x_{i}-\mathbb{E}\left(x_{i}\vert w_{i}^{T}\theta^{*}\right)$
for $i=1,\,...,\, n$) would automatically satisfy these previously
established conditions for estimation consistency. For a broad class
of sub-Gaussian matrices formed by the true residuals $v_{i}=x_{i}-\mathbb{E}\left(x_{i}\vert w_{i}^{T}\theta^{*}\right)$
for $i=1,\,...,\, n$ whenever $y_{1i}=1$, this paper provides results
that imply the RE condition (\ref{eq:6}) holds for $\hat{v}^{T}\hat{v}$
with high probability provided Assumption 2.3 is satisfied. Verifications
of the RE condition provide a finite-sample guarantee of Assumption
2.3 when the unknown residuals $v$ are replaced with their estimate
$\hat{v}$ and the expectation is replaced with a sample average.
\\
\\
While the RE assumption is a natural sufficient condition for analyzing
$l_{2}-$consistency of the Lasso estimator $\hat{\beta}_{Las}$,
$l_{2}-$consistency of the Dantzig selector $\hat{\beta}_{Dan}$
can be related to a different sufficient condition, the \textit{sensitivity
characteristics}, on the term $|X^{T}Xv|_{\infty}$. These sensitivity
characteristics were originally introduced in Ye and Zhang (2010)
as the cone invertibility factors and used in Gautier and Tsybakov
(2011) for high-dimensional instrumental variable regressions. Gautier
and Tsybakov (2011) shows that the sensitivity characteristics can
be larger than the usual RE condition of Bickel, et. al (2009) and
therefore the Dantzig-type estimators may lead to better results in
certain cases%
\footnote{Recently, another weaker version of the RE condition tailored to the
square-root Lasso is developed in Belloni, Chernozhukov, and Wang
(2014).%
}. The analysis of a pivotal Dantzig selector in this paper for estimating
the high-dimensional linear coefficients relies on the following definition
based on Gautier and Tsybakov (2011):\\
\\
\textbf{Definition 2.3} ($l_{2}-$\textit{sensitivity}). The matrix
$X\in\mathbb{R}^{n\times p}$ satisfies the $l_{2}-$\textit{sensitivity}
condition over a subset $S\subseteq\{1,\,2,...,p\}$ with parameters
$(q,\,\delta,\,\kappa^{'},\,\varphi)$ if 
\begin{equation}
\frac{\frac{1}{n}|X^{T}X\Delta|_{\infty}^{2}}{|v|_{2}^{2}}\geq\kappa^{'}>0\qquad\textrm{for all nonzero}\;\Delta\in\mathbb{C}(S;\, q,\,\varphi)\cap\mathbb{S}_{\delta}\label{eq:7}
\end{equation}
where $\mathbb{C}(S;\, q,\,\varphi)\cap\mathbb{S}_{\delta}$ is defined
in Assumption 2.3.\\
\\
When $y_{i}$ in (\ref{eq:4}) is a latent variable with only an observable
sign, other models such as the high-dimensional binary response models
may be considered. In analyzing these models, the RE condition can
be replaced with a similar notion, the \textit{restricted strong convexity}
(RSC) condition, originally formalized by Negahban, et. al. (2012)
in the context of the regularized \textit{M}-estimation with a general,
convex and differentiable loss function. The following definition
from Negahban, et. al. (2012) is adopted in this paper to analyze
the theoretical properties of an estimator for the high-dimensional
logit and probit model: \\
\\
\textbf{Definition 2.4} (RSC condition). A convex and differentiable
loss function $\mathcal{L}(\theta;\, z_{1}^{n})$ satisfies the RSC
condition over a subset $S\subseteq\{1,\,2,...,p\}$ with parameters
$(q,\,\delta,\,\kappa^{''},\,\varphi)$ where $\kappa^{''}>0$ if
\[
\mathcal{L}(\theta^{*}+\Delta;\, z_{1}^{n})-\mathcal{L}(\theta^{*};\, z_{1}^{n})-\left\langle \nabla\mathcal{L}(\theta^{*};\, z_{1}^{n}),\,\Delta\right\rangle \geq\kappa^{''}\left|\Delta\right|_{2}^{2}\quad\textrm{for all nonzero}\;\Delta\in\mathbb{C}(S;\, q,\,\varphi)\cap\mathbb{S}_{\delta}.
\]
where $\nabla\mathcal{L}(\theta^{*};\, z_{1}^{n})$ denotes the derivative
of $\mathcal{L}(\theta;\, z_{1}^{n})$ evaluated at $\theta=\theta^{*}$,
and $\mathbb{C}(S;\, q,\,\varphi)\cap\mathbb{S}_{\delta}$ is defined
in Assumption 2.3.

\section{Estimation procedures }

This section presents a 3-stage estimation procedure for the high-dimensional
linear coefficients in the main equation and two estimators of the
selection bias function. In terms of applicability, the proposed estimators
enjoy many computational advantages and can be easily implemented
using existing software packages.

\subsection{The multi-stage estimator of the high-dimensional linear coefficients}

To facilitate the presentation of the multi-stage estimator, we reverse
the order of the three stages when discussing the estimation procedure;
in particular, we will introduce the third-stage estimator and then
followed by the second-stage and the first-stage estimators. Note
that the second- and third-stage estimations concern only the selected
sample (observations with $y_{1i}=1$) and the first-stage estimation
concerns the entire sample. For the third-stage estimation, this paper
considers a non-pivotal Lasso procedure whose regularization parameter
depends on the unknown variance of $\eta_{i}$, and a pivotal Dantzig
selector (Gautier and Tsybakov, 2011) whose regularization parameter
does not involve the unknown variance of $\eta_{i}$. \\
\\
\textbf{Non-pivotal third-stage estimation}\\
Revisiting equation (\ref{eq:3}) in Section 1 suggests that if an
estimate of 
\[
\left(\mathbb{E}\left(x_{i1}\vert w_{i}^{T}\theta^{*}\right),\,...\,,\mathbb{E}\left(x_{ip}\vert w_{i}^{T}\theta^{*}\right)\right)
\]
is available to us, then we can form estimates 
\begin{eqnarray*}
\hat{v}_{i} & = & \left(x_{i1}-\hat{\mathbb{E}}\left(x_{i1}\vert w_{i}^{T}\hat{\theta}\right),\,...\,,x_{ip}-\hat{\mathbb{E}}\left(x_{ip}\vert w_{i}^{T}\hat{\theta}\right)\right),\\
\hat{v}_{i0} & = & y_{2i}-\hat{\mathbb{E}}\left(y_{2i}\vert w_{i}^{T}\hat{\theta}\right).
\end{eqnarray*}
of the nonparametric residuals 
\begin{eqnarray*}
v_{i} & = & \left(x_{i1}-\mathbb{E}\left(x_{i1}\vert w_{i}^{T}\theta^{*}\right),\,...\,,x_{ip}-\mathbb{E}\left(x_{ip}\vert w_{i}^{T}\theta^{*}\right)\right),\\
v_{i0} & = & y_{2i}-\mathbb{E}\left(y_{2i}\vert w_{i}^{T}\theta^{*}\right).
\end{eqnarray*}
Then, an estimator of the high-dimensional linear coefficients in
the main equation (the \textbf{third-stage} estimator) can be obtained
by performing the following Lasso program: 
\begin{equation}
\hat{\beta}_{HSEL}\in\textrm{argmin}_{\beta\in\mathbb{R}^{p}}:\:\frac{1}{2n}|\hat{v}_{0}-\hat{v}\beta|_{2}^{2}+\lambda_{n,3}|\beta|_{1},\label{eq:8}
\end{equation}
where $\lambda_{n,3}>0$ is some regularization parameter whose choice
is to be discussed in Section 4. In general, the choice of $\lambda_{n,3}$
depends on $\mathbb{E}(v_{ij}^{2})$ and $\mathbb{E}(\eta_{i}^{2})$.
To make $\lambda_{n,3}$ and the estimate $\hat{\beta}_{HSEL}$ independent
of the effect from $\mathbb{E}(v_{ij}^{2})$, we can impose weights
on the penalty term as follows 
\begin{equation}
\textrm{min}_{\beta\in\mathbb{R}^{p}}:\:\frac{1}{2n}|\hat{v}_{0}-\hat{v}\beta|_{2}^{2}+\lambda_{n,3}\sum_{j=1}^{p}\hat{\sigma}_{v_{j}}|\beta_{j}|,\label{eq:9}
\end{equation}
where $\hat{\sigma}_{v_{j}}:=\sqrt{\frac{1}{n}\sum_{i=1}^{n}\hat{v}_{ij}^{2}}$.
To make $\lambda_{n,3}$ not depend on the unknown variance of $\eta_{i}$,
we can consider the pivotal version of the Dantzig selector as in
Gautier and Tsybakov (2011).\\
\\
\textbf{Pivotal third-stage estimation}\\
Set $v_{j*}:=\max_{i\in\{1,...,n\}}\left(\max\left\{ \left|2x_{ij}\right|,\,\left|\hat{v}_{ij}\right|\right\} \right)$
for $j=1,\,...\,,p$ and denote $D$ the diagonal $p\times p$ matrix
with diagonal entries $v_{j*}^{-1}$. Consider the following optimization
problem:
\begin{equation}
\min_{(\beta,\,\sigma)\in A}:\:\left(\left|D^{-1}\beta\right|_{1}+C\sigma\right)\label{eq:10}
\end{equation}
where 
\[
A=\left\{ (\beta,\,\sigma):\,\beta\in\mathbb{R}^{p},\,\sigma>0,\,\frac{1}{n}\left|D\hat{v}^{T}(\hat{v}_{0}-\hat{v}\beta)\right|_{\infty}\leq\sigma\xi,\,\frac{1}{n}\left|\hat{v}_{0}-\hat{v}\beta\right|_{2}^{2}\leq\sigma^{2}\right\} 
\]
for some tuning parameter $\xi>0$ (to be specified in Section 4).
The computational aspect of this pivotal estimator is detailed in
Gautier and Tsybakov (2011).\\
\\
\textbf{Remark}. The third-stage estimation needs not to be restricted
to the Lasso or the Dantzig selector. Other methods with different
loss functions (such as the square-root Lasso in Belloni, et. al 2011,
2014) or with different penalty functions (such as the SCAD in Fan
and Li, 2001, or the MCP in Zhang, 2010) can be used. This paper focuses
on the analysis of the non-pivotal Lasso and the pivotal Dantzig selector
laid out above for the third-stage estimation.\\
\\
\textbf{Second-stage estimation}\\
To simplify the notations in the following, write $\mathbb{E}\left(x_{ij}\vert w_{i}^{T}\theta\right):=m_{j}(w_{i}^{T}\theta)$,
$\hat{\mathbb{E}}\left(x_{ij}\vert w_{i}^{T}\theta\right):=\hat{m}_{j}(w_{i}^{T}\theta)$,
$\mathbb{E}\left(y_{2i}\vert w_{i}^{T}\theta\right):=m_{0}(w_{i}^{T}\theta)$,
and $\hat{\mathbb{E}}\left(y_{2i}\vert w_{i}^{T}\theta\right):=\hat{m}_{0}(w_{i}^{T}\theta)$.
To estimate $m_{j}(w_{i}^{T}\theta^{*})$ for each $j=0,\,...,\, p$,
we first need some estimate $\hat{\theta}$ of $\theta^{*}$in the
selection equation. Supposing such an estimate is available, to obtain
a (\textbf{second-stage}) estimator of $m_{j}(w_{i}^{T}\theta^{*})$,
we consider the following least squares estimator 
\begin{equation}
\hat{m}_{j}\in\arg\min_{\tilde{m}_{j}\in\mathcal{F}_{j}}\left\{ \frac{1}{n}\sum_{i=1}^{n}\left(z_{ij}-\tilde{m}_{j}(w_{i}\hat{\theta})\right)^{2}\right\} ,\label{eq:11}
\end{equation}
or the regularized least-squares estimator 
\begin{equation}
\hat{m}_{j}\in\arg\min_{\tilde{m}_{j}\in\mathcal{F}_{j}}\left\{ \frac{1}{n}\sum_{i=1}^{n}\left(z_{ij}-\tilde{m}_{j}(w_{i}\hat{\theta})\right)^{2}+\lambda_{nj,2}\left|\tilde{m}_{j}\right|_{\mathcal{F}_{j}}^{2}\right\} ,\label{eq:12}
\end{equation}
where $\left|\cdot\right|_{\mathcal{F}_{j}}$ is a norm associated
with the function class $\mathcal{F}_{j}$ and $\lambda_{nj,2}\geq0$
is a regularization parameter and $z_{i0}=y_{2i}$ and $z_{ij}=x_{ij}$
for each $j=1,...,p$. The choice of $\lambda_{nj,2}$ is specified
in Section 4. A nonparametric regression problem based on (\ref{eq:11})
or (\ref{eq:12}) is a standard setup in many modern statistics books
(e.g., van der Vaart and Wellner, 1996; van de Geer, 2000; Wainwright,
2015, etc). 

In words, the solutions to program (\ref{eq:11}) are least-squares
estimators based on imposing explicit constraints on the function
class $\mathcal{F}_{j}$. The function $\hat{m}_{j}$ is chosen such
that the vector 
\[
\left(\hat{m}_{j}(w_{1}\hat{\theta}),\,...,\,\hat{m}_{j}(w_{i}\hat{\theta}),\,...,\,\hat{m}_{j}(w_{n}\hat{\theta})\right)
\]
is closest in $l_{2}-$norm to the observation $\left(z_{1j},\,...,\, z_{ij},\,...,\, z_{nj}\right)$
for $j=0,...,p$ in terms of the {}``selected'' sample. Examples
of (\ref{eq:11}) include the linear regression as the simplest case,
sparse linear regressions, convex regression where $\mathcal{F}_{j}$
is the class of convex functions (e.g., Guntuboyina and Sen, 2013),
Lipschitz and Isotonic regression where $\mathcal{F}_{j}$ is the
class of monotone Lipschitz functions (e.g., Kakade, Kalai, Kanade,
and Shamir, 2011), etc. In general, this optimization problem defining
the non-parametric least squares estimator $\hat{m}_{j}$ is infinite-dimensional
in nature, since $\hat{m}_{j}$ ranges over the function class $\mathcal{F}_{j}$.
If the function class is {}``too large'', the solution may not exist,
in which case $\mathcal{F}_{j}$ is chosen to be a compact subset
of some larger function class by introducing a ball radius in some
norm. From the computational point of view, it is sometimes more convenient
to implement estimators based on explicit penalization or regularization
terms as in (\ref{eq:12}). Examples of (\ref{eq:12}) include kernel
ridge regression where $\left|\cdot\right|_{\mathcal{F}_{j}}$ is
the norm associated with a reproducing kernel Hilbert space (see e.g.,
Gu, 2002; Berlinet and Thomas-Agnan, 2004; Wainwright, 2015), estimators
based on series expansion (e.g., Cencov, 1962; Andrews, 1991; Newey,
1994, 1997), as well as sieves (e.g., van de Geer, 2000; Chen, 2008)
and spline methods (e.g., Wahba, 1980, 1990). A procedure based on
Lipschitz regression for the second-stage nonparametric estimation
is illustrated in Section 4 for a leading case. 

It is worth mentioning that although the theoretical guarantees of
the multi-stage procedure provided by this paper requires the second-stage
estimation to fit into either (\ref{eq:11}) or (\ref{eq:12}), other
nonparametric methods including kernel density estimators, local polynomials,
etc., could also be a valid second-stage estimator for the multi-stage
procedure in the context of high-dimensional semiparametric selection
models and verifying those methods both theoretically and empirically
is an open question for future research.\textbf{}\\
\textbf{}\\
\textbf{First-stage estimation}\\
Note that in the second-stage estimation of $m_{j}(w_{i}^{T}\theta^{*})$
for each $j=0,...,p$, the coefficient vector $\theta^{*}$ is unknown
and needs to be replaced by some {}``consistent'' \textbf{first-stage}
estimate $\hat{\theta}$. Parametric and semiparametric estimation
of $\theta^{*}$ in the classical low-dimensional settings when the
dimension of $\theta^{*}$ is small relative to the sample size $n$
is well-studied (see, e.g., Powell 1994; Pagan and Ullah, 1999). In
the high-dimensional settings where the dimension of $\theta^{*}$
grows with and exceeds $n$, estimation of $\theta^{*}$ in recent
development of high-dimensional statistics has been focused on the
case where $\theta^{*}$ is either exactly sparse or approximately
sparse, and a distributional assumption is imposed on the error term
in the linear latent utility models in (\ref{eq:1}). Theoretical
guarantees have been established for the high-dimensional binary logit
models in the context of Generalized Linear Models (GLM) and \textit{M}-estimation
(e.g., van de Geer, 2008; Bühlmann and van de Geer, 2011; Negahban,
et. al, 2012; Loh and Wainwright, 2013). While the main theoretical
results of this paper concern estimators of the high-dimensional linear
coefficients $\beta^{*}$ in the main equation and estimators of the
selection bias function $g(\cdot)$, we illustrate here and also in
later sections the high-dimensional parametric estimation procedure
for the binary logit and probit models as they are considered the
work-horse of many empirical literatures and probit models are widely
applied to study selection problems. 

As for the high-dimensional sparse linear models, it is natural to
consider the estimator based on the $l_{1}-$\textit{regularized maximum
likelihood} for the binary logit and probit models, namely, 
\begin{equation}
\hat{\theta}\in\arg\min_{\theta\in\mathbb{R}^{d}}\left\{ -\frac{1}{n}\sum_{i=1}^{n}y_{1i}\phi_{1}(w_{i}^{T}\theta)+\frac{1}{n}\sum_{i=1}^{n}\phi_{2}(w_{i}^{T}\theta)+\lambda_{n,1}\left|\theta\right|_{1}\right\} \label{eq:13}
\end{equation}
where $n$ is the sample size of all observations. One can easily
verify that when $\phi_{1}(w_{i}^{T}\theta)=w_{i}^{T}\theta$ and
$\phi_{2}(w_{i}^{T}\theta)=\log(1+\exp(w_{i}^{T}\theta))$, the loss
function in the above program corresponds to a binary logit model;
when $\phi_{1}(w_{i}^{T}\theta)=\log\frac{\Phi(w_{i}^{T}\theta)}{1-\Phi(w_{i}^{T}\theta)}$
and $\phi_{2}(w_{i}^{T}\theta)=-\log\left[1-\Phi(w_{i}^{T}\theta)\right]$
($\Phi(\cdot)$ is the standard normal c.d.f.), the loss function
corresponds to a binary probit model. The loss function in (\ref{eq:13})
is written in terms of the negative of the likelihood and hence the
optimization program is a convex minimization problem. This paper
extends the analysis of the theoretical properties of these estimators
from the high-dimensional binary logit models to the high-dimensional
probit models, and focuses on the semiparametric estimation of $\beta^{*}$
instead of $\theta^{*}$. Developing semiparametric estimation techniques
for the high-dimensional sparse discrete choice models based upon
weak restrictions on the error distribution is left for future research.\\
\\
\textbf{Remark}. Upon solving (\ref{eq:13}), strategies such as the
thresholded-Lasso or the post-Lasso may be used before the second-stage
estimation, which might boost the performance of the multi-stage estimator
in certain situations.

\subsection{Estimators of the selection bias function}

Given the estimates $\hat{\theta}$ and $\hat{\beta}$ of $\theta^{*}$
and $\beta^{*}$, there are two ways to estimate the selection bias
function $g(w_{i}^{T}\theta^{*})$. Recalling (\ref{eq:2}) from Section
1, 
\[
y_{2i}=x_{i}^{T}\beta^{*}+g(w_{i}^{T}\theta^{*})+\eta_{i}.
\]
where by construction $\mathbb{E}\left[\eta_{i}\vert w_{i},\, x_{i},\, y_{1i}=1\right]=0$.
Taking the conditional expectation of the above leads to 
\[
\mathbb{E}\left(y_{2i}\vert w_{i}^{T}\theta^{*}\right)=\mathbb{E}\left(x_{i}\,\vert\, w_{i}^{T}\theta^{*}\right)\beta^{*}+g(w_{i}^{T}\theta^{*}),
\]
and as a result, 
\begin{equation}
g(w_{i}^{T}\theta^{*})=\mathbb{E}\left(y_{2i}\vert w_{i}^{T}\theta^{*}\right)-\mathbb{E}\left(x_{i}\,\vert\, w_{i}^{T}\theta^{*}\right)\beta^{*},\label{eq:14}
\end{equation}
where $\mathbb{E}\left(x_{i}\,\vert\, w_{i}^{T}\theta^{*}\right):=\left(\mathbb{E}\left(x_{i1}\vert w_{i}^{T}\theta^{*}\right),\,...\,,\mathbb{E}\left(x_{ip}\vert w_{i}^{T}\theta^{*}\right)\right)$.
Replacing $\mathbb{E}\left(y_{2i}\vert w_{i}^{T}\theta^{*}\right)$,
$\mathbb{E}\left(x_{i}\,\vert\, w_{i}^{T}\theta^{*}\right)$, and
$\beta^{*}$ with their estimates from Section 3.1 yields the estimator
$\hat{g}(w_{i}^{T}\hat{\theta})$ of $g(w_{i}^{T}\theta^{*})$:

\begin{equation}
\hat{g}(w_{i}^{T}\hat{\theta}):=\hat{\mathbb{E}}\left(y_{2i}\,\vert\, w_{i}^{T}\hat{\theta}\right)-\hat{\mathbb{E}}\left(x_{i}\,\vert\, w_{i}^{T}\hat{\theta}\right)\hat{\beta},\label{eq:15}
\end{equation}
where $\hat{\mathbb{E}}\left(x_{i}\,\vert\, w_{i}^{T}\hat{\theta}\right):=\left(\hat{\mathbb{E}}\left(x_{i1}\vert w_{i}^{T}\hat{\theta}\right),\,...\,,\hat{\mathbb{E}}\left(x_{ip}\vert w_{i}^{T}\hat{\theta}\right)\right)$
is the second-stage estimate of $\mathbb{E}\left(x_{i}\,\vert\, w_{i}^{T}\theta^{*}\right):=\left(\mathbb{E}\left(x_{i1}\vert w_{i}^{T}\theta^{*}\right),\,...\,,\mathbb{E}\left(x_{ip}\vert w_{i}^{T}\theta^{*}\right)\right)$. 

Alternatively, like how we obtain the second-stage estimates in Section
3.1, one can estimate $g(w_{i}^{T}\theta^{*})$ by solving the following
least-squares estimator 
\begin{equation}
\tilde{g}\in\arg\min_{f\in\mathcal{F}}\frac{1}{n}\sum_{i=1}^{n}\left(y_{2i}-x_{i}^{T}\hat{\beta}-f(w_{i}^{T}\hat{\theta})\right)^{2},\label{eq:16}
\end{equation}
or the regularized least-squares estimator 
\[
\tilde{g}\in\arg\min_{f\in\mathcal{F}}\left\{ \frac{1}{n}\sum_{i=1}^{n}\left(y_{2i}-x_{i}^{T}\hat{\beta}-f(w_{i}^{T}\hat{\theta})\right)^{2}+\lambda_{n}^{*}\left|f\right|_{\mathcal{F}}^{2}\right\} ,
\]
where $\lambda_{n}^{*}\geq0$ is a regularization parameter. Although
this alternative estimator $\tilde{g}(w_{i}^{T}\hat{\theta})$ of
$g(w_{i}^{T}\theta^{*})$ is computationally more involved relative
to the closed-form estimator $\hat{g}(w_{i}^{T}\hat{\theta})$, its
rate of convergence turns out to be faster in a leading case as we
will see in Section 4.

\section{Main theoretical results }

For notational simplicity, in the main theoretical results presented
below, we assume the regime of interest is $p\geq n$ and $d\geq n$
(i.e., the number of regressors grows with and exceed the sample size
$n$). The modification to allow $p<n$ or $d<n$ is trivial. Also,
as a general rule for this paper, all the $b$ constants denote positive
constants that are independent of $n$, $p$, $d$, $R_{q_{1}}$ and
$R_{q_{2}}$ but possibly depending on model specific parameters;
all the $c$ constants denote universal positive constants that are
independent of $n$, $p$, $d$, $R_{q_{1}}$ and $R_{q_{2}}$ as
well as model specific parameters. The specific values of these constants
may change from place to place. \\
\textbf{}\\
Recall from programs (\ref{eq:11}) and (\ref{eq:12}), $\tilde{m}_{j}(\cdot)\in\mathcal{F}_{j}$.
Suppose $m_{j}(\cdot)\in\mathcal{F}_{j}^{*}$, which may be different
from $\mathcal{F}_{j}$. Define the shifted version of the function
class $\mathcal{F}_{j}$ 
\[
\bar{\mathcal{F}}_{j}:=\left\{ f=f^{'}-f^{''}\,:\, f^{'},\, f^{''}\in\mathcal{F}_{j}\right\} .
\]
The following assumptions are imposed to obtain the theoretical results
in this section.\\
\textbf{}\\
\textbf{Assumption 4.1}: For any $j=0,...,p$, $\bar{\mathcal{F}}_{j}$
is a \textit{star-shaped} function class; i.e., for any $f\in\bar{\mathcal{F}}_{j}$,
the entire line $\left\{ \alpha f,\:\alpha\in[0,\,1]\right\} $ is
also contained within $\bar{\mathcal{F}}_{j}$. \\
\\
\textbf{Remark}. The star-shaped condition is often seen in literature
of nonparametric statistics (see e.g., van der Vaart and Wellner,
1996; Wainwright, 2015; and other textbooks on mathematical statistics).
It is relatively mild; for instance, it is satisfied whenever the
set $\bar{\mathcal{F}}_{j}$ is convex and contains the function $f=0$.
It is also satisfied by various non-convex sets of functions, such
as in the case of sparse linear regression.\\
\\
\textbf{Assumption 4.2}: The random vector $v_{j}$ for $j=0,...,p$
is sub-Gaussian with parameter at most $\sigma_{v_{j}}$. The matrix
$v\in\mathbb{R}^{n\times p}$ is sub-Gaussian with parameters $(\Sigma_{v},\,\sigma_{v}^{2})$
where the\textbf{ $j$}th column of $v$ is $v_{j}$ and $\sigma_{v}:=\max_{j=0,...,p}\sigma_{v_{j}}$.\\
\\
\textbf{Assumption 4.3}: The random vector $\eta$ is sub-Gaussian
with parameter at most $\sigma_{\eta}$. \\
\\
\textbf{Remark.} In the literature of nonparametric estimation, common
measures of function complexities associated with sub-Gaussian variables
can be controlled with \textit{standard} maximal inequalities as in
van der Vaart and Wellner (1996) and van de Geer (2000), etc. There
are some special cases of Assumptions 4.2 are 4.3 where other concentration
results (e.g., Maurey, 1991; Ledoux, 1996; Bobkov, 1999; Bobkov and
Ledoux, 2000) may provide sharper constants in the tail probability
when we relax the identicalness of \textbf{$\{w_{i},\, x_{i}\}$}
in Assumption 2.1. These special cases include: $v_{j}$ for $j=0,...,p$
and $\eta$ are (i) sub-Gaussian with \textit{strongly log-concave}
distribution (defined below) for some $\gamma_{v_{j}}>0$ and $\gamma_{\eta}>0$,
respectively; or, (ii) a bounded vector%
\footnote{A random vector with bounded elements is sub-Gaussian.%
} such that for every $i=1,...,n$, $v_{ij}$ and $\eta$ are supported
on the interval $(a_{v_{j}}^{'},\, a_{v_{j}}^{''})$ with $B_{v_{j}}:=a_{v_{j}}^{''}-a_{v_{j}}^{'}$,
and on $(a_{\eta}^{'},\, a_{\eta}^{''})$ with $B_{\eta}:=a_{\eta}^{''}-a_{\eta}^{'}$;
or, (iii) a mixture of (i) and (ii) in terms of its probability measure.
\\
\\
\textbf{Definition 4.1} (Strongly log-concave distributions). A distribution
$\mathbb{P}$ with density $\mathtt{p}$ (with respect to the Lebesgue
measure) is a \textit{strongly log-concave distribution} if the function
$\log\mathtt{p}$ is strongly concave. Equivalently stated, the density
can be written in the form $\mathtt{p}(x)=\exp\left(-\psi(x)\right)$,
where the function $\psi:\,\mathbb{R}^{n}\rightarrow\mathbb{R}$ is
strongly convex, meaning that there is some $\gamma>0$ such that
\[
\lambda\psi(x)+(1-\lambda)\psi(y)-\psi(\lambda x+(1-\lambda)y)\geq\frac{\gamma}{2}\lambda(1-\lambda)\left|x-y\right|_{2}^{2}
\]
for all $\lambda\in[0,\,1]$, and $x,\, y\in\mathbb{R}^{n}$.\\
\\
\textbf{Remark.} It is easy to verify that the distribution of a standard
Gaussian vector in $n$ dimensions is strongly log-concave with parameter
$\gamma=1$. More generally, any Gaussian distribution with covariance
matrix $\Sigma\succ0$ is strongly log-concave with parameter $\gamma=\lambda_{\min}(\Sigma^{-1})$.
In addition, there are a variety of non-Gaussian distributions that
are also strongly log-concave.

\subsection{Properties of the non-pivotal Lasso estimator of the high-dimensional
linear coefficients }

\subsubsection{General upper bounds and $l_{2}-$consistency }

The following theorem (Theorem 4.1) provides a general upper bound
on the error $\left|\hat{\beta}_{HSEL}-\beta^{*}\right|_{2}$ when
the second-stage estimation concerns a program as in (\ref{eq:11}).
This result is an {}``oracle-inequality'' type which does not assume
the unknown function $m(\cdot)$ belongs to the function class over
which the nonparametric estimator from (\ref{eq:11}) is defined.
In such settings, the performance of the estimator involves both the
estimation error and an approximation error, arising from the fact
that $m_{j}\notin\mathcal{F}_{j}$.

To state Theorem 4.1, we need to introduce a set of definitions. First,
we define a quantity that measures the complexity of the function
class $\mathcal{F}_{j}$ (a notion often used in nonparametric literature;
e.g., van der Vaart and Wellner, 1996; van de Geer, 2000; Barlett
and Mendelson, 2002; Koltchinski, 2006; Wainwright, 2015, etc.). For
any radius $r_{j}>0$, define the conditional \textit{local complexity}
\[
\mathcal{G}_{n}(r_{j};\,\mathcal{F}_{j}):=\mathbb{E}_{v_{j}}\left[\sup_{f\in\Omega(r_{j};\,\mathcal{F}_{j})}\left|\frac{1}{n}\sum_{i=1}^{n}v_{ij}f(w_{i}^{T}\theta^{*})\right|\vert w_{i}^{T}\theta^{*}\right],
\]
where variables $\left\{ v_{ij}\right\} _{i=1}^{n}$ for $j=0,...,p$
are \textit{i.i.d.} variates that satisfy Assumption 4.2, and 
\[
\Omega(r_{j};\,\mathcal{F}_{j})=\left\{ f\,:\, f\in\bar{\mathcal{F}}_{j}\,\left|f_{\theta^{*}}\right|_{n}\leq r_{j}\right\} ,
\]
where $\left|f_{\theta^{*}}\right|_{n}:=\sqrt{\frac{1}{n}\sum_{i=1}^{n}\left[f(w_{i}^{T}\theta^{*})\right]^{2}}$.
For any star-shaped shifted function class $\bar{\mathcal{F}}_{j}$,
the function $t\mapsto\frac{\mathcal{G}_{n}(t;\,\mathcal{F}_{j})}{t}$
is non-decreasing on the interval $(0,\,\infty)$. Second, let $T_{j}^{*}:=\sup_{f\in\mathcal{F}_{j}^{*}}\frac{1}{n}\sum_{i=1}^{n}\left[f(w_{i}^{T}\hat{\theta})-f(w_{i}^{T}\theta^{*})\right]^{2}$,
$T_{j}:=\sup_{f\in\mathcal{F}_{j}}\frac{1}{n}\sum_{i=1}^{n}\left[f(w_{i}^{T}\hat{\theta})-f(w_{i}^{T}\theta^{*})\right]^{2}$,
$T_{j}^{'}:=T_{j}^{*}\vee T_{j}$, and 
\begin{eqnarray*}
\mathcal{T}_{1} & = & \max_{j\in\{0,...,p\}}\left(T_{j}^{'}\vee\sqrt{T_{j}^{'}}\right)\\
\mathcal{T}_{2} & = & \max_{j\in\{0,...,p\}}t_{nj}^{2}\\
\mathcal{T}_{3} & = & \max_{j\in\{0,...,p\}}\inf_{\tilde{m}_{j}\in\mathcal{F}_{j}}\left(\frac{1}{n}\sum_{i=1}^{n}\left[\tilde{m}_{j}(w_{i}^{T}\hat{\theta})-m_{j}(w_{i}^{T}\hat{\theta})\right]^{2}+\sqrt{\frac{1}{n}\sum_{i=1}^{n}\left[\tilde{m}_{j}(w_{i}^{T}\theta^{*})-m_{j}(w_{i}^{T}\theta^{*})\right]^{2}}\right)\\
\mathcal{T}_{4} & = & \sigma_{v}\sigma_{\eta}\sqrt{\frac{\log p}{n}}.
\end{eqnarray*}
Third, recall in Section 2 the set we introduced:
\[
\mathbb{C}(S;\, q_{2},\,3):=\left\{ \Delta\in\mathbb{R}^{p}\,:\,|\Delta_{S^{c}}|_{1}\leq3|\Delta_{S}|_{1}+4|\beta_{S^{c}}^{*}|_{1}\right\} ,
\]
and the spherical set 
\[
\mathbb{S}_{\delta}:=\left\{ \Delta\in\mathbb{R}^{p}\,:\,\left|\Delta\right|_{2}\geq\delta\right\} ,
\]
and the intersection of these two sets $\mathbb{C}(S;\, q_{2},\,3)\cap\mathbb{S}_{\delta}$.
When $\beta^{*}$ is approximately sparse (namely, $q_{2}\in(0,\,1]$),
we choose $S$ in $\mathbb{C}(S;\, q_{2},\,3)$ to be the following
thresholded subset 
\[
S_{\underbar{\ensuremath{\tau}}}:=\left\{ j\in\left\{ 1,\,2,\,...,\, p\right\} \,:\,\left|\beta_{j}^{*}\right|>\underbar{\ensuremath{\tau}}\right\} 
\]
with the threshold parameter $\underbar{\ensuremath{\tau}}=\frac{\lambda_{n,3}}{\kappa_{L}}$
(recall $\lambda_{n,3}$ is the third-stage regularization parameter
whose choice is specified in the theorems and the parameter $\kappa_{L}$
is defined in Assumption 2.3, Section 2). When $\beta^{*}$ is exactly
sparse (namely, $q_{2}=0$), we set $\delta=\underbar{\ensuremath{\tau}}=0$
and choose $S=J(\beta^{*})$, which reduces the set $\mathbb{C}(S;\, q_{2},\,3)\cap\mathbb{S}_{\delta}$
to the following cone: 
\[
\mathbb{C}(J(\beta^{*});\,0,\,3):=\left\{ \Delta\in\mathbb{R}^{p}\,:\,|\Delta_{J(\beta^{*})^{c}}|_{1}\leq3|\Delta_{J(\beta^{*})}|_{1}\right\} .
\]
\textbf{Theorem 4.1}: Let the \textit{critical radius} $r_{nj}>0$
be the smallest positive quantity satisfying the \textit{critical
inequality} 
\[
\mathcal{G}_{n}\left(r_{nj};\,\mathcal{F}_{j}\right)\leq\frac{r_{nj}^{2}}{\sigma_{v_{j}}}.
\]
Suppose the second-stage estimator solves program (\ref{eq:11}) and
Assumptions 2.1, 2.2, 4.1-4.3 hold. Additionally, let Assumption 2.3
hold over the restricted set $\mathbb{C}(J(\beta^{*});\,0,\,3)$ for
the exact sparsity case ($q_{2}=0$ with $R_{q_{2}}=k_{2}$), and
over $\mathbb{C}(S_{\underbar{\ensuremath{\tau}}};\, q_{2},\,3)\cap\mathbb{S}_{\delta}$
where $\delta\asymp R_{q_{2}}^{\frac{1}{2}}\left(\lambda_{n,3}\right)^{1-\frac{q_{2}}{2}}$
and $\underbar{\ensuremath{\tau}}=\frac{\lambda_{n,3}}{\kappa_{L}}$
for the approximate sparsity case ($q_{2}\in(0,\,1]$), respectively.
For any $t_{nj}\geq r_{nj}$, if the third-stage regularization parameter
$\lambda_{n,3}$ satisfies 
\begin{equation}
\lambda_{n,3}\geq b(\sigma_{v},\,\sigma_{\eta})\left|\beta^{*}\right|_{1}\left(\mathcal{T}_{1}+\mathcal{T}_{2}+\mathcal{T}_{3}\right)+\mathcal{T}_{4}:=\bar{M},\label{eq:17}
\end{equation}
where $b(\sigma_{v},\,\sigma_{\eta})$ is a known function that only
depend on the parameters $\sigma_{v}$ and $\sigma_{\eta}$ (and independent
of $n$, $d$, $p$, $R_{q_{2}}$), and the condition 
\begin{equation}
R_{q_{2}}\underbar{\ensuremath{\tau}}^{-q_{2}}\left(\frac{\log p}{n}+\mathcal{T}_{1}+\mathcal{T}_{2}+\mathcal{T}_{3}\right)=O(\kappa_{L}),\label{eq:18}
\end{equation}
holds, then, 
\begin{equation}
|\hat{\beta}_{HSEL}-\beta^{*}|_{2}\leq\frac{c^{''}R_{q_{2}}^{\frac{1}{2}}}{\kappa_{L}^{1-\frac{q_{2}}{2}}}\left[\bar{M}\vee\lambda_{n,3}\right]^{1-\frac{q_{2}}{2}}\label{eq:19}
\end{equation}
with probability at least $1-c_{1}\exp\left(-c_{2}\log p\right)-c_{3}\sum_{j=0}^{p}\exp\left(-nC_{j}^{*}t_{nj}^{2}\right)$
for some $C_{j}^{*}$ independent of $n$, $d$, $p$, $R_{q_{2}}$.\\
\\
The following theorem (Theorem 4.2) provides a general upper bound
on the error $\left|\hat{\beta}_{HSEL}-\beta^{*}\right|_{2}$ when
the second-stage estimation concerns a regularized program as in (\ref{eq:12}).
As in Theorem 4.1, this result is an {}``oracle-inequality'' type
which does not assume the unknown function $m(\cdot)$ belongs to
the function class over which the nonparametric estimator from (\ref{eq:12})
is defined. For Theorem 4.2, let the \textit{local complexity} measure
$\mathcal{G}_{n}\left(r_{j};\,\mathcal{F}_{j}\right)$ be defined
over the set 
\[
\Omega(r_{j};\,\mathcal{F}_{j})=\left\{ f\,:\, f\in\bar{\mathcal{F}}_{j}\,\left|f_{\theta^{*}}\right|_{n}\leq r_{j},\,\left|f\right|_{\mathcal{F}_{j}}\leq1\right\} 
\]
where $\left|f_{\theta^{*}}\right|_{n}:=\sqrt{\frac{1}{n}\sum_{i=1}^{n}\left[f(w_{i}^{T}\theta^{*})\right]^{2}}$
and $j=0,...,p$. Also define the following quantities: 
\begin{eqnarray*}
\mathcal{T}_{1} & = & \max_{j\in\{0,...,p\}}\left(T_{j}^{'}\vee\sqrt{T_{j}^{'}}\right)\\
\mathcal{T}_{2} & = & \max_{j\in\{0,...,p\}}\bar{R}_{j}^{2}t_{nj}^{2}\\
\mathcal{T}_{3} & = & \max_{j\in\{0,...,p\}}\inf_{\tilde{m}_{j}\in\mathcal{F}_{j},\,\left|\tilde{m}_{j}\right|_{\mathcal{F}_{j}}\leq\bar{R}_{j}}\left(\frac{1}{n}\sum_{i=1}^{n}\left[\tilde{m}_{j}(w_{i}^{T}\hat{\theta})-m_{j}(w_{i}^{T}\hat{\theta})\right]^{2}+\sqrt{\frac{1}{n}\sum_{i=1}^{n}\left[\tilde{m}_{j}(w_{i}^{T}\theta^{*})-m_{j}(w_{i}^{T}\theta^{*})\right]^{2}}\right)\\
\mathcal{T}_{4} & = & \sigma_{v}\sigma_{\eta}\sqrt{\frac{\log p}{n}},
\end{eqnarray*}
where $T_{j}^{'}$ is defined prior to the presentation of Theorem
4.1. \\
\textbf{}\\
\textbf{Theorem 4.2}: Let the \textit{critical radius} $r_{nj}>0$
be the smallest positive quantity satisfying the \textit{critical
inequality} 
\[
\mathcal{G}_{n}\left(r_{nj};\,\mathcal{F}_{j}\right)\leq\frac{\bar{R}_{j}r_{nj}^{2}}{\sigma_{v_{j}}},
\]
where $\bar{R}_{j}>0$ is a user-defined radius. Suppose the second-stage
estimator solves the regularized program (\ref{eq:12}) and Assumptions
2.1, 2.2, 4.1-4.3 hold. Additionally, let Assumption 2.3 hold over
the restricted set $\mathbb{C}(J(\beta^{*});\,0,\,3)$ for the exact
sparsity case ($q_{2}=0$ with $R_{q_{2}}=k_{2}$), and over $\mathbb{C}(S_{\underbar{\ensuremath{\tau}}};\, q_{2},\,3)\cap\mathbb{S}_{\delta}$
where $\delta\asymp R_{q_{2}}^{\frac{1}{2}}\left(\lambda_{n,3}\right)^{1-\frac{q_{2}}{2}}$
and $\underbar{\ensuremath{\tau}}=\frac{\lambda_{n,3}}{\kappa_{L}}$
for the approximate sparsity case ($q_{2}\in(0,\,1]$), respectively.
For any $t_{nj}\geq r_{nj}$, if the second-stage regularization parameter
$\lambda_{nj,2}=2t_{nj}^{2}+\varsigma$ for any small positive constant
$\varsigma>0$ and the third-stage regularization parameter $\lambda_{n,3}$
satisfies (\ref{eq:17}), and condition (\ref{eq:18}) holds, then,
the upper bound (\ref{eq:19}) holds (where the terms $\mathcal{T}_{k}$,
$k=1,...,4$ correspond to the ones defined for Theorem 4.2) with
probability at least 
\[
1-c_{1}\exp\left(-c_{2}\log p\right)-c_{3}\sum_{j=0}^{p}\exp\left(-nC_{j}^{*}\bar{R}_{j}^{2}t_{nj}^{2}\right)
\]
for some $C_{j}^{*}$ independent of $n$, $d$, $p$, $R_{q_{2}}$.\\
\\
\textbf{Comments}: \\
(a) For the probability guarantees in Theorems 4.1 and 4.2, the constant
\[
C_{j}^{*}=c\frac{\gamma_{v_{j}}\wedge(B_{v_{j}}^{2}\vee B_{\eta}^{2})^{-1}}{\sigma_{v_{j}}^{2}\vee\sigma_{\eta}^{2}}
\]
when $v_{j}$ for $j=0,...,p$ and $\eta$ are (i) sub-Gaussian with
\textit{strongly log-concave} distribution for some $\gamma_{v_{j}}>0$
and $\gamma_{\eta}>0$, respectively; or, (ii) a bounded vector such
that for every $i=1,...,n$, $v_{ij}$ and $\eta$ are supported on
the interval $(a_{v_{j}}^{'},\, a_{v_{j}}^{''})$ with $B_{v_{j}}:=a_{v_{j}}^{''}-a_{v_{j}}^{'}$,
and on $(a_{\eta}^{'},\, a_{\eta}^{''})$ with $B_{\eta}:=a_{\eta}^{''}-a_{\eta}^{'}$;
or, (iii) a mixture of (i) and (ii) in terms of its probability measure.
\\
\\
(b) Condition (\ref{eq:18}) in Theorems 4.1 and 4.2 ensures that
with high probability, $\frac{\hat{v}^{T}\hat{v}}{n}$ satisfies the
RE condition (\ref{eq:6}) over $\mathbb{C}(J(\beta^{*});\,0,\,3)$
for the exact sparsity case ($q_{2}=0$ with $R_{q_{2}}=k_{2}$),
and over $\mathbb{C}(S_{\underbar{\ensuremath{\tau}}};\, q_{2},\,3)\cap\mathbb{S}_{\delta}$
where $\delta\asymp R_{q_{2}}^{\frac{1}{2}}\left(\lambda_{n,3}\right)^{1-\frac{q_{2}}{2}}$
and $\underbar{\ensuremath{\tau}}=\frac{\lambda_{n,3}}{\kappa_{L}}$
for the approximate sparsity case ($q_{2}\in(0,\,1]$), respectively.
An implication of this scaling condition is that it provides a finite-sample
guarantee of the population identification condition (Assumption 2.3)
subject to the underlying restricted sets. This result is formalized
in the following corollary.\\
\\
\textbf{Corollary 4.3}: Under the assumptions in Theorem 4.1 (respectively,
the assumptions in Theorem 4.2), we have, with the same probability
guarantees in Theorem 4.1 (respectively, in Theorem 4.2), 
\[
\frac{1}{n}\sum_{i=1}^{n}y_{1i}\left(x_{i}-\hat{\mathbb{E}}\left[x_{i}\,\vert\, w_{i}^{T}\hat{\theta},\, y_{1i}=1\right]\right)\left(x_{i}-\hat{\mathbb{E}}\left[x_{i}\,\vert\, w_{i}^{T}\hat{\theta},\, y_{1i}=1\right]\right)^{T}
\]
is nonsingular on the restricted sets subject to those in Theorem
4.1 (respectively, Theorem 4.2).\\
\textbf{}\\
\textbf{Remarks} \textbf{on Theorems 4.1 and 4.2}

The main proofs for Theorem 4.1, Theorem 4.2, and Corollary 4.3 are
provided in Sections A.1-A.4. 

These theorems imply that if $\lambda_{n,3}\asymp\bar{M}$ and 
\[
\frac{c^{''}R_{q_{2}}^{\frac{1}{2}}}{\kappa_{L}^{1-\frac{q_{2}}{2}}}\left[b(\sigma_{v},\,\sigma_{\eta})\left|\beta^{*}\right|_{1}\left(\mathcal{T}_{1}+\mathcal{T}_{2}+\mathcal{T}_{3}\right)+\mathcal{T}_{4}\right]^{1-\frac{q_{2}}{2}}\rightarrow0,
\]
as $n\rightarrow\infty$, then the two-stage estimator $\hat{\beta}_{HSEL}$
is $l_{2}-$consistent for \textbf{$\beta^{*}$}. From Theorems 4.1
and 4.2, it can be seen that the general upper bounds on $|\hat{\beta}_{HSEL}-\beta^{*}|_{2}$
depend on four sources of errors, $\mathcal{T}_{k}$, $k=1,...,4$.
The terms $\mathcal{T}_{1}$, $\mathcal{T}_{2}$, $\mathcal{T}_{3}$,
and $\mathcal{T}_{4}$ are related to the statistical error of the
first-stage estimation, the statistical error of the second-stage
nonparametric regression, the approximation error arising from the
fact that $m_{j}\notin\mathcal{F}_{j}$, and the statistical error
of the third-stage estimation, respectively. Inspecting the error
term $\mathcal{T}_{1}$ suggests that, given appropriate identification
assumptions, the upper bounds on $|\hat{\beta}_{HSEL}-\beta^{*}|_{2}$
in Theorems 4.1 and 4.2 also hold for the more general structure where
$\mathbb{E}\left(\epsilon_{2i}\vert w_{i},\, y_{1i}=1\right)=g(h(w_{i}^{T},\,\theta^{*}))$
and $h(w_{i}^{T},\,\theta^{*})$ is a scalar unobservable index. 

The extra factor $\left|\beta^{*}\right|_{1}$ (in the case of exact
sparsity,$\left|\beta^{*}\right|_{1}\asymp k_{2}$) in front of $\mathcal{T}_{1}$,
$\mathcal{T}_{2}$, and $\mathcal{T}_{3}$ in the upper bounds on
$|\hat{\beta}_{HSEL}-\beta^{*}|_{2}$ as well as in the choice of
$\lambda_{n,3}$ is unimprovable and arising from the fact that the
estimator is a sequential multi-stage procedure%
\footnote{Other plug-in type Lasso estimators for the exactly sparse case such
as the ones in Rosenbaum and Tsybakov (2011) and the high-dimensional
two-stage least-squares estimator in Zhu (2013), also involve the
extra factor $\left|\beta^{*}\right|_{1}$.%
}based on plugging the first-stage estimator $\hat{\theta}$ in the
place of $\theta^{*}$. When $q_{2}=1$, the extra factor $\left|\beta^{*}\right|_{1}$
in front of $\mathcal{T}_{1}$, $\mathcal{T}_{2}$, and $\mathcal{T}_{3}$
in the upper bounds on $|\hat{\beta}_{HSEL}-\beta^{*}|_{2}$ as well
as in the in the choice of $\lambda_{n,3}$ is crucial in order for
the argument in our analysis to go through. To see this, suppose $\sqrt{\frac{\log p}{n}}$
is small relative to $\mathcal{T}_{1}$, $\mathcal{T}_{2}$, and $\mathcal{T}_{3}$,
in which case, condition (\ref{eq:18}) can be reduced to 
\[
R_{q_{2}}\left(\left|\beta^{*}\right|_{1}\right)^{-q_{2}}\max_{k}\left\{ \mathcal{T}_{k}^{1-q_{2}}:\, k=1,\,2,\,3\right\} =O\left(\kappa_{L}^{1-q_{2}}\right).
\]
When $q_{2}=1$, $R_{q_{2}}=\left|\beta^{*}\right|_{1}$ and $\sqrt{R_{q_{2}}}\left(\left|\beta^{*}\right|_{1}\right)^{-q_{2}}=1$,
so the above condition holds. On the other hand, when $q_{2}\in[0,\,1)$,
condition (\ref{eq:18}) is easier to be satisfied.

When $\mathcal{T}_{1}$, $\mathcal{T}_{2}$, and $\mathcal{T}_{3}$
are small relative to $\mathcal{T}_{4}$ and $\beta^{*}$ is exactly
sparse with at most $k_{2}$ non-zero coefficients, if we set $\kappa_{L}=\lambda_{\min}(\Sigma_{v})$,
the upper bounds in Theorems 4.1 and 4.2 reduce to $|\hat{\beta}_{HSEL}-\beta^{*}|_{2}\precsim\frac{\sigma_{v}\sigma_{\eta}}{\lambda_{\min}(\Sigma_{v})}\sqrt{\frac{k_{2}\log p}{n}}$.
Note that the scaling $\sqrt{\frac{k_{2}\log p}{n}}$ is the optimal
rate of the Lasso for the usual high-dimensional linear regression
model (\ref{eq:4}) with exact sparsity, and the factor $\frac{\sigma_{v}\sigma_{\eta}}{\lambda_{\min}(\Sigma_{v})}$
has a natural interpretation of an inverse signal-to-noise ratio when
$v_{i}$ is a zero-mean Gaussian matrix with covariance $\Sigma_{v}=\sigma_{v}^{2}I_{p\times p}$:
one has $\lambda_{\min}(\Sigma_{v})=\sigma_{v}^{2}$, so $\frac{\sigma_{v}\sigma_{\eta}}{\lambda_{\min}(\Sigma_{v})}=\frac{\sigma_{\eta}}{\sigma_{v}}$,
which measures the inverse signal-to-noise ratio of the regressors. 

For the case of approximately sparse $\beta^{*}$ with $q_{1},\, q_{2}\in(0,\,1]$,
the rate 
\begin{equation}
\frac{c^{''}R_{q_{2}}^{\frac{1}{2}}}{\kappa_{L}^{1-\frac{q_{2}}{2}}}\left[b(\sigma_{v},\,\sigma_{\eta})\left|\beta^{*}\right|_{1}\left(\mathcal{T}_{1}+\mathcal{T}_{2}+\mathcal{T}_{3}\right)+\mathcal{T}_{4}\right]^{1-\frac{q_{2}}{2}}\label{eq:20}
\end{equation}
can be interpreted with the following heuristic. Suppose we choose
to the top $s_{2}$ coefficients of $\beta^{*}$ in absolute values
to estimate, then the fast decay imposed by the $l_{q_{2}}-$balls
condition on $\beta^{*}$ would mean that the remaining $p-s_{2}$
coefficients would have relatively little impact. With this intuition,
the rate for $q_{2}>0$ can be viewed as the rate that would be achieved
by choosing 
\[
s_{2}=\frac{c^{''}R_{q_{2}}}{\kappa_{L}^{-q_{2}}}\left[b(\sigma_{v},\,\sigma_{\eta})\left|\beta^{*}\right|_{1}\left(\mathcal{T}_{1}+\mathcal{T}_{2}+\mathcal{T}_{3}\right)+\mathcal{T}_{4}\right]^{-q_{2}}
\]
and then proceeding as if the problem were an instance of an exactly
sparse problem $q_{2}=0$ with $k_{2}=s_{2}$. For such a problem,
we would expect to obtain the rate 
\[
\frac{c^{''}\sqrt{s_{2}}}{\kappa_{L}}\left[b(\sigma_{v},\,\sigma_{\eta})\left|\beta^{*}\right|_{1}\left(\mathcal{T}_{1}+\mathcal{T}_{2}+\mathcal{T}_{3}\right)+\mathcal{T}_{4}\right],
\]
which is exactly equal to (\ref{eq:20}). 

Notice that the choice of $t_{nj}$ incurs a trade-off between $\mathcal{T}_{2}$
and the term $O\left(\sum_{j=0}^{p}\exp\left(-nt_{nj}^{2}\right)\right)$
in the probability guarantees in Theorems 4.1 and 4.2. This is a general
phenomenon for these tail bounds. For the problems considered in this
paper, $t_{nj}$ may be chosen in the way that $\mathcal{T}_{2}$
is dominated by $\mathcal{T}_{1}$, $\mathcal{T}_{3}$, and $\mathcal{T}_{4}$
while the probability guarantees are maximized to allow for the least
restrictive requirement on the sample size for $l_{2}-$consistency.
Section 4.1.2 provides a specific example in terms of the choice of
$t_{nj}$. When we set $t_{nj}=r_{nj}$, note that the statistical
error related to the second-stage nonparametric regression, $\mathcal{T}_{2}$,
is on the order of $O\left(\max_{j}r_{nj}^{2}\right)$ instead of
the optimal rate $O\left(\max_{j}r_{nj}\right)$ that one would expect
from a nonparametric regression as (\ref{eq:11}) or (\ref{eq:12}).
As long as $\max_{j}r_{nj}<1$, we have: $\max_{j}r_{nj}^{2}<\max_{j}r_{nj}$,
and provided $\max_{j}r_{nj}^{2}$ is small relative to $\mathcal{T}_{1}$,
$\mathcal{T}_{3}$ and $\mathcal{T}_{4}$, the convergence rate of
the estimator of the high-dimensional linear coefficients in the main
equation behaves as if the unknown nonparametric selection bias were
known. 

This result establishes the non-asymptotic counterpart of the familiar
asymptotic {}``oracle'' type of results from previous literature.
One of the drivers behind this oracle result lies on carefully controlling
for the term $\left|\frac{1}{n}\sum_{i=1}^{n}\hat{v}_{ij}\left[\hat{m}_{j}(w_{i}^{T}\hat{\theta})-\tilde{m}_{j}(w_{i}^{T}\hat{\theta})\right]\right|$
utilizing the fact that $\hat{v}_{ij}$ estimates the true residual
$v_{ij}$ which is obtained by projecting $x_{ij}$ or $y_{2i}$ onto
$w_{i}^{T}\theta^{*}$, namely, $v_{ij}=x_{ij}-\mathbb{E}\left(x_{ij}\vert w_{i}^{T}\theta^{*}\right)$
or $v_{i0}=y_{2i}-\mathbb{E}\left(y_{2i}\vert w_{i}^{T}\theta^{*}\right)$.
When this projection procedure is applied to classical low-dimensional
semilinear models with fixed $p$ and $d$ (in which case, there is
no first-stage related error $\mathcal{T}_{1}$), our general upper
bounds would imply that the nonparametric component needs to be estimated
at a rate \textit{no slower} than $O\left((\frac{1}{n})^{\frac{1}{4}}\right)$
in order for the estimator of the parametric component to achieve
the rate of $O\left(\sqrt{\frac{1}{n}}\right)$. In contrast to the
semilinear models, low-dimensional selection models require the rate
of the nonparametric component to be at least $O\left((\frac{1}{n})^{\frac{1}{3}}\right)$
because the nonparametric component in the selection model involves
an unknown single index that also needs to be estimated.

\subsubsection{Upper bounds and $l_{2}-$consistency for a leading case example}

An important consequence of Theorems 4.1 and 4.2 is when $m_{j}(\cdot)\in\mathcal{F}_{j}$
for every $j=0,...,p$ and $\mathcal{F}_{j}$ in
\[
T_{j}^{'}:=\sup_{f\in\mathcal{F}_{j}}\frac{1}{n}\sum_{i=1}^{n}\left[f(w_{i}^{T}\hat{\theta})-f(w_{i}^{T}\theta^{*})\right]^{2}
\]
can be restricted to the class of Lipschitz functions, as a result,
$\mathcal{T}_{3}=0$ and $T_{j}^{'}=\frac{1}{n}\sum_{i=1}^{n}L^{2}\left[w_{i}^{T}\hat{\theta}-w_{i}^{T}\theta^{*}\right]^{2}:=L^{2}B^{'}$.
Results regarding this leading case are provided in the following
corollaries (Corollaries 4.4 and 4.5). Before stating these results,
a procedure based on Lipschitz regression for the second-stage estimation
is presented and its theoretical guarantees are provided in Corollaries
4.4 and 4.5. 

We say that a function $f:\:\mathbb{R}\rightarrow\mathbb{R}$ is \textit{L-Lipschitz
}if 
\begin{equation}
\left|f(t)-f(t^{'})\right|\leq L\left|t-t^{'}\right|\label{eq:21}
\end{equation}
for all $t,\, t^{'}\in\mathbb{R}$. When $\mathcal{F}_{j}$ satisfies
the Lipschitz assumption, we restrict $\mathcal{F}_{j}$ in (\ref{eq:11})
to be the class of Lipschitz functions and consider $\tilde{m}_{j}$
in this class only, namely, 
\[
\hat{m}_{j}\in\textrm{arg }\min_{\begin{array}{c}
\tilde{m}_{j}:\:\mathbb{R}\rightarrow\mathbb{R}\\
\tilde{m}_{j}\,\textrm{is }L\textrm{-Lipschitz}
\end{array}}\left\{ \frac{1}{n}\sum_{i=1}^{n}\left(z_{ij}-\tilde{m}_{j}(w_{i}\hat{\theta})\right)^{2}\right\} \qquad\textrm{for }j=0,...,p.
\]
It can be easily verified that $\bar{\mathcal{F}}_{j}$, the shifted
class of Lipschitz functions is also Lipschitz and satisfies Assumption
4.1; i.e., it is star-shaped. By exploiting the structure of Lipschitz
functions, the program above can be converted to an equivalent finite-dimensional
problem by applying the constraint (\ref{eq:21}) to each of the sampled
points $w_{i}\hat{\theta}$ so that there must exist a real-valued
vector $(\tilde{z}_{1j},...,\tilde{z}_{ij},...,\tilde{z}_{nj})$ which
satisfies the constraints in the following convex program 
\begin{eqnarray}
(\hat{z}_{1j},...,\hat{z}_{ij},...,\hat{z}_{nj}) & \in & \textrm{arg}\min_{(\tilde{z}_{1j},...,\tilde{z}_{ij},...,\tilde{z}_{nj})}\left\{ \frac{1}{n}\sum_{i=1}^{n}\left(z_{ij}-\tilde{z}_{ij}\right)^{2}\right\} \nonumber \\
\textrm{s.t.}\quad\tilde{z}_{ij}-\tilde{z}_{i^{'}j} & \leq & L\left(w_{i}-w_{i^{'}}\right)^{T}\hat{\theta}\;\textrm{for all }i,\, i^{'}=1,...,n.\label{eq:22}
\end{eqnarray}
Given an optimal solution $(\hat{z}_{1j},...,\hat{z}_{ij},...,\hat{z}_{nj})$,
a Lipschitz function $\hat{m}_{j}$ can be constructed by interpolating
linearly between $\hat{z}_{ij}$s and the resulting function $\hat{m}_{j}$
is an estimate of $m_{j}$ (namely, the second-stage estimator). Moreover,
one can easily see that $\hat{m}_{j}(w_{i}^{T}\hat{\theta})=\hat{z}_{ij}$.
Note that the optimization problem above is a convex program with
a quadratic cost function and a total of $\left(\begin{array}{c}
n\\
2
\end{array}\right)$ linear constraints and $n$ variables ($n$ here denotes the sample
size of the observations for the main equation). There are many computationally
efficient algorithms for solving programs like this (e.g., the interior
point method). When $m_{j}(\cdot)$ is a monotonic Lipschitz function,
we can impose additional monotonicity constraints together with the
Lipschitz constraints in the above convex program. Kakade, et. al
(2011) provides an algorithm with provable guarantees for this type
of minimization problems. 

In the case where the Lipschitz constant \textit{$L$ }is unknown,
cross-validation methods can be used to determine \textit{$L$.} For
example, we can first solve the optimization problem (\ref{eq:22})
on a subsample of observations by imposing an additional constraint
$0\leq L\leq L^{(0)}$ for a chosen constant $L^{(0)}$ and obtain
$(\hat{z}_{1j},...,\hat{z}_{ij},...,\hat{z}_{nj},\, L):=\varpi^{0}$.
We then test for the prediction quality of this optimal solution $\varpi^{0}$
by comparing its predicted values (from interpolating linearly between
$\hat{z}_{ij}$s) for the remaining subsample with the actual observed
values. If the optimal solution $\varpi^{0}$ returns $L\approx L^{(0)}$,
we can iterate the process by imposing $0\leq L\leq L^{(1)}=2L^{(0)}$
in (\ref{eq:22}) and comparing the new optimal solution $\varpi^{1}$
with the previous one $\varpi^{0}$ and also testing for the prediction
quality of $\varpi^{1}$. \textbf{}\\
\textbf{}\\
\textbf{Assumption 4.4}: The matrix $w$ consists of bounded elements%
\footnote{A random matrix with bounded elements is sub-Gaussian.%
}.\\
\\
The following proposition (Proposition 4.1) regarding the \textit{critical
radius} $r_{nj}$ in Theorems 4.1 and 4.2 is based on results from
van der Vaart and Wellner (1996), van de Geer (2000), and Wainwright
(2015).\\
\\
\textbf{Proposition 4.1}: Let Assumptions 2.1 and 4.4 hold and $m_{j}(\cdot)\in\mathcal{F}_{j}$
for $j=0,...,p$. Suppose $\mathcal{F}_{j}$ belongs to the class
of $L-$Lipschitz functions and the Lipschitz regression procedure
(\ref{eq:22}) is applied. Then, for every $j=0,...,p$, $\mathcal{T}_{3}=0$
and $T_{j}^{'}=\frac{1}{n}\sum_{i=1}^{n}L^{2}\left[w_{i}^{T}\hat{\theta}-w_{i}^{T}\theta^{*}\right]^{2}:=L^{2}B^{'}$,
and the \textit{critical radius} $r_{nj}=O\left((\frac{|\theta^{*}|_{1}}{n})^{\frac{1}{3}}\right)$,
in Theorems 4.1 and 4.2.\\
\\
The following corollaries (Corollaries 4.4 and 4.5) provide results
regarding the leading case where for every $j=0,...,p$, $T_{j}^{'}=\frac{1}{n}\sum_{i=1}^{n}L^{2}\left[w_{i}^{T}\hat{\theta}-w_{i}^{T}\theta^{*}\right]^{2}:=L^{2}B^{'}$,$\mathcal{T}_{3}=0$,
and the \textit{critical radius} $r_{nj}=O\left((\frac{|\theta^{*}|_{1}}{n})^{\frac{1}{3}}\right)$,
in Theorems 4.1 and 4.2. These conditions are ensured by Proposition
4.1. The two corollaries differ by the upper bounds on the quantity
$B^{'}$. Justifications of these upper bounds on $B^{'}$ are given
by Propositions 4.2 and 4.3. Let $\Upsilon_{w,\theta^{*}}$ be a known
function depending only on $w$ and $\theta^{*}$. The quantity $\Upsilon_{w,\theta^{*}}$
changes according to the assumptions on $w$, which is to be made
clear by Propositions 4.2 and 4.3. To facilitate the discussion and
a later comparison with the minimax lower bounds in Section 4.4, the
results in Corollaries 4.4 and 4.5 are presented for the case of exact
sparsity on $\beta^{*}$ and $\theta^{*}$ ($q_{1}=q_{2}=0$). The
case of general sparsity on $\theta^{*}$ and $\beta^{*}$ ($q_{1},\, q_{2}\in[0,\,1]$)
is presented in Corollary 4.6 (which contains Corollary 4.4 as a special
case). \textbf{}\\
\textbf{}\\
\textbf{Corollary 4.4} ($q_{1}=q_{2}=0$): Suppose $\theta^{*}$ is
exactly sparse with at most $k_{1}$ non-zero coefficients. Suppose
for every $j=0,...,p$, $T_{j}^{'}=\frac{1}{n}\sum_{i=1}^{n}L^{2}\left[w_{i}^{T}\hat{\theta}-w_{i}^{T}\theta^{*}\right]^{2}:=L^{2}B^{'}$,$\mathcal{T}_{3}=0$,
and the \textit{critical radius} $r_{nj}=O\left((\frac{k_{1}}{n})^{\frac{1}{3}}\right)$,
and 
\[
B^{'}=\frac{1}{n}\sum_{i=1}^{n}\left[w_{i}^{T}\hat{\theta}-w_{i}^{T}\theta^{*}\right]^{2}\leq c\Upsilon_{w,\theta^{*}}\frac{k_{1}\log d}{n}
\]
with probability at least $1-O\left(\frac{1}{d}\right)$. Assume $t_{nj}^{2}$
in $\mathcal{T}_{2}$ is chosen such that $\left|\beta^{*}\right|_{1}\mathcal{T}_{2}$
is at most 
\[
O\left(\sqrt{\frac{\log p}{n}}\vee\left(\left|\beta^{*}\right|_{1}\sqrt{\frac{k_{1}\log d}{n}}\right)\right)
\]
and $nt_{nj}^{2}\succsim\log p$. Suppose Assumptions 2.1, 4.2-4.4
hold. Additionally, let $\beta^{*}$ satisfy the exact sparsity in
Assumption 2.2 ($q_{2}=0$ with $R_{q_{2}}=k_{2}$) and Assumption
2.3 hold over the restricted set $\mathbb{C}(J(\beta^{*});\,0,\,3)$.
Assume 
\[
\kappa_{2}\frac{k_{2}\log p}{n}+k_{2}\sqrt{\frac{k_{1}\log d}{n}}=O(\kappa_{1}),
\]
for some strictly positive constants $(\kappa_{1},\,\kappa_{2})$
depending only on $\kappa_{L}$, $\sigma_{v}$, $\Upsilon_{w,\theta^{*}}$,
and $L$. If the third-stage regularization parameter $\lambda_{n,3}$
satisfies 
\[
\lambda_{n,3}\geq c\left(\sigma_{v}\sigma_{\eta}\sqrt{\frac{\log p}{n}}\right)\vee\left(Lb\left(\sigma_{v},\,\sigma_{\eta}\right)\left|\beta^{*}\right|_{1}\sqrt{\Upsilon_{w,\theta^{*}}}\sqrt{\frac{k_{1}\log d}{n}}\right):=\bar{M}
\]
then, with probability at least $1-O\left(\frac{1}{p\wedge d}\right)$,
we have 
\[
|\hat{\beta}_{HSEL}-\beta^{*}|_{2}\leq\frac{c_{1}\sqrt{k_{2}}}{\kappa_{L}}\left[\bar{M}\vee\lambda_{n,3}\right]
\]
where $b\left(\sigma_{v},\,\sigma_{\eta}\right)$ is some known function
depending only on $\sigma_{v}$ and $\sigma_{\eta}$ (and independent
of $n$, $d$, $p$, $k_{1}$, and $k_{2}$). \\
\\
The following assumptions and proposition provide an example in which
the upper bound on $B^{'}$ in Corollary 4.4 is achieved. In particular,
it requires the eigenvalues of $\Sigma_{w}$ to be well-behaved over
some restricted set. \\
\textbf{}\\
\textbf{Assumption 4.5}: In program (\ref{eq:13}), we have: either
(a) $\phi_{1}(w_{i}^{T}\theta)=w_{i}^{T}\theta$ and $\phi_{2}(w_{i}^{T}\theta)=\log(1+\exp(w_{i}^{T}\theta))$;
namely, the loss function corresponds to a binary logit model. Or,
(b) $\phi_{1}(w_{i}^{T}\theta)=\log\frac{\Phi(w_{i}^{T}\theta)}{1-\Phi(w_{i}^{T}\theta)}$
and $\phi_{2}(w_{i}^{T}\theta)=-\log\left[1-\Phi(w_{i}^{T}\theta)\right]$;
namely, the loss function corresponds to a binary probit model. \\
\\
\textbf{Assumption 4.6}: The random matrix $w$ is sub-Gaussian with
parameters $(\Sigma_{w},\,\sigma_{w}^{2})$. For all $\Delta\in\mathbb{C}(J(\theta^{*});\,0,\,3)\backslash\{\mathbf{0}\}$,
the matrix $\Sigma_{w}$ satisfies 
\[
0<\kappa_{L}^{w}\leq\frac{\Delta^{T}\Sigma_{w}\Delta}{\left|\Delta\right|_{2}^{2}}\leq\kappa_{U}^{w}<\infty
\]
\textbf{Proposition 4.2}: Suppose the number of regressors $d(=d_{n})$
can grow with and exceed the sample size $n$ and\textbf{ }the number
of non-zero components in $\theta^{*}$ is at most $k_{1}(=k_{1n})$
and $k_{1}$ can increase to infinity with $n$ but slowly compared
to $n$. Let Assumptions 2.1, 4.5-4.6 hold. If $\hat{\theta}$ solves
program (\ref{eq:13}) with $\lambda_{n,1}\geq c\sigma_{w}\sqrt{\alpha_{u}}\sqrt{\frac{\log d}{n}}$
and $n\succsim k_{1}\log d$, then, with probability at least $1-O\left(\frac{1}{d}\right)$,
\[
\frac{1}{n}\sum_{i=1}^{n}\left[w_{i}^{T}(\hat{\theta}-\theta^{*})\right]^{2}\leq c^{'}\frac{\kappa_{U}^{w}}{\left(\kappa_{L}^{w}\right)^{2}}k_{1}\left(\left(\lambda_{n,1}\right)^{2}\vee\left(\sigma_{w}^{2}\alpha_{u}\frac{\log d}{n}\right)\right),
\]
where $\alpha_{u}>0$ is a scalar such that $\phi_{2}^{''}(u)\leq\alpha_{u}$
for all $u\in\mathbb{R}$. \\
\textbf{}\\
\textbf{Remark}. From Proposition 4.2, we can set $\Upsilon_{w,\theta^{*}}:=\frac{\kappa_{U}^{w}\sigma_{w}^{2}\alpha_{u}}{\left(\kappa_{L}^{w}\right)^{2}}$
in Corollary 4.4. The boundedness on $\phi_{2}^{''}(u)$ holds automatically
for the binary logit model and binary probit model. For the logit
model, we have $\phi_{2}^{''}(u_{i})=\frac{\exp(u_{i})}{1+\exp(u_{i})}\left(1-\frac{\exp(u_{i})}{1+\exp(u_{i})}\right)$.
For the probit model, note that $\phi_{2}^{''}(u_{i})$ is $1-Var\left(\epsilon_{1i}\,\vert\,\epsilon_{1i}\leq u_{i}\right)$
when $y_{1i}=1$ and $1-Var\left(\epsilon_{1i}\,\vert\,\epsilon_{1i}\geq-u_{i}\right)$
when $y_{1i}=0$ and the unconditional variance is normalized to $1$.
Since truncation always reduces variances (Greene, 2003), $\phi_{2}^{''}(u)$
is bounded from above. If $\lambda_{n,1}\asymp\sigma_{w}\sqrt{\alpha_{u}}\sqrt{\frac{\log d}{n}}$,
then 
\[
\frac{1}{n}\sum_{i=1}^{n}\left[w_{i}^{T}(\hat{\theta}-\theta^{*})\right]^{2}\leq c^{'}\frac{\kappa_{U}^{w}}{\left(\kappa_{L}^{w}\right)^{2}}\sigma_{w}^{2}\alpha_{u}\frac{k_{1}\log d}{n}.
\]
\textbf{Corollary 4.5}: Suppose $\theta^{*}$ is exactly sparse with
at most $k_{1}$ non-zero coefficients. Suppose for every $j=0,...,p$,
$T_{j}^{'}=\frac{1}{n}\sum_{i=1}^{n}L^{2}\left[w_{i}^{T}\hat{\theta}-w_{i}^{T}\theta^{*}\right]^{2}:=L^{2}B^{'}$,$\mathcal{T}_{3}=0$,
the \textit{critical radius} $r_{nj}=O\left((\frac{k_{1}}{n})^{\frac{1}{3}}\right)$,
and 
\[
B^{'}=\frac{1}{n}\sum_{i=1}^{n}\left[w_{i}^{T}\hat{\theta}-w_{i}^{T}\theta^{*}\right]^{2}\leq c\Upsilon_{w,\theta^{*}}\left|\theta^{*}\right|_{1}\sqrt{\frac{\log d}{n}}
\]
with probability at least $1-O\left(\frac{1}{d}\right)$. Assume $t_{nj}^{2}$
in $\mathcal{T}_{2}$ is chosen such that $\left|\beta^{*}\right|_{1}\mathcal{T}_{2}$
is at most 
\[
O\left(\sqrt{\frac{\log p}{n}}\vee\left(\left|\beta^{*}\right|_{1}\left(\frac{k_{1}^{2}\log d}{n}\right)^{\frac{1}{4}}\right)\right)
\]
and $nt_{nj}^{2}\succsim\log p$. Suppose Assumptions 2.1, 4.2-4.4
hold. Additionally, let $\beta^{*}$ satisfy the exact sparsity in
Assumption 2.2 ($q_{2}=0$ with $R_{q_{2}}=k_{2}$) and Assumption
2.3 hold over the restricted set $\mathbb{C}(J(\beta^{*});\,0,\,3)$.
Assume 
\[
\kappa_{2}\frac{k_{2}\log p}{n}+k_{2}\left(\frac{k_{1}^{2}\log d}{n}\right)^{\frac{1}{4}}=O(\kappa_{1}),
\]
for some strictly positive constants $(\kappa_{1},\,\kappa_{2})$
depending only on $\kappa_{L}$, $\sigma_{v}$, $\Upsilon_{w,\theta^{*}}$,
and $L$, if the third-stage regularization parameter $\lambda_{n,3}$
satisfies 
\[
\lambda_{n,3}\geq c^{'}\left(\sigma_{v}\sigma_{\eta}\sqrt{\frac{\log p}{n}}\right)\vee\left(Lb\left(\sigma_{v},\,\sigma_{\eta}\right)\left|\beta^{*}\right|_{1}\sqrt{\Upsilon_{w,\theta^{*}}}\left(\frac{\left|\theta^{*}\right|_{1}^{2}\log d}{n}\right)^{\frac{1}{4}}\right):=\bar{M}
\]
then, with probability at least $1-O\left(\frac{1}{p\wedge d}\right)$,
we have 
\[
|\hat{\beta}_{HSEL}-\beta^{*}|_{2}\leq\frac{c_{2}\sqrt{k_{2}}}{\kappa_{L}}\left[\bar{M}\vee\lambda_{n,3}\right]
\]
where $b\left(\sigma_{v},\,\sigma_{\eta}\right)$ is some known function
depending only on $\sigma_{v}$ and $\sigma_{\eta}$ (and independent
of $n$, $d$, $p$, $k_{1}$, and $k_{2}$). \\
\textbf{}\\
The following proposition provides an example in which the upper bound
on $B^{'}$ in Corollary 4.5 is achieved. Let $\rho_{i,\,\theta}:=-y_{1i}\phi_{1}(w_{i}^{T}\theta)+\phi_{2}(w_{i}^{T}\theta)$
and $\rho_{i,\,\bar{\theta}}^{''}$ be the second derivative of $\rho_{i,\,\theta}$,
evaluated at $\theta=\bar{\theta}$, where $\bar{\theta}$ is some
intermediate value between $\theta^{*}$ and $\hat{\theta}$, the
solution to program (\ref{eq:13}). Assumption 4.4 implies that there
is some $\alpha_{l}>0$ such that $\rho_{i,\,\bar{\theta}}^{''}\geq\alpha_{l}$
for all $i=1,...,n$.\\
\textbf{}\\
\textbf{Proposition 4.3}: Let Assumptions 2.1, 4.4 and 4.5 hold. Suppose
the number of regressors $d(=d_{n})$ can grow with and exceed the
sample size $n$ and\textbf{ }the number of non-zero components in
$\theta^{*}$ is at most $k_{1}(=k_{1n})$ and $k_{1}$ can increase
to infinity with $n$ but slowly compared to $n$. If $\hat{\theta}$
solves program (\ref{eq:13}) with the regularization parameter%
\footnote{The choice of $\lambda_{n,1}$ is detailed in Theorems 2.1 or 2.2
in Van de Geer (2008).%
} $\lambda_{n,1}\geq c\sqrt{\frac{\log d}{n}}$, then, 
\[
\frac{1}{n}\sum_{i=1}^{n}\left[w_{i}^{T}(\hat{\theta}-\theta^{*})\right]^{2}\leq c_{1}\Upsilon_{w,\theta^{*}}\left|\theta^{*}\right|_{1}\left(\sqrt{\frac{\log d}{n}}\vee\lambda_{n,1}\right)
\]
with probability at least $1-O\left(\frac{1}{d}\right)$, where $\Upsilon_{w,\theta^{*}}:=\alpha_{l}^{-1}$.
\\
\\
\textbf{Remarks on Corollary 4.4-4.5}

The proofs for Corollaries 4.4-4.5 and Propositions 4.1-4.3 are provided
in Sections A.5-A.7. 

Corollaries 4.4 and 4.5 imply that if $\lambda_{n,3}\asymp\bar{M}$
and the upper bounds on $|\hat{\beta}_{HSEL}-\beta^{*}|_{2}$ tend
to $0$ as $n\rightarrow\infty$, then the two-stage estimator $\hat{\beta}_{HSEL}$
is $l_{2}-$consistent for \textbf{$\beta^{*}$}. The difference between
Corollary 4.4 and Corollary 4.5 lies in that the statistical error
from the first-stage estimation is smaller in Corollary 4.4 relative
to Corollary 4.5 and as a result, the estimator $\hat{\beta}_{HSEL}$
has a faster rate of convergence in Corollary 4.4. The smaller first-stage
statistical error in Corollary 4.4 is at the expense of imposing conditions
on the eigenvalues of $\Sigma_{w}$, as shown in Proposition 4.2.
Consistency of $\hat{\beta}_{HSEL}$ \textit{per se} does not require
restrictions on the eigenvalues of $\Sigma_{w}$, which could be useful
in certain applications. Proposition 4.3 provides an example where
a slower rate of convergence is obtained by the first-stage estimator
upon relaxing the assumptions on the eigenvalues of $\Sigma_{w}$.

By Proposition 4.1, $\max_{j}r_{nj}^{2}=O\left(\left(\frac{k_{1}}{n}\right)^{\frac{2}{3}}\right)$.
Let us examine various choices of $t_{nj}^{2}\geq r_{nj}^{2}$ in
Corollary 4.4 (the analysis for Corollary 4.5 is similar). Setting
$t_{nj}^{2}\asymp\sqrt{\frac{\log p}{n\left|\beta^{*}\right|_{1}^{2}}}\vee\sqrt{\frac{k_{1}\log d}{n}}\succsim r_{nj}^{2}$
makes the second-stage error $\left|\beta^{*}\right|_{1}\mathcal{T}_{2}$
on the same order of $\sqrt{\frac{\log p}{n}}\vee\left(\left|\beta^{*}\right|_{1}\sqrt{\frac{k_{1}\log d}{n}}\right)$.
Under this choice of $t_{nj}^{2}$, we require $\sqrt{\frac{n\log p}{\left|\beta^{*}\right|_{1}^{2}}}\vee\sqrt{nk_{1}\log d}\succsim\log p$
in order for the upper bound on $|\hat{\beta}_{HSEL}-\beta^{*}|_{2}$
to hold with probability at least $1-O\left(\frac{1}{p\wedge d}\right)$.
Setting $t_{nj}^{2}\asymp\left(\frac{\log p\vee(k_{1}\log d)}{n}\right)^{\frac{2}{3}}\succsim r_{nj}^{2}$
requires $n^{\frac{1}{3}}\left(\log p\vee(k_{1}\log d)\right)^{\frac{2}{3}}\succsim\log p$
for the upper bound on $|\hat{\beta}_{HSEL}-\beta^{*}|_{2}$ to hold
with probability at least $1-O\left(\frac{1}{p\wedge d}\right)$.
If instead, we set $t_{nj}^{2}=r_{nj}^{2}$, then the probability
guarantee of $1-O\left(\frac{1}{p\wedge d}\right)$ would require
$k_{1}^{\frac{2}{3}}n^{\frac{1}{3}}\succsim\log p$. Given the exact
sparsity of $\beta^{*}$ (so$\left|\beta^{*}\right|_{1}^{2}\asymp k_{2}^{2}$),
if $k_{2}^{2}$ is sufficiently small relative to $n\log p$, the
first choice of $t_{nj}^{2}$ would provide the least restrictive
requirement on the sample size. A later result that concerns with
the selection consistency of $\hat{\beta}_{HSEL}$ assumes this choice
for $t_{nj}^{2}$ and the scaling condition $\sqrt{\frac{n\log p}{\left|\beta^{*}\right|_{1}^{2}}}\vee\sqrt{nk_{1}\log d}\succsim\log p$
on the sample size. When $p$ and $d$ are fixed and small relative
to $n$, the analysis above generalizes existing asymptotic {}``oracle''
results in semiparametric estimation of low-dimensional selection
models from specific estimators (such as a series estimator) to a
unified framework of nonparametric least squares estimators and regularized
nonparametric least squares estimators. 

More generally, when $\mathcal{F}_{j}$ belongs to a Hölder class
of order $\nu>0$, we have $\max_{j}r_{nj}^{2}=O\left(\left(\frac{k_{1}}{n}\right)^{\frac{2\nu}{2\nu+1}}\right)$.
When $v\geq1$, $\mathcal{T}_{1}\asymp\sqrt{\frac{k_{1}\log d}{n}}$
and as long as we choose $t_{nj}^{2}\asymp\sqrt{\frac{\log p}{n\left|\beta^{*}\right|_{1}^{2}}}\vee\sqrt{\frac{k_{1}\log d}{n}}$,
the second-stage error $\left|\beta^{*}\right|_{1}\mathcal{T}_{2}$
and consequently the upper bound on $|\hat{\beta}_{HSEL}-\beta^{*}|_{2}$
would be on the same order of $\sqrt{\frac{\log p}{n}}\vee\left(\left|\beta^{*}\right|_{1}\sqrt{\frac{k_{1}\log d}{n}}\right)$.
On the other hand, when $v\in(0,\,1)$, we have $\frac{2\nu}{2\nu+1}>\frac{\nu}{2}$
and $\mathcal{T}_{1}\asymp\left(\sqrt{B^{'}}\right)^{\nu}$. Provided
$B^{'}\geq O(\frac{1}{n})$ (which is indeed the case for Corollaries
4.4 and 4.5) and the choice of $t_{nj}^{2}\asymp\frac{\mathcal{T}_{4}}{\left|\beta^{*}\right|_{1}}\vee\left(\sqrt{B^{'}}\right)^{\nu}$,
then $|\hat{\beta}_{HSEL}-\beta^{*}|_{2}$ is bounded above by $\mathcal{T}_{4}\vee\left(\left|\beta^{*}\right|_{1}\left(\sqrt{B^{'}}\right)^{\nu}\right)$.
However, note in the simple example where $B^{'}=\frac{1}{n}$ and
$\mathcal{T}_{4}=\sqrt{\frac{1}{n}}$, we have $\mathcal{T}_{1}=\left(\frac{1}{n}\right)^{\frac{\nu}{2}}>\sqrt{\frac{1}{n}}$
for any $\nu\in(0,\,1)$ and therefore $|\hat{\beta}_{HSEL}-\beta^{*}|_{2}$
is bounded above by $\left(\frac{1}{n}\right)^{\frac{\nu}{2}}$. Consequently,
the minimum requirement for the {}``oracle'' result to hold in the
low-dimensional semiparametric selection models with fixed $p$ and
$d$ is to have $v=1$. For the high-dimensional selection models
considered in Corollary 4.4, the minimum requirement is to have $O\left(\frac{k_{1}\log d}{n}\right)^{\frac{\nu}{2}}=O\left(\sqrt{\frac{\log p}{n}}\right)$.
In sharp contrast to the low-dimensional semilinear model, the fact
that the nonparametric component in the selection model involves an
unknown single index that also needed to be estimated increases the
requirement on the rate of nonparametric estimation\textit{ per se.} 

Note that the regularization parameter $\lambda_{n,3}$ and the upper
bounds on $|\hat{\beta}_{HSEL}-\beta^{*}|_{2}$ depend on $\sigma_{v}$
and $\sigma_{\eta}$, which is intuitive. It is possible to {}``remove''
the dependence on $\sigma_{v}$ from the choice of $\lambda_{n,3}$
by imposing weights $\hat{\sigma}_{v_{j}}:=\sqrt{\frac{1}{n}\sum_{i=1}^{n}\hat{v}_{ij}^{2}}$
, $j=1,...,p$ on the penalty term as in (\ref{eq:9}). An application
of Lemmas A.11 and A.12 yields $\max_{j=1,...,p}\hat{\sigma}_{v_{j}}\leq2\sigma_{v}$
with probability at least $1-O\left(\frac{1}{p\wedge d}\right)$.
The first-stage estimator $\hat{\theta}$ in Corollaries 4.4 and 4.5
may be replaced with a post-Lasso estimator where a usual low-dimensional
estimation procedure is performed on the regressors selected by $\hat{\theta}$
(in a spirit similar to Belloni and Chernozhukov, 2011b, for example);
and upon perfect selection or near-perfect selection%
\footnote{Ravikumar, Wainwright, and Lafferty (2010) studies selection of a
$l_{1}-$regularized logistic regression in the high-dimensional setting.%
} of $\hat{\theta}$, the term $\sqrt{\frac{k_{1}\log d}{n}}$ from
Corollary 4.4 and the term $\left(\frac{k_{1}^{2}\log d}{n}\right)^{\frac{1}{4}}$
from Corollary 4.5 in the upper bounds on $|\hat{\beta}_{HSEL}-\beta^{*}|_{2}$
can be reduced to $\sqrt{\frac{k_{1}}{n}}$ and $\left(\frac{k_{1}^{2}}{n}\right)^{\frac{1}{4}}$,
respectively.

We now present a result for the general sparsity case where $\theta^{*}$
and $\beta^{*}$ belong to the general $l_{q_{1}}-$ and $l_{q_{2}}-$
{}``balls'' with $q_{1},\, q_{2}\in[0,\,1]$.\\
\textbf{}\\
\textbf{Corollary 4.6} ($q_{1},\, q_{2}\in[0,\,1]$): Suppose for
every $j=0,...,p$, $T_{j}^{'}=\frac{1}{n}\sum_{i=1}^{n}L^{2}\left[w_{i}^{T}\hat{\theta}-w_{i}^{T}\theta^{*}\right]^{2}:=L^{2}B^{'}$,$\mathcal{T}_{3}=0$,
and the \textit{critical radius} $r_{nj}=O\left((\frac{|\theta^{*}|_{1}}{n})^{\frac{1}{3}}\right)$.
Also, assume $\theta^{*}\in\mathcal{B}_{q_{1}}^{d}(R_{q_{1}})$ for
$q_{1}\in[0,\,1]$ with {}``radius'' $R_{q_{1}}$, and 
\[
B^{'}=\frac{1}{n}\sum_{i=1}^{n}\left[w_{i}^{T}\hat{\theta}-w_{i}^{T}\theta^{*}\right]^{2}\leq c\Upsilon_{w,\theta^{*}}R_{q_{1}}\left(\frac{\log d}{n}\right)^{1-\frac{q_{1}}{2}}
\]
with probability at least $1-O\left(\frac{1}{d}\right)$. Assume $t_{nj}^{2}$
in $\mathcal{T}_{2}$ is chosen such that $\left|\beta^{*}\right|_{1}\mathcal{T}_{2}$
is at most $O(\bar{M})$, where
\[
\bar{M}:=\max\left\{ \sqrt{\frac{\log p}{n}},\:\left|\beta^{*}\right|_{1}\left(\frac{\left|\theta^{*}\right|_{1}}{n}\right)^{\frac{2}{3}},\:\left|\beta^{*}\right|_{1}R_{q_{1}}^{\frac{1}{2}}\left(\sqrt{\frac{\log d}{n}}\right)^{1-\frac{q_{1}}{2}}\right\} 
\]
and $nt_{nj}^{2}\succsim\log p$. Moreover, condition (\ref{eq:18})
and Assumptions 2.1, 2.2, 4.2-4.4 hold. Additionally, let Assumption
2.3 hold over the restricted set $\mathbb{C}(J(\beta^{*});\,0,\,3)$
for the exact sparsity case ($q_{2}=0$ with $R_{q_{2}}\asymp k_{2}$),
and over $\mathbb{C}(S_{\underbar{\ensuremath{\tau}}};\, q_{2},\,3)\cap\mathbb{S}_{\delta}$
where $\delta\asymp R_{q_{2}}^{\frac{1}{2}}\left(\lambda_{n,3}\right)^{1-\frac{q_{2}}{2}}$
and $\underbar{\ensuremath{\tau}}=\frac{\lambda_{n,3}}{\kappa_{L}}$
for the approximate sparsity case ($q_{2}\in(0,\,1]$), respectively.
If the third-stage regularization parameter $\lambda_{n,3}\geq b_{0}\bar{M}$,
then, with probability at least $1-O\left(\frac{1}{p\wedge d}\right)$,
we have 
\[
|\hat{\beta}_{HSEL}-\beta^{*}|_{2}\leq\frac{b_{1}\sqrt{R_{q_{2}}}}{\kappa_{L}^{1-\frac{q_{2}}{2}}}\left(\bar{M}\vee\lambda_{n,3}\right)^{1-\frac{q_{2}}{2}}
\]
where $b_{0}$ and $b_{1}$ are some known functions depending only
on $\sigma_{v}$, $\sigma_{\eta}$, $\Upsilon_{w,\theta^{*}}$, and
$L$ (and independent of $n$, $d$, $p$, $R_{q_{1}}$, and $R_{q_{2}}$).
\\
\\
\textbf{Comment on Corollary 4.6}. Corollary 4.6 contains Corollary
4.4 as a special case with $q_{2}=0$, $R_{q_{2}}=k_{2}$ and $q_{1}=0$,
$R_{q_{1}}=k_{1}$. When $q_{1}=0$ so that $R_{q_{1}}=k_{1}$ and
$\left|\theta^{*}\right|_{1}\asymp k_{1}$, the second term in $\bar{M}$
is on the order of $O\left(\left|\beta^{*}\right|_{1}\left(\frac{k_{1}}{n}\right)^{\frac{2}{3}}\right)$
and therefore dominated by $\sqrt{\frac{\log p}{n}}\vee\left|\beta^{*}\right|_{1}\sqrt{\frac{k_{1}\log d}{n}}$,
as we have seen previously. For more general sparsity of $\theta^{*}$
($q_{1}\in(0,\,1]$), the second term in $\bar{M}$ may still be small
relative to the first and third terms and therefore the {}``oracle''
result is likely to hold for a range of scaling conditions on $n$,
$p$, $d$, $R_{q_{1}}$, and $\left|\theta^{*}\right|_{1}$.

\subsubsection{Variable-selection consistency of a leading case example with exact
sparsity}

The following theorem (Theorem 4.7) addresses the question:\textbf{
}given $\hat{\beta}_{HSEL}$, when does $\hat{\beta}_{HSEL}$ correctly
select the non-zero coefficients in the main equation with high probability?
This property is referred to as \textbf{variable-selection consistency},
which is relevant to case of exactly sparse $\beta^{*}$ (and therefore
this section assumes $\beta^{*}$ is exactly sparse with at most $k_{2}$
non-zero coefficients). We say $\hat{\beta}_{HSEL}$ achieves \textit{perfect
selection} if $\mathbb{P}[J(\hat{\beta}_{HSEL})=J(\beta^{*})]\rightarrow1$
and \textit{near-perfect selection} if $\mathbb{P}[J(\hat{\beta}_{HSEL})\supseteq J(\beta^{*})]\rightarrow1$
and the number of wrong components selected is on the order of $o_{p}(k_{2})$.
Upon perfect selection or near-perfect selection of the regressors,
we can then apply low-dimensional techniques to estimate and conduct
inference on the important coefficients. 

In order for the number of wrong components selected by the standard
Lasso to be on the order of $o_{p}(k_{2})$ in the context of standard
linear regression models, it is known that the so-called {}``neighborhood
stability condition'' (Meinshausen and Bühlmann, 2006) on the design
matrix, re-formulated in a nicer form as the {}``irrepresentable
condition'' by Zhao and Yu, 2006, or the {}``mutual incoherence
condition'' by Wainwright (2009), is sufficient and necessary. Furthermore,
it can be shown that the {}``irrepresentable condition'' implies
the RE condition (see, e.g., Bühlmann and van de Geer, 2011).\\
\textbf{}\\
\textbf{Assumption 4.7}:\textbf{ $\left\Vert \mathbb{E}\left[v_{1,J(\beta^{*})^{c}}^{T}v_{1,J(\beta^{*})}\right]\left[\mathbb{E}(v_{1,J(\beta^{*})}^{T}v_{1,J(\beta^{*})})\right]^{-1}\right\Vert _{\infty}\leq1-\phi$
}for some constant $\phi\in(0,\,1]$.\textbf{}\\
\\
Assumption 4.7, the so-called {}``mutual incoherence condition''
originally formalized by Wainwright (2009), captures the intuition
that the large number of irrelevant covariates cannot exert an overly
strong effect on the subset of relevant covariates. In the most desirable
case, the columns indexed by $j\in J(\beta^{*})^{c}$ would all be
orthogonal to the columns indexed by $j\in J(\beta^{*})$\textbf{
}and then we would have $\phi=1$.\textbf{ }In the high-dimensional
setting, this perfect orthogonality is hard to achieve, but one can
still hope for a type of {}``near orthogonality'' to hold. 

Assumptions 2.1 and 4.2 ensure that the left-hand-side of the inequality
in Assumption 4.7 always falls in $[0,\,1)$. To see this, note that
under Assumptions 2.1 and 4.2, each column $v_{j}$, $j=1,...,p$
is consisted of \textit{i.i.d.} sub-Gaussian variables. Without loss
of generality, we can assume $\mathbb{E}(v_{1j})=0$ for all $j=1,...,p$.
Consequently, the normalization $\max_{j=1,...,p}\frac{|v_{j}|_{2}}{\sqrt{n}}\leq\kappa_{c}$
where $0<\kappa_{c}<\infty$ follows from a standard bound for the
norms of zero-mean sub-Gaussian vectors and a union bound 
\[
\mathbb{P}\left[\max_{j=1,...,p}\frac{|v_{j}|_{2}}{\sqrt{n}}\leq\kappa_{c}\right]\geq1-2\exp(-cn+\log p)\geq1-2\exp(-c^{'}n),
\]
where the last inequality follows from $n>\log p$. For example, if
$v_{j}$ has a Gaussian design, then we have 
\[
\max_{j=1,...,p}\frac{|v_{j}|_{2}}{\sqrt{n}}\leq\max_{j=1,...,p}\Sigma_{jj}\left(1+\sqrt{\frac{32\log p}{n}}\right),
\]
where $\max_{j=1,..,p}\Sigma_{jj}$ corresponds to the maximal variance
of any element of $v$ (see Raskutti, et. al, 2011). \\
\\
\textbf{Theorem 4.7}: Under the assumptions in Corollary 4.4 and Assumption
4.7, if $n\succsim\left(k_{2}^{3}\log p\right)\vee\left(k_{2}^{2}k_{1}\log d\right)$,
$\sqrt{\frac{n\log p}{\left|\beta^{*}\right|_{1}^{2}}}\vee\sqrt{nk_{1}\log d}\succsim\log p$,
$\sqrt{\frac{k_{1}\log d}{n}}=o(1)$, and $\lambda_{n,3}$ satisfies
\[
\lambda_{n,3}\geq c\frac{8(2-\frac{\phi}{4})}{\phi}\left[\left(\sigma_{v}\sigma_{\eta}\sqrt{\frac{\log p}{n}}\right)\vee\left(L\left|\beta^{*}\right|_{1}\sqrt{\Upsilon_{w,\theta^{*}}}b\left(\sigma_{v},\,\sigma_{\eta}\right)\sqrt{\frac{k_{1}\log d}{n}}\right)\right],
\]
then, we have: (a) the support $J(\hat{\beta}_{HSEL})\subseteq J(\beta^{*})$;
(b) if $\min_{j\in J(\beta^{*})}|\beta_{j}^{*}|>\bar{B}$, where 
\[
\bar{B}:=\frac{c\sqrt{k_{2}}}{\lambda_{\min}\left(\mathbb{E}\left[v_{1,J(\beta^{*})}^{T}v_{1,J(\beta^{*})}\right]\right)}\left[\left(\sigma_{v}\sigma_{\eta}\sqrt{\frac{\log p}{n}}\right)\vee\left(L\left|\beta^{*}\right|_{1}\sqrt{\Upsilon_{w,\theta^{*}}}b\left(\sigma_{v},\,\sigma_{\eta}\right)\sqrt{\frac{k_{1}\log d}{n}}\right)\right]
\]
then $J(\hat{\beta}_{HSEL})\supseteq J(\beta^{*})$ and hence $\hat{\beta}_{HSEL}$
is variable-selection consistent, i.e., $J(\hat{\beta}_{HSEL})=J(\beta^{*})$,
with probability at least $1-O\left(\frac{1}{p\wedge d}\right)$.\\
\\
\textbf{Remark}. The proof for Theorem 4.7 is provided in Section
A.8. Part (a) of Theorem 4.7 guarantees that the Lasso does not falsely
include elements that are not in the support of $\beta^{*}$. This
result hinges on Assumption 4.7, namely, the mutual incoherence condition.
Part (b) implies that as long as the minimum value of $|\beta_{j}^{*}|$
over $j\in J(\beta^{*})$ is not too small, then the two-stage Lasso
does not falsely exclude elements that are in the support of $\beta^{*}$
with high probability. Combining the claims from (a) and (b), the
multi-stage estimator is variable-selection consistent with high probability.
\\
\\
\textbf{Inference with perfect or near perfect selection}

When the mutual incoherence condition and the assumption that the
true parameters $\beta_{j}^{*}$ over $j\in J(\beta^{*})$ is well
separated from $0$ are plausible for the empirical problem of interest,
conditioning on the perfect selection or near-perfect selection result
from Theorem 4.7, we can then apply low-dimensional techniques to
conduct inference on the important coefficients. In the following
discussion, we consider the simple case where $k_{1}$ and $k_{2}$
are fixed. Then, for example, one can apply the estimator 
\[
\tilde{\beta}:=\left(\hat{v}_{\hat{J}}^{T}\hat{v}_{\hat{J}}\right)^{-1}\left(\hat{v}_{\hat{J}}^{T}\hat{v}_{0}\right)
\]
where $\hat{J}:=J(\hat{\beta}_{HSEL})$. In the multi-stage procedure
proposed by this paper, if the second-stage nonparametric estimation
uses the series estimator from Newey (1991), then the post-selection
estimator $\tilde{\beta}$ can be shown to be algebraically equivalent
to the two-stage estimator of Newey (1991) for the semiparametric
selection models when the linear coefficients in the main equation
is low-dimensional. In deriving the $\sqrt{n}-$consistency and the
asymptotic normality of the two-stage estimator, Newey requires $\sqrt{n}-$consistency
on the first-stage estimator of the coefficients in the selection
equation. This suggests that in order for the results from Newey (1991)
to be applied on the estimator $\tilde{\beta}$, perfect selection
or near-perfect selection of $\hat{\theta}$ defined in (\ref{eq:15})
may be required. We may consider a variant of $\hat{\beta}_{HSEL}$.
This variant differs from $\hat{\beta}_{HSEL}$ in that, before the
second-stage estimation, a post-Lasso procedure is performed on the
regressors selected by the first-stage estimator $\hat{\theta}$ to
obtain $\tilde{\theta}$, which is then used to form the single index.
Rather than imposing perfect selection or near-perfect selection of
$\hat{\theta}$, another option is to use the procedure proposed by
Ahn and Powell (1993), which does not require $\sqrt{n}-$consistency
on the first-stage estimator and may allow imperfect selection of
$\hat{\theta}$. For all these post-selection estimators discussed
here, the asymptotic covariance matrix is rather complicated as it
involves the derivative of the unobservable selection function. Ahn
and Powell (1993) proposes a plug-in estimator for the asymptotic
covariance matrix. Alternatively, bootstrap variance estimation can
be used to obtain the standard errors of these post-selection estimators.

It is worth noting that while selection-consistency is a desirable
property of the Lasso that allows us to conduct post-selection inference,
it requires assumptions such as the mutual incoherence condition or
the irrepresentable condition which might not hold in economic problems
where the design matrices exhibit strong (empirical) correlations.
When selection consistency is not achieved by the Lasso procedure,
other inference procedures may be useful. While it is possible to
construct confidence intervals for individual coefficients (e.g.,
Belloni, Chernozhukov, and Hansen, 2014) and linear combinations of
several of them in certain high-dimensional models using a low dimensional
projection approach (e.g., El Karoui, 2013; Zhang and Zhang, 2013;
Javanmard and Montanari, 2014), general inference theory with high-dimensional
data is still underexplored owing to the complexity of the sampling
distributions of existing estimators (see e.g., Efron, 2010). Rather
than relying on distributional theory to conduct inference, The following
section (Section 4.2) provides an alternative way of constructing
confidence sets based on the pivotal Dantzig selector (\ref{eq:10})
from Section 3. Although developing inference and asymptotic theory
for low-dimensional parameters in the high-dimensional selection models
is not the focus of this paper, it makes an interesting topic for
future research.

\subsection{The pivotal Dantzig selector of the high-dimensional linear coefficients
and confidence sets}

The pivotal Dantzig selector (\ref{eq:10}) was originally proposed
by Gautier and Tsybakov (2011) in the context of high-dimensional
IV regression. For the particular case of this paper where the instruments
are the fitted regressors $\hat{v}$ themselves, this pivotal estimator
is an extension of the Dantzig selector to accommodate for the fact
that the variance of the noise $\eta$ is unknown. It can be related
to the square-root Lasso of Belloni, Chernozhukov, and Wang (2010)
and Belloni, Chernozhukov, and Wang (2014). The non-asymptotic bounds
derived in this section only apply to the case of exactly sparse $\beta^{*}$.
However, all these results can be extended to the case of approximately
sparse $\beta^{*}$ by applying analysis similar to those from previous
sections. The confidence sets are the by-products of the non-asymptotic
bounds on the pivotal estimator. Construction of confidence sets is
based on the following theorem (Theorem 4.8), which uses a bound for
moderate deviations of self-normalized sums of random variables established
by Jing, Shao and Wang (2003). This tool was first applied by Belloni,
Chen, and Chernozhukov (2010) and later by Gautier and Tsybakov (2011)
as well as Belloni, Chernozhukov, and Wang (2014). The following assumption
is needed for this deviation bound to be applied in obtaining Theorem
4.8.\\
\textbf{}\\
\textbf{Assumption 4.8}: For all $i=1,...,n$, $j=1,...,p$ and some
constant $\delta^{'}>0$, $\mathbb{E}\left[\left|v_{ij}\eta_{i}\right|^{2+\delta^{'}}\right]<\infty$
and neither of $v_{ij}\eta_{i}$ is almost surely equal to $0$. \\
Define 
\[
b_{n,\delta^{'}}:=\min_{j=1,...,p}\frac{\sqrt{\sum_{i=1}^{n}\mathbb{E}\left[v_{ij}^{2}\eta_{i}^{2}\right]}}{\left(\sum_{i=1}^{n}\mathbb{E}\left[\left|v_{ij}\eta_{i}\right|^{2+\delta^{'}}\right]\right)^{1/(2+\delta^{'})}}.
\]
Given, for $j=1,...,p$, the variables $v_{ij}\eta_{i}$ are \textit{i.i.d.},
we have 
\begin{equation}
b_{n,\delta^{'}}:=n^{\frac{\delta^{'}}{4+2\delta^{'}}}\min_{j=1,...,p}\frac{\sqrt{\mathbb{E}\left[v_{ij}^{2}\eta_{i}^{2}\right]}}{\left(\mathbb{E}\left[\left|v_{ij}\eta_{i}\right|^{2+\delta^{'}}\right]\right)^{1/(2+\delta^{'})}}.\label{eq:23}
\end{equation}
For $a\geq1$, set 
\begin{equation}
\alpha=2L\left(1-\Phi\left(a\sqrt{2\log p}\right)\right)+2a_{0}\frac{\left(1+a\sqrt{2\log p}\right)^{1+\delta^{'}}}{p^{a^{2}-1}b_{n,\,\delta^{'}}^{2+\delta^{'}}},\label{eq:24}
\end{equation}
where $a_{0}>0$ is the absolute constant from the formula (2.11)
in Jing, Shao and Wang (2003), and $\Phi(\cdot)$ is the standard
normal c.d.f.\\
\\
\textbf{Notation}. For Theorem 4.8, define the quantities $\hat{Q}(\beta):=\frac{1}{n}\left|\hat{v}_{0}-\hat{v}\beta\right|_{2}^{2}$,
and the $l_{2}-$sensitivity 
\[
\kappa_{J(\beta^{*})}^{*}=\inf_{\Delta\in\mathbb{C}(J(\beta^{*});0,\,\varphi)}\frac{\frac{1}{n}|\hat{v}^{T}\hat{v}\Delta|_{\infty}}{|\Delta|_{2}}
\]
for some $\varphi>1$. Recall from Section 4.1 the notation 
\[
B^{'}:=\frac{1}{n}\sum_{i=1}^{n}\left[w_{i}^{T}\hat{\theta}-w_{i}^{T}\theta^{*}\right]^{2}\leq c\Upsilon_{w,\theta^{*}}\frac{k_{1}\log d}{n}
\]
where $\Upsilon_{w,\theta^{*}}$ is a known function depending only
on $w$ and $\theta^{*}$, and from Section 3 the notations 
\[
v_{j*}:=\max_{i\in\{1,...,n\}}\left\{ \left|2x_{ij}\right|\vee\left|\hat{v}_{ij}\right|\right\} 
\]
for $j=1,\,...\,,p$, and $D$ the diagonal $p\times p$ matrix with
diagonal entries $v_{j*}^{-1}$, $j=1,\,...\,,p$.\\
\\
\textbf{Remark}. Under Assumptions 4.2 and 4.3, the condition $\mathbb{E}\left[\left|v_{ij}\eta_{i}\right|^{2+\delta^{'}}\right]<\infty$
is implied by the fact that $v_{ij}$ (for all $j=1,...,p$) and $\eta$
are sub-Gaussian. To see this, note that the random variable $v_{ij}\eta_{i}$
is sub-Exponential (using the fact that the product of two sub-Gaussian
variables is sub-Exponential) and one of the characterizations of
sub-Exponential variables says a zero-mean random variable $X$ is
sub-Exponential if and only if the quantity $\sup_{k\geq2}\left[\frac{\mathbb{E}(X^{k})}{k!}\right]^{1/k}$
is finite (see, e.g., Wainwright, 2015).\\
\\
\textbf{Theorem 4.8}: Suppose the assumptions in Corollary 4.4 and
Assumption 4.8 hold. For $a\geq1$, choose $\alpha$ as in (\ref{eq:24})
and set the tuning parameter 
\begin{equation}
\xi\geq a\max\left\{ c_{0}\sqrt{\frac{\log p}{n}},\:\left(\hat{Q}(\beta^{*})\right)^{-\frac{1}{2}}\left|\beta^{*}\right|_{1}\frac{Lb(\sigma_{v})\sqrt{B^{'}}}{\min_{j=1,...,p}v_{j*}}\right\} \label{eq:25}
\end{equation}
where $c_{0}>1$ and $b(\sigma_{v})$ is some known function depending
only on $\sigma_{v}$. If $p\leq\exp\left(\frac{b_{n,\delta^{'}}^{2}}{2a^{2}}\right)$,
then with probability at least $1-\alpha-O\left(\frac{1}{p\wedge d}\right)$,
for any solution $(\hat{\beta},\,\hat{\sigma})$ of program (\ref{eq:10}),
we have 
\begin{eqnarray}
\left|D^{-1}(\hat{\beta}-\beta^{*})\right|_{2} & \leq & \frac{1}{\kappa_{J(\beta^{*})}^{*}}\left[\frac{Lb(\sigma_{v})\sqrt{B^{'}}}{\min_{j=1,...,p}v_{j*}}\left|\hat{\beta}\right|_{1}+2\xi\hat{\sigma}\right]\left[1-\frac{\xi^{2}}{\kappa_{J(\beta^{*})}^{*}}\right]^{-1}\nonumber \\
 &  & \cdot\left[1-\frac{1}{\kappa_{J(\beta^{*})}^{*}}\left[\frac{Lb(\sigma_{v})\sqrt{k_{2}B^{'}}}{\left(\min_{j=1,...,p}v_{j*}\right)^{2}}\right]\left[1-\frac{\xi^{2}}{\kappa_{J(\beta^{*})}^{*}}\right]^{-1}\right]^{-1}.\label{eq:26}
\end{eqnarray}
and, for all $j=1,...,p$, 
\begin{eqnarray}
\left|\hat{\beta}_{j}-\beta_{j}^{*}\right| & \leq & \frac{1}{v_{j*}\kappa_{J(\beta^{*})}^{*}}\left[\frac{Lb(\sigma_{v})\sqrt{B^{'}}}{\min_{j=1,...,p}v_{j*}}\left|\hat{\beta}\right|_{1}+2\xi\hat{\sigma}\right]\left[1-\frac{\xi^{2}}{\kappa_{J(\beta^{*})}^{*}}\right]^{-1}\nonumber \\
 &  & \cdot\left[1-\frac{1}{\kappa_{J(\beta^{*})}^{*}}\left[\frac{Lb(\sigma_{v})\sqrt{k_{2}B^{'}}}{\left(\min_{j=1,...,p}v_{j*}\right)^{2}}\right]\left[1-\frac{\xi^{2}}{\kappa_{J(\beta^{*})}^{*}}\right]^{-1}\right]^{-1}.\label{eq:27}
\end{eqnarray}
Furthermore, 
\begin{eqnarray}
C\hat{\sigma} & \leq & \left|\Delta_{J(\beta^{*})}\right|_{1}+C\sqrt{\hat{Q}(\beta^{*})}\nonumber \\
 & \leq & \frac{\left|\Psi_{n}\Delta\right|_{\infty}}{\kappa_{J(\beta^{*}),J(\beta^{*})}^{*}}+C\sqrt{\hat{Q}(\beta^{*})}.\label{eq:28}
\end{eqnarray}

The proof for Theorem 4.8 is provided in Section A.9. 

To construct confidence sets based on Theorem 4.8, notice that the
bounds in (\ref{eq:26})-(\ref{eq:28}) are meaningful if $\kappa_{J(\beta^{*})}^{*}\geq\bar{\kappa}>0$
(i.e., the $l_{2}-$sensitivity is strictly positive and bounded away
from $0$). In spite of the appearance, bound (\ref{eq:26}) has the
same scaling as the bound in Corollary 4.4. This can be verified by
Proposition 9.3 in Gautier and Tsybakov (2011) which shows that apart
from some positive universal constant, the $l_{2}-$sensitivity is
no smaller than the restricted eigenvalue multiplied by $k_{2}^{-\frac{1}{2}}$.
However, in cases where the $l_{2}-$sensitivity is strictly larger,
bound (\ref{eq:26}) would be sharper than the bound in Corollary
4.4. Gautier and Tsybakov (2011) provides a data-driven approach of
computing $\kappa_{J(\beta^{*})}^{*}$ without knowing $J(\beta^{*})$.
As long as the tuning parameter $\xi$ is sufficiently small, i.e.,
$\left\{ \sqrt{\frac{\log p}{n}},\:\left|\beta^{*}\right|_{1}\sqrt{\frac{k_{1}\log d}{n}}\right\} \longrightarrow0$,
then the term $1-\frac{\xi^{2}}{\kappa_{J(\beta^{*})}^{*}}$ in (\ref{eq:26})-(\ref{eq:27})
is close to $1$ provided $\kappa_{J(\beta^{*})}^{*}\geq\bar{\kappa}>0$.
The choice of $\xi$ specified by (\ref{eq:25}) has the same scaling
as the choice of $\lambda_{n,3}$ for the non-pivotal Lasso estimator
in Corollary 4.4 (in either case, the scaling of the tuning parameter
needs to match the scaling of the maximum of the first-stage related
error and the third-stage related error) except that the choice of
$\xi$ does not involve the unknown variance of $\eta_{i}$ (and hence
\textit{pivotal}). In addition, notice that the upper bounds (\ref{eq:26})-(\ref{eq:27})
are also \textit{pivotal} to the unknown variance of $\eta_{i}$ .
The only terms that can involve unknown parameters in the choice of
$\xi$ and therefore the upper bounds (\ref{eq:26})-(\ref{eq:27})
are: $Lb(\sigma_{v})\sqrt{B^{'}}$ and $\left(\hat{Q}(\beta^{*})\right)^{-\frac{1}{2}}\left|\beta^{*}\right|_{1}$. 

The term $Lb(\sigma_{v})\sqrt{B^{'}}$ is relatively easy to deal
with: $b(\sigma_{v})$ can be replaced with 
\[
b(\hat{\sigma}_{v}):=b\left(\max_{j}\sqrt{\frac{1}{n}\sum_{i=1}^{n}\hat{v}_{ij}^{2}}\right)
\]
and an application of Lemma A.11 yields $\hat{\sigma}_{v}\leq2\sigma_{v}$
with probability at least $1-O\left(\frac{1}{p\wedge d}\right)$;
construction of confidence intervals (that do not contain any unknown
parameters) for $B^{'}$ has been considered in the context of several
Generalized Linear models (see, e.g., Van de Geer, 2008) and we will
assume in this discussion that these confidence sets $\hat{B}^{'}$
for $B^{'}$ are available. Consequently, whenever the term $Lb(\sigma_{v})\sqrt{B^{'}}$
shows up in the bounds (\ref{eq:26})-(\ref{eq:27}), we will replace
it with $Lb(\hat{\sigma}_{v})\sqrt{\hat{B}^{'}}$. In the case where
the\textit{ }constant \textit{$L$ }is unknown, Section 4.1.2 discusses
methods to determine this constant.

The term $\left(\hat{Q}(\beta^{*})\right)^{-\frac{1}{2}}\left|\beta^{*}\right|_{1}$
is the harder one here as $\beta^{*}$ is unknown and in fact the
parameters we want to estimate. One possibility is to consider the
following heuristic: 
\begin{enumerate}
\item In Step $k=0$ (initialization), solve program (\ref{eq:10}) with
$\xi^{k}=c_{0}\sqrt{\frac{\log p}{n}}$ to obtain $\hat{\beta}^{k}$
for some $c_{0}>1$; update $\xi^{k}$ with 
\begin{equation}
\xi^{k+1}\geq a\max\left\{ c_{0}\sqrt{\frac{\log p}{n}},\:\left(\hat{Q}(\hat{\beta}^{k})\right)^{-\frac{1}{2}}\left|\hat{\beta}^{k}\right|_{1}\frac{\hat{L}b(\hat{\sigma}_{v})\sqrt{\hat{B}^{'}}}{\min_{j=1,...,p}v_{j*}}\right\} .\label{eq:29}
\end{equation}

\item In Step $k+1$, solve program (\ref{eq:10}) with $\xi^{k+1}$ to
obtain $\hat{\beta}^{k+1}$ and update $\xi^{k+1}$ with $\xi^{k+2}$
as in (\ref{eq:29}). Repeat this step till a pre-specified tolerance
level on $\left|\hat{\beta}^{k+1}-\hat{\beta}^{k}\right|_{2}$ is
reached. 
\end{enumerate}
Establishing theoretical guarantees for the heuristic provided above
is pursued in a separate ongoing project. In special cases, we may
be able to circumvent the fact that $\beta^{*}$ is unknown. For example,
when $p$ is large relative to $d$ so that $\left|\beta^{*}\right|_{1}\sqrt{\frac{k_{1}\log d}{n}}\ll\sqrt{\frac{\log p}{n}}$,
then the result in Theorem 4.8 is essentially reduced to the case
where the pivotal Dantzig selector is applied to the standard high-dimensional
linear models with exact sparsity. In a related scenario where a post-Lasso
procedure is performed on the regressors selected by the first-stage
estimator $\hat{\theta}$ defined in (\ref{eq:13}), upon perfect
selection or near-perfect selection of $\hat{\theta}$, the factor
$\sqrt{\frac{k_{1}\log d}{n}}$ is reduced to $\sqrt{\frac{k_{1}}{n}}$
which may be smaller relative to $\sqrt{\frac{\log p}{n}}$.

\subsection{Properties of the estimators of the selection bias function}

Given the availability of estimates $\hat{\theta}$ and $\hat{\beta}$
of the high-dimensional linear coefficients from either the non-pivotal
procedure or the pivotal procedure, two different estimation strategies
for the nonparametric selection bias are considered: one is the closed
form estimator (\ref{eq:15}) and the other is the plug-in nonparametric
least squares estimator (\ref{eq:16}) which can be obtained from
the Lipschitz regression described in Section 4.1.2 if we assume $g(\cdot)$
belongs to the class $\mathcal{F}$ of Lipschitz functions. Despite
the fact that (\ref{eq:16}) is computationally more involved relative
to (\ref{eq:15}), its rate of convergence turns out to be faster
as shown in the following. To facilitate the discussion and a later
comparison with the minimax lower bounds in Section 4.4, we break
down the presentations of the results into the case of exact sparsity
on $\beta^{*}$ and $\theta^{*}$ ($q_{1}=q_{2}=0$) in Theorems 4.9
and 4.10, and the case of general sparsity on $\beta^{*}$ and $\theta^{*}$
($q_{1},\, q_{2}\in[0,\,1]$) in Theorems 4.11 and 4.12 (which contain
Theorems 4.9 and 4.10 as special cases, respectively).\\
\textbf{}\\
\textbf{Theorem 4.9} ($q_{1}=q_{2}=0$): Let the assumptions in Corollary
4.4 hold. Suppose $g(\cdot)$ belongs to the class $\mathcal{F}$
of Lipschitz functions. For the estimator $\hat{g}(\cdot)$ of $g(\cdot)$
obtained by (\ref{eq:15}), 
\[
\left(\mathbb{E}\left[\hat{g}(w_{i}^{T}\hat{\theta})-g(w_{i}^{T}\theta^{*})\right]^{2}\right)^{\frac{1}{2}}\leq cb\max\left\{ k_{2}\sqrt{\frac{\log p}{n}},\:\left|\beta^{*}\right|_{1}\left(\frac{k_{1}\log d}{n}\right)^{\frac{1}{4}},\: k_{2}\left|\beta^{*}\right|_{1}\sqrt{\frac{k_{1}\log d}{n}}\right\} 
\]
where $b$ is some constant depending only on the model-specific structure
(and independent of $n$, $d$, $p$, $k_{1}$, and $k_{2}$).\textbf{}\\
\textbf{}\\
\textbf{Theorem 4.10} ($q_{1},\, q_{2}\in[0,\,1]$): Let the assumptions
in Corollary 4.6 hold. Suppose $g(\cdot)$ belongs to the class $\mathcal{F}$
of Lipschitz functions. For the estimator $\hat{g}(\cdot)$ of $g(\cdot)$
obtained by (\ref{eq:15}), 
\[
\left(\mathbb{E}\left[\hat{g}(w_{i}^{T}\hat{\theta})-g(w_{i}^{T}\theta^{*})\right]^{2}\right)^{\frac{1}{2}}\leq cb\max\left\{ R_{q_{2}}\bar{M}^{1-q_{2}},\,\:\left|\beta^{*}\right|_{1}R_{q_{1}}^{\frac{1}{4}}\left(\sqrt{\frac{\log d}{n}}\right)^{\frac{1}{2}-\frac{q_{2}}{4}}\,\left|\beta^{*}\right|_{1}\left(\frac{\left|\theta^{*}\right|_{1}}{n}\right)^{\frac{1}{3}}\right\} ,
\]
where $\bar{M}$ is defined in Corollary 4.6 and $b$ is some constant
depending only on the model-specific structure (and independent of
$n$, $d$, $p$, $R_{q_{1}}$, and $R_{q_{2}}$).\textbf{}\\
\textbf{}\\
\textbf{Theorem 4.11} ($q_{1}=q_{2}=0$): Let the assumptions in Corollary
4.4 hold. Suppose $g(\cdot)$ belongs to the class $\mathcal{F}$
of Lipschitz functions and the random matrix $x$ is sub-Gaussian
with parameters $(\Sigma_{x},\,\sigma_{x}^{2})$. For all $\Delta\in\mathbb{C}(J(\beta^{*});\,0,\,3)\backslash\{\mathbf{0}\}$,
the matrix $\Sigma_{x}$ satisfies $\frac{\Delta^{T}\Sigma_{x}\Delta}{\left|\Delta\right|_{2}^{2}}\leq\kappa_{U}^{x}<\infty$.
For the estimator $\tilde{g}(\cdot)$ of $g(\cdot)$ obtained by (\ref{eq:16}),
\[
\left(\mathbb{E}\left[\tilde{g}(w_{i}^{T}\hat{\theta})-g(w_{i}^{T}\theta^{*})\right]^{2}\right)^{\frac{1}{2}}\leq c^{'}b^{'}\max\left\{ \sqrt{\frac{k_{2}\log p}{n}},\;\left|\beta^{*}\right|_{1}\sqrt{\frac{k_{1}k_{2}\log d}{n}},\;\left(\frac{k_{1}}{n}\right)^{\frac{1}{3}}\right\} 
\]
where $b^{'}$ is some constant depending only on the model-specific
structure (and independent of $n$, $d$, $p$, $k_{1}$, and $k_{2}$).
\textbf{}\\
\textbf{}\\
\textbf{Theorem 4.12} ($q_{1},\, q_{2}\in[0,\,1]$): Let the assumptions
in Corollary 4.6. Suppose $g(\cdot)$ belongs to the class $\mathcal{F}$
of Lipschitz functions and the random matrix $x$ is sub-Gaussian
with parameters $(\Sigma_{x},\,\sigma_{x}^{2})$. For all non-zero
$\Delta\in\mathbb{C}(S_{\underbar{\ensuremath{\tau}}};\, q_{2},\,3)\cap\mathbb{S}_{\delta}$
where $\mathbb{C}(S_{\underbar{\ensuremath{\tau}}};\, q_{2},\,3)\cap\mathbb{S}_{\delta}$
is defined in Corollary 4.6, the matrix $\Sigma_{x}$ satisfies $\frac{\Delta^{T}\Sigma_{x}\Delta}{\left|\Delta\right|_{2}^{2}}\leq\kappa_{U}^{x}<\infty$.
For the estimator $\tilde{g}(\cdot)$ of $g(\cdot)$ obtained by (\ref{eq:16}),
\[
\left(\mathbb{E}\left[\tilde{g}(w_{i}^{T}\hat{\theta})-g(w_{i}^{T}\theta^{*})\right]^{2}\right)^{\frac{1}{2}}\leq c^{'}b^{'}\max\left\{ \sqrt{R_{q_{2}}}\bar{M}^{1-\frac{q_{2}}{2}},\;\left(\frac{\left|\theta^{*}\right|_{1}}{n}\right)^{\frac{1}{3}}\right\} ,
\]
where $\bar{M}$ is defined in Corollary 4.6 and $b^{'}$ is some
constant depending only on the model-specific structure (and independent
of $n$, $d$, $p$, $R_{q_{1}}$, and $R_{q_{2}}$). \textbf{}\\
\textbf{}\\
\textbf{Remark}. The proofs for Theorems 4.9-4.12 are provided in
Sections A.10 and A.11. First let us look at the case of exactly sparse
$\beta^{*}$ and $\theta^{*}$ ($q_{1}=q_{2}=0$). From Theorem 4.11,
notice that the terms $\sqrt{\frac{k_{2}\log p}{n}}$ and $\left|\beta^{*}\right|_{1}\sqrt{\frac{k_{1}k_{2}\log d}{n}}$
are expected from the statistical error of $\hat{\beta}$ that we
plug into the nonparametric regression (\ref{eq:16}); and the term
$\left(\frac{k_{1}}{n}\right)^{\frac{1}{3}}$ is expected from the
fact that $g(\cdot)$ belongs to the class $\mathcal{F}$ of Lipschitz
functions%
\footnote{Note that when $\epsilon_{1i}$ and $\epsilon_{2i}$ in (\ref{eq:1})
are bivariate normal, the selection bias characterized by the Inverse
Mills Ratio is a $1-$Lipschitz function (see, e.g., Ruud, 2000).
Furthermore, if $m_{j}(\cdot)\in\mathcal{F}_{j}$ and $\mathcal{F}_{j}$
is the class of $L-$Lipschitz functions for every $j=1,...,p$, then
$g(\cdot)$ is a Lipschitz function.%
}. On the other hand, the term $\left(\frac{k_{1}}{n}\right)^{\frac{1}{3}}$
is suppressed by $\left(\frac{k_{1}\log d}{n}\right)^{\frac{1}{4}}$
in Theorem 4.9 for the closed-form estimator (\ref{eq:15}). When
$\beta^{*}$ is approximately sparse with $q_{2}=1$, Theorem 4.10
implies that the \textit{$\sqrt{MSE}$} of the closed-form estimator
(\ref{eq:15}) is bounded above by $R_{q_{2}}\bar{M}^{1-q_{2}}=\left|\beta^{*}\right|_{1}$.
This upper bound is unimprovable and as a result, it is not possible
for (\ref{eq:15}) to achieve \textit{MSE}-consistency even if $n\rightarrow\infty$
when $q_{2}=1$. In contrast to (\ref{eq:15}), the nonparametric
least squares estimator (\ref{eq:16}) is consistent in \textit{MSE}
as $n\rightarrow\infty$ when $q_{2}=1$. The key behind the sharp
rate achieved by the plug-in nonparametric least squares estimator
(\ref{eq:16}) in Theorems 4.11 and 4.12 lies on the random variables
\[
\frac{1}{n}\sum_{i=1}^{n}\eta_{i}\left[\tilde{g}(w_{i}^{T}\theta^{*})-g(w_{i}^{T}\theta^{*})\right],
\]
and 
\[
U_{n}:=\sup_{\delta\in\mathcal{S}(r_{1},\, r_{2})}\frac{1}{n}\left|\eta^{T}w\delta\right|,
\]
where 
\[
\mathcal{S}(r_{1},\, r_{2}):=\left\{ \delta\in\mathbb{R}^{d}\,\vert\,\left|\delta\right|_{1}\leq r_{1},\:\left|\delta\right|_{2}\leq r_{2}\right\} .
\]
The analysis for controlling the first term uses a {}``local function
complexity'' argument similar to what is done in the proofs for Theorems
4.1 and 4.2. To upper bound the second term $U_{n}$, we can apply
a discretization argument over the set $\mathcal{S}(r_{1},\, r_{2})$
together with results on metric entropy and the fact $\mathbb{E}\left[\eta_{i}\vert w_{i}\right]=0$.
{}``Small'' values of $r_{1}$ and $r_{2}$ are guaranteed by the
upper bounds on $|\hat{\theta}-\theta^{*}|_{2}$ from Lemma A.7 and
as a result we only need to work with a {}``small'' $\mathcal{S}(r_{1},\, r_{2})$.
The sharp rates provided by these types of analysis seem to be driven
by the projection nature of the underlying nonparametric least-squares
estimators%
\footnote{In fact, a recent paper by Chatterjee (2014) shows that the least
squares estimators are always \textit{admissible} up to a universal
constant in many modern statistics problems.%
}. As we will see in the following section, the overall convergence
rate of the estimator $\hat{\beta}$ (obtained by either the non-pivotal
procedure or the pivotal procedure) and the plug-in nonparametric
least squares estimator (\ref{eq:16}) is \textit{minimax optimal}
in terms of the $(n,\, d,\, p)-$scaling for the case of exactly sparse
$\beta^{*}$. However, we will also see that this minimax optimality
result does not apply to the case of approximately sparse $\beta^{*}$
because of the first-stage related estimation error.

\subsection{Statistical efficiency via lower bounds on minimax risks}

This section studies efficiency of the proposed estimators by deriving
lower bounds on minimax rate for the case of $l_{2}-$loss. Complementary
to the understanding of computationally efficient procedures are the
information-theoretic limitations of statistical estimation, applicable
to any procedure regardless of its computational cost. There is a
rich literature on such information-theoretic limits, which can have
two types of consequences. First, they can reveal gaps between the
performance of an optimal procedure in theory and known computationally
efficient methods. Second, they can demonstrate regimes in which practical
procedures achieve these information-theoretic limits. While one way
of interpreting minimax lower bounds is to view the choice of unknown
parameters in an adversarial manner, and to compare the estimators
based on their worst-case performance, many techniques for deriving
minimax lower bounds can be motivated by the Bayesian approach which
views the unknown parameters as random variables (e.g., Guntuboyina
2011).

The minimax lower bounds in this section are derived for model (\ref{eq:2}),
implied by the original selection model (\ref{eq:1}). As a consequence,
these lower bounds provide information-theoretic limits for any procedure
constructed based on model (\ref{eq:2}) for estimating model (\ref{eq:1}),
regardless of its computational cost. For $q_{1},\, q_{2}\in[0,\,1]$,
define $\mathcal{H}=\mathcal{B}_{q_{2}}^{p}(R_{q_{2}})\times\mathcal{F}\circ\mathcal{B}_{q_{1}}^{d}(R_{q_{1}})$,
where the $l_{q}-$ {}``ball'' is defined in Section 2 and $\mathcal{F}$
is the class of functions such that $g\in\mathcal{F}\,:\,\mathbb{R}\rightarrow\mathbb{R}$.
When $\beta^{*}\in\mathcal{B}_{0}^{p}(k_{2})$ and $\theta^{*}\in\mathcal{B}_{0}^{d}(k_{1})$,
model (\ref{eq:2}) corresponds to the case of exact sparsity on $\beta^{*}$
and $\theta^{*}$. When $\beta^{*}\in\mathcal{B}_{q_{2}}^{p}(R_{q_{2}})$
and $\theta^{*}\in\mathcal{B}_{q_{1}}^{d}(R_{q_{1}})$ for $q\in(0,\,1]$,
model (\ref{eq:2}) corresponds to the case of approximate sparsity
based upon imposing a certain decay rate on the ordered entries of
$\beta^{*}$ and $\theta^{*}$. Theorem 4.13 (Theorem 4.14) presents
a minimix lower bound for the case of exact sparsity $q_{1}=q_{2}=0$
(respectively, the case of approximate sparsity $q_{1},\, q_{2}\in(0,\,1]$).\\
\\
\textbf{Assumption 4.9}: There exists a constant $\underbar{\ensuremath{\kappa}}_{x}>0$
and a function $f_{l}(R_{q_{2}},\, q_{2},\, n,\, p)$ such that 
\[
\frac{1}{\sqrt{n}}\left|x\beta\right|_{2}\geq\underbar{\ensuremath{\kappa}}_{x}\left|\beta\right|_{2}-f_{l}(R_{q_{2}},\, q_{2},\, n,\, p)\quad\textrm{for all}\:\beta\in\mathcal{B}_{q_{2}}^{p}(R_{q_{2}}).
\]
\\
\textbf{Assumption 4.10}: There is no measurable function $f(w_{i}^{T}\theta)$
such that $x_{i}^{T}\lambda=f(w_{i}^{T}\theta)$ when $y_{1i}=1$
for $\lambda\in\mathcal{B}_{q_{2}}^{p}(R_{q_{2}})\backslash\{\mathbf{0}\}$.\\
\\
\textbf{Remark}. Assumptions 4.9 and 4.10 ensure the identifiability
of model (\ref{eq:2}), without which, lower bounds for high-dimensional
linear models usually involve a maximum of two quantities: a term
involving the diameter of the null-space restricted to the $l_{q}-$ball,
measuring the degree of non-identifiability of the model, and a term
arising from the metric entropy structure for $l_{q}-$balls (see
Raskutti, et. al, 2011). Assumption 4.9 together with Assumption 4.10
incurs an upper bound on the $\mathcal{B}_{q}^{p}(R)-$kernel diameter
in $l_{2}-$norm (this result is formalized in Lemma A.10 and proved
in Section A.12), and consequently the identifiability of model (\ref{eq:2}).
\\
\\
\textbf{Theorem 4.13} ($q_{1}=q_{2}=0$): Let $\mathcal{F}$ be the
class of $L-$Lipschitz functions and Assumptions 4.9-4.10 hold with
$f_{l}(R_{q_{2}},\, q_{2},\, n,\, p)=0$ and $\underbar{\ensuremath{\kappa}}_{x}>0$.
Define the parameter space $\Theta$ as 
\[
\left\{ \theta\in\mathcal{B}_{0}^{d}(k_{1}):\,\textrm{for any}\:\lambda\in\mathcal{B}_{0}^{p}(k_{2})\backslash\{\mathbf{0}\},\:\nexists\textrm{ a measurable }g(\cdot)\in\mathcal{F}\textrm{ s.t. }x_{i}^{T}\lambda=g(w_{i}^{T}\theta)\textrm{ when }y_{1i}=1\right\} .
\]
Moreover, $\frac{\left|w\theta\right|_{2}}{\sqrt{n}\left|\theta\right|_{2}}\leq\kappa_{u}$
for all $\theta\in\Theta$ and $\frac{\left|x\beta\right|_{2}}{\sqrt{n}\left|\beta\right|_{2}}\leq\kappa_{u}^{'}$
for all $\beta\in\mathcal{B}_{0}^{p}(k_{2})$. If the vector $\eta\sim N(0,\,\sigma_{\eta}I_{n\times n})$,
then, for some constant $b$ depending only on the model-specific
structure (and independent of $n$, $d$, $p$, $k_{1}$, and $k_{2}$),
\begin{eqnarray*}
 &  & \min_{\tilde{\beta},\,\tilde{f},\,\tilde{\theta}}\max_{\begin{array}{c}
\theta\in\Theta\\
f(\cdot)\in\mathcal{F}\\
\beta\in\mathcal{B}_{0}^{p}(k_{2})
\end{array}}\left(\mathbb{E}\left|\tilde{\beta}-\beta\right|_{2}^{2}\right)^{\frac{1}{2}}+\left(\mathbb{E}\left[\tilde{f}(w_{i}^{T}\tilde{\theta})-f(w_{i}^{T}\theta)\right]^{2}\right)^{\frac{1}{2}}\\
 & \geq & b\max\left\{ \sqrt{\frac{k_{1}\log d}{n}},\,\left(\frac{k_{1}}{n}\right)^{\frac{1}{3}},\,\sqrt{\frac{k_{2}\log p}{n}}\right\} .
\end{eqnarray*}
\textbf{Theorem 4.14} ($q_{1},\, q_{2}\in(0,\,1]$): Let $\mathcal{F}$
be the class of $L-$Lipschitz functions and Assumptions 4.9-4.10
hold with $f_{l}(R_{q_{2}},\, q_{2},\, n,\, p)=o\left(R_{q_{2}}^{\frac{1}{2}}\left(\frac{\log p}{n}\right)^{\frac{1}{2}-\frac{q_{2}}{4}}\right)$
and $\underbar{\ensuremath{\kappa}}_{x}>0$. Moreover, $\frac{1}{\sqrt{n}}\max_{j=1,...,d}\left|w_{j}\right|_{2}\leq\kappa_{w}<\infty$
and $\frac{1}{\sqrt{n}}\max_{j=1,...,p}\left|x_{j}\right|_{2}\leq\kappa_{x}<\infty$.
If the vector $\eta\sim N(0,\,\sigma_{\eta}I_{n\times n})$, then,
for some constant $b^{'}$ depending only on the model-specific structure
(and independent of $n$, $d$, $p$, $R_{q_{1}}$, and $R_{q_{2}}$),
\begin{eqnarray*}
 &  & \min_{\tilde{\beta},\,\tilde{f},\,\tilde{\theta}}\max_{\begin{array}{c}
\theta\in\Theta\\
f(\cdot)\in\mathcal{F}\\
\beta\in\mathcal{B}_{q_{2}}^{p}(R_{q_{2}})
\end{array}}\left(\mathbb{E}\left|\tilde{\beta}-\beta\right|_{2}^{2}\right)^{\frac{1}{2}}+\left(\mathbb{E}\left[\tilde{f}(w_{i}^{T}\tilde{\theta})-f(w_{i}^{T}\theta)\right]^{2}\right)^{\frac{1}{2}}\\
 & \geq & b^{'}\max\left\{ R_{q_{1}}^{\frac{1}{2}}\left(\frac{\log d}{n}\right)^{\frac{2-q_{1}}{4}},\,\left(\frac{R^{*}}{n}\right)^{\frac{1}{3}},\, R_{q_{2}}^{\frac{1}{2}}\left(\frac{\log p}{n}\right)^{\frac{2-q_{1}}{4}}\right\} ,
\end{eqnarray*}
where the parameter space $\Theta$ is defined in Theorem 4.13 with
$\mathcal{B}_{0}^{d}(k_{1})$ replaced by $\mathcal{B}_{q_{1}}^{d}(R_{q_{1}})$
and $\mathcal{B}_{0}^{p}(k_{2})$ replaced by $\mathcal{B}_{q_{2}}^{p}(R_{q_{2}})$,
and $R^{*}$ is the {}``radius'' $R_{q_{1}}$ when $q_{1}=1$.\\
\textbf{}\\
\textbf{Remark}. The proofs for Theorem 4.13 and Theorem 4.14 are
provided in Sections A.12 and A.13, respectively. These proofs are
information-theoretic in nature and based on Fano's inequality (see,
e.g., Guntuboyina, 2011; Wainwright, 2015) and results on the metric
entropy of the $l_{q}-$balls. By Lemma A.10, the conditions on $f_{l}(R_{q_{2}},\, q_{2},\, n,\, p)$
in Theorems 4.13 and 4.14 together with Assumption 4.10 ensure that
the kernel diameter for the nullspace of $\mathcal{B}_{q_{2}}^{p}(R_{q_{2}})$
is dominated by the term related to the metric entropy of $\mathcal{B}_{q_{2}}^{p}(R_{q_{2}})$.
In Theorem 4.13, we require Assumption 4.9 to hold with $f_{l}(R_{q_{2}},\, q_{2},\, n,\, p)=0$
and $\underbar{\ensuremath{\kappa}}_{x}>0$, which is closely related
to the restricted eigenvalue condition on the matrix $\frac{x^{T}x}{n}$
over the set $\mathcal{B}_{0}^{p}(k_{2})$. When $x$ is a sub-Gaussian
matrix with parameters $(\Sigma_{x},\,\sigma_{x}^{2})$ and for all
$\Delta\in\mathcal{B}_{0}^{p}(k_{2})\backslash\{\mathbf{0}\}$, the
matrix $\Sigma_{x}$ satisfies $\frac{\Delta^{T}\Sigma_{x}\Delta}{\left|\Delta\right|_{2}^{2}}\geq\underbar{\ensuremath{\kappa}}_{x}>0$,
then Lemma B.2 guarantees Assumption 4.9 to hold for $f_{l}(R_{q_{2}},\, q_{2},\, n,\, p)=0$
with high probability. Additionally, if $w$ is a sub-Gaussian matrix
with parameters $(\Sigma_{w},\,\sigma_{w}^{2})$, $\frac{\Delta^{T}\Sigma_{w}\Delta}{\left|\Delta\right|_{2}^{2}}\leq\bar{\kappa}^{w}<\infty$
for all $\theta\in\Theta$, and $\frac{\Delta^{T}\Sigma_{x}\Delta}{\left|\Delta\right|_{2}^{2}}\leq\bar{\kappa}^{x}<\infty$
for all $\beta\in\mathcal{B}_{0}^{p}(k_{2})$, Lemma B.2 also guarantees
that $\frac{\left|w\theta\right|_{2}^{2}}{n\left|\theta\right|_{2}^{2}}\leq c\bar{\kappa}^{w}$
for all $\theta\in\Theta$ and $\frac{\left|x\beta\right|_{2}^{2}}{n\left|\beta\right|_{2}^{2}}\leq c\bar{\kappa}^{x}$
for all $\beta\in\mathcal{B}_{0}^{p}(k_{2})$ hold with high probability
(in Theorem 4.13). Similarly, the conditions $\frac{1}{\sqrt{n}}\max_{j=1,...,d}\left|w_{j}\right|_{2}\leq\kappa_{w}<\infty$
and $\frac{1}{\sqrt{n}}\max_{j=1,...,p}\left|x_{j}\right|_{2}\leq\kappa_{x}<\infty$
(in Theorem 4.14) are also implied by Lemma B.2 with high probability
given $w_{j}$ ($j=1,...,d$) and $x_{j}$ ($j=1,...,p$) are sub-Gaussian. 

Compare the scaling of the lower bound in Theorem 4.13 with the upper
bounds in Corollary 4.4 and Theorem 4.11. In particular, from the
previous upper bounds, we have 
\[
\left(\mathbb{E}\left|\hat{\beta}_{HSEL}-\beta^{*}\right|_{2}^{2}\right)^{\frac{1}{2}}+\left(\mathbb{E}\left[\tilde{g}(w_{i}^{T}\hat{\theta})-g(w_{i}^{T}\theta^{*})\right]^{2}\right)^{\frac{1}{2}}\precsim\max\left\{ \left|\beta^{*}\right|_{1}\sqrt{\frac{k_{2}k_{1}\log d}{n}},\,\left(\frac{k_{1}}{n}\right)^{\frac{1}{3}},\,\sqrt{\frac{k_{2}\log p}{n}}\right\} 
\]
(The upper bound on $\mathbb{E}\left|\hat{\beta}_{HSEL}-\beta^{*}\right|_{2}^{2}$
is obtained by converting $\left|\hat{\beta}_{HSEL}-\beta^{*}\right|_{2}^{2}$
with a standard integration over the tail probability in the exponential
form). Notice that the scaling in the upper bound above matches the
lower bound in Theorem 4.13 in terms of $(n,\, d,\, p)-$factors.
The only difference between these bounds is that the upper bound exceeds
the lower bound by a factor of $\left|\beta^{*}\right|_{1}\sqrt{k_{2}}\asymp k_{2}^{\frac{3}{2}}$
in the term related to the complexity of the set $\Theta$, which
is likely due to the fact that the estimator $\hat{\beta}_{HSEL}$
is a sequential multi-stage procedure based on plugging in the first-stage
estimator $\hat{\theta}$ in the place of the unknown coefficient
$\theta^{*}$ in the selection equation. In a different but somewhat
related context which concerns with the high-dimensional sparse linear
regression models with many endogenous regressors and instruments
(see Zhu 2014), it is found that the upper bound on the \textit{$\sqrt{MSE}$}
of the $l_{1}-$regularized two-stage estimator exceeds the minimax
lower bound in Zhu (2014) by a factor of $k_{2}$ (where $k_{2}$
is the sparsity parameter for the second-stage model).

On the other hand, this minimax optimality result does not apply to
the case of approximately sparse $\beta^{*}$ when we compare the
scaling of the lower bound in Theorem 4.14 with the upper bounds in
Corollary 4.6 and Theorem 4.12 for the case $q_{2}\in(0,\,1]$. In
particular, from the previous upper bounds, we have, 
\begin{eqnarray*}
 & \left(\mathbb{E}\left|\hat{\beta}_{HSEL}-\beta^{*}\right|_{2}^{2}\right)^{\frac{1}{2}}+\left(\mathbb{E}\left[\tilde{g}(w_{i}^{T}\hat{\theta})-g(w_{i}^{T}\theta^{*})\right]^{2}\right)^{\frac{1}{2}}\\
\precsim & \left\{ \sqrt{R_{q_{2}}}\left[\max\left\{ \sqrt{\frac{\log p}{n}},\:\left|\beta^{*}\right|_{1}\left(\frac{\left|\theta^{*}\right|_{1}}{n}\right)^{\frac{2}{3}},\:\left|\beta^{*}\right|_{1}R_{q_{1}}^{\frac{1}{2}}\left(\sqrt{\frac{\log d}{n}}\right)^{1-\frac{q_{1}}{2}}\right\} \right]^{1-\frac{q_{2}}{2}}\right\} \vee\left\{ \left(\frac{\left|\theta^{*}\right|_{1}}{n}\right)^{\frac{1}{3}}\right\}  & .
\end{eqnarray*}
As in the case of exactly sparse $\beta^{*},$ the terms $\sqrt{R_{q_{2}}}\left(\sqrt{\frac{\log p}{n}}\right)^{1-\frac{q_{2}}{2}}$
and $\left(\frac{\left|\theta^{*}\right|_{1}}{n}\right)^{\frac{1}{3}}$
in the above upper bound match the scalings of the term related to
the complexity of the set $\mathcal{B}_{q_{2}}^{p}(R_{q_{2}})$ and
the term related to the complexity of the set $\mathcal{F}$, respectively.
In sharp contrast to the case of exactly sparse $\beta^{*}$ where
our sequential multi-stage procedure based on plugging in the first-stage
estimator $\hat{\theta}$ only exceeds the minimax optimal result
by a factor of $k_{2}^{\frac{3}{2}}$ in the term related to the complexity
of the set $\Theta$, the term $\sqrt{R_{q_{2}}}\left[\left|\beta^{*}\right|_{1}R_{q_{1}}^{\frac{1}{2}}\left(\sqrt{\frac{\log d}{n}}\right)^{1-\frac{q_{1}}{2}}\right]^{1-\frac{q_{2}}{2}}$
in the upper bound above is now worsened by an exponent of $1-\frac{q_{2}}{2}$
and a factor of $\sqrt{R_{q_{2}}}\left(\left|\beta^{*}\right|_{1}\right)^{1-\frac{q_{2}}{2}}$
when compared to the term related to $\Theta$, $R_{q_{2}}^{\frac{1}{2}}\left(\sqrt{\frac{\log p}{n}}\right)^{1-\frac{q_{1}}{2}}$,
in the lower bound of Theorem 4.14. When $q_{2}\in[0,\,1]$, note
that 
\[
\left[R_{q_{1}}^{\frac{1}{2}}\left(\sqrt{\frac{\log d}{n}}\right)^{1-\frac{q_{1}}{2}}\right]^{1-\frac{q_{2}}{2}}\geq R_{q_{1}}^{\frac{1}{2}}\left(\sqrt{\frac{\log d}{n}}\right)^{1-\frac{q_{1}}{2}}
\]
with {}``$=$'' holds only if $q_{2}=0$ (the case of exactly sparse
$\beta^{*}$). 

The lower bound in either Theorem 4.13 or Theorem 4.14 is a {}``point''
result. Even if the main equation in the original selection model
(\ref{eq:1}) has a normal error, the normality of $\eta$ is plausible
in model (\ref{eq:2}) only if $g(w_{i}^{T}\theta^{*})=0$, i.e.,
when there is no selection activity. Nevertheless, these {}``point''
results provided by Theorems 4.13 and 4.14 are still useful because
whether $g(w_{i}^{T}\theta^{*})$ equals $0$ or not would be unknown
in general and the error from having to estimate $g(\cdot)$ still
appears in the lower bounds. Moreover, even if the {}``point'' result
does not hold {}``globally'', given that the lower bounds are derived
for the minimax risks of the high-dimensional linear coefficients
together with the nonparametric selection bias function, at least
the second and third terms in the lower bounds of Theorems 4.13 and
4.14 should be unimprovable in any {}``global'' result. It is possible
to impose distributional assumptions other than normality on $\eta$
but the derivation of the lower bounds in the proofs may involve more
difficult computations related to the Kullback-Leibler divergence
or the more general $f-$divergence where $f$ is a convex function
with $f(1)=0$ (see Guntuboyina 2011 for a unified treatment of existing
techniques for obtaining lower bounds). For this reason, existing
literature on minimax lower bounds almost exclusively focuses on the
case of normal errors and lower bounds with less restrictive distributional
assumptions other than normality (e.g., sub-Gaussianity) on a random
vector are in general impossible to obtain. Recent work of efficiency
bounds (e.g., Hansen B., 2014) that proposes a shrinking neighborhood
analysis may provide a promising direction for extending these {}``point''
results to the case where $g(w_{i}^{T}\theta^{*})$ is in a shrinking
neighborhood of $0$.

\subsection{Estimation of high-dimensional semilinear models with a two-stage
projection strategy}

In this section, we discuss how the theory developed in this paper
can be applied to the semilinear models in high-dimensional settings.
The multi-stage estimator proposed in this paper is also useful for
estimating the linear coefficients of the following semilinear model:
\begin{equation}
y_{i}=x_{i}^{T}\beta^{*}+g(w_{i})+\eta_{i}\label{eq:30}
\end{equation}
where $x_{i}$ is a $p-$dimensional vector of regressors (and $p$
can grow with and exceed the sample size $n$). Furthermore, $g(\cdot):\,\mathbb{R}^{d}\rightarrow\mathbb{R}$
is an unknown function and $w_{i}$ is a $d-$dimensional vector of
regressors. Our multi-stage projection strategy is now reduced to
a two-stage procedure. Based on the analysis in this paper, it is
straightforward to see that Theorems 4.1 and 4.2 remain valid except
that there is no first-stage related error $\mathcal{T}_{1}$ in the
upper bounds on the estimator of $\beta^{*}$. As we have mentioned
before in Section 4.1, when $p$ and $d$ are fixed and small relative
to $n$, as long as $\mathcal{F}_{j}$ ($j=0,...,p$) in Theorems
4.1 and 4.2 are sufficiently smooth so that $\mathbb{E}(x_{ij}\,\vert\, w_{i})$
and $\mathbb{E}(y_{i}\,\vert\, w_{i})$ can be estimated at a rate
\textit{no slower} than $O\left(n^{-\frac{1}{4}}\right)$, the {}``oracle''
property will be achieved. 

The case where the dimensions $p$ and $d$ are both large relative
to $n$ (namely, $p\geq n$ and $d\geq n$) generalizes the semilinear
model considered in Belloni, et. al (2014) in which $d\geq n$ and
$p$ remains finite. When $p\geq n$, it is unclear whether the procedure
proposed by Belloni, et. al (2014)%
\footnote{When $d\geq n$ and $p$ remains finite, the procedure from Belloni,
et. al (2014) includes three steps: First, apply the Lasso to the
regression of $y_{i}$ on $w_{i}$; second, apply the Lasso to the
regression of $x_{ij}$ on $w_{i}$ for every $j=1,...,p$, respectively;
and third, regress $y_{i}$ on $x_{i}$ and the components of $w_{i}$
selected by the first and second step.%
} can be easily extended because the effect from imperfect selection
in the second step of Belloni, et. al (2014) may not be negligible
anymore when the number of components in $x_{i}$ is also large relative
to $n$. Instead, the projection strategy proposed in this paper can
be used to estimate $\beta^{*}$ in the semilinear model (\ref{eq:30})
when $p\geq n$ and $d\geq n$. One way to reduce the curse of dimensionality
in the joint multivariate nonparametric component $\mathbb{E}(z_{ij}\,\vert\, w_{i})$
(recall $z_{j}=x_{j}$ for $j=1,...,p$ and $z_{0}=y$) is to consider
the class of additive models of the form (Hastie and Tibshirani, 1999):
\[
\mathbb{E}(z_{ij}\,\vert\, w_{i}):=f_{j}(w_{i})=\sum_{l=1}^{d}f_{jl}(w_{il})
\]
where $f_{jl}(\cdot)\in\mathcal{F}_{jl}$ for $l=1,...,d$ and $j=0,...,p$. 

Let us consider the simplest case of $f_{jl}(w_{il})=w_{il}\theta_{l}^{*}$
where $\theta_{l}^{*}$ is a scalar and $\theta^{*}\in\mathcal{B}_{q_{1}}^{d}(R_{q_{1}})$
and $\beta^{*}\in\mathcal{B}_{q_{2}}^{p}(R_{q_{2}})$ for $q_{1},\, q_{2}\in[0,\,1]$.
Theorem 4.1 implies that the $l_{2}-$error of the two-stage estimator
is bounded above by 
\[
O\left(\left[R_{q_{2}}^{\frac{1}{2}}\left(\sqrt{\frac{\log p}{n}}\right)^{1-\frac{q_{2}}{2}}\right]\vee\left[R_{q_{2}}^{\frac{1}{2}}\left(\left|\beta^{*}\right|_{1}R_{q_{1}}\left(\frac{\log d}{n}\right)^{1-\frac{q_{1}}{2}}\right)^{1-\frac{q_{2}}{2}}\right]\right).
\]
For a more general structure on $f_{j}$, suppose $J(f_{j}):=\left\{ l\,:\, f_{jl}\neq0\right\} $
and $k_{1j}=\left|J(f_{j})\right|$, the cardinality of $J(f_{j})$,
can increase to infinity with $n$ but slowly compared to $n$ (i.e.,
$f_{j}$ is exactly sparse) and 
\[
\mathbb{E}(z_{ij}\,\vert\, w_{il})=\sum_{k=1}^{\infty}\vartheta_{jlk}\phi_{jlk}(w_{il})
\]
where $B_{jl}=\left(\phi_{jlk}\right)_{k=1}^{\infty}$ is an orthonormal
basis for $\mathcal{F}_{jl}$. For a truncation parameter $M$, also
define 
\[
\mathbb{E}^{M}(z_{ij}\,\vert\, w_{il})=\sum_{k=1}^{M}\vartheta_{jlk}\phi_{jlk}(w_{il}).
\]
Let $\Psi_{jl}$ denote the $n\times M$ matrix with $\Psi_{jl}(i,\, k)=\phi_{jlk}(w_{il})$.
For the first-stage estimation in our two-stage procedure, consider
the following minimization problem: 
\begin{equation}
\min_{\vartheta_{jl}\in\mathbb{R}^{M}}:\;\frac{1}{2n}\left|z_{j}-\sum_{l=1}^{d}\Psi_{jl}\vartheta_{jl}\right|_{2}^{2}+\lambda_{n}\sum_{l=1}^{d}\sqrt{\frac{1}{n}\vartheta_{jl}^{T}\Psi_{jl}^{T}\Psi_{jl}\vartheta_{jl}}\label{eq:31}
\end{equation}
for some regularization parameter $\lambda_{n}>0$. Program (\ref{eq:31})
is the sample version (with truncation) of the following: 
\[
\min_{f_{jl}\in\mathcal{F}_{jl}}:\;\frac{1}{2}\mathbb{E}\left(z_{ij}-\sum_{l=1}^{d}f_{jl}(w_{il})\right)^{2}+\lambda_{n}\sum_{l=1}^{d}\sqrt{\mathbb{E}\left(f_{jl}^{2}(w_{il})\right)}.
\]
The optimization program (\ref{eq:31}) is considered in Ravikumar,
Lafferty, and Wasserman (2009) and can be viewed as a functional version
of the grouped Lasso (Yuan and Lin, 2006). It can be solved with a
coordinate descent algorithm proposed by Ravikumar, et. al (2009).
Theoretical properties of the two-stage estimator with the first-stage
estimation based on (\ref{eq:31}) are being analyzed in a chapter
of the PhD thesis by Zhu (2015).

\section{Monte-Carlo simulation}

In this section, simulations are conducted to gain preliminary understanding
of the small-sample performance of the non-pivotal multi-stage estimator
$\hat{\beta}_{HSEL}$; ongoing work involves implementation of the
pivotal procedure described in Section 4.2. We consider model (\ref{eq:1})
where $w\in\mathbb{R}^{n\times d}$ is a matrix consisted of independent
uniform zero-mean random variables on $[-2,\,2]$ with variance $\sigma_{w}\approx1.33$
and $x$ takes on the first $p$ columns of $w$. The \textit{i.i.d.}
errors $\epsilon_{1i}\sim\mathcal{N}\left(0,\,1\right)$ for $i=1,...,n$
where $n$ denotes the number of observations generated for the selection
equation. We consider two scenarios where $n=88$ and $n=200$. Given
the setup here, on average $44$ (when $n=88$) and $100$ (when $n=200$)
observations, respectively, will be used for estimating the main equation.
Conditional on the observations $i$s with $y_{1i}=1$, the \textit{i.i.d.}
errors $(\epsilon_{1i},\,\epsilon_{2i})$ have the following joint
normal distribution 
\[
(\epsilon_{1i},\,\epsilon_{2i})\thicksim\mathcal{N}\left(\left(\begin{array}{c}
0\\
0
\end{array}\right),\;\left(\begin{array}{cc}
1 & \rho\sigma_{2}\\
\rho\sigma_{2} & \sigma_{2}
\end{array}\right)\right),
\]
where $\rho\in\{0,\,0.9\}$ and $\sigma_{2}\in\{0.3,\,1,\,2\}$. We
set $d=90$, $p=45$, $k_{1}=4$, and $k_{2}=2$. When $n=88$, this
setup of dimensionality represents a selection model where the number
of regressors in the selection equation and the main equation, respectively,
exceeds the number of observations used to estimate the corresponding
equation, while the number of relevant regressors (ones with nonzero
coefficients) is small relative to the sample size. We set $\theta_{j}^{*}=0.5$
for $j=1$, $2$, $3$, $46$ and the rest of components in $\theta^{*}$
take on values of $0$; set $\beta_{1}^{*}=\beta_{45}^{*}=1$ and
the rest of components in $\beta^{*}$ take on values of $0$. This
set up ensures that there is at least one component $w_{ij}$ with
$\theta_{j}^{*}$ in the support set of $\theta^{*}$ such that $w_{ij}$
is excluded from $x_{i}$. 

We consider four sets of experiments. The first experiment (Experiment
1) concerns the multi-stage estimator $\hat{\beta}_{HSEL}$. As a
benchmark for Experiment 1, Experiment 2 applies a one-step Lasso
procedure (without correcting selection bias) to the same main equation.
Experiments 3 and 4 are benchmarks concerning classical low-dimensional
settings. Experiment 3 applies the Heckman's 2-step procedure to model
(\ref{eq:1}) where the selection equation and the main equation are
in the low-dimensional setting and the supports of the true parameters
in both equations are known \textit{a priori}; Experiment 4 applies
the OLS to the same low-dimensional model as Experiment 3. We simulate
100 sets of data following the process described above. For each set
$t=1,...,100$, we compute the estimates $\hat{\beta}^{t}$ of the
main-equation parameters $\beta^{*}$, $l_{2}-$errors of these estimates,
$|\hat{\beta}^{t}-\beta^{*}|_{2}$, and selection percentages of $\hat{\beta}^{t}$
(computed by the number of the elements in $\hat{\beta}^{t}$ sharing
the same sign as their corresponding elements in $\beta^{*}$, divided
by the total number of elements in $\beta^{*}$). Results reported
in this section include: \\

(\textbf{a}) the mean of the relevant estimates $\frac{1}{100}\sum_{t=1}^{100}\hat{\beta}_{1}^{t}$; 

(\textbf{b}) the mean of the relevant estimates $\frac{1}{100}\sum_{t=1}^{100}\hat{\beta}_{45}^{t}$; 

(\textbf{c}) the mean of the averaged irrelevant estimates $\frac{1}{43}\sum_{j\neq1,\,45}\frac{1}{100}\sum_{t=1}^{100}\hat{\beta}_{j}^{t}$; 

(\textbf{d}) the mean of the $l_{2}-$errors of the estimates $\hat{\beta}^{t}$
computed as $\frac{1}{100}\sum_{t=1}^{100}|\hat{\beta}^{t}-\beta^{*}|_{2}$; 

(\textbf{e}) the mean of the selection percentages (computed in a
similar fashion as the mean of the $l_{2}-$errors of the estimates); 

(\textbf{f}) the mean of the squared $l_{2}-$errors (i.e., the \textit{sample}
\textit{mean squared error}, SMSE, computed as $\frac{1}{100}\sum_{t=1}^{100}|\hat{\beta}^{t}-\beta^{*}|_{2}^{2}$); 

(\textbf{g}) the sample squared bias $\sum_{j=1}^{45}(\bar{\hat{\beta}}_{j}-\beta_{j}^{*})^{2}$
(where $\bar{\hat{\beta}}_{j}=\frac{1}{100}\sum_{t=1}^{100}\hat{\beta}_{j}^{t}$
for $j=1,...,45$). \\
\\
The results in this section regarding Experiment 1 are based on the
choices of the regularization parameter $\lambda_{n,1}=0.5\sqrt{\frac{\log d}{n}}$
for the first-stage estimation problem (\ref{eq:13}) and the regularization
parameter $\lambda_{n,3}=0.2k_{2}\sqrt{\frac{k_{1}\log d}{n_{s}}}$
for the third-stage estimation problem (\ref{eq:8}), where $n_{s}$
denotes the number of observations with $y_{1i}=1$. The scalings
of $\lambda_{n,1}$ and $\lambda_{n,3}$ are chosen according to Proposition
4.2 and Corollary 4.4, respectively. The choice of $0.2k_{2}\sqrt{\frac{k_{1}\log d}{n_{s}}}$
is also used in Experiment 2 for comparing the performance of the
proposed procedure and the Lasso without corrective measures. Note
that for $\sigma_{2}=1$ and $\sigma_{x}=\sigma_{w}\approx1.33$,
$0.2k_{2}\sqrt{\frac{k_{1}\log d}{n_{s}}}$ is slightly greater than
$2\sigma_{2}\cdot\sigma_{x}\sqrt{\frac{\log p}{n_{s}}}$, the smallest
value required for the Lasso estimation of the standard sparse high-dimensional
linear models (e.g., Bickel, et. al, 2009). The second-stage estimation
in Experiment 1 is based on solving (\ref{eq:22}) with $L=1$. Ongoing
work involves implementing the cross validation procedure described
in Section 4.1.2 to determine $L$ when $\epsilon_{2i}$ has a non-normal
distribution. 

\begin{spacing}{1.5}
From Table 5.1, we see that the direct Lasso estimator without correcting
selection bias outperforms the multi-stage estimator $\hat{\beta}_{HSEL}$
when $\rho=0$, and \textit{vice versa} when $\rho=0.9$. For the
design considered here, in the presence of substantial selection activity
($\rho=0.9$), the mean of the $l_{2}-$errors (\textbf{row d}) and
the sample squared bias of the estimates (\textbf{row g}) by the direct
Lasso procedure without corrective measures are exacerbated in the
high-dimensional setting and this exacerbation mainly comes from the
poorer estimates of the relevant regressors as the mean of the averaged
irrelevant estimates varies little from the case $\rho=0$ to the
case $\rho=0.9$. Other simulation results (not included here due
to space limit) show that when $\sigma_{2}$ is increased (decreased)
from $1$ to $2$ (respectively, from $1$ to $0.3$), $\hat{\beta}_{HSEL}$
performs worse (respectively, better) relative to the case $\sigma_{2}=1$,
and similar patterns are observed when $w_{ij}$s are drawn from independent
uniform zero-mean random variables on $[-1,\,1]$ (respectively, on
$[-4,\,4]$). Also, as $n$ increases from $88$ to $200$, $\hat{\beta}_{HSEL}$
performs substantially better. These findings are intuitive and expected.
It is worth noting that for the design considered here, in terms of
the mean of the selection percentages (\textbf{row e}), the direct
Lasso procedure without corrective measures is comparable to $\hat{\beta}_{HSEL}$
even in the case $\rho=0.9$. Ongoing work is exploring situations
where variable selection by $\hat{\beta}_{HSEL}$ substantially outperforms
variable selection by the direct Lasso procedure.\\
\\
{\scriptsize }%
\begin{tabular}{llllll|l|c|c|c|c|c|c|c|c|}
 &  &  &  &  & \multicolumn{1}{l}{} & \multicolumn{1}{l}{} & \multicolumn{8}{c}{{\small Table 5.1: Monte-Carlo simulation results for $n=88$ }}\tabularnewline
\cline{7-15} 
 &  &  &  &  &  & \multirow{2}{*}{\begin{sideways}
\end{sideways}} & \multicolumn{4}{c}{{\scriptsize $\rho=0$}} & \multicolumn{4}{c|}{{\scriptsize $\rho=0.9$}}\tabularnewline
\cline{8-15} 
 &  &  &  &  &  &  & \textbf{\scriptsize Exp 1} & \textbf{\scriptsize Exp 2} & \textbf{\scriptsize Exp 3} & \textbf{\scriptsize Exp 4} & \textbf{\scriptsize Exp 1} & \textbf{\scriptsize Exp 2} & \textbf{\scriptsize Exp 3} & \textbf{\scriptsize Exp 4}\tabularnewline
\cline{7-15} 
 &  &  &  &  &  & \textbf{\scriptsize a} & {\scriptsize 0.703} & {\scriptsize 0.730} & {\scriptsize 1.005} & {\scriptsize 1.010} & {\scriptsize 0.627} & {\scriptsize 0.605} & {\scriptsize 1.007 } & {\scriptsize 1.002}\tabularnewline
\cline{7-15} 
 &  &  &  &  &  & \textbf{\scriptsize b} & {\scriptsize 0.742} & {\scriptsize 0.736} & {\scriptsize 0.996} & {\scriptsize 0.994} & {\scriptsize 0.762} & {\scriptsize 0.757} & {\scriptsize 1.006} & {\scriptsize 1.020}\tabularnewline
\cline{7-15} 
 &  &  &  &  &  & \textbf{\scriptsize c} & {\scriptsize -0.001} & {\scriptsize -0.001 } & {\scriptsize NA} & {\scriptsize NA} & {\scriptsize -0.001} & {\scriptsize -0.001} & {\scriptsize NA} & {\scriptsize NA}\tabularnewline
\cline{7-15} 
 &  &  &  &  &  & \textbf{\scriptsize d} & {\scriptsize 0.446} & {\scriptsize 0.430} & {\scriptsize 0.159} & {\scriptsize 0.161} & {\scriptsize 0.474} & {\scriptsize 0.495} & {\scriptsize 0.145} & {\scriptsize 0.164}\tabularnewline
\cline{7-15} 
 &  &  &  &  &  & \textbf{\scriptsize e} & {\scriptsize 0.969} & {\scriptsize 0.965} & {\scriptsize NA} & {\scriptsize NA} & {\scriptsize 0.981} & {\scriptsize 0.979} & {\scriptsize NA} & {\scriptsize NA}\tabularnewline
\cline{7-15} 
 &  &  &  &  &  & \textbf{\scriptsize f} & {\scriptsize 0.227} & {\scriptsize 0.209} & {\scriptsize 0.032} & {\scriptsize 0.033} & {\scriptsize 0.249} & {\scriptsize 0.262} & {\scriptsize 0.025} & {\scriptsize 0.034}\tabularnewline
\cline{7-15} 
 &  &  &  &  &  & \textbf{\scriptsize g} & {\scriptsize 0.155} & {\scriptsize 0.143} & {\scriptsize $4\times10^{-5}$} & {\scriptsize $9\times10^{-5}$} & {\scriptsize 0.197} & {\scriptsize 0.217} & {\scriptsize $9\times10^{-5}$} & {\scriptsize $4\times10^{-4}$}\tabularnewline
\cline{7-15} 
\end{tabular}{\scriptsize \par}
\end{spacing}

\section{An empirical application to the retail gasoline market }

Having established the theoretical properties of the 3-step estimators,
we now apply one of these estimators to an empirical example of price-discrimination
in the retail gasoline market. When consumers have different valuations
for a good, a firm can increase profits by developing a pricing scheme
that distinguish consumers with different valuations. In most cases,
a firm knows the distribution of consumer valuations in the market
but not the exact valuation of any specific consumer prior to the
sale. In these cases, a firm can offer a menu of different prices,
appropriately bundled with other aspects of the product (such as product
quality), and force consumers to choose bundles consistent with their
preferences. 

When differences in costs incurred to produce various bundles in the
menu are small compared to the differences in prices, this menu-based
offering is a price discrimination mechanism. Retail gasoline markets
present a good context to study price discrimination since different
gasoline stations in a market typically face similar costs of procuring
gasoline. Therefore, any price differences across gasoline stations
are likely due to reasons unrelated to the cost of procuring gasoline.
Gasoline retailers can choose to be either a two-product station offering
both self-service and full-service gasoline or a single-product station
offering only full-service or self-service gasoline. A two-product
station, by charging different prices for full- and self-service gasoline,
induces consumers with different valuations to choose the products
consistent with their preferences, namely, a two-product station engages
in price-discrimination. A single-product station, on the other hand,
is unable to price discriminate. 

Shepard (1991) estimates pricing decisions of gasoline stations without
endogenizing their decisions to price discriminate, i.e., their choice
to be single versus multi-product. Iyer and Seetharaman (2003) explicitly
examines a firm\textquoteright{}s incentive to price-discriminate.
In doing so, they highlight the importance of accounting for self
selectivity considerations in empirical analysis of price discrimination
based on market data. Specifically, Iyer and Seetharaman employ a
binary probit framework to model a gasoline station\textquoteright{}s
decision to be single-product or multi-product as a function of market
and station characteristics, and then model the prices chosen by the
gasoline station for its product(s) by estimating linear regressions
with Heckman's self-selectivity correction conditional on the station\textquoteright{}s
decision to offer a single- or multi-product. They show that incorrect
inferences about the incentive to price discriminate and about the
differences in the prices charged between single-product and multi-product
stations would result if the endogeneity in the choice of the station-type
were ignored in the estimation. Their empirical analysis also shows
that a larger income spread in the market implies a greater likelihood
of the gasoline station being multi-product. However, Iyer and Seetharaman
(2003) did not account for interactions between the gas stations in
their empirical analysis. Studies show that pricing decisions of retail
gasoline stations may depend on the degree of competitive intensity
in the {}``market'' (e.g., Slade, 1992). In the empirical literature
on competitive gasoline markets, there have been various ways of defining
a {}``market'' (see, e.g., Slade, 1986; Pinkse, Slade and Brett,
2002; Iyer and Seetharaman, 2008). For example, Iyer and Seetharaman
(2008) defines mutually exclusive census tracts as local markets,
and treat each market as the unit of observation in their empirical
analysis. In previous research, markets have been defined based on
stations that fall within a circle of half a mile or one mile radius. 

One common feature of the previous definitions of competitive markets
is that they are subjective heuristics. It would be ideal if one can
control for the interactions between different stations without requiring
\textit{a priori }knowledge of the structure of the competitive market.
Recent work including Manresa (2014) and Bonaldi, Hortacsu, Kastl
(2014) develop econometric models to recover the underlying networks
in different applications. Both papers hinge on the availability of
panel data for each observation in the cross section. In particular,
Manresa considers settings where outcomes depend on an agent's own
characteristics and on the characteristics of other agents in the
data. She applies a Lasso type estimator to identify individuals generating
spillovers and their strength using panel data on outcomes and characteristics.
Bonaldi, et. al proposes a new measure of systemic risk based on estimating
spillovers between funding costs of individual banks with a Lasso
type procedure, which is applied to the panel of each individual bank
to recover the financial network. However, for the empirical application
considered in this paper, panel data of each gas station is not available
and as a consequence, the econometric model by either Manresa or Bonaldi
is not suitable. Instead, we use geographic information and spatial
data to create a set of measures that are high-dimensional to control
for the interactions between the gas stations and employ one of our
proposed estimators to identify the competitive market structure.
The following subsection describes the data followed by the empirical
model.

\subsection*{Data and the empirical model}

This paper uses the data set from Iyer and Seetharaman (2003). It
was collected during July 1998 from a cross-section of 249 gasoline
stations in the Greater Saint Louis metropolitan area. Among the 249
stations, 65 are multi-product stations and 172 are single-product
self-service stations. In addition, there were 12 single-product full-service
stations. In the United States, the low incidence of full-service
single product stations is typical and in certain regions full service
is required by law. As in Iyer and Seetharaman (2003), we exclude
them from the empirical analysis. The survey data include the prices
of three grades - 87, 89 and 93 octane levels - of gasoline, along
with station-specific characteristics, i.e., number of gasoline pumps,
special advertising for cigarettes and soda, presence of convenience
store, pay-at-pump facility, car wash, service station, and the number
of stations with prices that are visible to a given station. This
data set also contains demographic information including income, population
density, age distribution, home value, and education levels. This
information comes from 1990 U.S. census data, which contain demographic
information at the level of each census tract.

The data also records addresses of each station, from which {}``Bing
Maps REST Services'' is used to obtain geographic information including
longitude and latitude, travel distance in driving mode between any
pair of stations, etc. This information can be used to create variables
for partially controlling for the interactions between stations. In
particular, given any station, we can count the number of stations
and/or stations under one of the three national brands (namely, Amoco,
Shell and Mobil), that fall within 1km, 1km and 2km, and so on, from
this station. Each of the numbers is then divided by the area (in
$\textrm{km}^{2}$) of the corresponding layer. The number of stations
with prices that are visible to a given station is another useful
measure of interaction between stations and this information is available
in the data. Using the number of competitors as a measure of interaction
between firms has been seen in previous literature (e.g., Bresnahan
and Reiss 1991; Iyer and Seetharaman, 2008). The novelty introduced
by this paper lies in the data-driven nature of the approach: rather
than assuming \textit{a priori }knowledge of the structure of the
competitive network, it relies on the data to determine the geographic
pattern of interaction between stations. If panel data on prices and
time-varying instrumental variables for prices are available, we can
include prices of other stations in the main equation. Some variants
of the 3-step estimators in this paper combined with the high-dimensional
IV estimator in Gautier and Tsybakov (2011) or the high-dimensional
2SLS estimator in Zhu (2013) may be considered as an alternative to
identify the sets of competitive markets engaged in pricing. However,
in the retail gasoline market, it may be difficult to obtain valid
time-varying instrumental variables for prices.

As in Iyer and Seetharaman (2003), we use a binary probit model for
the selection of service types where $y_{1i}=1$ in (\ref{eq:1})
indicates that station $i$ offers multi-service and $y_{1i}=0$ indicates
that station $i$ offers single-self-service. The same set of explanatory
variables included in the binary probit model of Iyer and Seetharaman
is used here: average income (\textit{AVG}), income spread (\textit{SPREAD}),
brand (\textit{BRAND}), pay-at-pump facility (\textit{PAP}), presence
of convenience store (\textit{CONV}), car wash (\textit{WASH}), and
service station (\textit{SERV}). This paper differs from Iyer and
Seetharaman mainly in terms of the specifications of the linear pricing
model (the main equation): First, while Iyer and Seetharaman assume
the selection bias takes on the functional form of the Inverse Mills
Ratio, we assume the selection bias function to obey the more general
nonparametric single index restriction in (\ref{eq:1}); second, we
add a set of measures that are high-dimensional to partially control
for the competition effects from other stations. In particular, the
following explanatory variables are included in the pricing model:
\textit{AVG}, \textit{BRAND}, special advertising for cigarettes and
soda (\textit{ADSCC}), the number of stations with visible prices
(\textit{VISP}), the total number of stations and the number of stations
under one of the three national brands within 1km (\textit{TOT\_1}
and \textit{BRND\_1}), 1km and 2km (\textit{TOT\_2} and \textit{BRND\_2}),
$\cdots$, 34km and 35km (\textit{TOT\_35} and \textit{BRND\_35})
from a given station. In summary, for the pricing equation, we have
$n=172$, $p=74$ for stations that serve single-self-service grade-87
gasoline, $n=168$, $p=74$ for stations that serve single-self-service
grade-93 gasoline, $n=65$, $p=74$ for stations that serve multi-service
grade-87 gasoline, and $n=65$, $p=74$ for stations that serve multi-service
grade-93 gasoline. While Iyer and Seetharaman include only average
income and brand in their pricing model, they suggest that special
advertising for cigarettes and soda might be correlated with the retail
gasoline prices and hence we include this information in our pricing
model. The last group of variables are measures added to partially
control for the competition effects from other stations. Iyer and
Seetharaman found the indicators of the presence of pay-at-pump facilities
and service stations statistically significant and therefore, the
exclusion restriction required by the selection model considered in
this paper is likely to be satisfied given the setup. Moreover, they
justify the exclusion restriction by arguing that a station\textquoteright{}s
decision pertaining to the configuration of its station characteristics
- pay-at-pump, convenience store, car-wash, and service station -
involves costly investments that the station owner has made along
with the station-type decision while setting up the retail facility.
In contrast, the pricing decisions may vary on a daily basis. 

The following summarizes the estimation procedure and empirical findings.
We briefly discuss the results pertaining to the effects of the service-type
decisions and focus mainly on the empirical findings from the pricing
regression because the main difference between Iyer and Seetharaman
and the empirical analysis in this paper lies in the latter.

\subsection*{Estimation and empirical findings}

A standard maximum likelihood procedure for estimating low-dimensional
binary probit models is used to obtain estimates of the selection
equation of service-type decisions. The estimation is performed for
grade-87 stations and grade-93 stations, respectively, and the results
are reported in Table 6.1. \\
\\
{\scriptsize }%
\begin{tabular}{ll|c|c|c|c|c|c|c|c|c|}
\multicolumn{2}{l}{} & \multicolumn{9}{c}{{\small Table 6.1: Results of the binary probit model}}\tabularnewline
\cline{3-11} 
 &  &  & \textit{Intercept} & \textit{AVG} & \textit{SPREAD} & \textit{BRAND} & \textit{PAP} & \textit{CONV} & \textit{WASH} & \textit{SERV}\tabularnewline
\cline{3-11} 
 &  & \multirow{2}{*}{\textbf{\scriptsize Grade-87}} & {\scriptsize $-1.525^{***}$} & {\scriptsize $0.008^{*}$} & {\scriptsize $-2.591^{**}$} & {\scriptsize $0.962^{**}$} & {\scriptsize $-0.851^{**}$} & {\scriptsize $-0.341$} & {\scriptsize $-0.012$} & {\scriptsize $2.341^{***}$}\tabularnewline
 &  &  & {\scriptsize (0.453)} & {\scriptsize (0.005)} & {\scriptsize (1.227)} & {\scriptsize (0.411)} & {\scriptsize (0.414)} & {\scriptsize (0.321)} & {\scriptsize (0.346)} & {\scriptsize (0.292)}\tabularnewline
\cline{4-11} 
 &  & \multirow{2}{*}{\textbf{\scriptsize Grade-93}} & {\scriptsize $-1.522^{***}$} & {\scriptsize $0.008^{*}$} & {\scriptsize $-2.583^{**}$} & {\scriptsize $0.958^{**}$} & {\scriptsize $-0.848^{**}$} & {\scriptsize $-0.338$} & {\scriptsize $-0.004$} & {\scriptsize $2.335^{***}$}\tabularnewline
 &  &  & {\scriptsize (0.452)} & {\scriptsize (0.005)} & {\scriptsize (1.229)} & {\scriptsize (0.410)} & {\scriptsize (0.413)} & {\scriptsize (0.321)} & {\scriptsize (0.347)} & {\scriptsize (0.292)}\tabularnewline
\cline{3-11} 
\end{tabular}\\
\\
\\
Because individual-level income is not available for this data set,
it is not possible to compute the sample standard deviation in income
for each tract. Instead, we use two measures to approximate income
spread: one is the absolute difference between the percentages of
median-level income group and the low-level income group for each
tract; the other is the absolute difference between the percentages
of median-level income group and the high-level income group for each
tract. A smaller value in the first (second) absolute difference indicates
a more evenly distributed population in the low-income (respectively,
high-income) group and the median income group. It turns out that
the second measure is not statistically significant and hence we drop
this measure from the probit model. As a consequence, the negative
sign of the estimate for \textit{SPREAD} suggests that more heterogeneous
income levels below the $50^{th}-$percentile in the market implies
a greater likelihood of the station being multi-product.

For the linear pricing regression model conditional on the service
type, the non-pivotal estimator $\hat{\beta}_{HSEL}$ based on (\ref{eq:8})
is used to select the variables with non-zero coefficients and then
$\tilde{\beta}:=\left(\hat{v}_{\hat{J}}^{T}\hat{v}_{\hat{J}}\right)^{-1}\left(\hat{v}_{\hat{J}}^{T}\hat{v}_{0}\right)$
is computed with $\hat{J}:=J(\hat{\beta}_{HSEL})$ (this is the Post-Lasso
procedure discussed in Section 4.1.3). For the second-stage estimation,
program (\ref{eq:22}) is solved where the Lipschitz constant $L$
is determined by the cross-validation procedure described in Section
4.1.2 and the choice of $L=1$ turns out to be robust. For the third-stage
estimation, given the setup of our empirical model, the first-stage
related estimation error is likely to be dominated by the third-stage
related error in the choice of $\lambda_{n,3}$ from Corollary 4.4
and hence we choose $\lambda_{n,3}$ based on the third-stage related
error. Program (\ref{eq:8}) is first solved with the choice of $\lambda_{n,3}(t)=2.001\cdot\hat{\sigma}_{v}\hat{\sigma}_{\eta}^{t}\sqrt{\frac{\log p}{n_{s}}}$
for $t=0$ (initialization), where $\hat{\sigma}_{v}:=\max_{j=1,...,p}\sqrt{\frac{1}{n}\sum_{i=1}^{n}\hat{v}_{ij}^{2}}$,
$\hat{\sigma}_{\eta}^{0}=1$, and $n_{s}$ denotes the number of observations
used for the pricing regression. Let $\hat{\beta}_{HSEL}^{t}$ denote
the resulting estimate based on $\lambda_{n,3}(t)$ and $\hat{\sigma}_{\eta}^{t+1}$
denote the updated sample standard deviation of the fitted residuals
$\hat{\eta}_{i}^{t}:=\hat{v}_{i0}-\hat{v}_{i}\hat{\beta}_{HSEL}^{t}$
for $i=1,...,n_{s}$. Program (\ref{eq:8}) is then solved with the
updated $\lambda_{n,3}(t+1)=2.001\cdot\hat{\sigma}_{v}\hat{\sigma}_{\eta}^{t+1}\sqrt{\frac{\log p}{n_{s}}}$.
Repeat this process until a pre-specified tolerance level on $\left|\hat{\sigma}_{\eta}^{t+1}-\hat{\sigma}_{\eta}^{t}\right|$
is reached. The result shows that the choice of $\hat{\sigma}_{\eta}^{t}\approx0.04$
is robust. After experimenting with a range of values around the final
choice of $\lambda_{n,3}$ determined according to the described procedure,
the set of variables selected by $\hat{\beta}_{HSEL}$ for the range
of $\lambda_{n,3}$ and the post Lasso estimates $\tilde{\beta}$
based on these selected variables are reported in Table 6.2 for the
following groups:\\

SSL: single-self-service grade-87 gasoline; 

SSH: single-self-service grade-93 gasoline; 

MSL: multi-self-service grade-87 gasoline; 

MSH: multi-self-service grade-93 gasoline. \\
\\
The variables with blanks in Table 6.2 correspond to those that are
not selected by $\hat{\beta}_{HSEL}$ in a particular group. The numerical
values within parentheses are bootstrapped standard errors for $\tilde{\beta}$.
Estimates with three asterisks, two asterisks, and a single asterisk
are statistically significant at level $\alpha=0.01$, $\alpha=0.05$,
and $\alpha=0.1$, respectively. Note that our second-stage estimation
and third-stage estimation use the demeaned explanatory variables
and demeaned prices, so the intercept term is excluded from the pricing
regression model.\\
\\
{\scriptsize }%
\begin{tabular}{llllllllll|l|l|c|c|c|}
 &  &  &  &  &  &  &  &  & \multicolumn{1}{l}{} & \multicolumn{5}{c}{{\small Table 6.2: Results of pricing regression}}\tabularnewline
\cline{11-15} 
 &  &  &  &  &  &  &  &  &  &  & \textit{AVG} & \textit{BRAND} & \textit{TOT\_2} & \textit{TOT\_4}\tabularnewline
\cline{11-15} 
 &  &  &  &  &  &  &  &  &  & \multirow{2}{*}{SSL} &  & {\scriptsize $0.030^{***}$} &  & {\scriptsize $-0.052^{***}$}\tabularnewline
 &  &  &  &  &  &  &  &  &  &  &  & {\scriptsize (0.005)} &  & {\scriptsize (0.016)}\tabularnewline
\cline{12-15} 
 &  &  &  &  &  &  &  &  &  & \multirow{2}{*}{SSH} &  & {\scriptsize $0.049^{***}$} & {\scriptsize $-0.030^{**}$} & {\scriptsize $-0.030^{*}$}\tabularnewline
 &  &  &  &  &  &  &  &  &  &  &  & {\scriptsize (0.007)} & {\scriptsize (0.015)} & {\scriptsize (0.021)}\tabularnewline
\cline{12-15} 
 &  &  &  &  &  &  &  &  &  & \multirow{2}{*}{MSL} & {\scriptsize $2\times10^{-4**}$} &  &  & {\scriptsize $-0.070^{***}$}\tabularnewline
 &  &  &  &  &  &  &  &  &  &  & {\scriptsize $(1\times10^{-4})$} &  &  & {\scriptsize (0.029)}\tabularnewline
\cline{12-15} 
 &  &  &  &  &  &  &  &  &  & \multirow{2}{*}{MSH} &  &  &  & {\scriptsize $-0.096^{***}$}\tabularnewline
 &  &  &  &  &  &  &  &  &  &  &  &  &  & {\scriptsize (0.032)}\tabularnewline
\cline{11-15} 
\end{tabular}\\
\\
\\
Regarding the results of the pricing regression in Table 6.2, the
estimate of\textit{ BRAND }has a positive sign in SSL and SSH. Moreover,
\textit{AVG} (average income) has a positive effect on the pricing
decisions in MSL. In Iyer and Seetharaman\textit{ }(2003) which estimated
the low-dimensional linear regression counterpart by pooling observations
of the single-self-service and multi-self-service each with Heckman's
selectivity correction (for grade-87 and grade-93, respectively),\textit{
BRAND} and \textit{AVG }are the only two variables included in their
pricing model and found to be statistically significant. Based on
our empirical results which remove selection bias and partially control
for potential interactions between the stations simultaneously, we
see that \textit{TOT\_4} is selected by $\hat{\beta}_{HSEL}$ in \textit{all}
groups and statistically significant at level $0.01$ in SSL, MSL,
and MSH and at level $0.1$ in SSH. \textit{TOT\_2} is selected by
$\hat{\beta}_{HSEL}$ in SSH and statistically significant at level
$0.05$. The negative sign of the estimate for \textit{TOT\_4} in
all groups (\textit{TOT\_2} in SSH) suggests that the total number
of stations within 3km-4km (respectively, 1km-2km) of a given station
has a negative effect on its price. On the other hand, the variable
\textit{VISP} (the number of stations with visible prices) is not
selected by $\hat{\beta}_{HSEL}$, which is less intuitive. However,
it is possible that in the presence of several clusters of gas stations,
there is less competition from adjacent clusters relative to ones
that are somewhat further apart. For example, Iyer and Seetharaman\textit{
}(2008) analyzed a similar but richer data set on prices and station
characteristics gathered across stations in the Saint Louis metropolitan
area and found that closely located retailers who face sufficient
heterogeneity in preferences across consumers in a local market may
differentiate on product design and pricing strategies (also see,
e.g., Png and Reitman, 1994); in contrast, retailers that are farther
apart from each other may adopt similar product design and pricing
strategies if the market is relatively homogeneous (also see, e.g.,
Slade 1992). Another explanation for the finding where more competition
comes from somewhat intermediate retailers instead of closest ones
is that consumers of retail gasoline may travel from their suburban
homes located in a neighborhood of one cluster to their work places
or central shopping areas located in another cluster that may be somewhat
further away; the further located clusters of stations may be linked
by routes that are more convenient for commuting (these routes may
be more direct or less congested, etc.). This explanation may suggest
that retailers consider commuting behavior of their customers when
setting the retail price. Investigating this factor requires more
substantial empirical analysis and a data set that includes more detailed
information on the demographics and business environment, which will
be pursued in future research. Nevertheless, the main finding on \textit{TOT\_2}
and \textit{TOT\_4 }suggests that in modeling the pricing decisions
of retail gas stations, not only it is useful to account for the self
selectivity of service-type but also to take into considerations of
potential interactions between stations; in particular, competition
effects from retailers that are not in the same local market (e.g.,
the same census tract or neighborhood within a circle of half a mile
or one mile radius, etc.) should not be overlooked.

\section{Conclusion}

This paper provides estimation tools together with their theoretical
guarantees for the semiparametric sample selection model in high-dimensional
settings under a weak nonparametric restriction on the form of the
selection correction. In particular, the number of regressors in the
main equation, $p$, and the number of regressors in the selection
equation, $d$, can grow with and exceed the sample size $n$. The
main theoretical results of this paper are finite-sample bounds from
which sufficient scaling conditions on the sample size for estimation
consistency and variable-selection consistency (i.e., the multi-stage
high-dimensional estimation procedure correctly selects the non-zero
coefficients in the main equation with high probability) are established.\textcolor{black}{{}
}Statistical efficiency of the proposed estimators is studied via
lower bounds on minimax risks. Inference procedures for the coefficients
of the main equation, one based on a pivotal Dantzig selector to construct
non-asymptotic confidence sets and one based on a post-selection strategy
(when perfect or near-perfect selection of the high-dimensional coefficients
is achieved), are discussed. 

Small-sample performance of one of the proposed procedures is evaluated
by Monte-Carlo simulations and illustrated with an empirical application
to the retail gasoline market in the Greater Saint Louis area. The
preliminary simulation results show that the {}``bias'' from not
performing the selection correction is exacerbated in high-dimensional
settings. For the empirical application, this paper models a firm\textquoteright{}s
choice of either a single-product or multi-product service as a function
of market and station characteristics and then models the station\textquoteright{}s
pricing decision, conditional on the choice of the station type. Using
geographic information and spatial data, a set of variables that are
high-dimensional is introduced to control for interactions between
the gas stations. The empirical finding suggests that competition
effects from retailers that are not in the same local market should
not be overlooked.

\section*{Appendix I: Main proofs }

\section*{Appendix II: Technical lemmas and the proofs }

\textbf{Appendix I} and \textbf{Appendix II} can be found in \textbf{Section
A }and\textbf{ Section B} of the online supplementary material: Proofs
to {}``High-Dimensional Semiparametric Selection Models: Estimation
Theory with an Application to the Retail Gasoline Market''. \\
(\href{https://sites.google.com/site/yingzhu1215/home/JobMar_Proofs.pdf}{https://sites.google.com/site/yingzhu1215/home/JobMar\_{}Proofs.pdf}) 

\newpage{}

\end{document}


\appendix

\title{Supplementary Material: Proofs to {}``High-Dimensional Semiparametric
Selection Models: Estimation Theory with an Application to the Retail
Gasoline Market'' }

\author{Ying Zhu}

\maketitle

\lyxaddress{\textit{Haas School of Business, University of California, Berkeley,
2220 Piedmont Ave., Berkeley, CA 94720. Email: ying\_zhu@haas.berkeley.edu.
Tel: 406-465-0498. Fax: 510-643-4255}}

\section{Appendix I: Main proofs}

\subsection{Lemmas A.1-A.2}

\textbf{Lemma A.1}: Let $\hat{\Gamma}=\frac{1}{n}\hat{v}^{T}\hat{v}$
and 
\begin{eqnarray*}
e & = & \left[\hat{\mathbb{E}}\left(x\vert w\hat{\theta}\right)-\mathbb{E}\left(x\vert w\hat{\theta}\right)+\mathbb{E}\left(x\vert w\hat{\theta}\right)-\mathbb{E}\left(x\vert w\theta^{*}\right)\right]\beta^{*}\\
 & - & \left[\hat{\mathbb{E}}\left(y_{2}\vert w\hat{\theta}\right)-\mathbb{E}\left(y_{2}\vert w\hat{\theta}\right)+\mathbb{E}\left(y_{2}\vert w\hat{\theta}\right)-\mathbb{E}\left(y_{2}\vert w\theta^{*}\right)\right]+\eta.
\end{eqnarray*}
Suppose the assumptions in Theorems 4.1 or 4.2 hold. If $\lambda_{n,3}$
in program (8) satisfies 
\[
\lambda_{n,3}\geq2|\frac{1}{n}\hat{v}^{T}e|_{\infty}>0,
\]
and 
\[
R_{q_{2}}\underbar{\ensuremath{\tau}}^{-q_{2}}\left(\frac{\log p}{n}+\mathcal{T}_{1}+\mathcal{T}_{2}+\mathcal{T}_{3}\right)=O(\kappa_{L}),
\]
then there is a constant $c>0$ such that under Assumption 2.2 ($q_{2}\in[0,\,1]$
and when $q_{2}=0$, $R_{q_{2}}:=k_{2}$), 
\[
|\hat{\beta}_{HSEL}-\beta^{*}|_{2}\leq\frac{c}{\kappa_{L}^{1-\frac{q_{2}}{2}}}R_{q_{2}}^{\frac{1}{2}}\left(\lambda_{n,3}\vee|\frac{2}{n}\hat{v}^{T}e|_{\infty}\right)^{1-\frac{q_{2}}{2}}.
\]
\textbf{Proof}. First, write 
\begin{eqnarray*}
v_{0} & = & \hat{v}_{0}+\hat{\mathbb{E}}\left(y_{2}\vert w\hat{\theta}\right)-\mathbb{E}\left(y_{2}\vert w\hat{\theta}\right)+\mathbb{E}\left(y_{2}\vert w\hat{\theta}\right)-\mathbb{E}\left(y_{2}\vert w\theta^{*}\right)\\
 & = & v\beta^{*}+\eta=\left[\hat{v}+\hat{\mathbb{E}}\left(x\vert w\hat{\theta}\right)-\mathbb{E}\left(x\vert w\hat{\theta}\right)+\mathbb{E}\left(x\vert w\hat{\theta}\right)-\mathbb{E}\left(x\vert w\theta^{*}\right)\right]\beta^{*}+\eta,
\end{eqnarray*}
thus we have 
\begin{eqnarray*}
\hat{v}_{0} & = & \hat{v}\beta^{*}+\left[\hat{\mathbb{E}}\left(x\vert w\hat{\theta}\right)-\mathbb{E}\left(x\vert w\hat{\theta}\right)+\mathbb{E}\left(x\vert w\hat{\theta}\right)-\mathbb{E}\left(x\vert w\theta^{*}\right)\right]\beta^{*}\\
 &  & -\left[\hat{\mathbb{E}}\left(y_{2}\vert w\hat{\theta}\right)-\mathbb{E}\left(y_{2}\vert w\hat{\theta}\right)+\mathbb{E}\left(y_{2}\vert w\hat{\theta}\right)-\mathbb{E}\left(y_{2}\vert w\theta^{*}\right)\right]+\eta\\
 & := & \hat{v}\beta^{*}+e
\end{eqnarray*}
where 
\begin{eqnarray*}
e & = & \left[\hat{\mathbb{E}}\left(x\vert w\hat{\theta}\right)-\mathbb{E}\left(x\vert w\hat{\theta}\right)+\mathbb{E}\left(x\vert w\hat{\theta}\right)-\mathbb{E}\left(x\vert w\theta^{*}\right)\right]\beta^{*}\\
 & - & \left[\hat{\mathbb{E}}\left(y_{2}\vert w\hat{\theta}\right)-\mathbb{E}\left(y_{2}\vert w\hat{\theta}\right)+\mathbb{E}\left(y_{2}\vert w\hat{\theta}\right)-\mathbb{E}\left(y_{2}\vert w\theta^{*}\right)\right]+\eta.
\end{eqnarray*}
Define the thresholded subset 
\[
S_{\underbar{\ensuremath{\tau}}}:=\left\{ j\in\left\{ 1,\,2,\,...,\, p\right\} \,:\,\left|\beta_{j}^{*}\right|>\underbar{\ensuremath{\tau}}\right\} 
\]
where $\underbar{\ensuremath{\tau}}=\frac{\lambda_{n,3}}{\kappa_{L}}$
is the threshold parameter. Define $\hat{\Delta}=\hat{\beta}_{HSEL}-\beta^{*}$
and the Lagrangian $L(\beta;\,\lambda_{n,3})=\frac{1}{2n}|\hat{u}-\hat{v}\beta|_{2}^{2}+\lambda_{n,3}|\beta|_{1}$.
Since $\hat{\beta}_{HSEL}$ is optimal, we have 
\[
L(\hat{\beta}_{HSEL};\,\lambda_{n,3})\leq L(\beta^{*};\,\lambda_{n,3})=\frac{1}{2n}|e|_{2}^{2}+\lambda_{n,3}|\beta^{*}|_{1},
\]
Some algebraic manipulation of the \textit{basic inequality} above
yields 
\begin{eqnarray}
0 & \leq & \frac{1}{2n}|\hat{v}\hat{\Delta}|_{2}^{2}\leq\frac{1}{n}e^{T}\hat{v}\hat{\Delta}+\lambda_{n,3}\left\{ |\beta_{S_{\underbar{\ensuremath{\tau}}}}^{*}|_{1}+|\beta_{S_{\underbar{\ensuremath{\tau}}}^{c}}^{*}|_{1}-|(\beta_{S_{\underbar{\ensuremath{\tau}}}}^{*}+\hat{\Delta}_{S_{\underbar{\ensuremath{\tau}}}},\,\beta_{S_{\underbar{\ensuremath{\tau}}}^{c}}^{*}+\hat{\Delta}_{S_{\underbar{\ensuremath{\tau}}}^{c}})|_{1}\right\} \nonumber \\
 & \leq & |\hat{\Delta}|_{1}|\frac{1}{n}\hat{v}^{T}e|_{\infty}+\lambda_{n,3}\left\{ |\hat{\Delta}_{S_{\underbar{\ensuremath{\tau}}}}|_{1}-|\hat{\Delta}_{S_{\underbar{\ensuremath{\tau}}}^{c}}|_{1}+2|\beta_{S_{\underbar{\ensuremath{\tau}}}^{c}}^{*}|_{1}\right\} \nonumber \\
 & \leq & \frac{\lambda_{n,3}}{2}\left\{ 3|\hat{\Delta}_{S_{\underbar{\ensuremath{\tau}}}}|_{1}-|\hat{\Delta}_{S_{\underbar{\ensuremath{\tau}}}^{c}}|_{1}+4|\beta_{S_{\underbar{\ensuremath{\tau}}}^{c}}^{*}|_{1}\right\} ,\label{eq:1}
\end{eqnarray}
where the last inequality holds as long as $\lambda_{n,3}\geq2|\frac{1}{n}\hat{v}^{T}e|_{\infty}>0$.
Consequently, 
\[
|\hat{\Delta}|_{1}\leq4|\hat{\Delta}_{S_{\underbar{\ensuremath{\tau}}}}|_{1}+4|\beta_{S_{\underbar{\ensuremath{\tau}}}^{c}}^{*}|_{1}\leq4\sqrt{\left|S_{\underbar{\ensuremath{\tau}}}\right|}|\hat{\Delta}|_{2}+4|\beta_{S_{\underbar{\ensuremath{\tau}}}^{c}}^{*}|_{1}.
\]
We now upper bound the cardinality of $S_{\underbar{\ensuremath{\tau}}}$
in terms of the threshold $\underbar{\ensuremath{\tau}}$ and $l_{q}-$ball
radius of $R_{q_{2}}$. Note that we have 
\[
R_{q_{2}}\geq\sum_{j=1}^{p}\left|\beta_{j}^{*}\right|^{q_{2}}\geq\sum_{j\in S_{\underbar{\ensuremath{\tau}}}}\left|\beta_{j}^{*}\right|^{q_{2}}\geq\underbar{\ensuremath{\tau}}^{q_{2}}\left|S_{\underbar{\ensuremath{\tau}}}\right|
\]
and therefore $\left|S_{\underbar{\ensuremath{\tau}}}\right|\leq\underbar{\ensuremath{\tau}}^{-q_{2}}R_{q_{2}}$.
To upper bound the approximation error $|\beta_{S_{\underbar{\ensuremath{\tau}}}^{c}}^{*}|_{1}$,
we use the fact that $\beta^{*}\in\mathcal{B}_{q_{2}}^{p}(R_{q_{2}})$
and have 
\[
|\beta_{S_{\underbar{\ensuremath{\tau}}}^{c}}^{*}|_{1}=\sum_{j\in S_{\underbar{\ensuremath{\tau}}}^{c}}\left|\beta_{j}^{*}\right|=\sum_{j\in S_{\underbar{\ensuremath{\tau}}}^{c}}\left|\beta_{j}^{*}\right|^{q_{2}}\left|\beta_{j}^{*}\right|^{1-q_{2}}\leq R_{q_{2}}\underbar{\ensuremath{\tau}}^{1-q_{2}}.
\]
Putting the pieces together yields 
\begin{equation}
|\hat{\Delta}|_{1}\leq4\sqrt{\underbar{\ensuremath{\tau}}^{-q_{2}}R_{q_{2}}}|\hat{\Delta}|_{2}+4R_{q_{2}}\underbar{\ensuremath{\tau}}^{1-q_{2}}.\label{eq:2}
\end{equation}
Let us first prove the case of $q_{2}\in(0,\,1]$. Note that we also
have 
\begin{eqnarray*}
\frac{1}{2n}|\hat{v}\hat{\Delta}|_{2}^{2} & \leq & |\hat{\Delta}|_{1}|\frac{1}{n}\hat{v}^{T}e|_{\infty}+\lambda_{n,3}\left\{ |\hat{\Delta}_{S(\underbar{\ensuremath{\tau}})}|_{1}-|\hat{\Delta}_{S(\underbar{\ensuremath{\tau}})^{c}}|_{1}+|\beta_{S(\underbar{\ensuremath{\tau}})^{c}}^{*}|_{1}\right\} \\
 & \leq & \left(8\sqrt{\underbar{\ensuremath{\tau}}^{-q_{2}}R_{q_{2}}}|\hat{\Delta}|_{2}+4R_{q_{2}}\underbar{\ensuremath{\tau}}^{1-q_{2}}\right)\lambda_{n,3}\\
 & \leq & c_{1}\sqrt{\underbar{\ensuremath{\tau}}^{-q_{2}}R_{q_{2}}}|\hat{\Delta}|_{2}\lambda_{n,3}+c_{2}\underbar{\ensuremath{\delta}}\\
 & \leq & \max\left\{ c_{1}R_{q_{2}}^{\frac{1}{2}}\kappa_{L}^{\frac{q_{2}}{2}}\left(\lambda_{n,3}\right)^{1-\frac{q_{2}}{2}}|\hat{\Delta}|_{2},\: c_{2}\underbar{\ensuremath{\delta}}\right\} 
\end{eqnarray*}
where the third and fourth inequalities follow from our choices of
$\underbar{\ensuremath{\tau}}=\frac{\lambda_{n,3}}{\kappa_{L}}$ and
$\underbar{\ensuremath{\delta}}\asymp R_{q_{2}}\lambda_{n,3}\tau^{1-q_{2}}$.
Now we proceed by cases. If 
\[
\max\left\{ c_{1}R_{q_{2}}^{\frac{1}{2}}\kappa_{L}^{\frac{q_{2}}{2}}\left(\lambda_{n,3}\right)^{1-\frac{q_{2}}{2}}|\hat{\Delta}|_{2},\: c_{2}\underbar{\ensuremath{\delta}}\right\} =c_{1}R_{q_{2}}^{\frac{1}{2}}\kappa_{L}^{\frac{q_{2}}{2}}\left(\lambda_{n,3}\right)^{1-\frac{q_{2}}{2}}|\hat{\Delta}|_{2},
\]
so that under the condition 
\[
R_{q_{2}}\underbar{\ensuremath{\tau}}^{-q_{2}}\left(\frac{\log p}{n}+\mathcal{T}_{1}+\mathcal{T}_{2}+\mathcal{T}_{3}\right)=O(\kappa_{L}),
\]
and provided $c>0$ is sufficiently large, we have 
\begin{equation}
|\hat{\Delta}|_{2}\geq c\kappa_{L}^{-1+\frac{q_{2}}{2}}R_{q_{2}}^{\frac{1}{2}}\left(\lambda_{n,3}\right)^{1-\frac{q_{2}}{2}}\geq\delta^{*}\label{eq:3}
\end{equation}
where 
\[
\delta^{*}=O\left(\max\left\{ \mathcal{T}_{1}^{\frac{3}{2}-q_{2}},\,\mathcal{T}_{2}^{\frac{3}{2}-q_{2}},\,\mathcal{T}_{3}^{\frac{3}{2}-q_{2}},\,\mathcal{T}_{4}^{2-q_{2}}\right\} R_{q_{2}}\right).
\]
Consequently, (\ref{eq:1}) and (\ref{eq:3}) together imply that
\begin{equation}
\hat{\Delta}\in\mathbb{K}(\underbar{\ensuremath{\delta}},\, S_{\underbar{\ensuremath{\tau}}}):=\mathbb{C}(S_{\underbar{\ensuremath{\tau}}};q_{2},\,3)\cap\left\{ \Delta\in\mathbb{R}^{p}\,:\,\left|\Delta\right|_{2}\geq\delta^{*}\right\} \label{eq:4}
\end{equation}
where $\mathbb{C}(S_{\underbar{\ensuremath{\tau}}};\,3):=\left\{ \Delta\in\mathbb{R}^{p}\,:\,|\Delta_{S_{\underbar{\ensuremath{\tau}}}^{c}}|_{1}\leq3|\Delta_{S_{\underbar{\ensuremath{\tau}}}}|_{1}+4|\beta_{S_{\underbar{\ensuremath{\tau}}}^{c}}^{*}|_{1}\right\} $.
By Lemma A.2 and Lemma A.3 (or Lemma A.5), the random matrix $\hat{\Gamma}=\hat{v}^{T}\hat{v}$
satisfies the RE condition (2.1) over 
\[
\mathbb{C}(S_{\underbar{\ensuremath{\tau}}};\,3)\cap\left\{ \Delta\in\mathbb{R}^{p}\,:\,\left|\Delta\right|_{2}\geq\delta^{*}\right\} ,
\]
therefore, we have 
\[
\kappa_{L}|\hat{\Delta}|_{2}^{2}\leq\frac{1}{2n}|\hat{v}\hat{\Delta}|_{2}^{2}\leq c_{1}R_{q_{2}}^{\frac{1}{2}}\kappa_{L}^{\frac{q_{2}}{2}}\left(\lambda_{n,3}\right)^{1-\frac{q_{2}}{2}}|\hat{\Delta}|_{2}
\]
so the claim follows. On the other hand, if 
\[
\max\left\{ c_{1}R_{q_{2}}^{\frac{1}{2}}\kappa_{L}^{\frac{q_{2}}{2}}\left(\lambda_{n,3}\right)^{1-\frac{q_{2}}{2}}|\hat{\Delta}|_{2},\: c_{2}\underbar{\ensuremath{\delta}}\right\} =\underbar{\ensuremath{\delta}},
\]
then 
\[
|\hat{\Delta}|_{2}\leq c\kappa_{L}^{-1+\frac{q_{2}}{2}}R_{q_{2}}^{\frac{1}{2}}\left(\lambda_{n,3}\right)^{1-\frac{q_{2}}{2}}
\]
so again the claim follows.

To prove the case of $q_{2}=0$, simply choose $\underbar{\ensuremath{\delta}}=0$
and $S_{\underbar{\ensuremath{\tau}}}=J(\beta^{*})$ and the claim
follows trivially from the above argument.$\square$\\
\\
\textbf{Remark}. Inequality (\ref{eq:2}) implies that $|\hat{\Delta}|_{1}\precsim R_{q_{2}}\left(\lambda_{n,3}\right)^{1-q_{2}}$.\\
\\
\textbf{Lemma A.2}:\textbf{ }Define the thresholded subset 
\[
S_{\underbar{\ensuremath{\tau}}}:=\left\{ j\in\left\{ 1,\,2,\,...,\, p\right\} \,:\,\left|\beta_{j}^{*}\right|>\underbar{\ensuremath{\tau}}\right\} .
\]
Under the assumptions in Theorem 4.1 (or Theorem 4.2) and the choice
$\underbar{\ensuremath{\tau}}=\frac{\lambda_{n,3}}{\kappa_{L}}$,
if 
\[
R_{q_{2}}\underbar{\ensuremath{\tau}}^{-q_{2}}\left(\frac{\log p}{n}+\mathcal{T}_{1}+\mathcal{T}_{2}+\mathcal{T}_{3}\right)=O(\kappa_{L}),
\]
the RE condition (2.1) of $\frac{\hat{v}^{T}\hat{v}}{n}$ holds over
the set 
\[
\mathbb{C}(S_{\underbar{\ensuremath{\tau}}};q_{2},\,3)\cap\left\{ \Delta\in\mathbb{R}^{p}\,:\,\left|\Delta\right|_{2}\geq\delta^{*}\right\} 
\]
where 
\[
\delta^{*}=O\left(\max\left\{ \mathcal{T}_{1}^{\frac{3}{2}-q_{2}},\,\mathcal{T}_{2}^{\frac{3}{2}-q_{2}},\,\mathcal{T}_{3}^{\frac{3}{2}-q_{2}},\,\mathcal{T}_{4}^{2-q_{2}}\right\} R_{q_{2}}\right).
\]
\textbf{Proof}. The argument is similar to what is used in the proof
of Lemma 2 from Negahban, et. al (2010). For any $\Delta\in\mathbb{C}(S_{\underbar{\ensuremath{\tau}}};q_{2},\,3)$,
we have 
\[
|\Delta|_{1}\leq4|\Delta_{S_{\underbar{\ensuremath{\tau}}}}|_{1}+4|\beta_{S_{\underbar{\ensuremath{\tau}}}^{c}}^{*}|_{1}\leq4\sqrt{\left|S_{\underbar{\ensuremath{\tau}}}\right|}|\Delta|_{2}+4R_{q_{2}}\underbar{\ensuremath{\tau}}^{1-q_{2}}\leq4\sqrt{R_{q_{2}}}\underbar{\ensuremath{\tau}}^{-\frac{q_{2}}{2}}|\Delta|_{2}+4R_{q_{2}}\underbar{\ensuremath{\tau}}^{1-q_{2}},
\]
where we have used the bound in (\ref{eq:2}) from the proof of Lemma
A.1. Therefore, for any vector $\Delta\in\mathbb{C}(S_{\underbar{\ensuremath{\tau}}};q_{2},\,3)$
and the choice $\underbar{\ensuremath{\tau}}=\frac{\lambda_{n,3}}{\kappa_{L}}$,
substituting the upper bound $4\sqrt{R_{q_{2}}}\underbar{\ensuremath{\tau}}^{-\frac{q_{2}}{2}}|\Delta|_{2}+4R_{q_{2}}\underbar{\ensuremath{\tau}}^{1-q_{2}}$
on $|\Delta|_{1}$ into condition (\ref{eq:5}) from Lemma A.3 and
condition (\ref{eq:9}) from Lemma A.5 respectively yields 
\begin{eqnarray*}
\left|\Delta^{T}\frac{\hat{v}^{T}\hat{v}}{n}\Delta\right| & \geq & |\Delta|_{2}^{2}\left\{ c\kappa_{L}-b_{0}R_{q_{2}}\underbar{\ensuremath{\tau}}^{-q_{2}}\left(\frac{\log p}{n}+\mathcal{T}_{1}+\mathcal{T}_{2}+\mathcal{T}_{3}\right)\right\} \\
 &  & -b_{1}R_{q_{2}}^{2}\underbar{\ensuremath{\tau}}^{2-2q_{2}}\left(\frac{\log p}{n}+\mathcal{T}_{1}+\mathcal{T}_{2}+\mathcal{T}_{3}\right),
\end{eqnarray*}
for some positive constants $b_{0}$ and $b_{1}$. With the choice
of 
\[
\delta^{*}=O\left(\max\left\{ \mathcal{T}_{1}^{\frac{3}{2}-q_{2}},\,\mathcal{T}_{2}^{\frac{3}{2}-q_{2}},\,\mathcal{T}_{3}^{\frac{3}{2}-q_{2}},\,\mathcal{T}_{4}^{2-q_{2}}\right\} R_{q_{2}}\right),
\]
if 
\[
R_{q_{2}}\underbar{\ensuremath{\tau}}^{-q_{2}}\left(\frac{\log p}{n}+\mathcal{T}_{1}+\mathcal{T}_{2}+\mathcal{T}_{3}\right)=O(\kappa_{L}),
\]
we have 
\[
\left|\Delta^{T}\frac{\hat{v}^{T}\hat{v}}{n}\Delta\right|\geq c^{'}\kappa_{L}\left\{ |\Delta|_{2}^{2}-\frac{|\Delta|_{2}^{2}}{2}\right\} =c^{''}\kappa_{L}|\Delta|_{2}^{2}
\]
for any $\Delta$ such that $\left|\Delta\right|_{2}\geq\delta^{*}$.
$\square$

\subsection{Theorem 4.1 }

Lemma A.1 implies that the $l_{2}$-consistency of $\hat{\beta}_{HSEL}$
requires verifications of two conditions: (i) $\hat{\Gamma}=\hat{v}^{T}\hat{v}$
satisfies the RE condition (2.1) with $\gamma=3$, and (ii) the term
$|\frac{1}{n}\hat{v}^{T}e|_{\infty}\precsim\phi(k_{1},\, k_{2},\, d,\, p,\, n)$
with high probability. This is done via Lemmas A.3 and A.4. For the
proofs for Lemmas A.3 and A.4, let the \textit{local complexity} measure
$\mathcal{G}_{n}\left(r_{j};\,\mathcal{F}_{j}\right)$ be defined
over the set 
\[
\Omega(r_{j};\,\mathcal{F}_{j})=\left\{ f\,:\, f\in\bar{\mathcal{F}}_{j}\,\left|f_{\theta^{*}}\right|_{n}\leq r_{j}\right\} 
\]
where $\left|f_{\theta^{*}}\right|_{n}:=\sqrt{\frac{1}{n}\sum_{i=1}^{n}\left[f(w_{i}^{T}\theta^{*})\right]^{2}}$
and $j=0,...,p$. Recall the notations $\mathbb{E}\left(x_{ij}\vert w_{i}^{T}\theta\right):=m_{j}(w_{i}^{T}\theta)$,
$\hat{\mathbb{E}}\left(x_{ij}\vert w_{i}^{T}\theta\right):=\hat{m}_{j}(w_{i}^{T}\theta)$,
$\mathbb{E}\left(y_{2i}\vert w_{i}^{T}\theta\right):=m_{0}(w_{i}^{T}\theta)$,
$\hat{\mathbb{E}}\left(y_{2i}\vert w_{i}^{T}\theta\right):=\hat{m}_{0}(w_{i}^{T}\theta)$.
To avoid clustering, write $\left\{ \hat{m}_{j}(w_{i}^{T}\theta)\right\} _{i=1}^{n}:=\hat{m}_{j}(\theta)$,
$\left\{ \tilde{m}_{j}(w_{i}^{T}\theta)\right\} _{i=1}^{n}:=\tilde{m}_{j}(\theta)$,
and $\left\{ m_{j}(w_{i}^{T}\theta)\right\} _{i=1}^{n}:=m_{j}(\theta)$;
these definitions are somewhat abuse of notation but keep in mind
that the individual component $\hat{m}_{j}(\cdot)$ of $\hat{m}_{j}(\theta)$,
$\tilde{m}_{j}(\cdot)$ of $\tilde{m}_{j}(\theta)$, and $m_{j}(\cdot)$
of $m_{j}(\theta)$ are functions from $\mathbb{R}$ to $\mathbb{R}$.
Also recall $T_{j}^{'}:=T_{j}^{*}\vee T_{j}$, $T_{j}^{*}:=\sup_{f\in\mathcal{F}_{j}^{*}}\frac{1}{n}\sum_{i=1}^{n}\left[f(w_{i}^{T}\hat{\theta})-f(w_{i}^{T}\theta^{*})\right]^{2}$,
$T_{j}:=\sup_{f\in\mathcal{F}_{j}}\frac{1}{n}\sum_{i=1}^{n}\left[f(w_{i}^{T}\hat{\theta})-f(w_{i}^{T}\theta^{*})\right]^{2}$,
and the following definitions for Theorem 4.1: 
\begin{eqnarray*}
\mathcal{T}_{1} & = & \max_{j\in\{0,...,p\}}\left(T_{j}^{'}\vee\sqrt{T_{j}^{'}}\right)\\
\mathcal{T}_{2} & = & \max_{j\in\{0,...,p\}}t_{nj}^{2}\\
\mathcal{T}_{3} & = & \max_{j\in\{0,...,p\}}\inf_{\tilde{m}_{j}\in\mathcal{F}_{j}}\left(\left|\tilde{m}_{j}(\hat{\theta})-m_{j}(\hat{\theta})\right|_{n}^{2}+\left|\tilde{m}_{j}(\theta^{*})-m_{j}(\theta^{*})\right|_{n}\right)\\
\mathcal{T}_{4} & = & \sigma_{v}\sigma_{\eta}\sqrt{\frac{\log p}{n}}.
\end{eqnarray*}
To facilitate the presentation, the proofs for Lemmas A.3 and A.4
work through the case $t_{nj}=r_{nj}$. Inspecting these proofs suggests
that with minor notation changes, the case where $t_{nj}\geq r_{nj}$
can be analyzed with almost the same argument because for any star-shaped
function class $\bar{\mathcal{F}}_{j}:=\left\{ f=f^{'}-f^{''}\,:\, f^{'},\, f^{''}\in\mathcal{F}_{j}\right\} $,
the function $t\mapsto\frac{\mathcal{G}_{n}(t;\,\mathcal{F}_{j})}{t}$
is non-decreasing on $(0,\,\infty)$. \\
\\
\textbf{Lemma A.3} (RE condition): Define $\hat{\Delta}=\hat{\beta}_{HSEL}-\beta^{*}$.
Let $r_{nj}>0$ be the smallest positive quantity satisfying the \textit{critical
inequality} 
\[
\mathcal{G}_{n}\left(r_{nj};\,\mathcal{F}_{j}\right)\leq\frac{r_{nj}^{2}}{\sigma_{v_{j}}}.
\]
Under Assumptions 2.1, 2.3, 4.1-4.3, we have 
\begin{equation}
\left|\Delta^{T}\frac{\hat{v}^{T}\hat{v}}{n}\Delta\right|\geq\kappa_{1}|\Delta|_{2}^{2}-\kappa_{2}\frac{\log p}{n}|\Delta|_{1}^{2}-c\left(\mathcal{T}_{1}+\mathcal{T}_{2}+\mathcal{T}_{3}\right)\left|\Delta\right|_{1}^{2}\label{eq:5}
\end{equation}
with probability at least $1-c^{'}\sum_{j=1}^{p}\exp\left(-nC_{j}^{*}r_{nj}^{2}\right)$,
where $\kappa_{1}=\frac{\kappa_{L}}{2}$, $\kappa_{2}=c_{0}\kappa_{L}\max\left\{ \frac{\sigma_{v}^{4}}{\kappa_{L}^{2}},\,1\right\} $,
$\sigma_{v}=\max_{j}\sigma_{v_{j}}$.\\
\\
\textbf{Proof. }We have 
\[
\left|\Delta^{T}\frac{\hat{v}^{T}\hat{v}}{n}\Delta\right|+\left|\Delta^{T}\left(\frac{v^{T}v-\hat{v}^{T}\hat{v}}{n}\right)\Delta\right|\geq\left|\Delta^{T}\frac{v^{T}v}{n}\Delta\right|,
\]
which implies \textbf{
\begin{eqnarray*}
\left|\Delta^{T}\frac{\hat{v}^{T}\hat{v}}{n}\Delta\right| & \geq & \left|\Delta^{T}\frac{v^{T}v}{n}\Delta\right|-\left|\Delta^{T}\left(\frac{v^{T}v-\hat{v}^{T}\hat{v}}{n}\right)\Delta\right|\\
 & \geq & \left|\Delta^{T}\frac{v^{T}v}{n}\Delta\right|-\left|\frac{v^{T}v-\hat{v}^{T}\hat{v}}{n}\right|_{\infty}\left|\Delta\right|_{1}^{2}\\
 & \geq & \left|\Delta^{T}\frac{v^{T}v}{n}\Delta\right|-\left(\left|\frac{v^{T}(\hat{v}-v)}{n}\right|_{\infty}+\left|\frac{(\hat{v}-v)^{T}\hat{v}}{n}\right|_{\infty}\right)\left|\Delta\right|_{1}^{2}\\
 & \geq & \left|\Delta^{T}\frac{v^{T}v}{n}\Delta\right|-\left|\frac{v^{T}(\hat{v}-v)}{n}\right|_{\infty}\left|\Delta\right|_{1}^{2}\\
 &  & -\left|\frac{(\hat{v}-v)^{T}v}{n}\right|_{\infty}\left|\Delta\right|_{1}^{2}-\left|\frac{(\hat{v}-v)^{T}(\hat{v}-v)}{n}\right|_{\infty}\left|\Delta\right|_{1}^{2}.
\end{eqnarray*}
}To bound the term $\left|\frac{v^{T}(\hat{v}-v)}{n}\right|_{\infty}$,
let us first fix $(j,\, j^{'})$\textbf{ }and bound the $(j,\, j^{'})$
element of the matrix\textbf{ }$\frac{v^{T}(\hat{v}-v)}{n}$. Notice
that 
\begin{eqnarray*}
\left|\frac{1}{n}v_{j}^{T}(\hat{v}_{j^{'}}-v{}_{j^{'}})\right| & = & \left|\frac{1}{n}\sum_{i=1}^{n}v_{ij}(\hat{v}_{ij^{'}}-v{}_{ij^{'}})\right|\\
 & = & \left|\frac{1}{n}\sum_{i=1}^{n}v_{ij}\left(\hat{m}_{j^{'}}(w_{i}^{T}\hat{\theta})-m_{j^{'}}(w_{i}^{T}\hat{\theta})+m_{j^{'}}(w_{i}^{T}\hat{\theta})-m_{j^{'}}(w_{i}^{T}\theta^{*})\right)\right|\\
 & \leq & \left|\frac{1}{n}\sum_{i=1}^{n}v_{ij}\left[\hat{m}_{j^{'}}(w_{i}^{T}\hat{\theta})-m_{j^{'}}(w_{i}^{T}\hat{\theta})\right]\right|+\sqrt{\frac{1}{n}\sum_{i=1}^{n}v_{ij}^{2}}\sqrt{T_{j^{'}}^{*}},
\end{eqnarray*}
where $T_{j^{'}}^{*}:=\sup_{f\in\mathcal{F}_{j^{'}}^{*}}\frac{1}{n}\sum_{i=1}^{n}\left[f(w_{i}^{T}\hat{\theta})-f(w_{i}^{T}\theta^{*})\right]^{2}$. 

For the term $\sqrt{\frac{1}{n}\sum_{i=1}^{n}v_{ij}^{2}}$, by Lemma
B.1, we have 
\[
\mathbb{P}\left[\max_{j=1,...,p}\left\{ \frac{1}{n}\sum_{i=1}^{n}v_{ij}^{2}\right\} \geq\sigma_{v{}_{j}}^{2}+t\right]\leq2\exp(-c\min(\frac{nt^{2}}{4\sigma_{v_{j}}^{4}},\,\frac{nt}{2\sigma_{v_{j}}^{2}})+\log p).
\]
Hence, setting $t=\sigma_{v}^{2}:=\max_{j}\sigma_{v_{j}}^{2}$ and
under the condition $n\gtrsim\log p$, we have 
\[
\mathbb{P}\left[\max_{j^{'}}\left\{ \frac{1}{n}\sum_{i=1}^{n}v_{ij}^{2}\right\} \geq2\sigma_{v}^{2}\right]\leq2\exp(-c^{'}n).
\]
Let $\hat{\Delta}_{ij^{'}}(\theta):=\hat{m}_{j^{'}}(w_{i}^{T}\theta)-m_{j^{'}}(w_{i}^{T}\theta)$,
$\bar{\Delta}_{ij^{'}}(\theta):=\hat{m}_{j^{'}}(w_{i}^{T}\theta)-\tilde{m}_{j^{'}}(w_{i}^{T}\theta)$,
and $\tilde{\Delta}_{ij^{'}}(\theta):=\tilde{m}_{j^{'}}(w_{i}^{T}\theta)-m_{j^{'}}(w_{i}^{T}\theta)$.
Note that $\hat{m}_{j^{'}}\in\mathcal{F}_{j^{'}}^{*}$ and we choose
$\hat{m}_{j^{'}},\,\tilde{m}_{j^{'}}\in\mathcal{F}_{j^{'}}$. For
the term 
\[
\left|\frac{1}{n}\sum_{i=1}^{n}v_{ij}\left[\hat{m}_{j^{'}}(w_{i}^{T}\hat{\theta})-m_{j^{'}}(w_{i}^{T}\hat{\theta})\right]\right|,
\]
we have 
\begin{eqnarray*}
\left|\frac{1}{n}\sum_{i=1}^{n}v_{ij}\hat{\Delta}_{ij^{'}}(\hat{\theta})\right| & \leq & \left|\frac{1}{n}\sum_{i=1}^{n}v_{ij}\hat{\Delta}_{ij^{'}}(\theta^{*})\right|+\sqrt{\frac{1}{n}\sum_{i=1}^{n}v_{ij}^{2}}\sqrt{\frac{1}{n}\sum_{i=1}^{n}\left[\hat{\Delta}_{ij^{'}}(\hat{\theta})-\hat{\Delta}_{ij^{'}}(\theta^{*})\right]^{2}}\\
 & \leq & \left|\frac{1}{n}\sum_{i=1}^{n}v_{ij}\bar{\Delta}_{ij^{'}}(\theta^{*})\right|+\left|\frac{1}{n}\sum_{i=1}^{n}v_{ij}\tilde{\Delta}_{ij^{'}}(\theta^{*})\right|+2\sqrt{\frac{1}{n}\sum_{i=1}^{n}v_{ij}^{2}}\sqrt{T_{j^{'}}^{*}+T_{j^{'}}}\\
 & \leq & \left|\frac{1}{n}\sum_{i=1}^{n}v_{ij}\bar{\Delta}_{ij^{'}}(\theta^{*})\right|+\sqrt{\frac{1}{n}\sum_{i=1}^{n}v_{ij}^{2}}\left|\tilde{m}_{j^{'}}(\theta^{*})-m_{j^{'}}(\theta^{*})\right|_{n}+\sqrt{\frac{1}{n}\sum_{i=1}^{n}v_{ij}^{2}}\sqrt{8T_{j^{'}}^{'}}.
\end{eqnarray*}
To upper bound the term $\frac{1}{n}\sum_{i=1}^{n}v_{ij}\bar{\Delta}_{ij^{'}}(\theta^{*})$,
we argue similarly as in Lemma A.11. First, note that by Lemma A.11
and the triangle inequality, we have 
\begin{eqnarray*}
\left|\bar{\Delta}_{j^{'}}(\theta^{*})\right|_{n} & \leq & \left|\hat{m}_{j^{'}}(\hat{\theta})-m_{j^{'}}(\hat{\theta})\right|_{n}+\left|\tilde{m}_{j^{'}}(\hat{\theta})-m_{j^{'}}(\hat{\theta})\right|_{n}+2\sqrt{T_{j^{'}}}\\
 & \leq & c^{'}\left\{ \sqrt{T_{j^{'}}^{'}}+\left(\sigma_{v}^{2}T_{j^{'}}^{'}\right)^{\frac{1}{4}}+\left|\tilde{m}_{j^{'}}(\hat{\theta})-m_{j^{'}}(\hat{\theta})\right|_{n}+r_{nj^{'}}\right\} :=\bar{r}_{nj^{'}}.
\end{eqnarray*}
Setting $u=\tilde{r}_{nj^{'}}$ in $\Omega(u;\,\mathcal{F}_{j})$
from Lemma B.3 and following the argument in the proof for Lemma B.3,
we obtain

\[
\max_{j^{'}}\frac{1}{n}\sum_{i=1}^{n}v_{ij}\bar{\Delta}_{ij^{'}}(\theta^{*})\leq2\max_{j^{'}}\bar{r}_{nj^{'}}^{2}
\]
with probability at least $1-c\sum_{j=1}^{p}\exp\left(-nC_{j}^{*}\bar{r}_{nj}^{2}\right)$.
Hence, we obtain 
\[
\max_{j,\, j^{'}}\left|\frac{1}{n}\sum_{i=1}^{n}v_{ij}\left[\hat{m}_{j^{'}}(w_{i}^{T}\hat{\theta})-m_{j^{'}}(w_{i}^{T}\hat{\theta})\right]\right|\leq c\max_{j^{'}}\left\{ \bar{r}_{nj^{'}}^{2}+\sigma_{v}\left|\tilde{m}_{j^{'}}(\theta^{*})-m_{j^{'}}(\theta^{*})\right|_{n}+\sqrt{\sigma_{v}^{2}T_{j^{'}}^{'}}\right\} .
\]
Consequently, we have 
\[
\left|\frac{v^{T}(\hat{v}-v)}{n}\right|_{\infty}\leq c\max_{j}\left\{ \bar{r}_{nj}^{2}+\sigma_{v}\left|\tilde{m}_{j^{'}}(\theta^{*})-m_{j^{'}}(\theta^{*})\right|_{n}+\sqrt{\sigma_{v}^{2}T_{j}^{'}}\right\} .
\]
To bound the term $\left|\frac{(\hat{v}-v)^{T}(\hat{v}-v)}{n}\right|_{\infty}$,
note that we have 
\begin{eqnarray*}
\left|\frac{(\hat{v}-v)^{T}(\hat{v}-v)}{n}\right|_{\infty} & \leq & 2\max_{j,\, j^{'}}\sqrt{\frac{1}{n}\sum_{i=1}^{n}\left[\hat{m}_{j}(w_{i}^{T}\hat{\theta})-m_{j}(w_{i}^{T}\hat{\theta})\right]{}^{2}}\cdot\sqrt{\frac{1}{n}\sum_{i=1}^{n}\left[\hat{m}_{j^{'}}(w_{i}^{T}\hat{\theta})-m_{j^{'}}(w_{i}^{T}\hat{\theta})\right]{}^{2}}\\
 &  & +2\max_{j,\, j^{'}}\sqrt{\frac{1}{n}\sum_{i=1}^{n}\left[m_{j}(w_{i}^{T}\hat{\theta})-m_{j}(w_{i}^{T}\theta^{*})\right]^{2}}\cdot\sqrt{\frac{1}{n}\sum_{i=1}^{n}\left[m_{j^{'}}(w_{i}^{T}\hat{\theta})-m_{j^{'}}(w_{i}^{T}\theta^{*})\right]{}^{2}}\\
 & \leq & c\max_{j}\left[T_{j}^{'}+\sqrt{\sigma_{v}^{2}T_{j}^{'}}+\left|\tilde{m}_{j}(\hat{\theta})-m_{j}(\hat{\theta})\right|_{n}^{2}+r_{nj}^{2}\right]
\end{eqnarray*}
with probability at least $1-c^{'}\sum_{j=1}^{p}\exp\left(-nC_{j}^{*}r_{nj}^{2}\right)$. 

Putting everything together and applying Lemma B.2, we have 
\begin{eqnarray*}
\left|\Delta^{T}\frac{\hat{v}^{T}\hat{v}}{n}\Delta\right| & \geq & \kappa_{1}|\Delta|_{2}^{2}-\kappa_{2}\frac{\log p}{n}|\Delta|_{1}^{2}\\
 & - & b(\sigma_{v})\max_{j}\left\{ r_{nj}^{2}+\inf_{\tilde{m}_{j}\in\mathcal{F}_{j}}\left(\left|\tilde{m}_{j}(\hat{\theta})-m_{j}(\hat{\theta})\right|_{n}^{2}+\left|\tilde{m}_{j}(\theta^{*})-m_{j}(\theta^{*})\right|_{n}\right)+\sqrt{T_{j}^{'}\vee T_{j}^{'2}}\right\} \left|\Delta\right|_{1}^{2}
\end{eqnarray*}
where $\kappa_{1}=\frac{\kappa_{L}}{2}$, $\kappa_{2}=c_{0}\kappa_{L}\max\left\{ \frac{\sigma_{v}^{4}}{\kappa_{L}^{2}},\,1\right\} $,
$\sigma_{v}=\max_{j}\sigma_{v_{j}}$, and a known function $b(\sigma_{v})$
depending only on $\sigma_{v}$.$\square$\\
\\
\textbf{Lemma A.4} (Upper bound on $\left|\frac{1}{n}\hat{v}^{T}e\right|_{\infty}$):
Let $r_{nj}>0$ be the smallest positive quantity satisfying the \textit{critical
inequality} 
\[
\mathcal{G}_{n}\left(r_{nj};\,\mathcal{F}_{j}\right)\leq\frac{r_{nj}^{2}}{\sigma_{v_{j}}}.
\]
Under Assumptions 2.1 and 4.1-4.3, we have 
\[
\left|\frac{\hat{v}^{T}e}{n}\right|_{\infty}\leq\left|\beta^{*}\right|_{1}b(\sigma_{v},\,\sigma_{\eta})\left(\mathcal{T}_{1}+\mathcal{T}_{2}+\mathcal{T}_{3}\right)+\mathcal{T}_{4},
\]
with probability at least $1-c_{1}\exp\left(-c_{2}\log p\right)-c_{3}\sum_{j=0}^{p}\exp\left(-nC_{j}^{*}r_{nj}^{2}\right)$.\\
\\
\textbf{Proof.} Recall that 
\begin{eqnarray*}
e & = & \left[\hat{\mathbb{E}}\left(x\vert w\hat{\theta}\right)-\mathbb{E}\left(x\vert w\hat{\theta}\right)+\mathbb{E}\left(x\vert w\hat{\theta}\right)-\mathbb{E}\left(x\vert w\theta^{*}\right)\right]\beta^{*}\\
 & - & \left[\hat{\mathbb{E}}\left(y_{2}\vert w\hat{\theta}\right)-\mathbb{E}\left(y_{2}\vert w\hat{\theta}\right)+\mathbb{E}\left(y_{2}\vert w\hat{\theta}\right)-\mathbb{E}\left(y_{2}\vert w\theta^{*}\right)\right]+\eta.
\end{eqnarray*}
and 
\begin{eqnarray*}
\hat{v} & = & x-\hat{\mathbb{E}}\left(x\vert w\hat{\theta}\right)\\
 & = & x-\mathbb{E}\left(x\vert w\theta^{*}\right)+\mathbb{E}\left(x\vert w\theta^{*}\right)-\mathbb{E}\left(x\vert w\hat{\theta}\right)+\mathbb{E}\left(x\vert w\hat{\theta}\right)-\hat{\mathbb{E}}\left(x\vert w\hat{\theta}\right).
\end{eqnarray*}
Recall $v=x-\mathbb{E}\left(x\vert w\theta^{*}\right)$. Let us introduce
the following notations $m_{j}(w_{i}^{T}\hat{\theta})-m_{j}(w_{i}^{T}\theta^{*}):=T_{ij}^{''}$
and $\hat{m}_{j}(w_{i}^{T}\theta)-m_{j}(w_{i}^{T}\theta):=\hat{\Delta}_{j}(w_{i}^{T}\theta)$,
and 
\begin{eqnarray}
T_{1} & = & \max_{j,\, j^{'}=1,...,p}\left|\frac{1}{n}\sum_{i=1}^{n}\left[v_{ij}-T_{ij}^{''}-\hat{\Delta}_{j}(w_{i}^{T}\hat{\theta})\right]\left[T_{ij^{'}}^{''}+\hat{\Delta}_{j^{'}}(w_{i}^{T}\hat{\theta})\right]\right|\left|\beta^{*}\right|_{1}\label{eq:6}\\
T_{2} & = & \max_{j=1,...,p}\left|\frac{1}{n}\sum_{i=1}^{n}\left[v_{ij}-T_{ij}^{''}-\hat{\Delta}_{j}(w_{i}^{T}\hat{\theta})\right]\left[T_{i0}^{''}+\hat{\Delta}_{0}(w_{i}^{T}\hat{\theta})\right]\right|\label{eq:7}\\
T_{3} & = & \max_{j=1,...,p}\left|\frac{1}{n}\sum_{i=1}^{n}\left[v_{ij}-T_{ij}^{''}-\hat{\Delta}_{j}(w_{i}^{T}\hat{\theta})\right]\eta_{i}\right|.\label{eq:8}
\end{eqnarray}
Expanding the products in (\ref{eq:6})-(\ref{eq:8}), applying the
Cauchy-Schwarz inequality on some of the terms in the expansion, we
obtain the following inequalities 
\begin{eqnarray*}
\frac{1}{n}\sum_{i=1}^{n}v_{ij}T_{ij^{'}}^{''} & \leq & \sqrt{\frac{1}{n}\sum_{i=1}^{n}v_{ij}^{2}}\sqrt{\frac{1}{n}\sum_{i=1}^{n}T_{ij^{'}}^{''2}}\\
\frac{1}{n}\sum_{i=1}^{n}v_{ij}\hat{\Delta}_{ij^{'}}(\hat{\theta}) & \leq & \left|\frac{1}{n}\sum_{i=1}^{n}v_{ij}\hat{\Delta}_{j^{'}}(w_{i}^{T}\theta^{*})\right|+\sqrt{\frac{1}{n}\sum_{i=1}^{n}v_{ij}^{2}}\sqrt{\frac{1}{n}\sum_{i=1}^{n}\left[\hat{\Delta}_{j^{'}}(w_{i}^{T}\hat{\theta})-\hat{\Delta}_{j^{'}}(w_{i}^{T}\theta^{*})\right]^{2}}\\
\frac{1}{n}\sum_{i=1}^{n}T_{ij}^{''}T_{ij^{'}}^{''} & \leq & \sqrt{\frac{1}{n}\sum_{i=1}^{n}T_{ij}^{''2}}\sqrt{\frac{1}{n}\sum_{i=1}^{n}T_{ij^{'}}^{''2}}\\
\frac{1}{n}\sum_{i=1}^{n}T_{ij}^{''}\hat{\Delta}_{ij^{'}}(\hat{\theta}) & \leq & \sqrt{\frac{1}{n}\sum_{i=1}^{n}T_{ij}^{''2}}\sqrt{\frac{1}{n}\sum_{i=1}^{n}\hat{\Delta}_{j^{'}}^{2}(w_{i}^{T}\hat{\theta})}\\
\frac{1}{n}\sum_{i=1}^{n}\hat{\Delta}_{ij}(\hat{\theta})\hat{\Delta}_{ij^{'}}(\hat{\theta}) & \leq & \sqrt{\frac{1}{n}\sum_{i=1}^{n}\hat{\Delta}_{j}^{2}(w_{i}^{T}\hat{\theta})}\sqrt{\frac{1}{n}\sum_{i=1}^{n}\hat{\Delta}_{j^{'}}^{2}(w_{i}^{T}\hat{\theta})}\\
\frac{1}{n}\sum_{i=1}^{n}T_{ij}^{''}\eta_{i} & \leq & \sqrt{\frac{1}{n}\sum_{i=1}^{n}\eta_{i}^{2}}\sqrt{\frac{1}{n}\sum_{i=1}^{n}T_{ij^{'}}^{''2}}\\
\frac{1}{n}\sum_{i=1}^{n}\hat{\Delta}_{ij}(\hat{\theta})\eta_{i} & \leq & \left|\frac{1}{n}\sum_{i=1}^{n}\eta_{i}\hat{\Delta}_{j^{'}}(w_{i}^{T}\theta^{*})\right|+\sqrt{\frac{1}{n}\sum_{i=1}^{n}\eta_{i}^{2}}\sqrt{\frac{1}{n}\sum_{i=1}^{n}\left[\hat{\Delta}_{j^{'}}(w_{i}^{T}\hat{\theta})-\hat{\Delta}_{j^{'}}(w_{i}^{T}\theta^{*})\right]^{2}}.
\end{eqnarray*}
To upper bound $\left|\frac{\hat{v}^{T}e}{n}\right|_{\infty}$, we
need to control the RHS inequalities listed above as well as the additional
term $\max_{j}\left\{ \left|\frac{1}{n}\sum_{i=1}^{n}v_{ij}\eta_{i}\right|\right\} $.

Under Assumptions\uline{ }4.2 and 4.3, and the condition $n\gtrsim\log p$,
applying Lemma B.1 yields 
\begin{eqnarray*}
\mathbb{P}\left[\max_{j}\left\{ \frac{1}{n}\sum_{i=1}^{n}v_{ij}^{2}\right\} \geq2\sigma_{v}^{2}\right] & \leq & 2\exp(-c^{'}n),\\
\mathbb{P}\left[\frac{1}{n}\sum_{i=1}^{n}\eta_{i}^{2}\geq2\sigma_{\eta}^{2}\right] & \leq & 2\exp(-c^{''}n),\\
\mathbb{P}\left[\max_{j}\left\{ \left|\frac{1}{n}\sum_{i=1}^{n}v_{ij}\eta_{i}\right|\right\} \geq c_{0}\sigma_{v}\sigma_{\eta}\sqrt{\frac{\log p}{n}}\right] & \leq & 2\exp(-c^{'''}\log p).
\end{eqnarray*}
An upper bound on the term $\sqrt{\frac{1}{n}\sum_{i=1}^{n}\hat{\Delta}_{j}^{2}(w_{i}^{T}\hat{\theta})}$
is derived in Lemma A.11, which yields 
\[
\max_{j}\left\{ \sqrt{\frac{1}{n}\sum_{i=1}^{n}\hat{\Delta}_{j}^{2}(w_{i}^{T}\hat{\theta})}\right\} \leq c\max_{j}\left\{ \sqrt{T_{j}^{'}}+\left(\sigma_{v}^{2}T_{j}^{'}\right)^{\frac{1}{4}}+\left|\tilde{m}_{j}(\hat{\theta})-m_{j}(\hat{\theta})\right|_{n}+r_{nj}\right\} ,
\]
with probability at least $1-c^{'}\sum_{j=0}^{p}\exp\left(-nC_{j}^{*}r_{nj}^{2}\right)$.
An upper bound on the term $\sqrt{\frac{1}{n}\sum_{i=1}^{n}T_{ij}^{''2}}$
is simply $\max_{j}\sqrt{T_{j}^{'}}$ by recalling the definition
of $T_{j}^{'}$. 

To upper bound the term $\left|\frac{1}{n}\sum_{i=1}^{n}v_{ij}\hat{\Delta}_{j^{'}}(w_{i}^{T}\theta^{*})\right|$,
recall from the proof for Lemma A.3, we have

\[
\max_{j,\, j^{'}}\left|\frac{1}{n}\sum_{i=1}^{n}v_{ij}\hat{\Delta}_{j^{'}}(w_{i}^{T}\theta^{*})\right|\leq c\max_{j^{'}}\left\{ \bar{r}_{nj^{'}}^{2}+\sigma_{v}\left|\tilde{m}_{j^{'}}(\theta^{*})-m_{j^{'}}(\theta^{*})\right|_{n}+\sqrt{\sigma_{v}^{2}T_{j^{'}}^{'}}\right\} 
\]
with probability at least $1-c^{'}\sum_{j=1}^{p}\exp\left(-nC_{j}^{*}\bar{r}_{nj}^{2}\right)$.
Applying the same argument from above to control the term $\left|\frac{1}{n}\sum_{i=1}^{n}\eta_{i}\hat{\Delta}_{ij^{'}}(\theta^{*})\right|$
yields 
\[
\max_{j^{'}}\left|\frac{1}{n}\sum_{i=1}^{n}\eta_{i}\hat{\Delta}_{j^{'}}(w_{i}^{T}\theta^{*})\right|\leq c^{'}\max_{j^{'}}\left\{ \bar{r}_{nj}^{2}+\sigma_{\eta}\left|\tilde{m}_{j^{'}}(\theta^{*})-m_{j^{'}}(\theta^{*})\right|_{n}+\sqrt{\sigma_{\eta}^{2}T_{j}^{'}}\right\} 
\]
with probability at least $1-c^{''}\sum_{j=1}^{p}\exp\left(-nC_{j}^{*}\bar{r}_{nj}^{2}\right)$. 

Putting all the pieces together, we have 
\[
\left|\frac{\hat{v}^{T}e}{n}\right|_{\infty}\leq\left|\beta^{*}\right|_{1}b(\sigma_{v},\,\sigma_{\eta})\left(\mathcal{T}_{1}+\mathcal{T}_{2}+\mathcal{T}_{3}\right)+\mathcal{T}_{4}
\]
with probability at least $1-c_{1}\exp\left(-c_{2}\log p\right)-c_{3}\sum_{j=1}^{p}\exp\left(-nC_{j}^{*}r_{nj}^{2}\right)$.
$\square$

Now, by applying Lemma A.1 and setting $\lambda_{n,3}\geq\left|\beta^{*}\right|_{1}b(\sigma_{v},\,\sigma_{\eta})\left(\mathcal{T}_{1}+\mathcal{T}_{2}+\mathcal{T}_{3}\right)+\mathcal{T}_{4}$
and combining with Lemma A.3, we obtain 
\[
|\hat{\beta}_{HSEL}-\beta^{*}|_{2}\leq\frac{c^{''}R_{q_{2}}^{\frac{1}{2}}}{\kappa_{L}^{1-\frac{q_{2}}{2}}}\left[\left|\beta^{*}\right|_{1}b(\sigma_{v},\,\sigma_{\eta})\left(\mathcal{T}_{1}+\mathcal{T}_{2}+\mathcal{T}_{3}\right)+\mathcal{T}_{4}\right]^{1-\frac{q_{2}}{2}},\quad q_{2}\in[0,\,1]
\]
with probability at least $1-c_{1}\exp\left(-c_{2}\log p\right)-c_{3}\sum_{j=1}^{p}\exp\left(-nC_{j}^{*}r_{nj}^{2}\right)$.

\subsection{Theorem 4.2 }

As in the proofs for Theorem 4.1, the proofs for Lemmas A.5 and A.6
work through the case $t_{nj}=r_{nj}$. The case where $t_{nj}\geq r_{nj}$
can be analyzed with almost the same argument because for any star-shaped
function class $\bar{\mathcal{F}}_{j}:=\left\{ f=f^{'}-f^{''}\,:\, f^{'},\, f^{''}\in\mathcal{F}_{j}\right\} $,
the function $t\mapsto\frac{\mathcal{G}_{n}(t;\,\mathcal{F}_{j})}{t}$
is non-decreasing on $(0,\,\infty)$. For the proofs for Lemmas A.5
and A.6, let the \textit{local complexity} measure $\mathcal{G}_{n}\left(r_{j};\,\mathcal{F}_{j}\right)$
be defined over the set 
\[
\Omega(r_{j};\,\mathcal{F}_{j})=\left\{ f\,:\, f\in\bar{\mathcal{F}}_{j}\,\left|f_{\theta^{*}}\right|_{n}\leq r_{j},\,\left|f\right|_{\mathcal{F}_{j}}\leq1\right\} 
\]
where $\left|f_{\theta^{*}}\right|_{n}:=\sqrt{\frac{1}{n}\sum_{i=1}^{n}\left[f(w_{i}^{T}\theta^{*})\right]^{2}}$
and $j=0,...,p$. Recall the following definitions for Theorem 4.2:
\begin{eqnarray*}
\mathcal{T}_{1} & = & \max_{j\in\{0,...,p\}}\left(T_{j}^{'}\vee\sqrt{T_{j}^{'}}\right)\\
\mathcal{T}_{2} & = & \max_{j\in\{0,...,p\}}\bar{R}_{j}^{2}t_{nj}^{2}\\
\mathcal{T}_{3} & = & \max_{j\in\{0,...,p\}}\inf_{\tilde{m}_{j}\in\mathcal{F}_{j},\left|\tilde{m}_{j}\right|_{\mathcal{F}_{j}}\leq\bar{R}_{j}}\left(\left|\tilde{m}_{j}(\hat{\theta})-m_{j}(\hat{\theta})\right|_{n}^{2}+\left|\tilde{m}_{j}(\theta^{*})-m_{j}(\theta^{*})\right|_{n}\right)\\
\mathcal{T}_{4} & = & \sigma_{v}\sigma_{\eta}\sqrt{\frac{\log p}{n}}.
\end{eqnarray*}
\textbf{Lemma A.5} (RE condition): Let $r_{nj}>0$ be the smallest
positive quantity satisfying the \textit{critical inequality} 
\[
\mathcal{G}_{n}\left(r_{nj};\,\mathcal{F}_{j}\right)\leq\frac{\bar{R}_{j}r_{nj}^{2}}{\sigma_{v_{j}}},
\]
where $\bar{R}_{j}>0$ is a user-defined radius. Under Assumptions
2.1, 2.3, 4.1-4.3, if we choose the second-stage regularization parameter
$\lambda_{nj,2}=(2+\varsigma)r_{nj}^{2}$ for $j=1,...,p$ where $\varsigma$
is a small positive constant, then 
\begin{equation}
\left|\Delta^{T}\frac{\hat{v}^{T}\hat{v}}{n}\Delta\right|\geq\kappa_{1}|\Delta|_{2}^{2}-\kappa_{2}\frac{\log p}{n}|\Delta|_{1}^{2}-c\left(\mathcal{T}_{1}+\mathcal{T}_{2}+\mathcal{T}_{3}\right)\left|\Delta\right|_{1}^{2}\label{eq:9}
\end{equation}
with probability at least $1-c^{'}\sum_{j=1}^{p}\exp\left(-nC_{j}^{*}\bar{R}_{j}^{2}r_{nj}^{2}\right)$\\
\\
\textbf{Proof.} The proof is almost identical to that of Lemma A.3
except that in upper bounding the term $\frac{1}{n}\sum_{i=1}^{n}\left[\hat{m}_{j}(w_{i}^{T}\hat{\theta})-m_{j}(w_{i}^{T}\hat{\theta})\right]{}^{2}$,
we apply Lemma A.12 with the second-stage regularization parameter
$\lambda_{nj,2}=(2+\varsigma)r_{nj}^{2}$ where $\varsigma$ is a
small positive constant for $j=1,...,p$, and obtain 
\[
\frac{1}{n}\sum_{i=1}^{n}\left[\hat{m}_{j}(w_{i}^{T}\hat{\theta})-m_{j}(w_{i}^{T}\hat{\theta})\right]^{2}\leq c_{0}\left|\tilde{m}_{j}(\hat{\theta})-m_{j}(\hat{\theta})\right|_{n}^{2}+c_{1}\bar{R}_{j}^{2}r_{nj}^{2}+c_{2}T_{j}^{'}+c_{3}\sqrt{\sigma_{v}^{2}T_{j}^{'}}
\]
with probability at least $1-c\exp\left(-nC_{j}^{*}\bar{R}_{j}^{2}r_{nj}^{2}\right)$.
The second-stage regularization parameter $\lambda_{nj,2}$ is chosen
such that $\lambda_{nj,2}\geq2r_{nj}^{2}$ for $j=1,...,p$. Using
similar argument as in the proof for Lemma A.3 yields 
\begin{eqnarray*}
\left|\frac{v^{T}(\hat{v}-v)}{n}\right|_{\infty} & \leq\\
c\max_{j}\left\{ \inf_{\tilde{m}_{j}\in\mathcal{F}_{j},\,\left|\tilde{m}_{j}\right|_{\mathcal{F}_{j}}\leq\bar{R}_{j}}\left(\left|\tilde{m}_{j}(\hat{\theta})-m_{j}(\hat{\theta})\right|_{n}^{2}+\left|\tilde{m}_{j}(\theta^{*})-m_{j}(\theta^{*})\right|_{n}\right)+\bar{R}_{j}^{2}r_{nj}^{2}+T_{j}^{'}+\sqrt{\left(\sigma_{v}^{2}T_{j}^{'}\right)}\right\} 
\end{eqnarray*}
with probability at least $1-c_{1}\exp\left(-c_{2}\log p\right)-c_{3}\sum_{j=1}^{p}\exp\left(-nC_{j}^{*}\bar{R}_{j}^{2}r_{nj}^{2}\right)$. 

For the term $\left|\frac{(\hat{v}-v)^{T}(\hat{v}-v)}{n}\right|_{\infty}$,
we have 
\[
\left|\frac{(\hat{v}-v)^{T}(\hat{v}-v)}{n}\right|_{\infty}\leq c^{'}\max_{j}\left[\left|\tilde{m}_{j}(\hat{\theta})-m_{j}(\hat{\theta})\right|_{n}^{2}+\bar{R}_{j}^{2}r_{nj}^{2}+T_{j}^{'}+\sqrt{\sigma_{v}^{2}T_{j}^{'}}\right]
\]
with probability at least $1-c\sum_{j=1}^{p}\exp\left(-nC_{j}^{*}\bar{R}_{j}^{2}r_{nj}^{2}\right)$. 

Putting everything together and applying Lemma B.2, we have 
\begin{eqnarray}
\left|\Delta^{T}\frac{\hat{v}^{T}\hat{v}}{n}\Delta\right| & \geq & \kappa_{1}|\Delta|_{2}^{2}-\kappa_{2}\frac{\log p}{n}|\Delta|_{1}^{2}\nonumber \\
 & - & b(\sigma_{v})\max_{j}\left\{ \bar{R}_{j}^{2}r_{nj}^{2}+\inf_{\tilde{m}_{j}\in\mathcal{F}_{j}}\left(\left|\tilde{m}_{j}(\hat{\theta})-m_{j}(\hat{\theta})\right|_{n}^{2}+\left|\tilde{m}_{j}(\theta^{*})-m_{j}(\theta^{*})\right|_{n}\right)+\sqrt{T_{j}^{'}\vee T_{j}^{'2}}\right\} \left|\Delta\right|_{1}^{2}\nonumber \\
\label{eq:10}
\end{eqnarray}
where $\kappa_{1}=\frac{\kappa_{L}}{2}$, $\kappa_{2}=c_{0}\kappa_{L}\max\left\{ \frac{\sigma_{v}^{4}}{\kappa_{L}^{2}},\,1\right\} $,
and $\sigma_{v}=\max_{j}\sigma_{v_{j}}$, and a known function $b(\sigma_{v})$
depending only on $\sigma_{v}$. $\square$\\
\\
\textbf{Lemma A.6} (Upper bound on $\left|\frac{1}{n}\hat{v}^{T}e\right|_{\infty}$):
Let $r_{nj}>0$ be the smallest positive quantity satisfying the \textit{critical
inequality} 
\[
\mathcal{G}_{n}\left(r_{nj};\,\mathcal{F}_{j}\right)\leq\frac{\bar{R}_{j}r_{nj}^{2}}{\sigma_{v_{j}}},
\]
and $\bar{R}_{j}>0$ is a user-defined radius. Under Assumptions\uline{
}2.1 and 4.1-4.3 and if we choose the second-stage regularization
parameter $\lambda_{nj,2}=(2+\varsigma)r_{nj}^{2}$ where $\varsigma$
is a small positive constant for $j=1,...,p$, we have 
\[
\left|\frac{\hat{v}^{T}e}{n}\right|_{\infty}\leq\left|\beta^{*}\right|_{1}b(\sigma_{v},\,\sigma_{\eta})\left(\mathcal{T}_{1}+\mathcal{T}_{2}+\mathcal{T}_{3}\right)+\mathcal{T}_{4},
\]
with probability at least $1-c_{1}\exp\left(-c_{2}\log p\right)-c_{3}\sum_{j=1}^{p}\exp\left(-nC_{j}^{*}\bar{R}_{j}^{2}r_{nj}^{2}\right)$.\\
\\
\textbf{Proof.} The proof is almost identical to that of Lemma A.4
except that in upper bounding the term $\sqrt{\frac{1}{n}\sum_{i=1}^{n}\left[\hat{m}_{j}(w_{i}^{T}\hat{\theta})-m_{j}(w_{i}^{T}\hat{\theta})\right]{}^{2}}$,
we apply Lemma A.12 just as what we have done in the proof for Lemma
A.5. $\square$

Consequently, by applying Lemma A.1 and setting $\lambda_{n,3}\geq\left|\beta^{*}\right|_{1}b(\sigma_{v},\,\sigma_{\eta})\left(\mathcal{T}_{1}+\mathcal{T}_{2}+\mathcal{T}_{3}\right)+\mathcal{T}_{4}$
and combining with Lemma A.5, we obtain 
\[
|\hat{\beta}_{HSEL}-\beta^{*}|_{2}\leq\frac{c^{''}R_{q_{2}}^{\frac{1}{2}}}{\kappa_{L}^{1-\frac{q_{2}}{2}}}\left[\left|\beta^{*}\right|_{1}b(\sigma_{v},\,\sigma_{\eta})\left(\mathcal{T}_{1}+\mathcal{T}_{2}+\mathcal{T}_{3}\right)+\mathcal{T}_{4}\right]^{1-\frac{q_{2}}{2}},\quad q_{2}\in(0,\,1]
\]
with probability at least $1-c_{1}^{'}\exp\left(-c_{2}^{'}\log p\right)-c_{3}^{'}\sum_{j=1}^{p}\exp\left(-nC_{j}^{*}\bar{R}_{j}^{2}r_{nj}^{2}\right)$.

\subsection{Corollary 4.3}

\textbf{Proof.} This follows from the previous proofs for the RE condition
and Lemma A.1. $\square$

\subsection{Proposition 4.1 and Corollaries 4.4-4.6}

\textbf{Proof of Proposition 4.1}. Note that under the assumption
$m_{j}(\cdot)\in\mathcal{F}_{j}$ for $j=0,...,p$, we have $\inf\left|\tilde{m}_{j}(\hat{\theta})-m_{j}(\hat{\theta})\right|_{n}^{2}=0$.
When $\mathcal{F}_{j}$ is $L-$Lipschitz, under the Lipschitz regression
(22), 
\[
T_{j}^{'}=\frac{1}{n}\sum_{i=1}^{n}L^{2}\left[w_{i}^{T}\hat{\theta}-w_{i}^{T}\theta^{*}\right]^{2}:=L^{2}B^{'}
\]
for every $j=0,...,p$. Now let us compute the critical radius $r_{nj}$
that satisfies 
\[
\mathcal{G}_{n}\left(r_{nj};\,\mathcal{F}_{j}\right)\leq\frac{r_{nj}^{2}}{\sigma_{v_{j}}},
\]
where $\mathcal{F}_{j}$ is the space of $L-$Lipschitz functions.
The metric entropy of this class (in the sup-norm) scales as $\log N_{\infty}(t;\,\mathcal{\bar{F}}_{j})\asymp\frac{1}{t}$
by Lemma B.11. Consequently, applying Lemma B.9 yields 
\begin{eqnarray*}
\frac{1}{\sqrt{n}}\int_{0}^{r_{nj}}\sqrt{\log N_{n}(t;\,\mathbb{B}_{n}(r_{nj};\,\mathcal{H}_{j})}dt & \leq & \frac{1}{\sqrt{n}}\int_{0}^{r_{nj}}\sqrt{\log N_{\infty}(t;\,\textrm{sym}(\mathcal{H}_{j})}dt\\
 & \leq & \frac{c}{\sqrt{n}}\int_{0}^{r_{nj}}\left(\frac{|\theta^{*}|_{1}}{t}\right)^{\frac{1}{2}}dt\\
 & = & \frac{c}{\sqrt{n}}\sqrt{|\theta^{*}|_{1}r_{nj}}.
\end{eqnarray*}
Thus, it suffices to choose $r_{nj}>0$ such that 
\begin{eqnarray*}
\frac{\sqrt{r_{nj}}}{\sqrt{n}} & \asymp & \frac{r_{nj}^{2}}{\sigma_{v_{j}}}\\
\Longrightarrow r_{nj}^{2} & \asymp & \left(\frac{|\theta^{*}|_{1}\sigma_{v_{j}}^{2}}{n}\right)^{\frac{2}{3}}.
\end{eqnarray*}
\textbf{$\square$}\\
\textbf{}\\
\textbf{Proof of Corollaries 4.4-4.6}. In proving the two corollaries,
we follow the proof for Theorem 4.1. We prove the case where $\beta^{*}$
and $\theta^{*}$ are exactly sparse and the approximate sparse case
of $\beta^{*}$ and $\theta^{*}$ follows the same argument. We have
\[
\max_{j=0,...,p}T_{j}^{'}\leq\frac{1}{n}\sum_{i=1}^{n}L^{2}\left[w_{i}^{T}\hat{\theta}-w_{i}^{T}\theta^{*}\right]^{2}=L^{2}B^{'}
\]
Now, depending on the assumptions on $w_{i}$, we have 
\begin{eqnarray*}
B^{'} & = & \frac{1}{n}\sum_{i=1}^{n}\left[w_{i}^{T}\hat{\theta}-w_{i}^{T}\theta^{*}\right]^{2}\leq c\frac{\kappa_{U}^{w}\sigma_{w}^{2}\alpha_{u}}{(\kappa_{L}^{w})^{2}}\frac{k_{1}\log d}{n}\;\textrm{(Corollary 4.4)}\\
B^{'} & = & \frac{1}{n}\sum_{i=1}^{n}\left[w_{i}^{T}\hat{\theta}-w_{i}^{T}\theta^{*}\right]^{2}\leq c^{'}\frac{\left|\theta^{*}\right|_{1}}{\alpha_{l}}\sqrt{\frac{k_{1}\log d}{n}}\;\textrm{(Corollary 4.5)}
\end{eqnarray*}
with probability at least $1-O\left(\frac{1}{d}\right)$. Putting
the pieces together and applying Theorem 4.1, choosing $t_{nj}^{2}\geq r_{nj}^{2}$
in $\mathcal{T}_{2}$ such that $\left|\beta^{*}\right|_{1}\mathcal{T}_{2}$
is at most 
\[
O\left(\sqrt{\frac{\log p}{n}}\vee\left(\left|\beta^{*}\right|_{1}\sqrt{\frac{k_{1}\log d}{n}}\right)\right).
\]
and $\lambda_{n,3}$ in Corollaries 4.4 and 4.5, if $\sqrt{\frac{n\log p}{\left|\beta^{*}\right|_{1}^{2}}}\vee\sqrt{nk_{1}\log d}\succsim\log p$,
and 
\[
\kappa_{2}\frac{k_{2}\log p}{n}+ck_{2}\sigma_{v}\left[\left(\frac{k_{1}\sigma_{v}^{2}}{n}\right)^{\frac{2}{3}}\right]+k_{2}B^{'}+k_{2}\sqrt{\sigma_{v}^{2}B^{'}}=O(\kappa_{1}),
\]
then, with probability at least $1-O\left(\frac{1}{d\wedge p}\right)$,
the error satisfies the bound 
\begin{eqnarray*}
|\hat{\beta}_{HSEL}-\beta^{*}|_{2} & \leq & \frac{c_{1}\sqrt{k_{2}}}{\kappa_{L}}\left[\left(\sigma_{v}\sigma_{\eta}\sqrt{\frac{\log p}{n}}\right)\vee\left(L\sqrt{\frac{\kappa_{U}^{w}\sigma_{w}^{2}\alpha_{u}}{(\kappa_{L}^{w})^{2}}}b\left(\sigma_{v},\,\sigma_{\eta}\right)\left|\beta^{*}\right|_{1}\sqrt{\frac{k_{1}\log d}{n}}\right)\right]\\
 &  & \textrm{(Corollary 4.4)}
\end{eqnarray*}
and 
\begin{eqnarray*}
|\hat{\beta}_{HSEL}-\beta^{*}|_{2} & \leq & \frac{c_{2}\sqrt{k_{2}}}{\kappa_{L}}\left[\left(\sigma_{v}\sigma_{\eta}\sqrt{\frac{\log p}{n}}\right)\vee\left(L\sqrt{\frac{\left|\theta^{*}\right|_{1}}{\alpha_{l}}}b^{'}\left(\sigma_{v},\,\sigma_{\eta}\right)\left|\beta^{*}\right|_{1}\left(\frac{k_{1}\log d}{n}\right)^{\frac{1}{4}}\right)\right]\\
 &  & \textrm{(Corollary 4.5)}.
\end{eqnarray*}
\textbf{$\square$}

\subsection{Proposition 4.2}

\textbf{Proof.} First, by Corollary 5 in Negahban, et. al (2012) for
the logit model and Lemma A.7 for the probit model, we have, for $\lambda_{n,1}\asymp c\sigma_{w}\sqrt{\alpha_{u}}\sqrt{\frac{\log d}{n}}$,
\[
\left|\hat{\theta}-\theta^{*}\right|_{2}^{2}\leq c_{1}\frac{\sigma_{w}^{2}\alpha_{u}}{\left(\kappa_{L}^{w}\right)^{2}}\frac{k_{1}\log d}{n}
\]
with probability at least $1-O\left(\frac{1}{d}\right)$. Let $\Delta=\hat{\theta}-\theta^{*}$.
Under Assumption 4.6, Lemma B.2 implies 
\[
\frac{\left|w\Delta\right|{}_{2}^{2}}{n}\leq\frac{3\kappa_{U}^{w}}{2}\left|\Delta\right|{}_{2}^{2}+\alpha^{'}\frac{\log d}{n}\left|\Delta\right|{}_{1}^{2}.
\]
Combining the inequality above with Lemma 1 in Negahban, et. al which
shows that $\left|\Delta\right|_{1}\leq4\sqrt{k_{1}}\left|\Delta\right|_{2}$
yields the desired result. $\square$\\
\textbf{}\\
\textbf{Lemma A.7}: Suppose the number of regressors $d(=d_{n})$
can grow with and exceed the sample size $n$ and\textbf{ }the number
of non-zero components in $\theta^{*}$ is at most $k_{1}(=k_{1n})$
and $k_{1}$ can increase to infinity with $n$ but slowly compared
to $n$. Let Assumptions 2.1 and 4.6 (the lower bound part) hold.
If $\hat{\theta}$ solves program (13) with $\lambda_{n,1}\geq c\sigma_{w}\sqrt{\alpha_{u}}\sqrt{\frac{\log d}{n}}$
and $n\succsim k_{1}\log d$, then, 
\[
\left|\hat{\theta}-\theta^{*}\right|_{2}\leq c\frac{\sqrt{k_{1}}}{\kappa_{L}^{w}}\left(\sqrt{\alpha_{u}\sigma_{w}^{2}}\sqrt{\frac{\log d}{n}}\vee\lambda_{n,1}\right)
\]
with probability at least $1-O\left(\frac{1}{d}\right)$.\\
\textbf{}\\
\textbf{Proof.} Recall 
\[
L_{n}(\theta)=\left\{ -\frac{1}{n}\sum_{i=1}^{n}y_{1i}\phi_{1}(w_{i}^{T}\theta)+\frac{1}{n}\sum_{i=1}^{n}\phi_{2}(w_{i}^{T}\theta)+\lambda_{n,1}\left|\theta\right|_{1}\right\} ,
\]
where 
\begin{eqnarray*}
\phi_{1}(w_{i}^{T}\theta) & = & \log\frac{\Phi(w_{i}^{T}\theta)}{1-\Phi(w_{i}^{T}\theta)},\\
\phi_{2}(w_{i}^{T}\theta) & = & -\log\left[1-\Phi(w_{i}^{T}\theta)\right]=-\log\frac{1}{1+\exp\left(\phi_{1}(w_{i}^{T}\theta)\right)}
\end{eqnarray*}
for the probit model. Define the following quantity 
\[
\delta L_{n}(\theta^{*},\,\Delta):=L_{n}(\theta^{*}+\Delta)-L_{n}(\theta^{*})-\left\langle \nabla L_{n}(\theta^{*}),\,\Delta\right\rangle .
\]
To prove Lemma A.7, by Theorem 1 in Negahban, et. al (2012), it suffices
to show that Step 1: $\delta L_{n}(\theta^{*},\,\hat{\Delta})$ satisfies
the RSC condition (2.4) where $\hat{\Delta}=\hat{\theta}-\theta^{*}$;
and Step 2: with high probability, we have 
\[
\left|\nabla L_{n}(\theta^{*})\right|_{\infty}\leq\sqrt{\alpha_{u}\sigma_{w}^{2}\frac{\log d}{n}}.
\]
\textbf{Step 1}: We need to show that $\delta L_{n}(\theta^{*},\,\hat{\Delta})$
satisfies the RSC condition where $\hat{\Delta}=\hat{\theta}-\theta^{*}$.
Define $\rho_{i,\,\theta}=-y_{1i}\phi_{1}(w_{i}^{T}\theta)+\phi_{2}(w_{i}^{T}\theta)$,
$\hat{Q}_{1}(\theta)=\mathbb{E}_{n}\rho_{i,\,\theta}$, and $Q_{1}(\theta)=\mathbb{E}\rho_{i,\,\theta}$.
By the mean value theorem, we have 
\[
L_{n}(\theta^{*},\,\hat{\Delta})=\hat{Q}_{1}(\hat{\theta})-\hat{Q}_{1}(\theta^{*})=\frac{1}{n}\sum_{i=1}^{n}\rho_{i,\,\bar{\theta}}^{''}\left[w_{i}^{T}(\hat{\theta}-\theta^{*})\right]^{2},
\]
where $\rho_{i,\,\bar{\theta}}^{''}$ denotes the second derivative
of $\rho_{i,\,\bar{\theta}}$ and $\bar{\theta}$ is some intermediate
value. 

Proposition 1 in Loh and Wainwright (2013) implies that there exists
a constant $\underbar{\ensuremath{\alpha}}>0$, depending on the probit
model and $(\sigma_{w},\,\Sigma_{w})$, such that for all vectors
$\theta_{2}\in\mathcal{B}_{2}(3)\cap\mathcal{B}_{1}(R)$ 
\[
\hat{Q}_{1}(\theta_{2})-\hat{Q}_{1}(\theta_{1})\geq\begin{cases}
\begin{array}{c}
\frac{\underbar{\ensuremath{\alpha}}}{2}\left|\Delta\right|_{2}^{2}-\frac{c^{2}\sigma_{w}^{2}}{2\underbar{\ensuremath{\alpha}}}\frac{\log p}{n}\left|\Delta\right|_{1}^{2}\quad\textrm{for all}\:\left|\Delta\right|_{2}\leq3,\\
\frac{3\underbar{\ensuremath{\alpha}}}{2}\left|\Delta\right|_{2}-3c\sigma_{w}\sqrt{\frac{\log p}{n}}\left|\Delta\right|_{1}\quad\textrm{for all}\:\left|\Delta\right|_{2}\geq3,
\end{array}\end{cases}
\]
with probability at least $1-c_{1}\exp(-c_{2}n)$ and $\Delta=\theta_{2}-\theta_{1}$.
Combining this result with Lemma 1 in Negahban, et. al (2012) which
shows that $\left|\hat{\Delta}\right|_{1}\leq4\sqrt{k_{1}}\left|\hat{\Delta}\right|_{2}$
yields the RSC of $\delta L_{n}(\theta^{*},\,\hat{\Delta})$.\\
\\
\textbf{Step 2}: In the following, we show that $\left|\nabla L_{n}(\theta^{*})\right|_{\infty}\leq\sqrt{\alpha_{u}\sigma_{w}^{2}\frac{\log d}{n}}$
with high probability. For a fixed index $j=1,...,d$, we begin by
establishing an upper bound on 
\[
\frac{1}{n}\sum_{i=1}^{n}V_{ij}:=\frac{1}{n}\sum_{i=1}^{n}\left[w_{ij}y_{1i}\phi_{1}^{'}(w_{i}^{T}\theta^{*})-w_{ij}\phi_{2}^{'}(w_{i}^{T}\theta^{*})\right]
\]
where $\phi_{1}^{'}(u)=\frac{\partial\phi_{1}}{\partial u}$ and $\phi_{2}^{'}(u)=\frac{\partial\phi_{2}}{\partial\phi_{1}}\frac{\partial\phi_{1}}{\partial u}$.
Let us condition on $\{w_{i}\}_{i=1}^{n}$, so that $y_{1i}$ is drawn
from the exponential family with parameter $\phi_{1}(w_{i}^{T}\theta)=\log(\frac{\Phi(w_{i}^{T}\theta)}{1-\Phi(w_{i}^{T}\theta)})$
for the probit model. For any $t\in\mathbb{R}$, compute the cumulant
function 
\begin{eqnarray*}
 &  & \log\mathbb{E}\left[\exp(tV_{ij})\,\vert\, w_{i}\right]\\
 & = & \log\left\{ \mathbb{E}\left[\exp\left(tw_{ij}y_{1i}\phi_{1}^{'}(w_{i}^{T}\theta^{*})\right)\,|\, w_{i}\right]\exp\left(-tw_{ij}\phi_{2}^{'}(w_{i}^{T}\theta^{*})\right)\right\} \\
 & = & \phi_{2}\left(\phi_{1}^{'}(w_{i}^{T}\theta^{*})tw_{ij}+\phi_{1}(w_{i}^{T}\theta^{*})\right)-\phi_{2}(\phi_{1}(w_{i}^{T}\theta^{*}))-\phi_{2}^{'}(w_{i}^{T}\theta^{*}))(tw_{ij})\\
 & = & \frac{t^{2}}{2}w_{ij}^{2}\left[\phi_{2}^{''}\left(\phi_{1}(w_{i}^{T}\theta^{*})+\upsilon_{i}\phi_{1}^{'}(w_{i}^{T}\theta^{*})tw_{ij}\right)\right]
\end{eqnarray*}
for some $\upsilon_{i}\in[0,\,1]$. Since this upper bound holds for
each $i=1,...,n$, we obtain 
\begin{eqnarray*}
\frac{1}{n}\sum_{i=1}^{n}\log\mathbb{E}\left[\exp(tV_{ij})\,\vert\, w_{i}\right] & \leq & \frac{t^{2}}{2}\left\{ \frac{1}{n}\sum_{i=1}^{n}w_{ij}^{2}\left[\phi_{2}^{''}\left(\phi_{1}(w_{i}^{T}\theta^{*})+\upsilon_{i}\phi_{1}^{'}(w_{i}^{T}\theta^{*})tw_{ij}\right)\right]\right\} \\
 & \leq & \frac{t^{2}\alpha_{u}}{2}\left\{ \frac{1}{n}\sum_{i=1}^{n}w_{ij}^{2}\right\} 
\end{eqnarray*}
where the inequality follows since $\phi_{2}^{''}(u)\leq\alpha_{u}$
for some $\alpha_{u}>0$ and all $u\in\mathbb{R}$. Note that for
each $j=1,...,d$, the variables $\left\{ w_{ij}^{2}-\mathbb{E}\left[w_{ij}^{2}\right]\right\} _{i=1}^{n}$
are \textit{i.i.d.} zero-mean and sub-Gaussian with parameter at most
$\sigma_{w}$. Consequently, we have $\mathbb{E}\left[w_{ij}^{2}\right]\leq\sigma_{w}^{2}$.
Since the squared variables are sub-exponential, by Lemma B.1, we
have the tail bound 
\[
\mathbb{P}\left[\frac{1}{n}\sum_{i=1}^{n}w_{ij}^{2}\geq2\sigma_{w}^{2}\right]\leq2\exp\left(-c_{1}n\right).
\]
Define the event $\Xi=\left\{ \max_{j=1,...,d}\frac{1}{n}\sum_{i=1}^{n}w_{ij}^{2}\leq2\sigma_{w}^{2}\right\} $.
An application of union bound yields 
\[
\mathbb{P}\left[\Xi^{c}\right]\leq2\exp\left(-c_{1}n+\log d\right)\leq2\exp(-c_{1}n)
\]
under the condition $n\succsim\log d$. Conditioning on the event
$\Xi$, we obtain 
\[
\frac{1}{n}\sum_{i=1}^{n}\log\mathbb{E}\left[\exp(tV_{ij})\,\vert\, w_{i}\right]\leq t^{2}\alpha_{u}\sigma_{w}^{2},\quad\textrm{for each }j=1,...,d.
\]
An application of the Chernoff bound and the union bound yields 
\[
\mathbb{P}\left[\max_{j=1,...,d}\left|\frac{1}{n}\sum_{i=1}^{n}V_{ij}\right|\geq\varepsilon\,\vert\,\Xi\right]\leq2\exp\left(-c_{1}n\frac{\varepsilon^{2}}{\alpha_{u}\sigma_{w}^{2}}+\log d\right).
\]
Setting $\varepsilon=\sqrt{\alpha_{u}\sigma_{w}^{2}\frac{\log d}{n}}$
and combining with the bound on $\mathbb{P}\left[\Xi^{c}\right]$
yields 
\begin{eqnarray*}
\mathbb{P}\left[\max_{j=1,...,d}\left|\frac{1}{n}\sum_{i=1}^{n}V_{ij}\right|\geq\sqrt{\alpha_{u}\sigma_{w}^{2}\frac{\log d}{n}}\right] & \leq & \mathbb{P}\left[\Xi^{c}\right]+\mathbb{P}\left[\max_{j=1,...,d}\left|\frac{1}{n}\sum_{i=1}^{n}V_{ij}\right|\geq\varepsilon\,\vert\,\Xi\right]\\
 & \leq & c_{2}\exp\left(-c_{3}n\frac{\varepsilon^{2}}{\alpha_{u}\sigma_{w}^{2}}+\log d\right)\leq c_{4}\exp\left(-c_{5}\log d\right),
\end{eqnarray*}
where the last inequality follows since $p>n$ is the regime of our
interest. 

Now put the pieces together, we obtain 
\[
\left|\nabla L_{n}(\theta^{*})\right|_{\infty}\leq\sqrt{\alpha_{u}\sigma_{w}^{2}\frac{\log d}{n}}
\]
with probability at least $1-c^{'}\exp\left(-c^{''}\log d\right)$. 

Combining Step 1 and Step 2 and applying Theorem 1 in Negahban, et.
al (2012) yields the desired result. $\square$\\
\\
\textbf{Remark}. Theorem 1 in Negahban, et. al (2012) can be applied
to obtain the upper bound in Proposition 4.2 for the case of approximate
sparsity $q_{1}\in(0,\,1]$.

\subsection{Proposition 4.3}

\textbf{Proof. }Recall from the proof for Proposition 4.2: 
\[
\hat{Q}_{1}(\hat{\theta})-\hat{Q}_{1}(\theta^{*})=\frac{1}{n}\sum_{i=1}^{n}\rho_{i,\,\bar{\theta}}^{''}\left[w_{i}^{T}(\hat{\theta}-\theta^{*})\right]^{2},
\]
where $\rho_{i,\,\bar{\theta}}^{''}$ denotes the second derivative
of $\rho_{i,\,\bar{\theta}}$ and $\bar{\theta}$ is some intermediate
value. Assumption 4.4 implies $\rho_{i,\,\bar{\theta}}^{''}\geq\alpha_{l}>0$,
so we further have 
\begin{equation}
\frac{1}{n}\sum_{i=1}^{n}\left[w_{i}^{T}(\hat{\theta}-\theta^{*})\right]^{2}\leq\frac{1}{\alpha_{l}}\left[\hat{Q}_{1}(\hat{\theta})-\hat{Q}_{1}(\theta^{*})\right].
\end{equation}
Hence, it remains to upper bound $\hat{Q}_{1}(\hat{\theta})-\hat{Q}_{1}(\theta^{*})$.
Write 
\begin{eqnarray}
\hat{Q}_{1}(\hat{\theta})-\hat{Q}_{1}(\theta^{*}) & = & \hat{Q}_{1}(\hat{\theta})-Q_{1}(\hat{\theta})-\left[\hat{Q}_{1}(\theta^{*})-Q_{1}(\theta^{*})\right]+\left[Q_{1}(\hat{\theta})-Q_{1}(\theta^{*})\right]\nonumber \\
 & := & \Lambda_{\hat{\theta}}-\Lambda_{\theta^{*}}+E_{\hat{\theta}},\label{eq:12}
\end{eqnarray}
where we define $\hat{Q}_{1}(\hat{\theta})-Q_{1}(\hat{\theta})-\left[\hat{Q}_{1}(\theta^{*})-Q_{1}(\theta^{*})\right]:=\Lambda_{\hat{\theta}}-\Lambda_{\theta^{*}}$
and the excess risk $Q_{1}(\hat{\theta})-Q_{1}(\theta^{*}):=E_{\hat{\theta}}$.
Furthermore, by the definition of $\hat{\theta}$, we have $\hat{Q}_{1}(\hat{\theta})+\lambda_{n,1}\left|\hat{\theta}\right|_{1}\geq\hat{Q}_{1}(\theta^{*})+\lambda_{n,1}\left|\theta^{*}\right|_{1}$,
which yields the following \textit{basic inequality} 
\[
E_{\hat{\theta}}+\lambda_{n,1}\left|\hat{\theta}\right|_{1}\leq-\left[\Lambda_{\hat{\theta}}-\Lambda_{\theta^{*}}\right]+\lambda_{n,1}\left|\theta^{*}\right|_{1}+E_{\theta^{*}}
\]
Under our model setup, $E_{\theta^{*}}=0$. This implies 
\[
E_{\hat{\theta}}\leq\left|\Lambda_{\hat{\theta}}-\Lambda_{\theta^{*}}\right|+\lambda_{n,1}\left|\hat{\theta}-\theta^{*}\right|_{1},
\]
by the triangle inequality. To control for $\Lambda_{\hat{\theta}}-\Lambda_{\theta^{*}}$
and $E_{\hat{\theta}}$, we follow Lemma 6.7 in Bühlmann and Van de
Geer (2011) (alternatively, see Theorems 2.1 or 2.2 in Van de Geer,
2008) to obtain the desired upper bound on $\hat{Q}_{1}(\hat{\theta})-\hat{Q}_{1}(\theta^{*})$
in (\ref{eq:12}).$\square$

\subsection{Theorem 4.7}

The proof for Theorem 4.7 hinges on an intermediate result that shows
the {}``mutual incoherence'' assumption on $\mathbb{E}[v_{1}^{T}v_{1}]$
(the population version of $\frac{1}{n}v^{T}v$) guarantees that,
with high probability, analogous conditions hold for the estimated
quantity $\frac{1}{n}\hat{v}^{T}\hat{v}$. This result is established
in Lemma A.13. The main proof for Theorem 4.7 is based on a construction
called Primal-Dual Witness (PDW) method developed by Wainwright (2009).
This method constructs a pair $(\hat{\beta},\,\hat{\mu})$. When this
procedure succeeds, the constructed pair is primal-dual optimal, and
acts as a witness for the fact that the Lasso has a unique optimal
solution with the correct signed support. The procedure is described
in the following. 
\begin{enumerate}
\item Set $\hat{\beta}_{J(\beta^{*})^{c}}=0$.
\item Obtain $(\hat{\beta}_{J(\beta^{*})},\,\hat{\mu}_{J(\beta^{*})})$
by solving the oracle subproblem 
\[
\hat{\beta}_{J(\beta^{*})}\in\arg\min_{\beta_{J(\beta^{*})}\in\mathbb{R}^{k_{2}}}\{\frac{1}{2n}|\hat{v}_{0}-\hat{v}_{J(\beta^{*})}\beta_{J(\beta^{*})}|_{2}^{2}+\lambda_{n,3}|\beta_{J(\beta^{*})}|_{1}\},
\]
and choose $\hat{\mu}_{J(\beta^{*})}\in\partial|\hat{\beta}_{J(\beta^{*})}|{}_{1}$,
where $\partial|\hat{\beta}_{J(\beta^{*})}|_{1}$ denotes the set
of subgradients at $\hat{\beta}_{J(\beta^{*})}$ for the function
$|\cdot|_{1}:\,\mathbb{R}^{k_{2}}\rightarrow\mathbb{R}$.
\item Solve for $\hat{\mu}_{J(\beta^{*})^{c}}$ via the zero-subgradient
equation 
\[
\frac{1}{n}\hat{v}^{T}(\hat{v}_{0}-\hat{v}\hat{\beta})+\lambda_{n,3}\hat{\mu}=0,
\]
and check whether or not the \textit{strict dual feasibility} condition
$|\hat{\mu}_{J(\beta^{*})^{c}}|_{\infty}<1$ holds.
\end{enumerate}
The proof includes four parts. Lemma A.8 guarantees the uniqueness
of the optimal solution of the two-stage Lasso procedure, $\hat{\beta}_{HSEL}$,
and shows that $\hat{\beta}_{HSEL}=(\hat{\beta}_{J(\beta^{*})},\,\mathbf{0})$
where $\hat{\beta}_{J(\beta^{*})}$ is the solution obtained in step
2 of the PDW construction. Based on this uniqueness claim, one can
then talk unambiguously about the support of the two-stage Lasso estimate.
Lemma A.9 proves part (a) of Theorem 4.7 by verifying the \textit{strict
dual feasibility} condition in step 3 of PDW and ensures that $\hat{\beta}_{HSEL,\, J(\beta^{*})}$
is uniformly close to $\beta_{J(\beta^{*})}^{*}$ in the $l_{\infty}-$norm.
Part (b) of Theorem 4.7 is a consequence of this uniform norm bound
in Lemma A.9: as long as the minimum value of $|\beta_{j}^{*}|$ over
$j\in J(\beta^{*})$ is not too small, then the multi-stage estimator
does not falsely exclude elements that are in the support of $\beta^{*}$
with high probability. To simplify the notations in the following
analysis, let $J(\beta^{*})=K$, $J(\beta^{*})^{c}:=K^{c}$, $\Sigma_{K^{c}K}:=\mathbb{E}\left[v_{1,K^{c}}^{T}v_{1,K}\right]$,
$\hat{\Sigma}_{K^{c}K}:=\frac{1}{n}v_{K^{c}}^{T}v_{K}$, and $\check{\Sigma}_{K^{c}K}:=\frac{1}{n}\hat{v}_{K^{c}}^{T}\hat{v}_{K}$.
Similarly, let $\Sigma_{KK}:=\mathbb{E}\left[v_{1,K}^{T}v_{1,K}\right]$,
$\hat{\Sigma}_{KK}:=\frac{1}{n}v_{K}^{T}v_{K}$, and $\check{\Sigma}_{KK}:=\frac{1}{n}\hat{v}_{K}^{T}\hat{v}_{K}$.
\textbf{}\\
\textbf{}\\
\textbf{Lemma A.8}: If the PDW construction succeeds and if $\lambda_{\min}\left(\mathbb{E}\left[v_{1,J(\beta^{*})}^{T}v_{1,J(\beta^{*})}\right]\right)\geq C_{\min}>0$,
then the vector $(\hat{\beta}_{K},\,\mathbf{0})\in\mathbb{R}^{p}$
is the unique optimal solution of the Lasso.\\
\\
\textbf{Proof}. The proof for Lemma A.8 adopts the proof for Lemma
1 from Chapter 6.4.2 of Wainwright (2015). If the PDW construction
succeeds, then $\hat{\beta}=(\hat{\beta}_{K},\,\mathbf{0})$ is an
optimal solution with associated subgradient vector $\hat{\mu}\in\mathbb{R}^{p}$
satisfying $|\hat{\mu}_{K^{c}}|_{\infty}<1$, and $\left\langle \hat{\mu},\,\hat{\beta}\right\rangle =|\hat{\beta}|_{1}$.
Suppose $\tilde{\beta}$ is another optimal solution. Letting $F(\beta)=\frac{1}{2n}|\hat{v}_{0}-\hat{v}\beta|_{2}^{2}$,
then \textbf{$F(\hat{\beta})+\lambda_{n,3}\left\langle \hat{\mu},\,\hat{\beta}\right\rangle =F(\tilde{\beta})+\lambda_{n,3}|\tilde{\beta}|_{1}$},
and hence\textbf{ $F(\hat{\beta})-\lambda_{n,3}\left\langle \hat{\mu},\,\tilde{\beta}-\hat{\beta}\right\rangle =F(\tilde{\beta})+\lambda_{n,3}\left(|\tilde{\beta}|_{1}-\left\langle \hat{\mu},\,\tilde{\beta}\right\rangle \right)$}.
However, by the zero-subgradient%
\footnote{Given a convex function $g:\,\mathbb{R}^{p}\mapsto\mathbb{R}$, $\mu\in\mathbb{R}^{p}$
is a subgradient at $\beta$, denoted by $\mu\in\partial g(\beta)$,
if $g(\beta+\triangle)\geq g(\beta)+\left\langle \mu,\,\triangle\right\rangle $
for all $\triangle\in\mathbb{R}^{p}$. When $g(\beta)=|\beta|_{1}$,
notice that $\mu\in\partial|\beta|_{1}$ if and only if $\mu_{j}=\textrm{sign}(\beta_{j})$
for all $j=1,...,p$, where $\textrm{sign}(0)$ is allowed to be any
number in $[-1,\,1]$.%
} conditions for optimality, we have \textbf{$\lambda_{n,3}\hat{\mu}=-\nabla F(\hat{\beta})$},
which implies that \textbf{$F(\hat{\beta})+\left\langle \nabla F(\hat{\beta}),\,\tilde{\beta}-\hat{\beta}\right\rangle -F(\tilde{\beta})=\lambda_{n,3}\left(|\tilde{\beta}|_{1}-\left\langle \hat{\mu},\,\tilde{\beta}\right\rangle \right)$}.\textbf{
}By convexity of \textbf{$F$}, the left-hand side is non-positive,
which implies that $|\tilde{\beta}|_{1}\leq\left\langle \hat{\mu},\,\tilde{\beta}\right\rangle $.
But since we also have \textbf{$\left\langle \hat{\mu},\,\tilde{\beta}\right\rangle \leq|\hat{\mu}|_{\infty}|\tilde{\beta}|_{1}$},
we must have\textbf{ $|\tilde{\beta}|_{1}=\left\langle \hat{\mu},\,\tilde{\beta}\right\rangle $}.
Since\textbf{ $|\hat{\mu}_{K^{c}}|_{\infty}<1$}, this equality can
only occur if $\tilde{\beta}_{j}=0$ for all\textbf{ $j\in K^{c}$}.
Thus, all optimal solutions are supported only on $K$, and hence
can be obtained by solving the oracle subproblem in the PDW procedure.
Given $\lambda_{\min}\left(\mathbb{E}\left[v_{1,J(\beta^{*})}^{T}v_{1,J(\beta^{*})}\right]\right)\geq C_{\min}>0$,
this subproblem is strictly convex, and hence it has a unique minimizer.\textbf{
$\square$}\\
\\
\textbf{Lemma A.9}: Suppose the assumptions in Theorem 4.7 hold. With
the choice of the regularization parameter 
\[
\lambda_{n,3}\geq c\frac{8(2-\frac{\phi}{4})}{\phi}\left[\left(\sigma_{v}\sigma_{\eta}\sqrt{\frac{\log p}{n}}\right)\vee\left(L\left|\beta^{*}\right|_{1}\sqrt{\Upsilon_{w,\theta^{*}}}b\left(\sigma_{v},\,\sigma_{\eta}\right)\sqrt{\frac{k_{1}\log d}{n}}\right)\right]
\]
and under the condition $n\succsim\left(k_{2}^{3}\log p\right)\vee\left(k_{2}^{2}k_{1}\log d\right)$,
$\sqrt{\frac{n\log p}{\left|\beta^{*}\right|_{1}^{2}}}\vee\sqrt{nk_{1}\log d}\succsim\log p$,
and $\sqrt{\frac{k_{1}\log d}{n}}=o(1)$, we have $|\hat{\mu}_{K^{c}}|_{\infty}\leq1-\frac{\phi}{8}$
with probability at least $1-c_{1}\exp(-c_{2}\log(p\wedge d))$. Furthermore,
\[
|\hat{\beta}_{K}-\beta_{K}^{*}|_{\infty}\leq\frac{c\sqrt{k_{2}}}{\lambda_{\min}(\Sigma_{KK})}\left[\left(\sigma_{v}\sigma_{\eta}\sqrt{\frac{\log p}{n}}\right)\vee\left(L\left|\beta^{*}\right|_{1}\sqrt{\Upsilon_{w,\theta^{*}}}b\left(\sigma_{v},\,\sigma_{\eta}\right)\sqrt{\frac{k_{1}\log d}{n}}\right)\right],
\]
with probability at least $1-c_{1}\exp(-c_{2}\log(p\wedge d))$. \\
\\
\textbf{Proof}. By construction, the sub-vectors $\hat{\beta}_{K}$,
$\hat{\mu}_{K}$, and $\hat{\mu}_{K^{c}}$ satisfy the zero-subgradient
condition in the PDW construction. Recall the definition of $e$ from
Lemma A.1. With the fact that $\hat{\beta}_{K^{c}}=\beta_{K^{c}}^{*}=0$,
we have 
\begin{eqnarray*}
\frac{1}{n}\hat{v}_{K}^{T}\hat{v}_{K}\left(\hat{\beta}_{K}-\beta_{K}^{*}\right)+\frac{1}{n}\hat{v}_{K}^{T}e+\lambda_{n,3}\hat{\mu}_{K} & = & 0,\\
\frac{1}{n}\hat{v}_{K^{c}}^{T}\hat{v}_{K}\left(\hat{\beta}_{K}-\beta_{K}^{*}\right)+\frac{1}{n}\hat{v}_{K^{c}}^{T}e+\lambda_{n,3}\hat{\mu}_{K^{c}} & = & 0.
\end{eqnarray*}
From the equations above, by solving for the vector $\hat{\mu}_{K^{c}}\in\mathbb{R}^{p-k_{2}}$,
we obtain 
\begin{eqnarray*}
\hat{\mu}_{K^{c}} & = & -\frac{1}{n\lambda_{n,3}}\hat{v}_{K^{c}}^{T}\hat{v}_{K}\left(\hat{\beta}_{K}-\beta_{K}^{*}\right)-\hat{v}_{K^{c}}^{T}\frac{e}{n\lambda_{n,3}},\\
\hat{\beta}_{K}-\beta_{K}^{*} & = & -\left(\frac{1}{n}\hat{v}_{K}^{T}\hat{v}_{K}\right)^{-1}\frac{\hat{v}_{K}^{T}e}{n}-\lambda_{n,3}\left(\frac{\hat{v}_{K}^{T}\hat{v}_{K}}{n}\right)^{-1}\hat{\mu}_{K},
\end{eqnarray*}
which yields 
\[
\hat{\mu}_{K^{c}}=\left(\check{\Sigma}_{K^{c}K}\check{\Sigma}_{KK}^{-1}\right)\hat{\mu}_{K}+\left(\hat{v}_{K^{c}}^{T}\frac{e}{n\lambda_{n,3}}\right)-\left(\check{\Sigma}_{K^{c}K}\check{\Sigma}_{KK}^{-1}\right)\hat{v}_{K}^{T}\frac{e}{n\lambda_{n,3}}.
\]
By the triangle inequality, we have 
\[
|\hat{\mu}_{K^{c}}|_{\infty}\leq\left\Vert \check{\Sigma}_{K^{c}K}\check{\Sigma}_{KK}^{-1}\right\Vert _{\infty}+\left|\hat{v}_{K^{c}}^{T}\frac{e}{n\lambda_{n,3}}\right|_{\infty}+\left\Vert \check{\Sigma}_{K^{c}K}\check{\Sigma}_{KK}^{-1}\right\Vert _{\infty}\left|\hat{v}_{K}^{T}\frac{e}{n\lambda_{n,3}}\right|_{\infty},
\]
where the fact that $|\hat{\mu}_{K}|_{\infty}\leq1$ is used in the
inequality above. By Lemma A.13, we have $\left\Vert \check{\Sigma}_{K^{c}K}\check{\Sigma}_{KK}^{-1}\right\Vert _{\infty}\leq1-\frac{\phi}{4}$
with probability at least $1-c\exp(-\log(p\wedge d))$. Hence, 
\begin{eqnarray*}
|\hat{\mu}_{K^{c}}|_{\infty} & \leq & 1-\frac{\phi}{4}+\left|\hat{v}_{K^{c}}^{T}\frac{e}{n\lambda_{n,3}}\right|_{\infty}+\left\Vert \check{\Sigma}_{K^{c}K}\check{\Sigma}_{KK}^{-1}\right\Vert _{\infty}\left|\hat{v}_{K}^{T}\frac{e}{n\lambda_{n,3}}\right|_{\infty}\\
 & \leq & 1-\frac{\phi}{4}+\left(2-\frac{\phi}{4}\right)\left|\hat{v}^{T}\frac{e}{n\lambda_{n,3}}\right|_{\infty}.
\end{eqnarray*}
Therefore, it suffices to show that $\left(2-\frac{\phi}{4}\right)\left|\hat{v}^{T}\frac{e}{n\lambda_{n,3}}\right|_{\infty}\leq\frac{\phi}{8}$
with high probability. This result is established in Lemma A.16. Thus,
we have $|\hat{\mu}_{K^{c}}|_{\infty}\leq1-\frac{\phi}{8}$ with high
probability. 

It remains to establish a bound on the $l_{\infty}-$norm of the error
$\hat{\beta}_{K}-\beta_{K}^{*}$. By the triangle inequality, we have
\begin{eqnarray*}
|\hat{\beta}_{K}-\beta_{K}^{*}|_{\infty} & \leq & \left|\left(\frac{\hat{v}_{K}^{T}\hat{v}_{K}}{n}\right)^{-1}\frac{\hat{v}_{K}^{T}e}{n}\right|_{\infty}+\lambda_{n,3}\left\Vert \left(\frac{\hat{v}_{K}^{T}\hat{v}_{K}}{n}\right)^{-1}\right\Vert _{\infty}\\
 & \leq & \left\Vert \left(\frac{\hat{v}_{K}^{T}\hat{v}_{K}}{n}\right)^{-1}\right\Vert _{\infty}\left|\frac{\hat{v}_{K}^{T}e}{n}\right|_{\infty}+\lambda_{n,3}\left\Vert \left(\frac{\hat{v}_{K}^{T}\hat{v}_{K}}{n}\right)^{-1}\right\Vert _{\infty},
\end{eqnarray*}
Using bound (\ref{eq:35}) from Lemma A.15, we have 
\[
\left\Vert \left(\frac{\hat{v}_{K}^{T}\hat{v}_{K}}{n}\right)^{-1}\right\Vert _{\infty}\leq\frac{2\sqrt{k_{2}}}{\lambda_{\min}(\hat{\Sigma}_{KK})}\leq\frac{4\sqrt{k_{2}}}{\lambda_{\min}(\Sigma_{KK})}.
\]
From the proof for Corollary 4.4, we have, with probability at least
$1-c_{1}\exp(-c_{2}\log(p\wedge d))$, 
\[
|\frac{1}{n}\hat{v}^{T}e|_{\infty}\precsim\left(\sigma_{v}\sigma_{\eta}\sqrt{\frac{\log p}{n}}\right)\vee\left(L\left|\beta^{*}\right|_{1}\sqrt{\Upsilon_{w,\theta^{*}}}b\left(\sigma_{v},\,\sigma_{\eta}\right)\sqrt{\frac{k_{1}\log d}{n}}\right).
\]
Putting everything together, with the choice of $\lambda_{n,3}$,
we obtain 
\[
|\hat{\beta}_{K}-\beta_{K}^{*}|_{\infty}\leq\frac{c\sqrt{k_{2}}}{\lambda_{\min}(\Sigma_{KK})}\left[\left(\sigma_{v}\sigma_{\eta}\sqrt{\frac{\log p}{n}}\right)\vee\left(L\left|\beta^{*}\right|_{1}\sqrt{\Upsilon_{w,\theta^{*}}}b\left(\sigma_{v},\,\sigma_{\eta}\right)\sqrt{\frac{k_{1}\log d}{n}}\right)\right],
\]
with probability at least $1-c_{1}\exp(-c_{2}\log(p\wedge d))$),
as claimed. $\square$

\subsection{Theorem 4.8}

\textbf{Proof}. Let $\Delta:=D^{-1}(\hat{\beta}-\beta^{*})$ and $\Psi_{n}:=\frac{1}{n}D\hat{v}^{T}\hat{v}D$
and we have 
\[
\left|\Psi_{n}\Delta\right|=\left|\frac{1}{n}D\hat{v}^{T}(\hat{v}_{0}-\hat{v}\hat{\beta})-\frac{1}{n}D\hat{v}^{T}(\hat{v}_{0}-\hat{v}\beta^{*})\right|.
\]
Consequently, 
\begin{eqnarray*}
\left|\Psi_{n}\Delta\right|_{\infty} & \leq & \left|\frac{1}{n}D\hat{v}^{T}(\hat{v}_{0}-\hat{v}\hat{\beta})\right|_{\infty}+\left|\frac{1}{n}D\hat{v}^{T}(\hat{v}_{0}-\hat{v}\beta^{*})\right|_{\infty}\\
 & \leq & \hat{\sigma}\xi+\left|\frac{1}{n}D\hat{v}^{T}e\right|_{\infty},
\end{eqnarray*}
recalling $e:=\hat{v}_{0}-\hat{v}\beta^{*}$ defined in Lemma A.1.
As in the nonpivotal case, to upper bound $\left|D\frac{\hat{v}^{T}e}{n}\right|_{\infty}$,
we control the terms except $\max_{j}\left\{ \frac{1}{n}\sum_{i=1}^{n}v_{ij}\eta_{i}\right\} $
in (\ref{eq:6})-(\ref{eq:8}) with the argument from the proofs for
Lemma A.4 and Corollary 4.4. To control for the term $\max_{j}\left\{ \frac{1}{n}\sum_{i=1}^{n}v_{ij}\eta_{i}\right\} $,
we use an argument similar to Gautier and Tsybakov (2011). Let $Q(\beta):=\frac{1}{n}\left|v_{0}-v\beta\right|_{2}^{2}$
and define the event 
\[
\mathcal{E}=\left\{ \left|\frac{1}{n}Dv^{T}\eta\right|_{\infty}\leq\sqrt{Q(\beta^{*})}\xi\right\} .
\]
Since $Q(\beta^{*})=\frac{1}{n}\sum_{i=1}^{n}\eta^{2}$, we have 
\begin{eqnarray*}
\mathbb{P}(\mathcal{E}^{c}) & \leq & \sum_{j=1}^{p}\mathbb{P}\left(\frac{1}{n}\left|\frac{\sum_{i=1}^{n}v_{ij}\eta_{i}}{v_{j*}\sqrt{\frac{1}{n}\sum_{i=1}^{n}\eta_{i}^{2}}}\right|\geq\xi\right)\\
 & \leq & \sum_{j=1}^{p}\mathbb{P}\left(\left|\frac{\sum_{i=1}^{n}v_{ij}\eta_{i}}{\sqrt{\sum_{i=1}^{n}\left(v_{ij}\eta_{i}\right)^{2}}}\right|\geq\sqrt{n}\xi\right).
\end{eqnarray*}
By Lemma B.16, for all $j=1,\,...\,,p$, 
\[
\mathbb{P}\left(\left|\frac{\sum_{i=1}^{n}v_{ij}\eta_{i}}{\sqrt{\sum_{i=1}^{n}\left(v_{ij}\eta_{i}\right)^{2}}}\right|\geq\sqrt{n}\xi\right)\leq2\left(1-\Phi\left(\sqrt{n}\xi\right)\right)+2a_{0}\frac{\left(1+\sqrt{n}\xi\right)^{1+\delta^{'}}}{p^{a^{2}}b_{n,\,\delta^{'}}^{2+\delta^{'}}}.
\]
Thus, the event $\mathcal{E}$ holds with probability at least $1-\alpha$
where 
\[
\alpha=2L\left(1-\Phi\left(\sqrt{n}\xi\right)\right)+2a_{0}\frac{\left(1+\sqrt{n}\xi\right)^{1+\delta^{'}}}{p^{a^{2}-1}b_{n,\,\delta^{'}}^{2+\delta^{'}}}.
\]
Now, let $\hat{Q}(\beta):=\frac{1}{n}\left|\hat{v}_{0}-\hat{v}\beta\right|_{2}^{2}$.
The following shows that, conditioning on the event $\mathcal{E}$
and by appropriately choosing $\xi$, we have 
\[
\left|\frac{1}{n}D\hat{v}^{T}(\hat{v}_{0}-\hat{v}\beta^{*})\right|_{\infty}\leq\hat{Q}(\beta^{*})\xi,
\]
namely, $\left(\beta^{*},\,\hat{Q}(\beta^{*})\right)\in A$. First,
notice that we have 
\begin{eqnarray*}
\left|\frac{1}{n}D\hat{v}^{T}(\hat{v}_{0}-\hat{v}\beta^{*})\right|_{\infty} & \leq & \left|\frac{1}{n}Dv^{T}(v_{0}-v\beta^{*})\right|_{\infty}+\left|\frac{1}{n}D\hat{v}^{T}(\hat{v}_{0}-\hat{v}\beta^{*})-\frac{1}{n}Dv^{T}(v_{0}-v\beta^{*})\right|_{\infty}\\
 & \leq & \left|\frac{1}{n}Dv^{T}(v_{0}-v\beta^{*})\right|_{\infty}+2\left|\frac{1}{n}Dv^{T}(\hat{v}_{0}-v_{0})\right|_{\infty}+\left|\frac{1}{n}D(\hat{v}-v){}^{T}(\hat{v}_{0}-v_{0})\right|_{\infty}\\
 &  & +2\left|\frac{1}{n}Dv^{T}(\hat{v}-v)\beta^{*}\right|_{\infty}+\left|\frac{1}{n}D(\hat{v}-v){}^{T}(\hat{v}-v)\beta^{*}\right|_{\infty}\\
 & \leq & \sqrt{Q(\beta^{*})}\xi+\frac{Lb(\sigma_{v})\sqrt{B^{'}}}{\min_{j=1,...,p}v_{j*}}\left|\beta^{*}\right|_{1}\\
 & \leq & \sqrt{\hat{Q}(\beta^{*})}\xi+\left|\sqrt{\hat{Q}(\beta^{*})}-\sqrt{Q(\beta^{*})}\right|\xi+\frac{Lb(\sigma_{v})\sqrt{B^{'}}}{\min_{j=\{1,...,p\}}v_{j*}}\left|\beta^{*}\right|_{1}
\end{eqnarray*}
where the third inequality follows from argument similar to the proofs
of the upper bounds for Lemma A.3 and Corollary 4.4, and $b(\sigma_{v})$
is a positive constant only depending on $\sigma_{v}$. We now control
for the term $\left|\sqrt{\hat{Q}(\beta^{*})}-\sqrt{Q(\beta^{*})}\right|$.
Note that 
\begin{eqnarray*}
\left|\hat{Q}(\beta^{*})-Q(\beta^{*})\right| & \leq & \left|\frac{1}{n}\sum_{i=1}^{n}\left(\hat{v}_{i0}^{2}-v_{i0}^{2}\right)\right|+\left|\frac{2}{n}\sum_{i=1}^{n}\left(\hat{v}_{i0}\hat{v}_{i}-v_{i0}v_{i}\right)\beta^{*}\right|+\left|\frac{1}{n}\sum_{i=1}^{n}\beta^{*T}\left(\hat{v}_{i}^{T}\hat{v}_{i}-v_{i}^{T}v_{i}\right)\beta^{*}\right|\\
 & \leq & \left|\frac{1}{n}\sum_{i=1}^{n}\left(\hat{v}_{i0}^{2}-v_{i0}^{2}\right)\right|+\max_{j=1,...,p}\left|\frac{2}{n}\sum_{i=1}^{n}\left(\hat{v}_{i0}\hat{v}_{ij}-v_{i0}v_{ij}\right)\right|\left|\beta^{*}\right|_{1}\\
 &  & +\max_{j,\, j^{'}=1,...,p}\left|\frac{2}{n}\sum_{i=1}^{n}\left(\hat{v}_{ij}^{T}\hat{v}_{ij^{'}}-v_{ij}^{T}v_{ij^{'}}\right)\right|\left|\beta^{*}\right|_{1}^{2}.
\end{eqnarray*}
Following from argument similar to the proofs of the upper bounds
for Lemma A.3 and Corollary 4.4, we can upper bound the terms 
\begin{eqnarray*}
 & \left|\frac{1}{n}\sum_{i=1}^{n}\left(\hat{v}_{i0}^{2}-v_{i0}^{2}\right)\right|,\\
 & \max_{j=1,...,p}\left|\frac{2}{n}\sum_{i=1}^{n}\left(\hat{v}_{i0}\hat{v}_{ij}-v_{i0}v_{ij}\right)\right|,\\
 & \max_{j,\, j^{'}=1,...,p}\left|\frac{2}{n}\sum_{i=1}^{n}\left(\hat{v}_{ij}^{T}\hat{v}_{ij^{'}}-v_{ij}^{T}v_{ij^{'}}\right)\right|,
\end{eqnarray*}
which yields 
\[
\left|\hat{Q}(\beta^{*})-Q(\beta^{*})\right|=O\left(\left(\frac{1}{n}\right)^{\frac{2}{3}}\left|\beta^{*}\right|_{1}+\left(\frac{\log d}{n}\right)^{\frac{1}{2}}\left|\beta^{*}\right|_{1}\right).
\]
Now, notice that if $\left|a_{1}-a_{2}\right|\leq a_{3}$ , then $\left|\sqrt{a_{1}}-\sqrt{a_{2}}\right|\leq\sqrt{a_{3}}$,
we have 
\[
\left|\sqrt{\hat{Q}(\beta^{*})}-\sqrt{Q(\beta^{*})}\right|=O\left(\left(\frac{1}{n}\right)^{\frac{1}{3}}\left|\beta^{*}\right|_{1}^{\frac{1}{2}}+\left(\frac{\log d}{n}\right)^{\frac{1}{4}}\left|\beta^{*}\right|_{1}^{\frac{1}{2}}\right),
\]
and therefore, 
\[
\left|\sqrt{\hat{Q}(\beta^{*})}-\sqrt{Q(\beta^{*})}\right|\xi\leq c\sqrt{\hat{Q}(\beta^{*})}\xi,
\]
for $c>1$. Consequently, if we choose 
\[
\xi\geq a\max\left\{ c\sqrt{\frac{\log p}{n}},\:\frac{Lb(\sigma_{v})\sqrt{B^{'}}}{\min_{j=1,...,p}v_{j*}}\left(\hat{Q}(\beta^{*})\right)^{-\frac{1}{2}}\left|\beta^{*}\right|_{1}\right\} 
\]
we have 
\[
\left|\frac{1}{n}D\hat{v}^{T}(\hat{v}_{0}-\hat{v}\beta^{*})\right|_{\infty}<\sqrt{\hat{Q}(\beta^{*})}\xi
\]
with probability at least $1-O\left(\frac{1}{p\wedge d}\right)-\alpha$.
On the other hand, $(\hat{\beta},\,\hat{\sigma})$ minimizes the criterion
$\left|D^{-1}\beta\right|_{1}+C\sigma$ on the same set $A$. Thus,
we have 
\[
\left|D^{-1}\hat{\beta}\right|_{1}+C\hat{\sigma}\leq\left|D^{-1}\beta^{*}\right|_{1}+C\sqrt{\hat{Q}(\beta^{*})},
\]
which implies 
\begin{eqnarray*}
\left|\Delta_{J(\beta^{*})^{c}}\right|_{1} & = & \sum_{j\in J(\beta^{*})^{c}}\left|v_{j*}\hat{\beta}_{j}\right|\\
 & \leq & \sum_{j\in J(\beta^{*})}\left(\left|v_{j*}\beta_{j}^{*}\right|-\left|v_{j*}\hat{\beta}_{j}\right|\right)+C\left(\sqrt{\hat{Q}(\beta^{*})}-\sqrt{\hat{Q}(\hat{\beta})}\right)\\
 & \leq & \left|\Delta_{J(\beta^{*})}\right|_{1}+C\left(\sqrt{\hat{Q}(\beta^{*})}-\sqrt{\hat{Q}(\hat{\beta})}\right)\\
 & \leq & \left|\Delta_{J(\beta^{*})}\right|_{1}+C\left|\frac{\left(\frac{1}{n}\sum_{i=1}^{n}\hat{v}_{ij}e_{i}\right)D\Delta}{\sqrt{\frac{1}{n}\sum_{i=1}^{n}e_{i}^{2}}}\right|\\
 & \leq & \left|\Delta_{J(\beta^{*})}\right|_{1}+C\left|\frac{\left(\frac{1}{n}\sum_{i=1}^{n}\hat{v}_{ij}e_{i}\right)D}{\sqrt{\frac{1}{n}\sum_{i=1}^{n}e_{i}^{2}}}\right|_{\infty}\left|\Delta\right|_{1}\\
 & \leq & \left|\Delta_{J(\beta^{*})}\right|_{1}+C\left|\Delta\right|_{1},
\end{eqnarray*}
where the third inequality follows from the convexity of $\beta\mapsto\sqrt{\hat{Q}(\beta)}$
and the last inequality follows from the concavity of $\sqrt{\cdot}$
and the Cauchy-Schwarz inequality. Consequently, we have the following
cone condition: 
\[
\left|\Delta_{J(\beta^{*})^{c}}\right|_{1}\leq\frac{1+C}{1-C}\left|\Delta_{J(\beta^{*})}\right|_{1}.
\]
We now finish upper bounding $\left|\Psi_{n}\Delta\right|_{\infty}$.
Recall that 
\begin{eqnarray*}
\left|\Psi_{n}\Delta\right|_{\infty} & \leq & \hat{\sigma}\xi+\left|\frac{1}{n}D\hat{v}^{T}e\right|_{\infty}\\
 & \leq & \hat{\sigma}\xi+\sqrt{Q(\beta^{*})}\xi+\frac{Lb(\sigma_{v})\sqrt{B^{'}}}{\min_{j=1,...,p}v_{j*}}\left|\beta^{*}\right|_{1}\\
 & \leq & \xi\left(2\hat{\sigma}+\sqrt{Q(\beta^{*})}-\sqrt{Q(\hat{\beta})}\right)+\xi\left(\sqrt{Q(\hat{\beta})}-\sqrt{\hat{Q}(\hat{\beta})}\right)\\
 &  & +\frac{Lb(\sigma_{v})\sqrt{B^{'}}}{\min_{j=1,...,p}v_{j*}}\left|\beta^{*}\right|_{1}.
\end{eqnarray*}
We first upper bound the term $\sqrt{Q(\beta^{*})}-\sqrt{Q(\hat{\beta})}$.
Arguing as in the proof for the cone condition, we have 
\begin{eqnarray*}
\sqrt{Q(\beta^{*})}-\sqrt{Q(\hat{\beta})} & \leq & \left|\frac{\left(\frac{1}{n}\sum_{i=1}^{n}v_{ij}\eta_{i}\right)D\Delta}{\sqrt{\frac{1}{n}\sum_{i=1}^{n}\eta_{i}^{2}}}\right|\\
 & \leq & \max_{j\in1,...,p}\left|\frac{\left(\frac{1}{n}\sum_{i=1}^{n}v_{ij}\eta_{i}\right)}{\sqrt{\frac{1}{n}\sum_{i=1}^{n}v_{ij}^{2}\eta_{i}^{2}}}\right|\left|\Delta\right|_{1}\\
 & \leq & \xi\left|\Delta\right|_{1}.
\end{eqnarray*}
To upper bound the term $\sqrt{Q(\hat{\beta})}-\sqrt{\hat{Q}(\hat{\beta})}$,
we argue as in the proof for upper bounding $\left|\sqrt{\hat{Q}(\beta^{*})}-\sqrt{Q(\beta^{*})}\right|$
from the above and obtain 
\[
\sqrt{\hat{Q}(\hat{\beta})}-\sqrt{Q(\hat{\beta})}=O\left(\left(\frac{1}{n}\right)^{\frac{1}{3}}\left|\hat{\beta}\right|_{1}^{\frac{1}{2}}+\left(\frac{\log d}{n}\right)^{\frac{1}{4}}\left|\hat{\beta}\right|_{1}^{\frac{1}{2}}\right),
\]
and therefore, 
\[
\xi\left(\sqrt{\hat{Q}(\hat{\beta})}-\sqrt{Q(\hat{\beta})}\right)<\xi\hat{\sigma}.
\]
Combining the pieces together, we obtain 
\begin{eqnarray*}
\left|\Psi_{n}\Delta\right|_{\infty} & \leq & 2\xi\hat{\sigma}+\xi^{2}\left|\Delta\right|_{1}+\frac{Lb(\sigma_{v})\sqrt{B^{'}}}{\min_{j=1,...,p}v_{j*}}\left|\beta^{*}\right|_{1}\\
 & \leq & \xi\left(2\hat{\sigma}+\xi\frac{\left|\Psi_{n}\Delta\right|_{\infty}}{\kappa_{J(\beta^{*})}^{*}}\right)+\frac{Lb(\sigma_{v})\sqrt{B^{'}}}{\min_{j=1,...,p}v_{j*}}\left|\beta^{*}\right|_{1}
\end{eqnarray*}
where the second inequality uses the definition of the sensitivity
$\kappa_{J(\beta^{*})}^{*}$. The above yields 
\[
\left|\Psi_{n}\Delta\right|_{\infty}\leq\left[\frac{Lb(\sigma_{v})\sqrt{B^{'}}}{\min_{j=1,...,p}v_{j*}}\left|\beta^{*}\right|_{1}+2\xi\hat{\sigma}\right]\left[1-\frac{\xi^{2}}{\kappa_{J(\beta^{*})}^{*}}\right]^{-1}.
\]
Now, we can apply the cone condition to obtain 
\[
\left|D^{-1}(\hat{\beta}-\beta^{*})\right|_{2}\leq\frac{1}{\kappa_{J(\beta^{*})}^{*}}\left[\frac{Lb(\sigma_{v})\sqrt{B^{'}}}{\min_{j=1,...,p}v_{j*}}\left|\beta^{*}\right|_{1}+2\xi\hat{\sigma}\right]\left[1-\frac{\xi^{2}}{\kappa_{J(\beta^{*})}^{*}}\right]^{-1}.
\]
Some further algebraic manipulation shows 
\begin{eqnarray*}
\left|D^{-1}(\hat{\beta}-\beta^{*})\right|_{2} & \leq & \frac{1}{\kappa_{J(\beta^{*})}^{*}}\left[\frac{Lb(\sigma_{v})\sqrt{B^{'}}}{\min_{j=1,...,p}v_{j*}}\left|\hat{\beta}\right|_{1}+2\xi\hat{\sigma}\right]\left[1-\frac{\xi^{2}}{\kappa_{J(\beta^{*})}^{*}}\right]^{-1}\\
 & + & \frac{1}{\kappa_{J(\beta^{*})}^{*}}\left[\frac{Lb(\sigma_{v})\sqrt{B^{'}}}{\left(\min_{j=1,...,p}v_{j*}\right)^{2}}\sqrt{k_{2}}\left|D^{-1}\left(\hat{\beta}-\beta^{*}\right)\right|_{2}\right]\left[1-\frac{\xi^{2}}{\kappa_{J(\beta^{*})}^{*}}\right]^{-1}
\end{eqnarray*}
which implies 

\begin{eqnarray*}
\left|D^{-1}(\hat{\beta}-\beta^{*})\right|_{2} & \leq & \frac{1}{\kappa_{J(\beta^{*})}^{*}}\left[\frac{Lb(\sigma_{v})\sqrt{B^{'}}}{\min_{j=1,...,p}v_{j*}}\left|\hat{\beta}\right|_{1}+2\xi\hat{\sigma}\right]\left[1-\frac{\xi^{2}}{\kappa_{J(\beta^{*})}^{*}}\right]^{-1}\\
 &  & \cdot\left[1-\frac{1}{\kappa_{J(\beta^{*})}^{*}}\left[\frac{Lb(\sigma_{v})\sqrt{k_{2}B^{'}}}{\left(\min_{j=1,...,p}v_{j*}\right)^{2}}\right]\left[1-\frac{\xi^{2}}{\kappa_{J(\beta^{*})}^{*}}\right]^{-1}\right]^{-1}.
\end{eqnarray*}
Simple calculations yield, for all $j=1,...,p$, 
\begin{eqnarray*}
\left|\hat{\beta}_{j}-\beta_{j}^{*}\right| & \leq & \frac{1}{v_{j*}\kappa_{J(\beta^{*})}^{*}}\left[\frac{Lb(\sigma_{v})\sqrt{B^{'}}}{\min_{j=1,...,p}v_{j*}}\left|\hat{\beta}\right|_{1}+2\xi\hat{\sigma}\right]\left[1-\frac{\xi^{2}}{\kappa_{J(\beta^{*})}^{*}}\right]^{-1}\\
 &  & \cdot\left[1-\frac{1}{\kappa_{J(\beta^{*})}^{*}}\left[\frac{Lb(\sigma_{v})\sqrt{k_{2}B^{'}}}{\left(\min_{j=1,...,p}v_{j*}\right)^{2}}\right]\left[1-\frac{\xi^{2}}{\kappa_{J(\beta^{*})}^{*}}\right]^{-1}\right]^{-1}.
\end{eqnarray*}
To establish the upper bound on $\hat{\sigma}$, note that the optimality
of $(\hat{\beta},\,\hat{\sigma})$ implies 
\begin{eqnarray*}
C\hat{\sigma} & \leq & \left|\Delta_{J(\beta^{*})}\right|_{1}+C\sqrt{\hat{Q}(\beta^{*})}\\
 & \leq & \frac{\left|\Psi_{n}\Delta\right|_{\infty}}{\kappa_{J(\beta^{*}),J(\beta^{*})}^{*}}+C\sqrt{\hat{Q}(\beta^{*})}.
\end{eqnarray*}
$\square$

\subsection{Theorems 4.9 and 4.10}

\textbf{Proof. }The estimator of $g(w_{i}^{T}\theta^{*})$ in Theorems
4.9 and 4.10 takes on the following form 
\[
\hat{\mathbb{E}}\left(y_{2i}\,\vert\, w_{i}^{T}\hat{\theta}\right)-\hat{\mathbb{E}}\left(x_{i}\,\vert\, w_{i}^{T}\hat{\theta}\right)^{T}\hat{\beta}:=\hat{g}(w_{i}^{T}\hat{\theta}),
\]
where $\hat{\beta}=\hat{\beta}_{HSEL}$. We now derive an upper bound
for 
\[
\left|\hat{g}(\hat{\theta})-g(\theta^{*})\right|_{2}:=\sqrt{\frac{1}{n}\sum_{i=1}^{n}\left[\hat{g}(w_{i}^{T}\hat{\theta})-g(w_{i}^{T}\theta^{*})\right]^{2}}
\]
Note that we can write 
\begin{eqnarray*}
 & \hat{g}(w_{i}^{T}\hat{\theta})-g(w_{i}^{T}\theta^{*})\\
= & \left[\hat{\mathbb{E}}\left(y_{2i}\,\vert\, w_{i}^{T}\hat{\theta}\right)-\mathbb{E}\left(y_{2i}\,\vert\, w_{i}^{T}\theta^{*}\right)\right]-\left[\hat{\mathbb{E}}\left(x_{i}\,\vert\, w_{i}^{T}\hat{\theta}\right)^{T}\hat{\beta}-\mathbb{E}\left(x_{i}\,\vert\, w_{i}^{T}\theta^{*}\right)^{T}\beta^{*}\right]\\
= & \left[\hat{\mathbb{E}}\left(y_{2i}\,\vert\, w_{i}^{T}\hat{\theta}\right)-\mathbb{E}\left(y_{2i}\,\vert\, w_{i}^{T}\theta^{*}\right)\right]-\left[\hat{\mathbb{E}}\left(x_{i}\,\vert\, w_{i}^{T}\hat{\theta}\right)-\mathbb{E}\left(x_{i}\,\vert\, w_{i}^{T}\hat{\theta}\right)\right]^{T}\hat{\beta}\\
 & -\left[\mathbb{E}\left(x_{i}\,\vert\, w_{i}^{T}\hat{\theta}\right)-\mathbb{E}\left(x_{i}\,\vert\, w_{i}^{T}\theta^{*}\right)\right]^{T}\hat{\beta}-\mathbb{E}\left(x_{i}\,\vert\, w_{i}^{T}\theta^{*}\right)^{T}(\hat{\beta}-\beta^{*}).
\end{eqnarray*}
By an elementary inequality and Hölder's inequality, we have 
\begin{eqnarray*}
 & \left[\hat{g}(w_{i}^{T}\hat{\theta})-g(w_{i}^{T}\theta^{*})\right]^{2}\\
\leq & 4\left[\hat{\mathbb{E}}\left(y_{2i}\,\vert\, w_{i}^{T}\hat{\theta}\right)-\mathbb{E}\left(y_{2i}\,\vert\, w_{i}^{T}\theta^{*}\right)\right]^{2}+4\left[\left[\hat{\mathbb{E}}\left(x_{i}\,\vert\, w_{i}^{T}\hat{\theta}\right)-\mathbb{E}\left(x_{i}\,\vert\, w_{i}^{T}\hat{\theta}\right)\right]^{T}\hat{\beta}\right]^{2}\\
 & +4\left[\left[\mathbb{E}\left(x_{i}\,\vert\, w_{i}^{T}\hat{\theta}\right)-\mathbb{E}\left(x_{i}\,\vert\, w_{i}^{T}\theta^{*}\right)\right]^{T}\hat{\beta}\right]^{2}+4\left[\mathbb{E}\left(x_{i}\,\vert\, w_{i}^{T}\theta^{*}\right)^{T}(\hat{\beta}-\beta^{*})\right]^{2}\\
\leq & 8\left[\hat{\mathbb{E}}\left(y_{2i}\,\vert\, w_{i}^{T}\hat{\theta}\right)-\mathbb{E}\left(y_{2i}\,\vert\, w_{i}^{T}\hat{\theta}\right)\right]^{2}+8\left[\mathbb{E}\left(y_{2i}\,\vert\, w_{i}^{T}\hat{\theta}\right)-\mathbb{E}\left(y_{2i}\,\vert\, w_{i}^{T}\theta^{*}\right)\right]^{2}\\
 & +8\max_{j}\left[\hat{\mathbb{E}}\left(x_{ij}\,\vert\, w_{i}^{T}\hat{\theta}\right)-\mathbb{E}\left(x_{ij}\,\vert\, w_{i}^{T}\hat{\theta}\right)\right]^{2}\left|\beta^{*}\right|_{1}^{2}+8\max_{j}\left[\hat{\mathbb{E}}\left(x_{ij}\,\vert\, w_{i}^{T}\hat{\theta}\right)-\mathbb{E}\left(x_{ij}\,\vert\, w_{i}^{T}\hat{\theta}\right)\right]^{2}\left|\hat{\beta}-\beta^{*}\right|_{1}^{2}\\
 & +8\max_{j}\left[\mathbb{E}\left(x_{ij}\,\vert\, w_{i}^{T}\hat{\theta}\right)-\mathbb{E}\left(x_{ij}\,\vert\, w_{i}^{T}\theta^{*}\right)\right]^{2}\left|\beta^{*}\right|_{1}^{2}+8\max_{j}\left[\mathbb{E}\left(x_{ij}\,\vert\, w_{i}^{T}\hat{\theta}\right)-\mathbb{E}\left(x_{ij}\,\vert\, w_{i}^{T}\theta^{*}\right)\right]^{2}\left|\hat{\beta}-\beta^{*}\right|_{1}^{2}\\
 & +4\max_{j}\left(\mathbb{E}\left(x_{ij}\,\vert\, w_{i}^{T}\theta^{*}\right)\right)^{2}\left|\hat{\beta}-\beta^{*}\right|_{1}^{2}.
\end{eqnarray*}
Applying previous results yields the desired result. By a standard
integration over the tail probability (in the exponential form), the
bound above with high probability can be converted to a bound in expectation
as stated in Theorem 4.8.$\square$

\subsection{Theorems 4.11 and 4.12}

\textbf{Proof. }We prove the case where $\beta^{*}$ and $\theta^{*}$
are exactly sparse and the approximate sparse case of $\beta^{*}$
and $\theta^{*}$ follows the same argument (see the remark at the
end of this section). The estimator of $g(w_{i}^{T}\theta^{*})$ in
Theorem 4.11 is obtained by solving the following program 
\[
\tilde{g}\in\arg\min_{f\in\mathcal{F}}\frac{1}{n}\sum_{i=1}^{n}\left(y_{2i}-x_{i}^{T}\hat{\beta}-f(w_{i}^{T}\hat{\theta})\right)^{2}
\]
where $\hat{\beta}=\hat{\beta}_{HSEL}$. We derive an upper bound
for 
\[
\left|\tilde{g}(\hat{\theta})-g(\theta^{*})\right|_{2}:=\sqrt{\frac{1}{n}\sum_{i=1}^{n}\left[\tilde{g}(w_{i}^{T}\hat{\theta})-g(w_{i}^{T}\theta^{*})\right]^{2}}.
\]
Since $\tilde{g}$ is optimal and $g$ is feasible, we have the basic
inequality 
\[
\frac{1}{n}\sum_{i=1}^{n}\left(y_{2i}-x_{i}^{T}\hat{\beta}-\tilde{g}(w_{i}^{T}\hat{\theta})\right)^{2}\leq\frac{1}{n}\sum_{i=1}^{n}\left(y_{2i}-x_{i}^{T}\hat{\beta}-g(w_{i}^{T}\hat{\theta})\right)^{2}.
\]
Some algebra leads to the equivalent expression 
\begin{eqnarray*}
\frac{1}{2}\left|\tilde{g}(\hat{\theta})-g(\hat{\theta})\right|_{n}^{2} & \leq & \frac{1}{n}\sum_{i=1}^{n}\left[x_{i}^{T}(\beta^{*}-\hat{\beta})\right]\left[\tilde{g}(w_{i}^{T}\hat{\theta})-g(w_{i}^{T}\hat{\theta})\right]\\
 &  & +\frac{1}{n}\sum_{i=1}^{n}\left[g(w_{i}^{T}\theta^{*})-g(w_{i}^{T}\hat{\theta})\right]\left[\tilde{g}(w_{i}^{T}\hat{\theta})-g(w_{i}^{T}\hat{\theta})\right]\\
 &  & +\frac{1}{n}\sum_{i=1}^{n}\eta_{i}\left[\tilde{g}(w_{i}^{T}\theta^{*})-g(w_{i}^{T}\theta^{*})\right]\\
 &  & +\frac{1}{n}\sum_{i=1}^{n}\eta_{i}\left[\tilde{g}(w_{i}^{T}\hat{\theta})-\tilde{g}(w_{i}^{T}\theta^{*})\right]\\
 &  & -\frac{1}{n}\sum_{i=1}^{n}\eta_{i}\left[g(w_{i}^{T}\hat{\theta})-g(w_{i}^{T}\theta^{*})\right].
\end{eqnarray*}
Now, by the Fenchel-Young inequality, we have 
\begin{eqnarray*}
 & \frac{1}{n}\sum_{i=1}^{n}\left[x_{i}^{T}(\beta^{*}-\hat{\beta})\right]\left[\tilde{g}(w_{i}^{T}\hat{\theta})-g(w_{i}^{T}\hat{\theta})\right]+\frac{1}{n}\sum_{i=1}^{n}\left[g(w_{i}^{T}\theta^{*})-g(w_{i}^{T}\hat{\theta})\right]\left[\tilde{g}(w_{i}^{T}\hat{\theta})-g(w_{i}^{T}\hat{\theta})\right]\\
 & \leq\frac{2}{n}\sum_{i=1}^{n}\left[x_{i}^{T}(\beta^{*}-\hat{\beta})\right]^{2}+\frac{1}{8n}\sum_{i=1}^{n}\left[\tilde{g}(w_{i}^{T}\hat{\theta})-g(w_{i}^{T}\hat{\theta})\right]^{2}\\
 & +\frac{2}{n}\sum_{i=1}^{n}\left[g(w_{i}^{T}\theta^{*})-g(w_{i}^{T}\hat{\theta})\right]^{2}+\frac{1}{8n}\sum_{i=1}^{n}\left[\tilde{g}(w_{i}^{T}\hat{\theta})-g(w_{i}^{T}\hat{\theta})\right]^{2}.
\end{eqnarray*}
Consequently, we have 
\begin{eqnarray*}
\frac{1}{4}\left|\tilde{g}(\hat{\theta})-g(\hat{\theta})\right|_{n}^{2} & \leq & \frac{2}{n}\sum_{i=1}^{n}\left[x_{i}^{T}(\beta^{*}-\hat{\beta})\right]^{2}+\frac{2}{n}\sum_{i=1}^{n}\left[g(w_{i}^{T}\theta^{*})-g(w_{i}^{T}\hat{\theta})\right]^{2}\\
 &  & +\frac{1}{n}\sum_{i=1}^{n}\eta_{i}\left[\tilde{g}(w_{i}^{T}\theta^{*})-g(w_{i}^{T}\theta^{*})\right]\\
 &  & +\frac{1}{n}\sum_{i=1}^{n}\eta_{i}\left[\tilde{g}(w_{i}^{T}\hat{\theta})-\tilde{g}(w_{i}^{T}\theta^{*})\right]\\
 &  & -\frac{1}{n}\sum_{i=1}^{n}\eta_{i}\left[g(w_{i}^{T}\hat{\theta})-g(w_{i}^{T}\theta^{*})\right].
\end{eqnarray*}
First, recall the upper bound on $\left|\hat{\beta}-\beta^{*}\right|_{2}^{2}$
from Corollary 4.4. Let $\Delta=\hat{\beta}-\beta^{*}$. Applying
the inequality in Lemma B.2 yields 
\[
\frac{\left|x\Delta\right|{}_{2}^{2}}{n}\leq\frac{3\kappa_{U}^{x}}{2}\left|\Delta\right|{}_{2}^{2}+\alpha^{'}\frac{\log d}{n}\left|\Delta\right|{}_{1}^{2},
\]
and combining with the fact that $\left|\Delta\right|_{1}\leq4\sqrt{k_{2}}\left|\Delta\right|_{2}$
shown previously yields an upper bound on the term $\frac{1}{n}\sum_{i=1}^{n}\left[x_{i}^{T}(\beta^{*}-\hat{\beta})\right]^{2}$. 

Now, by the $\bar{L}-$Lipschitz condition on $g(\cdot)$, we have
\[
\frac{1}{n}\sum_{i=1}^{n}\left[g(w_{i}^{T}\theta^{*})-g(w_{i}^{T}\hat{\theta})\right]^{2}\leq\frac{\bar{L}^{2}}{n}\sum_{i=1}^{n}\left[w_{i}^{T}\hat{\theta}-w_{i}^{T}\theta^{*}\right]^{2}
\]
which is upper bounded by Proposition 4.2.

For the term $\frac{1}{n}\sum_{i=1}^{n}\eta_{i}\left[\tilde{g}(w_{i}^{T}\theta^{*})-g(w_{i}^{T}\theta^{*})\right]$,
we follow the argument as in the proof for Corollary 4.4. Using the
fact that $\mathbb{E}\left[\eta_{i}\vert w_{i}\right]=0$, we obtain
\[
\frac{1}{n}\sum_{i=1}^{n}\eta_{i}\left[\tilde{g}(w_{i}^{T}\theta^{*})-g(w_{i}^{T}\theta^{*})\right]\leq c\left(\frac{k_{1}\sigma_{\eta}^{2}}{n}\right)^{\frac{2}{3}}.
\]
It remains to upper bound the terms $\frac{1}{n}\sum_{i=1}^{n}\eta_{i}\left[\tilde{g}(w_{i}^{T}\hat{\theta})-\tilde{g}(w_{i}^{T}\theta^{*})\right]$
and $\frac{1}{n}\sum_{i=1}^{n}\eta_{i}\left[g(w_{i}^{T}\hat{\theta})-g(w_{i}^{T}\theta^{*})\right]$.
By the $\bar{L}-$Lipschitz assumption on $g(\cdot)$ and $\tilde{g}(\cdot)$,
we have 
\begin{eqnarray*}
\frac{1}{n}\sum_{i=1}^{n}\eta_{i}\left[g(w_{i}^{T}\theta^{*})-g(w_{i}^{T}\hat{\theta})\right] & \leq & \frac{\bar{L}}{n}\sum_{i=1}^{n}\eta_{i}\left[w_{i}^{T}\hat{\theta}-w_{i}^{T}\theta^{*}\right],\\
\frac{1}{n}\sum_{i=1}^{n}\eta_{i}\left[\tilde{g}(w_{i}^{T}\theta^{*})-\tilde{g}(w_{i}^{T}\hat{\theta})\right] & \leq & \frac{\bar{L}}{n}\sum_{i=1}^{n}\eta_{i}\left[w_{i}^{T}\hat{\theta}-w_{i}^{T}\theta^{*}\right].
\end{eqnarray*}
We provide upper bounds on the RHS terms in the inequalities above.
For $r_{1},\, r_{2}>0$, define the set 
\[
\mathcal{S}(r_{1},\, r_{2}):=\left\{ \delta\in\mathbb{R}^{d}\,\vert\,\left|\delta\right|_{1}\leq r_{1},\:\left|\delta\right|_{2}\leq r_{2}\right\} ,
\]
and the random variables $U_{n}=U_{n}(r_{1},\, r_{2})$ given by 
\[
U_{n}:=\sup_{\delta\in\mathcal{S}(r_{1},\, r_{2})}\frac{1}{n}\left|\eta^{T}w\delta\right|.
\]
For a given $t\in(0,\,1)$ to be chosen, let us upper bound the minimal
cardinality of a set that covers $\mathcal{S}(r_{1},\, r_{2})$ up
to $r_{2}t-$accuracy in $l_{2}-$norm. By Lemma B.10, we can find
such a covering set $\left\{ \delta^{1},\,...\,,\delta^{N}\right\} \subset\mathcal{S}(r_{1},\, r_{2})$
with cardinality $N=N(r_{1},\, r_{2},\, t)$ that is upper bounded
as 
\[
\log N(r_{1},\, r_{2},\, t)\leq c_{0}r_{1}^{2}\left(\frac{1}{t}\right)^{2}\log d.
\]
Consequently, for each $\delta\in\mathcal{S}(r_{1},\, r_{2})$, we
may find some $\delta^{i}$ such that $\left|\delta-\delta^{i}\right|_{2}\leq r_{2}t$.
By the triangle inequality, we then have 
\begin{eqnarray*}
\frac{1}{n}\left|\eta^{T}w\delta\right| & \leq & \frac{1}{n}\left|\eta^{T}w\delta^{i}\right|+\frac{1}{n}\left|\eta^{T}w(\delta-\delta^{i})\right|\\
 & \leq & \frac{1}{n}\left|\eta^{T}w\delta^{i}\right|+\frac{\left|\eta\right|_{2}}{\sqrt{n}}\frac{\left|w(\delta-\delta^{i})\right|_{2}}{\sqrt{n}}.
\end{eqnarray*}
Given the assumptions on $w$, we have $\frac{\left|w(\delta-\delta^{i})\right|_{2}}{\sqrt{n}}\leq\sqrt{\kappa_{U}^{w}}\left|\delta-\delta^{i}\right|_{2}\leq\sqrt{\kappa_{U}^{w}}r_{2}t$.
Moreover, by Lemma B.1, we have $\frac{\left|\eta\right|_{2}}{\sqrt{n}}\leq2\sigma_{\eta}$
with probability $1-c_{1}\exp(-c_{2}n)$. Putting together the pieces,
we obtain 
\[
\frac{1}{n}\left|\eta^{T}w\delta\right|\leq\frac{1}{n}\left|\eta^{T}w\delta^{i}\right|+2\sqrt{\kappa_{U}^{w}}\sigma_{\eta}r_{2}t
\]
with probability $1-c_{1}\exp(-c_{2}n)$. Taking the supremum over
$\delta$ on both sides yields 
\[
U_{n}\leq\max_{i=1,...,N}\frac{1}{n}\left|\eta^{T}w\delta^{i}\right|+2\sqrt{\kappa_{U}^{w}}\sigma_{\eta}r_{2}t.
\]
It remains to bound the finite maximum over the covering set. Note
that conditioning on $\left\{ w_{i}\right\} _{i=1}^{n}$, each variate
$\frac{1}{n}\left|\eta^{T}w\delta^{i}\right|$ is zero-mean (by the
fact that $\mathbb{E}\left[\eta_{i}\vert w_{i}\right]=0$) sub-Gaussian
with parameter $\frac{\sigma_{\eta}^{2}\kappa_{U}^{w}r_{2}^{2}}{n}$.
Under the assumptions on $\delta^{i}$ and $w$, by the property of
sub-Gaussian variables (Lemma B.17), we conclude that 
\[
U_{n}\leq\sigma_{\eta}r_{2}\sqrt{\kappa_{U}^{w}}\left\{ \sqrt{\frac{3\log N(r_{1},\, r_{2},\, t)}{n}}+2t\right\} .
\]
Now, by Proposition 4.2, we can set $r_{1}=c^{'}k_{1}\sqrt{\frac{\log d}{n}}$
and $r_{2}=c^{''}\sqrt{\frac{k_{1}\log d}{n}}$. By choosing $t\asymp\sqrt{\frac{k_{1}\log d}{n}}$
so that $\sqrt{\frac{3\log N(r_{1},\, r_{2},\, t)}{n}}=2t$, we obtain
\[
U_{n}\leq c\sigma_{\eta}\sqrt{\kappa_{U}^{w}}\frac{k_{1}\log d}{n}.
\]
Consequently, we conclude that, 
\begin{eqnarray*}
\frac{1}{n}\sum_{i=1}^{n}\eta_{i}\left[g(w_{i}^{T}\theta^{*})-g(w_{i}^{T}\hat{\theta})\right] & \leq & c\sigma_{\eta}\sqrt{\kappa_{U}^{w}}\frac{k_{1}\log d}{n},\\
\frac{1}{n}\sum_{i=1}^{n}\eta_{i}\left[\tilde{g}(w_{i}^{T}\theta^{*})-\tilde{g}(w_{i}^{T}\hat{\theta})\right] & \leq & c\sigma_{\eta}\sqrt{\kappa_{U}^{w}}\frac{k_{1}\log d}{n}.
\end{eqnarray*}
Putting the pieces together, we have, 
\begin{eqnarray*}
\left|\tilde{g}(\hat{\theta})-g(\theta^{*})\right|_{n} & \leq & \left|\tilde{g}(\hat{\theta})-g(\hat{\theta})\right|_{n}+\left|g(\hat{\theta})-g(\theta^{*})\right|_{n}\\
 & \leq & c^{'}b^{'}\max\left\{ \sqrt{\frac{k_{2}\log p}{n}},\;\left|\beta^{*}\right|_{1}\sqrt{\frac{k_{1}k_{2}\log d}{n}},\;\left(\frac{k_{1}}{n}\right)^{\frac{1}{3}}\right\} .
\end{eqnarray*}
where $b^{'}$ is some constant depending only on the model-specific
structure and independent of $n$, $d$, $p$, $k_{1}$, and $k_{2}$.
By a standard integration over the tail probability, the bound above
with high probability can be converted to a bound in expectation as
stated in Theorem 4.11.$\square$\\
\\
\textbf{Remark}. To prove for the general sparsity case on $\beta^{*}$
and $\theta^{*}$($q_{1},\, q_{2}\in[0,\,1]$), we replace $r_{1}=c^{'}k_{1}\sqrt{\frac{\log d}{n}}$
and $r_{2}=c^{''}\sqrt{\frac{k_{1}\log d}{n}}$ with $r_{1}=c^{'}R_{q_{2}}\left(\sqrt{\frac{\log d}{n}}\right)^{1-q_{2}}$
and $r_{2}=c^{''}R_{q_{2}}^{\frac{1}{2}}\left(\sqrt{\frac{\log d}{n}}\right)^{1-\frac{q_{2}}{2}}$.

\subsection{Theorem 4.13}

Define $\mathcal{H}=\mathcal{B}_{q_{2}}^{p}(R_{q_{2}})\times\mathcal{F}\circ\mathcal{B}_{q_{1}}^{d}(R_{q_{1}})$
where $q_{1},\, q_{2}\in[0,\,1]$. Let $\tilde{x}:=(x,\,\overrightarrow{f_{\theta}})\in\mathbb{R}^{n\times(p+1)}$
where $\overrightarrow{f_{\theta}}:=\left(f(w_{1}^{T}\theta),\,...\,,f(w_{n}^{T}\theta)\right)^{T}$
with $\theta\in\mathcal{B}_{q_{1}}^{d}(R_{q_{1}})$ and $f(\cdot)\in\mathcal{F}$.
The lower bounds derived in the following involve the set defined
by intersecting the kernel of $\tilde{x}$ with $\mathcal{B}_{q_{2}}^{p}(R_{q_{2}})$,
which we denote $\mathcal{N}_{q_{2}}(\tilde{x}):=\textrm{Ker}(\tilde{x})\cap\mathcal{B}_{q_{2}}^{p}(R_{q_{2}})$
where 
\[
\textrm{Ker}(\tilde{x}):=\left\{ \beta\,:\,\tilde{x}\left(\begin{array}{c}
\beta\\
1
\end{array}\right)=0\right\} ,
\]
and define the $\mathcal{B}_{q_{2}}^{p}(R_{q_{2}})-$kernel diameter
in the $l_{2}-$norm 
\[
\textrm{diam}\left(\mathcal{N}_{q_{2}}(\tilde{x})\right):=\max_{\beta\in\mathcal{N}_{q_{2}}(\tilde{x})}\left|\beta\right|_{2}.
\]
The following Lemma controls for the term $\textrm{diam}\left(\mathcal{N}_{q_{2}}(\tilde{x})\right)$.
\\
\\
\textbf{Lemma A.10}: If Assumptions 4.9-4.10 hold for any $q_{2}\in[0,\,1]$,
then the $\mathcal{B}_{q_{2}}^{p}(R_{q_{2}})-$kernel diameter in
$l_{2}-$norm is upper bounded as 
\[
\textrm{diam}\left(\mathcal{N}_{q_{2}}(\tilde{x})\right)\leq\frac{f_{l}(R_{q_{2}},\, q_{2},\, n,\, p)}{\kappa_{l}}.
\]
\textbf{Proof.} This proof is similar to the proof for Lemma 1 in
Raskutti, et. al (2011). Under Assumption 4.10, we are guaranteed
that there is no measurable function $f(w_{i}^{T}\theta)$ such that
$x_{i}^{T}\lambda=f(w_{i}^{T}\theta)$ when $y_{1i}=1$ for $\lambda\in\mathcal{B}_{q_{2}}^{p}(R_{q_{2}})$.
On the other hand, if 
\[
\textrm{diam}\left(\mathcal{N}_{q_{2}}(\tilde{x})\right)>\frac{f_{l}(R_{q_{2}},\, q_{2},\, n,\, p)}{\kappa_{l}},
\]
then there must exist some $\beta\in\mathcal{B}_{q_{2}}^{p}(R_{q_{2}})$
with $x\beta=0$ and $\left|\beta\right|_{2}>\frac{f_{l}(R_{q_{2}},\, q_{2},\, n,\, p)}{\kappa_{l}}$.
We then have 
\[
0=\frac{1}{\sqrt{n}}\left|x\beta\right|_{2}<\kappa_{l}\left|\beta\right|_{2}-f_{l}(R_{q_{2}},\, q_{2},\, n,\, p),
\]
which contradicts the lower bound condition on $\frac{1}{\sqrt{n}}\left|x\beta\right|_{2}$
in Assumption 4.9. $\square$\\
\textbf{}\\
\textbf{Proof.} Let $M=M_{2}(\delta_{n};\,\mathcal{H})$ be the cardinality
of a maximal packing of $\mathcal{H}$ in the $L^{2}(\mathbb{P})-$metric,
say with elements $\left\{ h^{1},\,...\,,h^{M}\right\} $. Also let
the random vector $H\in\mathbb{R}^{p+1}$ be uniformly distributed
over the packing set $\left\{ h^{1},\,...\,,h^{M}\right\} $. A standard
argument yields a lower bound on the minimax $L^{2}(\mathbb{P})-$risk
in terms of the error in a multi-way hypothesis testing problem: in
particular, we have 
\begin{equation}
\min_{\hat{h}}\max_{h\in\mathcal{H}}\mathbb{E}\left|\hat{h}-h\right|_{2}^{2}\geq\frac{1}{4}\delta_{n}^{2}\min_{\tilde{h}}\mathbb{P}\left[\tilde{h}\neq H\right]\label{eq:13}
\end{equation}
where the estimator $\tilde{h}$ takes values in the packing set.
Fano's inequality implies the lower bound 
\[
\mathbb{P}\left[\tilde{h}\neq H\right]\geq1-\frac{I(H;\, y_{2})+\log2}{\log M_{2}(\delta_{n};\,\mathcal{H})}
\]
where $I(H;\, y_{2})$ is the mutual information between random parameters
$H$ in the packing set and the observation vector $y_{2}\in\mathbb{R}^{n}$. 

Let us first upper bound the mutual information; this is done via
the procedure of Yang and Barron (1999), which is based on covering
the model space $\mathcal{P}(\mathcal{H}):=\left\{ \mathbb{P}_{h}\,:\, h\in\mathcal{H}\right\} $
under the square-root Kullback-Leibler divergence. Following the argument
of Yang and Barron (1999), the mutual information is upper bounded
as 
\[
I(H;\, y_{2})\leq\log N(\varepsilon_{n};\,\mathcal{P}(\mathcal{H}))+\frac{c^{2}n}{\sigma_{\eta}^{2}}(\kappa_{1u}\vee\kappa_{2u})\varepsilon_{n}^{2}
\]
where $N(\varepsilon_{n};\,\mathcal{P}(\mathcal{H}))$ is the covering
number of $\mathcal{P}(\mathcal{H})$ under the KL divergence. Combining
this upper bound with the Fano lower bound above yields 
\begin{equation}
\mathbb{P}\left[\tilde{h}\neq H\right]\geq1-\frac{\log N(\varepsilon_{n};\,\mathcal{P}(\mathcal{H}))+\frac{c^{2}n}{\sigma_{\eta}^{2}}(\kappa_{1u}\vee\kappa_{2u})\varepsilon_{n}^{2}+\log2}{\log M_{2}(\delta_{n};\,\mathcal{H})}.\label{eq:14}
\end{equation}
An upper bound on $\log N(\varepsilon_{n};\,\mathcal{P}(\mathcal{H}))$
can be obtained using the following argument. For the Gaussian models
considered here, the KL divergence takes the form 
\begin{eqnarray*}
D(\mathbb{P}_{h} & \parallel & \mathbb{P}_{h^{'}})=\frac{1}{2\sigma_{\eta}^{2}}\left|x(\beta-\beta^{'})+f_{\theta}-f_{\theta^{'}}^{'}\right|_{2}^{2}\\
 & \leq & \frac{1}{\sigma_{\eta}^{2}}\left|x(\beta-\beta^{'})\right|_{2}^{2}+\frac{n}{\sigma_{\eta}^{2}}\left|f_{\theta}-f_{\theta^{'}}^{'}\right|_{n}^{2}\\
 & \leq & \frac{1}{\sigma_{\eta}^{2}}\left|x(\beta-\beta^{'})\right|_{2}^{2}+\frac{2n}{\sigma_{\eta}^{2}}\left|f_{\theta}-f_{\theta}^{'}\right|_{n}^{2}+\frac{2n}{\sigma_{\eta}^{2}}\left|f_{\theta}^{'}-f_{\theta^{'}}^{'}\right|_{n}^{2}\\
 & \leq & \frac{n\kappa_{2u}}{\sigma_{\eta}^{2}}\left|\beta-\beta^{'}\right|_{2}^{2}+\frac{2n}{\sigma_{\eta}^{2}}\left|f_{\theta}-f_{\theta}^{'}\right|_{\infty}^{2}+\frac{2n\bar{L}^{2}\kappa_{1u}}{\sigma_{\eta}^{2}}\left|\theta-\theta^{'}\right|_{2}^{2}
\end{eqnarray*}
where the last inequality uses the $\bar{L}-$Lipschitz assumption
on $f^{'}\in\mathcal{F}$ and the upper-RE condition on $w$ and $x$
as well as the fact that $\left|f-g\right|_{n}\leq\left|f-g\right|_{\infty}$. 

Now, applying the upper bound on the covering number of $\mathcal{B}_{0}^{d}(k_{1})$
in the $l_{2}-$metric provided by Lemma B.14, we conclude that there
exists a set $\left\{ \theta^{1},\,...\,,\theta^{N_{1}}\right\} $
such that for all $\theta\in\mathcal{B}_{0}^{d}(k_{1})$, there exists
some index $i$ with$\left|\theta-\theta^{i}\right|_{2}^{2}\leq\varepsilon_{n}^{2}$
and $\log N_{2}\left(\varepsilon_{n};\,\mathcal{B}_{0}^{d}(k_{1})\right)\leq k_{1}\log d+k_{1}\log\frac{\sqrt{k_{1}}}{\varepsilon_{n}}$;
similarly, there also exists a set $\left\{ \beta^{1},\,...\,,\beta^{N_{2}}\right\} $
such that for all $\beta\in\mathcal{B}_{0}^{p}(k_{2})$, there exists
some index $j$ with $\left|\beta-\beta^{j}\right|_{2}^{2}\leq\varepsilon_{n}^{2}$
and $\log N_{2}\left(\varepsilon_{n};\,\mathcal{B}_{0}^{p}(k_{2})\right)\leq k_{2}\log p+k_{2}\log\frac{\sqrt{k_{2}}}{\varepsilon_{n}}$.
Applying the upper bound on the covering number of the $\bar{L}-$Lipschitz
function class $\mathcal{F}$ in the sup-norm provided by Lemma B.11,
we conclude that there exists a set $\left\{ f_{\theta}^{1},\,...\,,f_{\theta}^{N_{3}}\right\} $
such that for all $f\in\mathcal{F}$ and $\theta\in\left\{ \theta^{1},\,...\,,\theta^{N_{1}}\right\} $,
there exists some index $l$ with $\sup_{w_{i}\theta}\left|f(w_{i}^{T}\theta)-f^{l}(w_{i}^{T}\theta)\right|\leq\varepsilon_{n}$
and $\log N_{\infty}\left(\varepsilon_{n};\,\mathcal{F}_{\theta}\right)\leq\frac{2k_{1}}{\varepsilon_{n}}$
for all $\theta\in\left\{ \theta^{1},\,...\,,\theta^{N_{1}}\right\} $.
It is not hard to see that a covering set of $\mathcal{P}(\mathcal{H})$
under the KL divergence can be formed by 
\[
\left\{ \theta^{1},\,...\,,\theta^{N_{1}}\right\} \times\left\{ \beta^{1},\,...\,,\beta^{N_{2}}\right\} \times\left\{ f_{\theta}^{1},\,...\,,f_{\theta}^{N_{3}}\,:\,\theta\in\left\{ \theta^{1},\,...\,,\theta^{N_{1}}\right\} \right\} 
\]
and 
\begin{eqnarray}
\log N(\varepsilon_{n};\,\mathcal{P}(\mathcal{H})) & \leq & \log\left[\left(N_{2}\left(\varepsilon_{n};\,\mathcal{B}_{0}^{d}(k_{1})\right)\right)^{2}\cdot N_{2}\left(\varepsilon_{n};\,\mathcal{B}_{0}^{p}(k_{2})\right)\cdot N_{\infty}\left(\varepsilon_{n};\,\mathcal{F}\right)\right]\nonumber \\
 & \leq & c^{'}\left\{ k_{1}\log d+k_{1}\log\frac{\sqrt{k_{1}}}{\varepsilon_{n}}+k_{2}\log p+k_{2}\log\frac{\sqrt{k_{2}}}{\varepsilon_{n}}+\frac{k_{1}}{\varepsilon_{n}}\right\} .\label{eq:15}
\end{eqnarray}
It remains to lower bound $\log M_{2}(\delta_{n};\,\mathcal{H})$;
we do so by constructing a minimal covering set of $\mathcal{H}=\mathcal{B}_{0}^{p}(k_{2})\times\mathcal{F}\circ\mathcal{B}_{0}^{d}(k_{1})$
in the $L^{2}(\mathbb{P})-$metric namely, $\log N_{2}(\delta_{n};\,\mathcal{H})$.
Note that 
\begin{eqnarray*}
\left|h-h^{'}\right|_{L^{2}(\mathbb{P})}^{2} & = & \left|\beta-\beta^{'}\right|_{2}^{2}+\mathbb{E}\left[f(w_{i}\theta)-f^{'}(w_{i}\theta^{'})\right]^{2}\\
 & \leq & \left|\beta-\beta^{'}\right|_{2}^{2}+2\mathbb{E}\left[f(w_{i}\theta)-f^{'}(w_{i}\theta)\right]^{2}+2\mathbb{E}\left[f^{'}(w_{i}\theta)-f^{'}(w_{i}\theta^{'})\right]^{2}\\
 & \leq & \left|\beta-\beta^{'}\right|_{2}^{2}+2\mathbb{E}\left[f(w_{i}\theta)-f^{'}(w_{i}\theta)\right]^{2}+2\kappa_{1u}\bar{L}^{2}\left|\theta-\theta^{'}\right|_{2}^{2}
\end{eqnarray*}
where the last inequality uses the $\bar{L}-$Lipschitz assumption
on $f^{'}\in\mathcal{F}$ and the population upper-RE condition on
$w$. Using the similar argument from above, we conclude from the
first inequality of Lemma B.12 and Lemma B.15 that there exists a
set $\left\{ \theta^{1},\,...\,,\theta^{N_{1}}\right\} $ such that
for all $\theta\in\mathcal{B}_{0}^{d}(k_{1})$, there exists some
index $i$ with$\left|\theta-\theta^{i}\right|_{2}^{2}\leq\delta_{n}^{2}$
and $\log N_{2}\left(\delta_{n};\,\mathcal{B}_{0}^{d}(k_{1})\right)\geq\log M_{2}\left(2\delta_{n};\,\mathcal{B}_{0}^{d}(k_{1})\right)\geq\frac{k_{1}}{2}\log\frac{d-k_{1}}{k_{1}/2}$;
similarly, there also exists a set $\left\{ \beta^{1},\,...\,,\beta^{N_{2}}\right\} $
such that for all $\beta\in\mathcal{B}_{0}^{p}(k_{2})$, there exists
some index $j$ with $\left|\beta-\beta^{j}\right|_{2}^{2}\leq\delta_{n}^{2}$
and $\log N_{2}\left(\delta_{n};\,\mathcal{B}_{0}^{p}(k_{2})\right)\geq\log M_{2}\left(2\delta_{n};\,\mathcal{B}_{0}^{p}(k_{2})\right)\geq\frac{k_{2}}{2}\log\frac{d-k_{2}}{k_{2}/2}$.
Furthermore, by Lemma B.11 and the first inequality of Lemma B.12,
we conclude that there exists a set $\left\{ f_{\theta}^{1},\,...\,,f_{\theta}^{N_{3}}\right\} $
such that for all $f\in\mathcal{F}$ and $\theta\in\left\{ \theta^{1},\,...\,,\theta^{N_{1}}\right\} $,
there exists some index $l$ with $\mathbb{E}\left[f(w_{i}^{T}\theta)-f^{l}(w_{i}^{T}\theta)\right]^{2}\leq\delta_{n}^{2}$
and $\log N_{2}\left(\delta_{n};\,\mathcal{F}_{\theta}\right)\geq\log M_{2}\left(2\delta_{n};\,\mathcal{F}_{\theta}\right)\geq\frac{k_{1}}{\delta_{n}}$
for all $\theta\in\left\{ \theta^{1},\,...\,,\theta^{N_{1}}\right\} $.
It is not hard to see that a covering set of $\mathcal{H}$ under
the $L^{2}(\mathbb{P})-$metric can be formed by 
\[
\left\{ \theta^{1},\,...\,,\theta^{N_{1}}\right\} \times\left\{ \beta^{1},\,...\,,\beta^{N_{2}}\right\} \times\left\{ f_{\theta}^{1},\,...\,,f_{\theta}^{N_{3}}\,:\,\theta\in\left\{ \theta^{1},\,...\,,\theta^{N_{1}}\right\} \right\} 
\]
and 
\begin{eqnarray*}
\log N_{2}(\delta_{n};\,\mathcal{H}) & \geq & \log\left[\left(M_{2}\left(2\delta_{n};\,\mathcal{B}_{0}^{d}(k_{1})\right)\right)^{2}\cdot M_{2}\left(2\delta_{n};\,\mathcal{B}_{0}^{p}(k_{2})\right)\cdot M_{2}\left(2\delta_{n};\,\mathcal{F}\right)\right]\\
 & \geq & c^{''}\left\{ \frac{k_{1}}{2}\log\frac{d-k_{1}}{k_{1}/2}+\frac{k_{2}}{2}\log\frac{d-k_{2}}{k_{2}/2}+\frac{k_{1}}{\delta_{n}}\right\} .
\end{eqnarray*}
Now applying the second inequality of Lemma B.12 yields 
\begin{eqnarray}
\log M(\delta_{n};\,\mathcal{H}) & \geq & \log N(\delta_{n};\,\mathcal{H})\nonumber \\
 & \geq & c^{''}\left\{ \frac{k_{1}}{2}\log\frac{d-k_{1}}{k_{1}/2}+\frac{k_{2}}{2}\log\frac{d-k_{2}}{k_{2}/2}+\frac{k_{1}}{\delta_{n}}\right\} .\label{eq:16}
\end{eqnarray}
Substituting (\ref{eq:15}) and (\ref{eq:16}) into (\ref{eq:14})
yields 
\[
\mathbb{P}\left[\tilde{h}\neq H\right]\geq1-\frac{c^{'}\left\{ k_{1}\log d+k_{1}\log\frac{\sqrt{k_{1}}}{\varepsilon_{n}}+k_{2}\log p+k_{2}\log\frac{\sqrt{k_{2}}}{\varepsilon_{n}}+\frac{k_{1}}{\varepsilon_{n}}\right\} +\frac{c^{2}n}{\sigma_{\eta}^{2}}(\kappa_{1u}\vee\kappa_{2u})\varepsilon_{n}^{2}+\log2}{c^{''}\left\{ \frac{k_{1}}{2}\log\frac{d-k_{1}}{k_{1}/2}+\frac{k_{2}}{2}\log\frac{d-k_{2}}{k_{2}/2}+\frac{k_{1}}{\delta_{n}}\right\} }.
\]
Under the condition $\varepsilon\geq\frac{1}{d}$, setting 
\begin{eqnarray*}
\frac{c^{2}n}{\sigma_{\eta}^{2}}(\kappa_{1u}\vee\kappa_{2u})\varepsilon_{n}^{2} & = & c^{'}\left\{ k_{1}\log d+k_{1}\log\frac{\sqrt{k_{1}}}{\varepsilon_{n}}+k_{2}\log p+k_{2}\log\frac{\sqrt{k_{2}}}{\varepsilon_{n}}+\frac{k_{1}}{\varepsilon_{n}}\right\} ,\\
\textrm{and}\:4\frac{c^{2}n}{\sigma_{\eta}^{2}}(\kappa_{1u}\vee\kappa_{2u})\varepsilon_{n}^{2} & = & c^{''}\left\{ \frac{k_{1}}{2}\log\frac{d-k_{1}}{k_{1}/2}+\frac{k_{2}}{2}\log\frac{d-k_{2}}{k_{2}/2}+\frac{k_{1}}{\delta_{n}}\right\} 
\end{eqnarray*}
yields the following 
\begin{eqnarray*}
\mathbb{P}\left[\tilde{h}\neq H\right] & \geq & \frac{1}{4},\\
\textrm{and\,\ }\delta_{n}^{2} & \asymp & \varepsilon_{n}^{2}\asymp\max\left\{ \frac{k_{1}\log d}{n},\,\frac{k_{2}\log p}{n},\,\left(\frac{k_{1}}{n}\right)^{\frac{2}{3}}\right\} .
\end{eqnarray*}
It can be easily verified that under the scaling of $\varepsilon_{n}$,
we are guaranteed to have $\varepsilon_{n}\geq\frac{1}{d}$.

Consequently, substituting the pieces above into (\ref{eq:13}) yields
\[
\min_{\hat{h}}\max_{h\in\mathcal{H}}\mathbb{E}\left|\hat{h}-h\right|_{2}^{2}\succsim\max\left\{ \frac{k_{1}\log d}{n},\,\frac{k_{2}\log p}{n},\,\left(\frac{k_{1}}{n}\right)^{\frac{2}{3}}\right\} .
\]
$\square$

\subsection{Theorem 4.14}

\textbf{Proof}. The proof for Theorem 4.14 is almost identical to
the proof for Theorem 4.13 except that instead of Lemma B.15, we apply
Lemma B.10 to lower bound the packing number of $\mathcal{B}_{q_{1}}^{d}(R_{q_{1}})$
and $\mathcal{B}_{q_{2}}^{p}(R_{q_{2}})$ in the $l_{2}-$metric where
$q_{1},\, q_{2}\in(0,\,1]$. These steps give us 
\[
\mathbb{P}\left[\tilde{h}\neq H\right]\geq1-\frac{u_{q_{1}}\left[R_{q_{1}}^{\frac{2}{2-q_{1}}}\left(\frac{1}{\varepsilon_{n}}\right)^{\frac{2q_{1}}{2-q_{1}}}\log d\right]+u_{q_{2}}\left[R_{q_{2}}^{\frac{2}{2-q_{2}}}\left(\frac{1}{\varepsilon_{n}}\right)^{\frac{2q_{2}}{2-q_{2}}}\log p\right]+c^{'}\frac{R}{\varepsilon_{n}}+\frac{c^{2}n}{\sigma_{\eta}^{2}}(\kappa_{1u}\vee\kappa_{2u})\varepsilon_{n}^{2}+\log2}{l_{q_{1}}\left[R_{q_{1}}^{\frac{2}{2-q_{1}}}\left(\frac{1}{\delta_{n}}\right)^{\frac{2q_{1}}{2-q_{1}}}\log d\right]+l_{q_{2}}\left[R_{q_{2}}^{\frac{2}{2-q_{2}}}\left(\frac{1}{\delta_{n}}\right)^{\frac{2q_{2}}{2-q_{2}}}\log p\right]+c^{''}\frac{R}{\delta_{n}}}.
\]
Setting 
\begin{eqnarray*}
\frac{c^{2}n}{\sigma_{\eta}^{2}}(\kappa_{1u}\vee\kappa_{2u})\varepsilon_{n}^{2} & = & u_{q_{1}}\left[R_{q_{1}}^{\frac{2}{2-q_{1}}}\left(\frac{1}{\varepsilon_{n}}\right)^{\frac{2q_{1}}{2-q_{1}}}\log d\right]+u_{q_{2}}\left[R_{q_{2}}^{\frac{2}{2-q_{2}}}\left(\frac{1}{\varepsilon_{n}}\right)^{\frac{2q_{2}}{2-q_{2}}}\log p\right]+c^{'}\frac{R}{\varepsilon_{n}},\\
\textrm{and}\:4\frac{c^{2}n}{\sigma_{\eta}^{2}}(\kappa_{1u}\vee\kappa_{2u})\varepsilon_{n}^{2} & = & \left[l_{q_{1}}R_{q_{1}}^{\frac{2}{2-q_{1}}}\left(\frac{1}{\delta_{n}}\right)^{\frac{2q_{1}}{2-q_{1}}}\log d\right]+\left[l_{q_{2}}R_{q_{2}}^{\frac{2}{2-q_{2}}}\left(\frac{1}{\delta_{n}}\right)^{\frac{2q_{2}}{2-q_{2}}}\log p\right]+c^{''}\frac{R}{\delta_{n}}
\end{eqnarray*}
yields the following 
\begin{eqnarray*}
\mathbb{P}\left[\tilde{h}\neq H\right] & \geq & \frac{1}{4},\\
\textrm{and\,\ }\delta_{n}^{2} & \asymp & \varepsilon_{n}^{2}\asymp\max\left\{ R_{q_{1}}\left(\frac{\log d}{n}\right)^{\frac{2-q_{1}}{2}},\, R_{q_{2}}\left(\frac{\log p}{n}\right)^{\frac{2-q_{1}}{2}},\,\left(\frac{R}{n}\right)^{\frac{2}{3}}\right\} .
\end{eqnarray*}
Consequently, substituting the pieces above into (\ref{eq:13}) yields
\[
\min_{\hat{h}}\max_{h\in\mathcal{H}}\mathbb{E}\left|\hat{h}-h\right|_{2}^{2}\succsim\max\left\{ R_{q_{1}}\left(\frac{\log d}{n}\right)^{\frac{2-q_{1}}{2}},\, R_{q_{2}}\left(\frac{\log p}{n}\right)^{\frac{2-q_{1}}{2}},\,\left(\frac{R}{n}\right)^{\frac{2}{3}}\right\} .
\]
$\square$

\subsection{Upper bounds on the second-stage estimators}

The following lemma upper bounds the term $\frac{1}{n}\sum_{i=1}^{n}\left[\hat{m}_{j}(w_{i}^{T}\hat{\theta})-m_{j}(w_{i}^{T}\hat{\theta})\right]^{2}$
(for $j=0,...,p$). Recall $y_{i}=z_{i0}$ and $x_{ij}=z_{ij}$ for
$j=1,...,p$. From Section 3 we have 
\[
z_{ij}=\mathbb{E}\left(z_{ij}\vert w_{i}^{T}\theta^{*}\right)+v_{ij}:=m_{j}(w_{i}^{T}\theta^{*})+v_{ij}\qquad\textrm{for }j=0,...,p
\]
which can be rewritten as 
\[
z_{ij}=m_{j}(w_{i}^{T}\hat{\theta})+v_{ij}-\left(m_{j}(w_{i}^{T}\hat{\theta})-m_{j}(w_{i}^{T}\theta^{*})\right).
\]
To avoid notation clustering, write $\left\{ \hat{m}_{j}(w_{i}^{T}\theta)\right\} _{i=1}^{n}:=\hat{m}_{j}(\theta)$.
Recall the set 
\[
\bar{\mathcal{F}}_{j}:=\left\{ f=f^{'}-f^{''}\,:\, f^{'},\, f^{''}\in\mathcal{F}_{j}\right\} ,
\]
and 
\[
\Omega(r_{j};\,\mathcal{F}_{j})=\left\{ f\,:\, f\in\mathcal{\bar{F}}_{j},\,\left|f_{\theta^{*}}\right|_{n}\leq r_{j}\right\} ,
\]
where $\left|f_{\theta^{*}}\right|_{n}:=\sqrt{\frac{1}{n}\sum_{i=1}^{n}\left[f(w_{i}^{T}\theta^{*})\right]^{2}}$,
and the conditional \textit{local complexity} 
\begin{equation}
\mathcal{G}_{n}(r_{j};\,\mathcal{F}_{j}):=\mathbb{E}_{v_{j}}\left[\sup_{f\in\Omega(r_{j};\,\mathcal{F}_{j})}\left|\frac{1}{n}\sum_{i=1}^{n}v_{ij}f(w_{i}^{T}\theta^{*})\right|\vert w_{i}^{T}\theta^{*}\right],\label{eq:17}
\end{equation}
where variables $\left\{ v_{ij}\right\} _{i=1}^{n}$ for $j=0,...,p$
are \textit{i.i.d.} variates that satisfy Assumption 4.2. In terms
of this complexity measure (\ref{eq:17}), we have the following result.\\
\\
\textbf{Lemma A.11}:\textbf{ }Let $r_{nj}>0$ be the smallest positive
quantity satisfying the \textit{critical inequality} 
\begin{equation}
\mathcal{G}_{n}\left(r_{nj};\,\mathcal{F}_{j}\right)\leq\frac{r_{nj}^{2}}{\sigma_{v_{j}}}.\label{eq:18}
\end{equation}
Then there are universal positive constants $(c_{1},\, c_{2})$ such
that for any $t\geq r_{nj}$ and $\tilde{m}_{j}\in\mathcal{F}_{j}$,
the non-parametric least squares estimate $\hat{m}_{j}$ via program
(11) satisfies 
\[
\left|\hat{m}_{j}(\hat{\theta})-m_{j}(\hat{\theta})\right|_{n}^{2}\leq7\left|\tilde{m}_{j}(\hat{\theta})-m_{j}(\hat{\theta})\right|_{n}^{2}+128tr_{nj}+4\sqrt{32\sigma_{v_{j}}^{2}T_{j}^{'}}+16T_{j}^{'}
\]
with probability at least $1-c_{1}\exp\left(-nC_{j}^{*}tr_{nj}\right)$.\\
\\
\textbf{Remark.} The constants in the error bound is not optimal and
can be improved.\\
\\
\textbf{Proof. }Since $\hat{m}_{j}$ is optimal and $\tilde{m}_{j}$
is feasible for the program (11), we have 
\[
\frac{1}{2}\sum_{i=1}^{n}\left(z_{ij}-\hat{m}_{j}(w_{i}^{T}\hat{\theta})\right)^{2}\leq\frac{1}{2}\sum_{i=1}^{n}\left(z_{ij}-\tilde{m}_{j}(w_{i}^{T}\hat{\theta})\right)^{2}.
\]
Define $\hat{\Delta}_{ij}(\theta)=\left(\hat{m}_{j}-m_{j}\right)\circ(w_{i}^{T}\theta)$,
$\left\{ \hat{\Delta}_{ij}(\theta)\right\} _{i=1}^{n}=\hat{\Delta}_{j}(\theta)$,
$\bar{\Delta}_{ij}(\theta)=\left(\hat{m}_{j}-\tilde{m}_{j}\right)\circ(w_{i}^{T}\theta)$,
$\left\{ \bar{\Delta}_{ij}(\theta)\right\} _{i=1}^{n}=\bar{\Delta}_{j}(\theta)$,
$\tilde{\Delta}_{ij}(\theta)=\left(\tilde{m}_{xj}-m_{xj}\right)\circ(w_{i}^{T}\theta)$,
and $\left\{ \tilde{\Delta}_{ij}(\theta)\right\} _{i=1}^{n}=\tilde{\Delta}_{j}(\theta)$.
Performing some algebra yields the \textit{basic inequality 
\[
\left|\hat{\Delta}_{j}(\hat{\theta})\right|_{n}^{2}\leq\left|\tilde{\Delta}_{j}(\hat{\theta})\right|_{n}^{2}+2\left|\frac{1}{n}\sum_{i=1}^{n}\left(v_{ij}+T_{ij}^{''}\right)\bar{\Delta}_{ij}(\hat{\theta})\right|
\]
}Furthermore, we have 
\begin{eqnarray*}
\frac{1}{n}\left|\sum_{i=1}^{n}T_{ij}^{''}\bar{\Delta}_{ij}(\hat{\theta})\right| & \leq & \frac{1}{n}\left|\sum_{i=1}^{n}T_{ij}^{''}\hat{\Delta}_{ij}(\hat{\theta})\right|+\frac{1}{n}\left|\sum_{i=1}^{n}T_{ij}^{''}\tilde{\Delta}_{ij}(\hat{\theta})\right|\\
 & \leq & \frac{2}{n}\left|\sum_{i=1}^{n}T_{ij}^{''2}\right|+\frac{1}{4n}\left|\sum_{i=1}^{n}\hat{\Delta}_{ij}^{2}(\hat{\theta})\right|+\frac{1}{4n}\left|\sum_{i=1}^{n}\tilde{\Delta}_{ij}^{2}(\hat{\theta})\right|,
\end{eqnarray*}
where the first inequality follows from the triangle inequality and
the second inequality follows from the Fenchel-Young inequality. Consequently,
we have 
\begin{equation}
\frac{1}{2}\left|\hat{\Delta}_{j}(\hat{\theta})\right|_{n}^{2}\leq\frac{3}{2}\left|\tilde{\Delta}_{j}(\hat{\theta})\right|_{n}^{2}+2\left|\frac{1}{n}\sum_{i=1}^{n}v_{ij}\bar{\Delta}_{ij}(\hat{\theta})\right|+4T_{j}^{'},\label{eq:19}
\end{equation}
recalling the definition of $T_{j}^{'}$. 

We need to upper bound the stochastic process $\frac{1}{n}\left|\sum_{i=1}^{n}v_{ij}\bar{\Delta}_{ij}(\hat{\theta})\right|$.
Now, note that by the elementary inequalities, we have 
\begin{eqnarray*}
\frac{1}{n}\left|\sum_{i=1}^{n}v_{ij}\bar{\Delta}_{ij}(\hat{\theta})\right| & \leq & \frac{1}{n}\left|\sum_{i=1}^{n}v_{ij}\left[\bar{\Delta}_{ij}(\hat{\theta})-\bar{\Delta}_{ij}(\theta^{*})\right]\right|+\frac{1}{n}\left|\sum_{i=1}^{n}v_{ij}\bar{\Delta}_{ij}(\theta^{*})\right|\\
 & \leq & \sqrt{\frac{1}{n}\sum_{i=1}^{n}v_{ij}^{2}}\sqrt{\frac{1}{n}\sum_{i=1}^{n}\left[\bar{\Delta}_{ij}(\hat{\theta})-\bar{\Delta}_{ij}(\theta^{*})\right]^{2}}+\frac{1}{n}\left|\sum_{i=1}^{n}v_{ij}\bar{\Delta}_{ij}(\theta^{*})\right|\\
 & \leq & \sqrt{\frac{1}{n}\sum_{i=1}^{n}v_{ij}^{2}}\sqrt{4\sup_{f\in\mathcal{F}_{j}}\frac{1}{n}\sum_{i=1}^{n}\left[f(w_{i}^{T}\hat{\theta})-f(w_{i}^{T}\theta^{*})\right]^{2}}+\frac{1}{n}\left|\sum_{i=1}^{n}v_{ij}\bar{\Delta}_{ij}(\theta^{*})\right|\\
 & \leq & \sqrt{\frac{1}{n}\sum_{i=1}^{n}v_{ij}^{2}}\sqrt{4T_{j}^{'}}+\frac{1}{n}\left|\sum_{i=1}^{n}v_{ij}\bar{\Delta}_{ij}(\theta^{*})\right|
\end{eqnarray*}
For the term $\sqrt{\frac{1}{n}\sum_{i=1}^{n}v_{ij^{'}}^{2}}$, by
Lemma B.1, we have 
\[
\mathbb{P}\left[\frac{1}{n}\sum_{i=1}^{n}v_{ij^{'}}^{2}\geq2\sigma_{v{}_{j^{'}}}^{2}\right]\leq2\exp(-cn).
\]
For the term $\frac{1}{n}\left|\sum_{i=1}^{n}v_{ij}\bar{\Delta}_{ij}(\theta^{*})\right|$,
notice that the error function $\bar{\Delta}_{j}$ belongs to the
set $\mathcal{\bar{F}}_{j}$. If $\left|\bar{\Delta}_{j}(\theta^{*})\right|_{n}\geq u=\sqrt{tr_{nj}}$,
an upper bound on the term $\frac{1}{n}\left|\sum_{i=1}^{n}v_{ij}\bar{\Delta}_{ij}(\theta^{*})\right|$
can be obtained by Lemma B.3: 
\begin{eqnarray*}
\frac{1}{n}\left|\sum_{i=1}^{n}v_{ij}\bar{\Delta}_{ij}(\theta^{*})\right| & \leq & \sup_{f(\theta^{*})\in\Omega(\left|\bar{\Delta}\right|_{n};\,\mathcal{F}_{j})}\left|\frac{1}{n}\sum_{i=1}^{n}v_{ij}f(w_{i}^{T}\theta^{*})\right|\\
 & \leq & 2\left|\bar{\Delta}_{j}(\theta^{*})\right|_{n}\sqrt{tr_{nj}},
\end{eqnarray*}
with probability at least $1-\exp\left(-nC_{j}^{*}tr_{nj}\right)$,
where $r_{nj}>0$ is the smallest positive quantity satisfying the
\textit{critical inequality} 
\[
\mathcal{G}_{n}\left(r_{nj};\,\mathcal{F}_{j}\right)\leq\frac{r_{nj}^{2}}{\sigma_{v_{j}}}.
\]
Applying the triangle inequality and the Fenchel-Young inequality
twice yields 
\[
\left|\bar{\Delta}_{j}(\theta^{*})\right|_{n}\sqrt{tr_{nj}}\leq8tr_{nj}+\frac{1}{16}\left(\left|\hat{\Delta}_{j}(\theta^{*})\right|_{n}^{2}+\left|\tilde{\Delta}_{j}(\theta^{*})\right|_{n}^{2}\right).
\]
Putting the pieces together, we obtain 
\[
\frac{1}{4}\left|\hat{\Delta}_{j}(\hat{\theta})\right|_{n}^{2}\leq\frac{7}{4}\left|\tilde{\Delta}_{j}(\hat{\theta})\right|_{n}^{2}+32tr_{nj}+\sqrt{32\sigma_{v_{j}}^{2}T_{j}^{'}}+4T_{j}^{'}
\]
where 
\[
T_{j}^{'}=\sup_{f\in\mathcal{F}_{j}}\frac{1}{n}\sum_{i=1}^{n}\left[f(w_{i}^{T}\hat{\theta})-f(w_{i}^{T}\theta^{*})\right]^{2}.
\]
If $\left|\bar{\Delta}_{j}(\theta^{*})\right|_{n}<\sqrt{tr_{nj}}$,
then the claim follows trivially by triangle inequality and the Fenchel-Young
inequality.$\square$\\
\textbf{}\\
\textbf{Lemma A.12}: For a user-defined radius $\bar{R}_{j}>0$, let
$r_{nj}>0$ be the smallest positive quantity satisfying the \textit{critical
inequality} 
\begin{equation}
\mathcal{G}_{n}\left(r_{nj};\,\mathcal{F}_{j}\right)\leq\frac{\bar{R}_{j}r_{nj}^{2}}{\sigma_{v_{j}}},\label{eq:20}
\end{equation}
where $\Omega(r_{j};\,\mathcal{F}_{j})=\left\{ f\,:\, f\in\bar{\mathcal{F}}_{j}\,\left|f_{\theta^{*}}\right|_{n}\leq r_{j},\,\left|f\right|_{\mathcal{F}_{j}}\leq1\right\} $
for $j=0,...,p$. Suppose that we solve program (12) with $\lambda_{nj,2}\geq2r_{nj}^{2}$.
Let $\tilde{m}_{j}\in\mathcal{F}_{j}$ with $\left|\tilde{m}_{j}\right|_{\mathcal{F}_{j}}\leq\bar{R}_{j}$.
Then there are universal positive constants $(c_{0},\, c_{1},\, c_{2},\, c_{3})$
such that the non-parametric least squares estimate $\hat{m}_{j}$
via (12) satisfies 
\[
\left|\hat{m}_{j}(\hat{\theta})-m_{j}(\hat{\theta})\right|_{n}^{2}\leq c_{0}\left|\tilde{m}_{j}(\hat{\theta})-m_{j}(\hat{\theta})\right|_{n}^{2}+c_{1}\bar{R}_{j}^{2}\left\{ r_{nj}^{2}+\lambda_{nj,2}\right\} +c_{2}T_{j}^{'}+c_{3}\sqrt{\sigma_{v}^{2}T_{j}^{'}}
\]
with probability at least $1-c\exp\left(-nC_{j}^{*}\bar{R}_{j}^{2}r_{nj}^{2}\right)$.\\
\textbf{}\\
\textbf{Proof. }Introduce the shorthand $\tilde{\sigma}_{v_{j}}=\frac{\sigma_{v_{j}}}{\bar{R}_{j}}$
and we work with an equivalent model with noise variance $\left(\frac{\sigma_{v_{j}}}{\bar{R}_{j}}\right)^{2}$
and the rescaled approximation error $\left|\tilde{m}_{j}(\hat{\theta})-m_{j}(\hat{\theta})\right|_{n}^{2}$.
Note that the final error then should be multiplied by $\bar{R}_{j}^{2}$.
Now, since $\hat{m}_{j}$ is optimal and $\tilde{m}_{j}$ is feasible
for program (12), we have
\[
\frac{1}{2}\sum_{i=1}^{n}\left(z_{ij}-\hat{m}_{j}(w_{i}^{T}\hat{\theta})\right)^{2}+\lambda_{nj,2}\left|\hat{m}_{j}\right|_{\mathcal{F}_{j}}^{2}\leq\frac{1}{2}\sum_{i=1}^{n}\left(z_{ij}-\tilde{m}_{j}(w_{i}^{T}\hat{\theta})\right)^{2}+\lambda_{nj,2}\left|\tilde{m}_{j}\right|_{\mathcal{F}_{j}}^{2}.
\]
Recall the definitions $\hat{\Delta}_{ij}(\theta):=\left(\hat{m}_{j}-m_{j}\right)\circ(w_{i}^{T}\theta)$,
$\tilde{\Delta}_{ij}(\theta):=\left(\tilde{m}_{j}-m_{j}\right)\circ(w_{i}^{T}\theta)$,
$\bar{\Delta}_{ij}(\theta):=\left(\hat{m}_{j}-\tilde{m}_{j}\right)\circ(w_{i}^{T}\theta)$,
$\left\{ \hat{\Delta}_{ij}(\theta)\right\} _{i=1}^{n}:=\hat{\Delta}_{j}(\theta)$,
$\left\{ \tilde{\Delta}_{ij}(\theta)\right\} _{i=1}^{n}:=\tilde{\Delta}_{j}(\theta)$,
and $\left\{ \bar{\Delta}_{ij}(\theta)\right\} _{i=1}^{n}:=\bar{\Delta}_{j}(\theta)$.
Performing some algebra yields the modified basic inequality 
\begin{eqnarray}
\frac{1}{2}\left|\hat{\Delta}_{j}(\hat{\theta})\right|_{n}^{2} & \leq & \frac{1}{2}\left|\tilde{\Delta}_{j}(\hat{\theta})\right|_{n}^{2}+\frac{1}{n}\left|\sum_{i=1}^{n}\left(v_{ij}+T_{ij}\right)\bar{\Delta}_{ij}(\hat{\theta})\right|+\lambda_{nj,2}\left\{ \left|\tilde{m}_{j}\right|_{\mathcal{F}_{j}}^{2}-\left|\hat{m}_{j}\right|_{\mathcal{F}_{j}}^{2}\right\} \nonumber \\
 & \leq & \frac{1}{2}\left|\tilde{\Delta}_{j}(\hat{\theta})\right|_{n}^{2}+\frac{1}{n}\left|\sum_{i=1}^{n}v_{ij}\bar{\Delta}_{ij}(\hat{\theta})\right|+\frac{1}{n}\left|\sum_{i=1}^{n}T_{ij}\bar{\Delta}_{ij}(\hat{\theta})\right|\nonumber \\
 &  & +\lambda_{nj,2}\left\{ \left|\tilde{m}_{j}\right|_{\mathcal{F}_{j}}^{2}-\left|\hat{m}_{j}\right|_{\mathcal{F}_{j}}^{2}\right\} .\label{eq:21}
\end{eqnarray}
Furthermore, we have 
\begin{eqnarray*}
\frac{1}{n}\left|\sum_{i=1}^{n}T_{ij}\bar{\Delta}_{ij}(\hat{\theta})\right| & \leq & \frac{1}{n}\left|\sum_{i=1}^{n}T_{ij}\hat{\Delta}_{ij}(\hat{\theta})\right|+\frac{1}{n}\left|\sum_{i=1}^{n}T_{ij}\tilde{\Delta}_{ij}(\hat{\theta})\right|\\
 & \leq & \frac{2}{n}\left|\sum_{i=1}^{n}T_{ij}^{2}\right|+\frac{1}{4n}\left|\sum_{i=1}^{n}\hat{\Delta}_{ij}^{2}(\hat{\theta})\right|+\frac{1}{4n}\left|\sum_{i=1}^{n}\tilde{\Delta}_{ij}^{2}(\hat{\theta})\right|,
\end{eqnarray*}
where the first inequality uses the triangle inequality and the second
uses the Fenchel-Young inequality. Substituting the inequality into
(\ref{eq:21}), we obtain 
\begin{equation}
\frac{1}{4}\left|\hat{\Delta}_{j}(\hat{\theta})\right|_{n}^{2}\leq\frac{3}{4}\left|\tilde{\Delta}_{j}(\hat{\theta})\right|_{n}^{2}+\frac{1}{n}\left|\sum_{i=1}^{n}v_{ij}\bar{\Delta}_{ij}(\hat{\theta})\right|+\frac{2}{n}\left|\sum_{i=1}^{n}T_{ij}^{2}\right|+\lambda_{nj,2}\left\{ \left|\tilde{m}_{j}\right|_{\mathcal{F}_{j}}^{2}-\left|\hat{m}_{j}\right|_{\mathcal{F}_{j}}^{2}\right\} .\label{eq:22}
\end{equation}
The remainder of the proof can be divided into two cases.\\
\textbf{\textit{Case }}\textbf{(}\textbf{\textit{i}}\textbf{): }$\left|\hat{m}_{j}\right|_{\mathcal{F}_{j}}\leq2$.
In this case, (\ref{eq:22}) implies that 
\[
\frac{1}{4}\left|\hat{\Delta}_{j}\right|_{n}^{2}\leq\frac{3}{4}\left|\tilde{\Delta}_{j}\right|_{n}^{2}+\frac{1}{n}\left|\sum_{i=1}^{n}v_{ij}\bar{\Delta}_{ij}(\hat{\theta})\right|+\frac{2}{n}\left|\sum_{i=1}^{n}T_{ij}^{2}\right|+9\lambda_{nj,2},
\]
using the fact that $\left|\tilde{m}_{j}\right|_{\mathcal{F}_{j}}\leq1$.
Since \textbf{$\left|\hat{m}_{j}\right|_{\mathcal{F}_{j}}\leq2$},
we have $\left|\bar{\Delta}_{j}\right|_{\mathcal{F}_{j}}\leq3$\textbf{
}and hence
\[
\left|\tilde{m}_{j}\right|_{\mathcal{F}_{j}}^{2}-\left|\hat{m}_{j}\right|_{\mathcal{F}_{j}}^{2}=\left(\left|\tilde{m}_{j}\right|_{\mathcal{F}_{j}}+\left|\hat{m}_{j}\right|_{\mathcal{F}_{j}}\right)\left(\left|\tilde{m}_{j}\right|_{\mathcal{F}_{j}}-\left|\hat{m}_{j}\right|_{\mathcal{F}_{j}}\right)\leq9
\]
by the triangle inequality. We need to upper bound the stochastic
process $\frac{1}{n}\left|\sum_{i=1}^{n}v_{ij}\bar{\Delta}_{ij}(\hat{\theta})\right|$.
\begin{eqnarray*}
\frac{1}{n}\left|\sum_{i=1}^{n}v_{ij}\bar{\Delta}_{ij}(\hat{\theta})\right| & \leq & \frac{1}{n}\left|\sum_{i=1}^{n}v_{ij}\left[\bar{\Delta}_{ij}(\hat{\theta})-\bar{\Delta}_{ij}(\theta^{*})\right]\right|+\frac{1}{n}\left|\sum_{i=1}^{n}v_{ij}\bar{\Delta}_{ij}(\theta^{*})\right|\\
 & \leq & \sqrt{\frac{1}{n}\sum_{i=1}^{n}v_{ij}^{2}}\sqrt{\frac{1}{n}\sum_{i=1}^{n}\left[\bar{\Delta}_{ij}(\hat{\theta})-\bar{\Delta}_{ij}(\theta^{*})\right]^{2}}+\frac{1}{n}\left|\sum_{i=1}^{n}v_{ij}\bar{\Delta}_{ij}(\theta^{*})\right|\\
 & \leq & \sqrt{\frac{1}{n}\sum_{i=1}^{n}v_{ij}^{2}}\sqrt{4\sup_{f\in\mathcal{F}_{j}}\frac{1}{n}\sum_{i=1}^{n}\left[f(w_{i}^{T}\hat{\theta})-f(w_{i}^{T}\theta^{*})\right]^{2}}+\frac{1}{n}\left|\sum_{i=1}^{n}v_{ij}\bar{\Delta}_{ij}(\theta^{*})\right|\\
 & \leq & \sqrt{\frac{1}{n}\sum_{i=1}^{n}v_{ij}^{2}}\sqrt{4T_{j}^{'}}+\frac{1}{n}\left|\sum_{i=1}^{n}v_{ij}\bar{\Delta}_{ij}(\theta^{*})\right|.
\end{eqnarray*}
Now, applying Lemma B.4 yields 
\begin{equation}
\frac{1}{n}\left|\sum_{i=1}^{n}v_{ij}\bar{\Delta}_{ij}(\theta^{*})\right|\leq6r_{nj}\left|\bar{\Delta}_{j}(\theta^{*})\right|_{n}+\frac{1}{32}\left|\bar{\Delta}_{j}(\theta^{*})\right|_{n}^{2}\label{eq:23}
\end{equation}
with probability at least $1-c\exp\left(-nC_{j}^{*}\bar{R}_{j}^{2}r_{nj}^{2}\right)$.
Consequently, 
\begin{eqnarray*}
\frac{1}{4}\left|\hat{\Delta}_{j}(\hat{\theta})\right|_{n}^{2} & \leq & \frac{3}{4}\left|\tilde{\Delta}_{j}(\hat{\theta})\right|_{n}^{2}+6r_{nj}\left|\bar{\Delta}_{j}(\hat{\theta})\right|_{n}+\frac{1}{32}\left|\bar{\Delta}_{j}(\hat{\theta})\right|_{n}^{2}+\frac{2}{n}\left|\sum_{i=1}^{n}T_{ij}^{2}\right|+9\lambda_{nj,2}\\
 & \leq & \frac{3}{4}\left|\tilde{\Delta}_{j}(\hat{\theta})\right|_{n}^{2}+216r_{nj}^{2}+\frac{1}{12}\left|\hat{\Delta}_{j}(\hat{\theta})\right|_{n}^{2}+\frac{1}{12}\left|\tilde{\Delta}_{j}(\hat{\theta})\right|_{n}^{2}\\
 &  & +\frac{1}{16}\left|\hat{\Delta}_{j}(\hat{\theta})\right|_{n}^{2}+\frac{1}{16}\left|\tilde{\Delta}_{j}(\hat{\theta})\right|_{n}^{2}+\frac{2}{n}\left|\sum_{i=1}^{n}T_{ij}^{2}\right|+9\lambda_{nj,2},
\end{eqnarray*}
where the second inequality follows by the triangle inequality, and
the Fenchel-Young inequality. This inequality above implies 
\[
\frac{5}{48}\left|\hat{\Delta}_{j}(\hat{\theta})\right|_{n}^{2}\leq\frac{43}{48}\left|\tilde{\Delta}_{j}(\hat{\theta})\right|_{n}^{2}+216r_{nj}^{2}+\frac{2}{n}\left|\sum_{i=1}^{n}T_{ij}^{2}\right|+9\lambda_{nj,2}.
\]
\textbf{}\\
\textbf{\textit{Case }}\textbf{(}\textbf{\textit{ii}}\textbf{): }$\left|\hat{m}_{j}\right|_{\mathcal{F}_{j}}>2>1\geq\left|\tilde{m}_{j}\right|_{\mathcal{F}_{j}}$.
In this case, we have 
\begin{eqnarray*}
\left|\tilde{m}_{j}\right|_{\mathcal{F}_{j}}^{2}-\left|\hat{m}_{j}\right|_{\mathcal{F}_{j}}^{2} & = & \left(\left|\tilde{m}_{j}\right|_{\mathcal{F}_{j}}+\left|\hat{m}_{j}\right|_{\mathcal{F}_{j}}\right)\left(\left|\tilde{m}_{j}\right|_{\mathcal{F}_{j}}-\left|\hat{m}_{j}\right|_{\mathcal{F}_{j}}\right)\\
 & \leq & \left|\tilde{m}_{j}\right|_{\mathcal{F}_{j}}-\left|\hat{m}_{j}\right|_{\mathcal{F}_{j}},
\end{eqnarray*}
since $\left|\tilde{m}_{j}\right|_{\mathcal{F}_{j}}-\left|\hat{m}_{j}\right|_{\mathcal{F}_{j}}<0$.
Writing $\hat{m}_{j}=\tilde{m}_{j}+\bar{\Delta}_{j}$ and applying
the triangle inequality $\left|\hat{m}_{j}\right|_{\mathcal{F}_{j}}\geq\left|\bar{\Delta}_{j}\right|_{\mathcal{F}_{j}}-\left|\tilde{m}_{j}\right|_{\mathcal{F}_{j}}$
yields 
\begin{eqnarray*}
\lambda_{nj,2}\left\{ \left|\tilde{m}_{j}\right|_{\mathcal{F}_{j}}^{2}-\left|\hat{m}_{j}\right|_{\mathcal{F}_{j}}^{2}\right\}  & \leq & \lambda_{nj,2}\left\{ 2\left|\tilde{m}_{j}\right|_{\mathcal{F}_{j}}-\left|\bar{\Delta}_{j}\right|_{\mathcal{F}_{j}}\right\} \\
 & \leq & 2\lambda_{nj,2}-\lambda_{nj,2}\left|\bar{\Delta}_{j}\right|_{\mathcal{F}_{j}}.
\end{eqnarray*}
Substituting this upper bound into (\ref{eq:22}) yields 
\begin{eqnarray}
\frac{1}{4}\left|\hat{\Delta}_{j}(\hat{\theta})\right|_{n}^{2} & \leq & \frac{3}{4}\left|\tilde{\Delta}_{j}(\hat{\theta})\right|_{n}^{2}+\frac{1}{n}\left|\sum_{i=1}^{n}v_{ij}\bar{\Delta}_{ij}(\hat{\theta})\right|+\frac{2}{n}\left|\sum_{i=1}^{n}T_{ij}^{2}\right|+2\lambda_{nj,2}-\lambda_{nj,2}\left|\bar{\Delta}_{j}\right|_{\mathcal{F}_{j}}\nonumber \\
 & \leq & \frac{3}{4}\left|\tilde{\Delta}_{j}(\hat{\theta})\right|_{n}^{2}+\frac{1}{n}\left|\sum_{i=1}^{n}v_{ij}\bar{\Delta}_{ij}(\theta^{*})\right|+\sqrt{\frac{1}{n}\sum_{i=1}^{n}v_{ij}^{2}}\sqrt{4T_{j}^{'}}\label{eq:24}\\
 &  & +\frac{2}{n}\left|\sum_{i=1}^{n}T_{ij}^{2}\right|+2\lambda_{nj,2}-\lambda_{nj,2}\left|\bar{\Delta}_{j}\right|_{\mathcal{F}_{j}}\nonumber 
\end{eqnarray}
Lemma B.5 upper bounds the stochastic component $\frac{1}{n}\left|\sum_{i=1}^{n}v_{ij}\bar{\Delta}_{ij}(\theta^{*})\right|$
in inequality (\ref{eq:24}): There are universal positive constants
$(c_{1},\, c_{2})$ such that for all $\left|\Delta_{j}\right|_{\mathcal{F}_{j}}\geq1$,
\begin{equation}
\frac{1}{n}\left|\sum_{i=1}^{n}v_{ij}\Delta_{ij}(\theta^{*})\right|\leq2r_{nj}\left|\Delta_{j}(\theta^{*})\right|_{n}+2r_{nj}^{2}\left|\Delta_{j}\right|_{\mathcal{F}_{j}}+\frac{1}{16}\left|\Delta_{j}(\theta^{*})\right|_{n}^{2}\label{eq:25}
\end{equation}
with probability at least $1-c\exp\left(-nC_{j}^{*}\bar{R}_{j}^{2}r_{nj}^{2}\right)$.
We now complete the proof using inequality (\ref{eq:25}). Since $\left|\hat{m}_{j}\right|_{\mathcal{F}_{j}}>2>1\geq\left|\tilde{m}_{j}\right|_{\mathcal{F}_{j}}$,
the triangle inequality implies that $\left|\bar{\Delta}_{j}\right|_{\mathcal{F}_{j}}\geq\left|\hat{m}_{j}\right|_{\mathcal{F}_{j}}-\left|\tilde{m}_{j}\right|_{\mathcal{F}_{j}}>1$,
so that inequality (\ref{eq:25}) can be applied. Consequently, 
\begin{eqnarray*}
\frac{1}{4}\left|\hat{\Delta}_{j}(\hat{\theta})\right|_{n}^{2} & \leq & \frac{3}{4}\left|\tilde{\Delta}_{j}(\hat{\theta})\right|_{n}^{2}+2r_{nj}\left|\bar{\Delta}_{j}(\hat{\theta})\right|_{n}+\left(2r_{nj}^{2}-\lambda_{nj,2}\right)\left|\bar{\Delta}_{j}\right|_{\mathcal{F}_{j}}\\
 &  & +\frac{1}{16}\left|\bar{\Delta}_{j}(\hat{\theta})\right|_{n}^{2}+\sqrt{\frac{1}{n}\sum_{i=1}^{n}v_{ij}^{2}}\sqrt{4T_{j}^{'}}+\frac{2}{n}\left|\sum_{i=1}^{n}T_{ij}^{2}\right|+2\lambda_{nj,2}\\
 & \leq & \frac{3}{4}\left|\tilde{\Delta}_{j}(\hat{\theta})\right|_{n}^{2}+32r_{nj}^{2}+\frac{1}{16}\left|\hat{\Delta}_{j}(\hat{\theta})\right|_{n}^{2}+\frac{1}{16}\left|\tilde{\Delta}_{j}(\hat{\theta})\right|_{n}^{2}\\
 &  & +\frac{1}{8}\left|\hat{\Delta}_{j}(\hat{\theta})\right|_{n}^{2}+\frac{1}{8}\left|\tilde{\Delta}_{j}(\hat{\theta})\right|_{n}^{2}+\sqrt{\frac{1}{n}\sum_{i=1}^{n}v_{ij}^{2}}\sqrt{4T_{j}^{'}}+\frac{2}{n}\left|\sum_{i=1}^{n}T_{ij}^{2}\right|+2\lambda_{nj,2},
\end{eqnarray*}
where the second inequality follows by the choice of $\lambda_{nj,2}$,
the triangle inequality, and the Fenchel-Young inequality. This inequality
above implies 
\[
\frac{1}{16}\left|\hat{\Delta}_{j}(\hat{\theta})\right|_{n}^{2}\leq\frac{15}{16}\left|\tilde{\Delta}_{j}(\hat{\theta})\right|_{n}^{2}+32r_{nj}^{2}+\sqrt{\frac{1}{n}\sum_{i=1}^{n}v_{ij}^{2}}\sqrt{4T_{j}^{'}}+\frac{2}{n}\left|\sum_{i=1}^{n}T_{ij}^{2}\right|+2\lambda_{nj,2}.
\]
\textbf{$\square$}

\subsection{Lemmas for Theorem 4.7}

\textbf{Lemma A.13}: If the assumptions in Corollary 4.4 and Assumption
4.7 hold, then under the conditions $n\succsim\left(k_{2}^{3}\log p\right)\vee\left(k_{2}^{2}k_{1}\log d\right)$,
$\sqrt{\frac{n\log p}{\left|\beta^{*}\right|_{1}^{2}}}\vee\sqrt{nk_{1}\log d}\succsim\log p$,
and $\sqrt{\frac{k_{1}\log d}{n}}=o(1)$, the sample matrix $\frac{1}{n}\hat{v}^{T}\hat{v}$
satisfies an analogous version of the {}``mutual incoherence'' assumption
with high probability, i.e., 
\[
\mathbb{P}\left[\left\Vert \frac{1}{n}\hat{v}_{K^{c}}^{T}\hat{v}_{K}\left(\frac{1}{n}\hat{v}_{K}^{T}\hat{v}_{K}\right)^{-1}\right\Vert _{\infty}\geq1-\frac{\phi}{4}\right]\leq O\left(\frac{1}{p\wedge d}\right).
\]
\textbf{Proof}. We use the following decomposition similar to the
method used in Ravikumar, et. al. (2010)\textbf{ 
\[
\check{\Sigma}_{K^{c}K}\check{\Sigma}_{KK}^{-1}-\Sigma_{K^{c}K}\Sigma_{KK}^{-1}=R_{1}+R_{2}+R_{3}+R_{4}+R_{5}+R_{6},
\]
}where\textbf{ 
\begin{eqnarray*}
R_{1} & = & \Sigma_{K^{c}K}[\hat{\Sigma}_{KK}^{-1}-\Sigma_{KK}^{-1}],\\
R_{2} & = & [\hat{\Sigma}_{K^{c}K}-\Sigma_{K^{c}K}]\Sigma_{KK}^{-1},\\
R_{3} & = & [\hat{\Sigma}_{K^{c}K}-\Sigma_{K^{c}K}][\hat{\Sigma}_{KK}^{-1}-\Sigma_{KK}^{-1}],\\
R_{4} & = & \hat{\Sigma}_{K^{c}K}[\check{\Sigma}_{KK}^{-1}-\hat{\Sigma}_{KK}^{-1}],\\
R_{5} & = & [\check{\Sigma}_{K^{c}K}-\hat{\Sigma}_{K^{c}K}]\hat{\Sigma}_{KK}^{-1},\\
R_{6} & = & [\check{\Sigma}_{K^{c}K}-\hat{\Sigma}_{K^{c}K}][\check{\Sigma}_{KK}^{-1}-\hat{\Sigma}_{KK}^{-1}].
\end{eqnarray*}
}By Assumption 4.7, we have 
\[
\left\Vert \Sigma_{K^{c}K}\Sigma_{KK}^{-1}\right\Vert _{\infty}\leq1-\phi.
\]
It suffices to show that $||R_{i}||_{\infty}\leq\frac{\phi}{6}$\textbf{
}for $i=1,...,3$ and $||R_{i}||_{\infty}\leq\frac{\phi}{12}$\textbf{
}for $i=4,...,6$. 

For the first term $R_{1}$, we have 
\[
R_{1}=-\Sigma_{K^{c}K}\Sigma_{KK}^{-1}[\hat{\Sigma}_{KK}-\Sigma_{KK}]\hat{\Sigma}_{KK}^{-1},
\]
Using the sub-multiplicative property $||AB||_{\infty}\leq||A||_{\infty}||B||_{\infty}$
and the elementary inequality $||A||_{\infty}\leq\sqrt{m}||A||_{2}$
for any symmetric matrix $A\in\mathbb{R}^{m\times m}$, we can bound
$R_{1}$ as follows: 
\begin{eqnarray*}
||R_{1}||_{\infty} & \leq & \left\Vert \Sigma_{K^{c}K}\Sigma_{KK}^{-1}\right\Vert _{\infty}\left\Vert \hat{\Sigma}_{KK}-\Sigma_{KK}\right\Vert _{\infty}\left\Vert \hat{\Sigma}_{KK}^{-1}\right\Vert _{\infty}\\
 & \leq & (1-\phi)\left\Vert \hat{\Sigma}_{KK}-\Sigma_{KK}\right\Vert _{\infty}\sqrt{k_{2}}\left\Vert \hat{\Sigma}_{KK}^{-1}\right\Vert _{2},
\end{eqnarray*}
where the last inequality follows from Assumption 4.7. Using bound
(\ref{eq:32}) from the proof for Lemma A.14, we have 
\[
\left\Vert \hat{\Sigma}_{KK}^{-1}\right\Vert _{2}\leq\frac{2}{\lambda_{\min}(\Sigma_{KK})}
\]
with probability at least $1-c_{1}\exp(-c_{2}n)$. Next, applying
bound (\ref{eq:27}) from Lemma A.14 with $\varepsilon=\frac{\phi\lambda_{\min}(\Sigma_{KK})}{12(1-\phi)\sqrt{k_{2}}}$,
we have 
\[
\mathbb{P}\left[\left\Vert \hat{\Sigma}_{KK}-\Sigma_{KK}\right\Vert _{\infty}\geq\frac{\phi\lambda_{\min}(\Sigma_{KK})}{12(1-\phi)\sqrt{k_{2}}}\right]\leq2\exp(-b\frac{n}{k_{2}^{3}}+2\log k_{2}).
\]
Then, we are guaranteed that 
\[
\mathbb{P}[||R_{1}||_{\infty}\geq\frac{\phi}{6}]\leq2\exp(-b\frac{n}{k_{2}^{3}}+2\log k_{2}).
\]
For the second term $R_{2}$, we first write 
\begin{eqnarray*}
||R_{2}||_{\infty} & \leq & \sqrt{k_{2}}\left\Vert \Sigma_{KK}^{-1}\right\Vert _{2}\left\Vert \hat{\Sigma}_{K^{c}K}-\Sigma_{K^{c}K}\right\Vert _{\infty}\\
 & \leq & \frac{\sqrt{k_{2}}}{\lambda_{\min}(\Sigma_{KK})}\left\Vert \hat{\Sigma}_{K^{c}K}-\Sigma_{K^{c}K}\right\Vert _{\infty}.
\end{eqnarray*}
An application of bound (\ref{eq:26}) from Lemma A.14 with $\varepsilon=\frac{\phi}{6}\frac{\lambda_{\min}(\Sigma_{KK})}{\sqrt{k_{2}}}$
to bound the term $\left\Vert \hat{\Sigma}_{K^{c}K}-\Sigma_{K^{c}K}\right\Vert _{\infty}$
yields 
\[
\mathbb{P}[||R_{2}||_{\infty}\geq\frac{\phi}{6}]\leq2\exp(-b\frac{n}{k_{2}^{3}}+\log(p-k_{2})+\log k_{2}).
\]
For the third term $R_{3}$, by applying bounds (\ref{eq:26}) from
Lemma A.14 with $\varepsilon=\frac{\phi\lambda_{\min}(\Sigma_{KK})}{6}$
to bound the term $\left\Vert \hat{\Sigma}_{K^{c}K}-\Sigma_{K^{c}K}\right\Vert _{\infty}$
and (\ref{eq:28}) from Lemma A.14 to bound the term $\left\Vert \hat{\Sigma}_{KK}^{-1}-\Sigma_{KK}^{-1}\right\Vert _{\infty}$,
we have

\[
\mathbb{P}[||R_{3}||_{\infty}\geq\frac{\phi}{6}]\leq2\exp(-b\frac{n}{k_{2}^{2}}+\log(p-k_{2})+\log k_{2}).
\]
Putting everything together, we conclude that 
\[
\mathbb{P}[||\hat{\Sigma}_{K^{c}K}\hat{\Sigma}_{KK}^{-1}||_{\infty}\geq1-\frac{\phi}{2}]\leq c\exp(-b\frac{n}{k_{2}^{3}}+2\log p).
\]
For the fourth term $R_{4}$, we have, with probability at least $1-c\exp(-bn\min\{\frac{1}{k_{2}^{3}},\,\frac{1}{k_{2}^{3/2}}\}+2\log p)$,
\begin{eqnarray*}
||R_{4}||_{\infty} & \leq & \left\Vert \hat{\Sigma}_{K^{c}K}\hat{\Sigma}_{KK}^{-1}\right\Vert _{\infty}\left\Vert \check{\Sigma}_{KK}-\hat{\Sigma}_{KK}\right\Vert _{\infty}\left\Vert \check{\Sigma}_{KK}^{-1}\right\Vert _{\infty}\\
 & \leq & (1-\frac{\phi}{2})\left\Vert \check{\Sigma}_{KK}-\hat{\Sigma}_{KK}\right\Vert _{\infty}\sqrt{k_{2}}\left\Vert \check{\Sigma}_{KK}^{-1}\right\Vert _{2},
\end{eqnarray*}
where the last inequality follows from the bound on $||\hat{\Sigma}_{K^{c}K}\hat{\Sigma}_{KK}^{-1}||_{\infty}$
established previously. Using bounds (\ref{eq:35}) from the proof
for Lemma A.15, we have 
\[
\left\Vert \check{\Sigma}_{KK}^{-1}\right\Vert _{2}\leq\frac{4}{\lambda_{\min}(\Sigma_{KK})}
\]
with probability at least $1-c_{1}\exp(-c_{2}\log(p\vee d))$. Next,
applying Lemma A.15 with $\varepsilon=\frac{\phi\lambda_{\min}(\Sigma_{KK})}{48(1-\frac{\phi}{2})\sqrt{k_{2}}}$
to bound the term $\left\Vert \check{\Sigma}_{KK}-\hat{\Sigma}_{KK}\right\Vert _{\infty}$
yields, 
\[
\mathbb{P}[||R_{4}||_{\infty}\geq\frac{\phi}{12}]\leq ck_{2}^{2}\cdot\exp\left(-\frac{c_{1}n}{k_{2}^{\frac{3}{2}}}\right)+c_{2}\exp(-c_{3}\log d)
\]
For the fifth term $R_{5}$, using bound (\ref{eq:32}) from the proof
for Lemma A.14, we have 
\begin{eqnarray*}
||R_{5}||_{\infty} & \leq & \sqrt{k_{2}}\left\Vert \hat{\Sigma}_{KK}^{-1}\right\Vert _{2}\left\Vert \check{\Sigma}_{K^{c}K}-\hat{\Sigma}_{K^{c}K}\right\Vert _{\infty}\\
 & \leq & \frac{2\sqrt{k_{2}}}{\lambda_{\min}(\Sigma_{KK})}\left\Vert \check{\Sigma}_{K^{c}K}-\hat{\Sigma}_{K^{c}K}\right\Vert _{\infty}.
\end{eqnarray*}
Applying Lemma A.15 with $\varepsilon=\frac{\phi\lambda_{\min}(\Sigma_{KK})}{24\sqrt{k_{2}}}$
to bound the term $\left\Vert \check{\Sigma}_{K^{c}K}-\hat{\Sigma}_{K^{c}K}\right\Vert _{\infty}$
yields 
\[
\mathbb{P}[||R_{5}||_{\infty}\geq\frac{\phi}{12}]\leq c(p-k_{2})k_{2}\cdot\exp\left(-\frac{c_{1}n}{k_{2}^{\frac{3}{2}}}\right)+c_{2}\exp(-c_{3}\log d)
\]
For the sixth term $R_{6}$, by applying Lemma A.15 to bound the terms
$\left\Vert \check{\Sigma}_{K^{c}K}-\hat{\Sigma}_{K^{c}K}\right\Vert _{\infty}$
and $\left\Vert \check{\Sigma}_{KK}^{-1}-\hat{\Sigma}_{KK}^{-1}\right\Vert _{\infty}$
respectively, with $\varepsilon=\frac{\phi}{12}\frac{\lambda_{\min}(\Sigma_{KK})}{8}$,
we are guaranteed that 
\[
\mathbb{P}[||R_{6}||_{\infty}\geq\frac{\phi}{12}]\leq c(p-k_{2})k_{2}\cdot\exp\left(-\frac{c_{1}n}{k_{2}}\right)+c_{2}\exp(-c_{3}\log d)
\]
Under the conditions $n\succsim\left(k_{2}^{3}\log p\right)\vee\left(k_{2}^{2}k_{1}\log d\right)$,
$\sqrt{\frac{n\log p}{\left|\beta^{*}\right|_{1}^{2}}}\vee\sqrt{nk_{1}\log d}\succsim\log p$,
and $\sqrt{\frac{k_{1}\log d}{n}}=o(1)$, putting the bounds on $R_{1}-R_{6}$
together, we conclude that 
\[
\mathbb{P}[||\check{\Sigma}_{K^{c}K}\check{\Sigma}_{KK}^{-1}||_{\infty}\geq1-\frac{\phi}{4}]\leq O\left(\frac{1}{p\wedge d}\right).
\]
$\square$\textbf{}\\
\textbf{}\\
\textbf{Lemma A.14}: Suppose Assumptions 2.1 and 4.2 hold. For any
$\varepsilon>0$ and constant $c$, we have

\begin{equation}
\mathbb{P}\left\{ \left\Vert \hat{\Sigma}_{K^{c}K}-\Sigma_{K^{c}K}\right\Vert _{\infty}\geq\varepsilon\right\} \leq(p-k_{2})k_{2}\cdot2\exp(-cn\min\{\frac{\varepsilon^{2}}{4k_{2}^{2}\sigma_{v}^{4}},\,\frac{\varepsilon}{2k_{2}\sigma_{v}^{2}}\}),\label{eq:26}
\end{equation}
\begin{equation}
\mathbb{P}\left\{ \left\Vert \hat{\Sigma}_{KK}-\Sigma_{KK}\right\Vert _{\infty}\geq\varepsilon\right\} \leq k_{2}^{2}\cdot2\exp(-cn\min\{\frac{\varepsilon^{2}}{4k_{2}^{2}\sigma_{v}^{4}},\,\frac{\varepsilon}{2k_{2}\sigma_{v}^{2}}\}).\label{eq:27}
\end{equation}
Furthermore, under the scaling $n\succsim k_{2}\log p$, for constants
$b_{1}$ and $b_{2}$, we have

\begin{equation}
\left\Vert \hat{\Sigma}_{KK}^{-1}-\Sigma_{KK}^{-1}\right\Vert _{\infty}\leq\frac{1}{\lambda_{\min}(\Sigma_{KK})},\label{eq:28}
\end{equation}
with probability at least $1-c_{1}\exp(-c_{2}n\min\{\frac{\lambda_{\min}^{2}(\Sigma_{KK})}{4k_{2}\sigma_{v}^{4}},\,\frac{\lambda_{\min}(\Sigma_{KK})}{2\sqrt{k_{2}}\sigma_{v}^{2}}\})$.\textbf{}\\
\textbf{}\\
\textbf{Proof}. Denote the element $(j^{'},\, j)$ of the matrix difference
$\hat{\Sigma}_{K^{c}K}-\Sigma_{K^{c}K}$ by $\omega_{j^{'}j}$. By
the definition of the $l_{\infty}$matrix norm, we have 
\begin{eqnarray*}
\mathbb{P}\left\{ \left\Vert \hat{\Sigma}_{K^{c}K}-\Sigma_{K^{c}K}\right\Vert _{\infty}\geq\varepsilon\right\}  & = & \mathbb{P}\left\{ \max_{j^{'}\in K^{c}}\sum_{j\in K}|\omega_{j^{'}j}|\geq\varepsilon\right\} \\
 & \leq & (p-k_{2})\mathbb{P}\left\{ \sum_{j\in K}|\omega_{j^{'}j}|\geq\varepsilon\right\} \\
 & \leq & (p-k_{2})\mathbb{P}\left\{ \exists j\in K\,\vert\,|\omega_{j^{'}j}|\geq\frac{\varepsilon}{k_{2}}\right\} \\
 & \leq & (p-k_{2})k_{2}\mathbb{P}\left\{ |\omega_{j^{'}j}|\geq\frac{\varepsilon}{k_{2}}\right\} \\
 & \leq & (p-k_{2})k_{2}\cdot2\exp(-cn\min\{\frac{\varepsilon^{2}}{k_{2}^{2}\sigma_{v}^{4}},\,\frac{\varepsilon}{k_{2}\sigma_{v}^{2}}\}),
\end{eqnarray*}
where the last inequality follows the deviation bound for sub-exponential
random variables, i.e., Lemma B.1. Bound (\ref{eq:27}) can be obtained
in a similar way except that the pre-factor $(p-k_{2})$ is replaced
by $k_{2}$. To prove the last bound (\ref{eq:28}), write 
\begin{eqnarray}
\left\Vert \hat{\Sigma}_{KK}^{-1}-\Sigma_{KK}^{-1}\right\Vert _{\infty} & = & \left\Vert \Sigma_{KK}^{-1}\left[\Sigma_{KK}-\hat{\Sigma}_{KK}\right]\hat{\Sigma}_{KK}^{-1}\right\Vert _{\infty}\nonumber \\
 & = & \sqrt{k_{2}}\left\Vert \Sigma_{KK}^{-1}\left[\Sigma_{KK}-\hat{\Sigma}_{KK}\right]\hat{\Sigma}_{KK}^{-1}\right\Vert _{2}\nonumber \\
 & = & \sqrt{k_{2}}\left\Vert \Sigma_{KK}^{-1}\right\Vert _{2}\left\Vert \Sigma_{KK}-\hat{\Sigma}_{KK}\right\Vert _{2}\left\Vert \hat{\Sigma}_{KK}^{-1}\right\Vert _{2}\nonumber \\
 & \leq & \frac{\sqrt{k_{2}}}{\lambda_{\min}(\Sigma_{KK})}\left\Vert \Sigma_{KK}-\hat{\Sigma}_{KK}\right\Vert _{2}\left\Vert \hat{\Sigma}_{KK}^{-1}\right\Vert _{2}.\label{eq:29}
\end{eqnarray}
To bound the term $\left\Vert \Sigma_{KK}-\hat{\Sigma}_{KK}\right\Vert _{2}$
in (\ref{eq:29}), applying Lemma B.1 with $v^{T}v=\hat{\Sigma}_{KK}$
and $t=\frac{\lambda_{\min}(\Sigma_{KK})}{2\sqrt{k_{2}}}$ yields
\[
\left\Vert \hat{\Sigma}_{KK}-\Sigma_{KK}\right\Vert _{2}\leq\frac{\lambda_{\min}(\Sigma_{KK})}{2\sqrt{k_{2}}},
\]
with probability at least $1-c_{1}\exp(-c_{2}n\min\{\frac{\lambda_{\min}^{2}(\Sigma_{KK})}{4k_{2}\sigma_{v}^{4}},\,\frac{\lambda_{\min}(\Sigma_{KK})}{2\sqrt{k_{2}}\sigma_{v}^{2}}\})$. 

To bound the term\textbf{$\left\Vert \hat{\Sigma}_{KK}^{-1}\right\Vert _{2}$}
in (\ref{eq:29}), note that we can write 
\begin{eqnarray}
\lambda_{\min}(\Sigma_{KK}) & = & \min_{||h^{'}||_{2}=1}h^{'T}\Sigma_{KK}h^{'}\nonumber \\
 & = & \min_{||h^{'}||_{2}=1}\left[h^{'T}\hat{\Sigma}_{KK}h^{'}+h^{'T}(\Sigma_{KK}-\hat{\Sigma}_{KK})h^{'}\right]\nonumber \\
 & \leq & h^{T}\hat{\Sigma}_{KK}h+h^{T}(\Sigma_{KK}-\hat{\Sigma}_{KK})h\label{eq:30}
\end{eqnarray}
where $h\in\mathbb{R}^{k_{2}}$ is a unit-norm minimal eigenvector
of $\hat{\Sigma}_{KK}$. Applying Lemma B.1 yields 
\[
\left|h^{T}\left(\Sigma_{KK}-\hat{\Sigma}_{KK}\right)h\right|\leq\frac{\lambda_{\min}(\Sigma_{KK})}{2}
\]
with probability at least $1-c_{1}\exp(-c_{2}n)$. Therefore, 
\begin{eqnarray}
\lambda_{\min}(\Sigma_{KK}) & \leq & \lambda_{\min}(\hat{\Sigma}_{KK})+\frac{\lambda_{\min}(\Sigma_{KK})}{2}\nonumber \\
\Longrightarrow\lambda_{\min}(\hat{\Sigma}_{KK}) & \geq & \frac{\lambda_{\min}(\Sigma_{KK})}{2},\label{eq:31}
\end{eqnarray}
and consequently,
\begin{equation}
\left\Vert \hat{\Sigma}_{KK}^{-1}\right\Vert _{2}\leq\frac{2}{\lambda_{\min}(\Sigma_{KK})}.\label{eq:32}
\end{equation}
Putting everything together, we have 
\[
\left\Vert \hat{\Sigma}_{KK}^{-1}-\Sigma_{KK}^{-1}\right\Vert _{\infty}\leq\frac{\sqrt{k_{2}}}{\lambda_{\min}(\Sigma_{KK})}\frac{\lambda_{\min}(\Sigma_{KK})}{2\sqrt{k_{2}}}\frac{2}{\lambda_{\min}(\Sigma_{KK})}=\frac{1}{\lambda_{\min}(\Sigma_{KK})}.
\]
with probability at least $1-c_{1}\exp(-c_{2}n\min\{\frac{\lambda_{\min}^{2}(\Sigma_{KK})}{4k_{2}\sigma_{v}^{4}},\,\frac{\lambda_{\min}(\Sigma_{KK})}{2\sqrt{k_{2}}\sigma_{v}^{2}}\})$.\textbf{
$\square$}\\
\\
\textbf{Lemma A.15}:\textbf{ }Suppose the assumptions in Lemma A.13
hold. For any $\varepsilon>0$, under the condition $n\succsim\left(k_{2}^{3}\log p\right)\vee\left(k_{2}^{2}k_{1}\log d\right)$,
$\sqrt{\frac{n\log p}{\left|\beta^{*}\right|_{2}}}\vee\sqrt{nk_{1}\log d}\succsim\log p$,
and $\sqrt{\frac{k_{1}\log d}{n}}=o(1)$, we have 
\[
\mathbb{P}\left\{ \left\Vert \check{\Sigma}_{K^{c}K}-\hat{\Sigma}_{K^{c}K}\right\Vert _{\infty}\geq\varepsilon\right\} \leq c(p-k_{2})k_{2}\cdot\exp\left(-\frac{n\varepsilon}{c_{1}k_{2}}\right)+c_{2}\exp(-c_{3}\log d)
\]
\[
\mathbb{P}\left\{ \left\Vert \check{\Sigma}_{KK}-\hat{\Sigma}_{KK}\right\Vert _{\infty}\geq\varepsilon\right\} \leq c^{'}k_{2}^{2}\cdot\exp\left(-\frac{n\varepsilon}{c_{1}k_{2}}\right)+c_{2}\exp(-c_{3}\log d)
\]
Furthermore, we have

\begin{equation}
\left\Vert \check{\Sigma}_{KK}^{-1}-\hat{\Sigma}_{KK}^{-1}\right\Vert _{\infty}\leq\frac{8}{\lambda_{\min}(\Sigma_{KK})}\quad\textrm{with probability at least }1-c_{1}\exp(-c_{2}\log(p\wedge d)).\label{eq:33}
\end{equation}
\textbf{Proof}. Denote the element $(j^{'},\, j)$ of the matrix difference
$\check{\Sigma}{}_{K^{c}K}-\hat{\Sigma}_{K^{c}K}$ by $\omega_{j^{'}j}$.
Using the same argument as in the proof for Lemma A.14, we have 
\[
\mathbb{P}\left\{ \left\Vert \check{\Sigma}{}_{K^{c}K}-\hat{\Sigma}_{K^{c}K}\right\Vert _{\infty}\geq\varepsilon\right\} \leq(p-k_{2})k_{2}\mathbb{P}\left\{ |\omega_{j^{'}j}|\geq\frac{\varepsilon}{k_{2}}\right\} .
\]
Following the derivation of the upper bounds on $\left|\frac{(\hat{x}-\tilde{x})^{T}\tilde{x}}{n}\right|_{\infty}$
and $\left|\frac{(\hat{x}-\tilde{x})^{T}(\hat{x}-\tilde{x})}{n}\right|_{\infty}$
in the proofs for Lemma A.3 and Corollary 4.4 and the identity 
\[
\frac{1}{n}\left(\check{\Sigma}_{K^{c}K}-\hat{\Sigma}_{K^{c}K}\right)=\frac{1}{n}v_{K^{c}}^{T}(\hat{v}_{K}-v_{K})+\frac{1}{n}(\hat{v}_{K^{c}}-v_{K^{c}})^{T}v_{K}+\frac{1}{n}(\hat{v}_{K^{c}}-v_{K^{c}})^{T}(\hat{v}_{K}-v_{K}),
\]
we notice that to upper bound $|\omega_{j^{'}j}|$, it suffices to
upper bound $3\cdot\left|\frac{1}{n}v_{j^{'}}^{T}(\hat{v}_{j}-v_{j})\right|$.
From the proofs for Lemma A.3 and Corollary 4.4, we have 
\[
\left|\frac{1}{n}v_{j^{'}}^{T}(\hat{v}_{j}-v_{j})\right|\leq b\left(u^{2}+\sqrt{\frac{k_{1}\log d}{n}}\right),
\]
with probability at least $1-c_{1}\exp\left(-nC_{j}^{*}u^{2}\right)-O\left(\frac{1}{d}\right)$.
If $n\gtrsim k_{2}^{2}k_{1}\log d$, setting $u=\sqrt{\frac{\varepsilon}{bk_{2}}}$
for any $\varepsilon>0$ yields 
\[
\mathbb{P}\left[|\omega_{j^{'}j}|\geq\frac{\varepsilon}{k_{2}}\right]\leq c\exp\left(-\frac{n\varepsilon}{c_{1}k_{2}}\right)+c_{2}\exp(-c_{3}\log d).
\]
Therefore, 
\[
\mathbb{P}\left\{ \left\Vert \check{\Sigma}_{K^{c}K}-\hat{\Sigma}_{K^{c}K}\right\Vert _{\infty}\geq\varepsilon\right\} \leq c(p-k_{2})k_{2}\cdot\exp\left(-\frac{n\varepsilon}{c_{1}k_{2}}\right)+c_{2}\exp(-c_{3}\log d).
\]
The second bound in Lemma A.15 can be obtained in a similar way except
that the pre-factor $(p-k_{2})$ is replaced by $k_{2}$. 

To prove the third bound in Lemma A.15, by applying the same argument
as in the proof for Lemma A.14, we have 
\begin{eqnarray*}
\left\Vert \check{\Sigma}_{KK}^{-1}-\hat{\Sigma}_{KK}^{-1}\right\Vert _{\infty} & \leq & \frac{\sqrt{k_{2}}}{\lambda_{\min}(\hat{\Sigma}_{KK})}\left\Vert \hat{\Sigma}_{KK}-\check{\Sigma}_{KK}\right\Vert _{2}\left\Vert \check{\Sigma}_{KK}^{-1}\right\Vert _{2}\\
 & \leq & \frac{2\sqrt{k_{2}}}{\lambda_{\min}(\Sigma_{KK})}\left\Vert \hat{\Sigma}_{KK}-\check{\Sigma}_{KK}\right\Vert _{2}\left\Vert \check{\Sigma}_{KK}^{-1}\right\Vert _{2},
\end{eqnarray*}
where the last inequality comes from bound (\ref{eq:31}). 

For the term $\left\Vert \hat{\Sigma}_{KK}-\check{\Sigma}_{KK}\right\Vert _{2}$,
setting $\varepsilon=\frac{\lambda_{\min}(\Sigma_{KK})}{\sqrt{k_{2}}}$
in the previous upper bound on $\left\Vert \hat{\Sigma}_{KK}-\check{\Sigma}_{KK}\right\Vert _{\infty}$
yields 
\begin{eqnarray}
\left\Vert \hat{\Sigma}_{KK}-\check{\Sigma}_{KK}\right\Vert _{2} & \leq & \left\Vert \hat{\Sigma}_{KK}-\check{\Sigma}_{KK}\right\Vert _{\infty}\nonumber \\
 & \leq & \frac{\lambda_{\min}(\Sigma_{KK})}{\sqrt{k_{2}}},\label{eq:34}
\end{eqnarray}
with probability at least $1-ck_{2}^{2}\cdot\exp\left(-\frac{n}{c_{1}k_{2}^{3/2}}\right)-c_{2}\exp(-c_{3}\log d)$.

To bound the term\textbf{ $\left\Vert \check{\Sigma}_{KK}^{-1}\right\Vert _{2}$},
again we have, 
\begin{eqnarray*}
\lambda_{\min}(\hat{\Sigma}_{KK}) & \leq & h^{T}\check{\Sigma}_{KK}h+h^{T}(\hat{\Sigma}_{KK}-\check{\Sigma}_{KK})h\\
 & \leq & h^{T}\check{\Sigma}_{KK}h+k_{2}\left|\hat{\Sigma}_{KK}-\check{\Sigma}_{KK}\right|_{\infty}\\
 & \leq & h^{T}\check{\Sigma}_{KK}h+b^{'}\left\{ \left(\frac{1}{n}\right)^{\frac{2}{3}}+\sqrt{\frac{k_{1}\log d}{n}}\right\} ,
\end{eqnarray*}
where $h\in\mathbb{R}^{k_{2}}$ is a unit-norm minimal eigenvector
of $\check{\Sigma}_{KK}$. The last inequality follows from the bounds
on $\left|\frac{(\hat{v}-v)^{T}v}{n}\right|_{\infty}$ and $\left|\frac{(\hat{v}-v)^{T}(\hat{v}-v)}{n}\right|_{\infty}$
from the proofs for Lemma A.3 and Corollary 4.4. Therefore, if $\sqrt{\frac{k_{1}\log d}{n}}=o(1)$,
then we have 

\[
\lambda_{\min}(\check{\Sigma}_{KK})\geq\frac{\lambda_{\min}(\hat{\Sigma}_{KK})}{2}
\]

\begin{eqnarray}
\Longrightarrow\left\Vert \check{\Sigma}_{KK}^{-1}\right\Vert _{2} & \leq & \frac{2}{\lambda_{\min}(\hat{\Sigma}_{KK})}\nonumber \\
 & \leq & \frac{4}{\lambda_{\min}(\Sigma_{KK})},\label{eq:35}
\end{eqnarray}
where the last inequality follows from bound (\ref{eq:31}) in the
proof for Lemma A.14. Putting everything together, we have 
\[
\left\Vert \hat{\Sigma}_{KK}^{-1}-\tilde{\Sigma}_{KK}^{-1}\right\Vert _{\infty}\leq\frac{2\sqrt{k_{2}}}{\lambda_{\min}(\Sigma_{KK})}\frac{\lambda_{\min}(\Sigma_{KK})}{\sqrt{k_{2}}}\frac{4}{\lambda_{\min}(\Sigma_{KK})}=\frac{8}{\lambda_{\min}(\Sigma_{KK})}.
\]
with probability at least $1-O\left(\frac{1}{d\wedge p}\right)$.\textbf{
$\square$}\\
\textbf{}\\
\textbf{Lemma A.16}: Suppose the conditions in Corollary 4.4 hold.
With the choice of the tuning parameter 
\[
\lambda_{n,3}\geq c\frac{8(2-\frac{\phi}{4})}{\phi}\left[\left(\sigma_{v}\sigma_{\eta}\sqrt{\frac{\log p}{n}}\right)\vee\left(L\left|\beta^{*}\right|_{1}\sqrt{\Upsilon_{w,\theta^{*}}}b\left(\sigma_{v},\,\sigma_{\eta}\right)\sqrt{\frac{k_{1}\log d}{n}}\right)\right]
\]
for some sufficiently large constant $c>0$, under the conditions
$n\succsim\left(k_{2}\log p\right)\vee\left(k_{2}^{2}k_{1}\log d\right)$
and $\sqrt{\frac{n\log p}{\left|\beta^{*}\right|_{1}^{2}}}\vee\sqrt{nk_{1}\log d}\succsim\log p$,
then, we have 
\[
\left(2-\frac{\phi}{4}\right)\left|\hat{v}^{T}\frac{e}{n\lambda_{n,3}}\right|_{\infty}\leq\frac{\phi}{8},
\]
with probability at least $1-c_{1}\exp(-c_{2}\log(d\wedge p))$. \\
\textbf{}\\
\textbf{Proof}. Recall from the proofs for Lemma A.4 and Corollary
4.4 on $|\frac{1}{n\lambda_{n,3}}\hat{v}^{T}e|_{\infty}$, so as long
as 
\[
\lambda_{n,3}\geq c\frac{8(2-\frac{\phi}{4})}{\phi}\left[\left(\sigma_{v}\sigma_{\eta}\sqrt{\frac{\log p}{n}}\right)\vee\left(L\left|\beta^{*}\right|_{1}\sqrt{\Upsilon_{w,\theta^{*}}}b\left(\sigma_{v},\,\sigma_{\eta}\right)\sqrt{\frac{k_{1}\log d}{n}}\right)\right]
\]
for some sufficiently large constant $c>0$, under the conditions
$n\succsim\left(k_{2}\log p\right)\vee\left(k_{2}^{2}k_{1}\log d\right)$
and $\sqrt{\frac{n\log p}{\left|\beta^{*}\right|_{1}^{2}}}\vee\sqrt{nk_{1}\log d}\succsim\log p$,
we have

\[
\left(2-\frac{\phi}{4}\right)\left|\hat{v}^{T}\frac{e}{n\lambda_{n,3}}\right|_{\infty}\leq\frac{\phi}{8},
\]
with probability at least $1-c_{1}\exp(-c_{2}\log(d\wedge p))$.$\square$

\section{Appendix II: Technical lemmas and the proofs }

\textbf{Lemma B.1}: Let $X\in\mathbb{R}^{n\times p_{1}}$ be a zero-mean
sub-Gaussian matrix with parameters $(\Sigma_{X},\,\sigma_{X}^{2})$.
For any fixed (unit) vector $\Delta\in\mathbb{R}^{p_{1}}$, we have
\[
\mathbb{P}\left[\left|\frac{\left|X\Delta\right|{}_{2}^{2}}{n}-\mathbb{E}[\frac{\left|X\Delta\right|{}_{2}^{2}}{n}]\right|\geq\varepsilon\right]\leq2\exp\left(-c_{0}n\min\left\{ \frac{\varepsilon^{2}}{\sigma_{X}^{4}},\,\frac{\varepsilon}{\sigma_{X}^{2}}\right\} \right).
\]
Moreover, if $Y\in\mathbb{R}^{n\times p_{2}}$ is a zero-mean sub-Gaussian
matrix with parameters $(\Sigma_{Y},\,\sigma_{Y}^{2})$, then 
\[
\mathbb{P}\left[\left|\frac{Y^{T}X}{n}-\textrm{cov}(Y_{i},\, X_{i})\right|{}_{\infty}\geq\varepsilon\right]\leq6\exp\left(-c_{1}n\min\left\{ \frac{\varepsilon^{2}}{\sigma_{X}^{2}\sigma_{Y}^{2}},\,\frac{\varepsilon}{\sigma_{X}\sigma_{Y}}\right\} +\log p_{1}+\log p_{2}\right),
\]
where $X_{i}$ and $Y_{i}$ are the $i^{th}$ rows of $X$ and $Y$,
respectively. \\
\\
\textbf{Remark. }This lemma is Lemma 14 in Loh and Wainwright (2012).
$\square$\\
\textbf{}\\
\textbf{Lemma B.2}: Let $X\in\mathbb{R}^{n\times p_{1}}$ be a sub-Gaussian
matrix with parameters $(\Sigma_{X},\,\sigma_{X}^{2})$. We have 
\begin{eqnarray*}
\frac{\left|X\Delta\right|{}_{2}^{2}}{n} & \geq & \frac{\underbar{\ensuremath{\alpha}}}{2}\left|\Delta\right|{}_{2}^{2}-\alpha^{'}\frac{\log p_{1}}{n}\left|\Delta\right|{}_{1}^{2},\qquad\textrm{for all }\Delta\in\mathbb{R}^{p_{1}}\\
\frac{\left|X\Delta\right|{}_{2}^{2}}{n} & \leq & \frac{3\bar{\alpha}}{2}\left|\Delta\right|{}_{2}^{2}+\alpha^{'}\frac{\log p_{1}}{n}\left|\Delta\right|{}_{1}^{2},\qquad\textrm{for all }\Delta\in\mathbb{R}^{p_{1}}
\end{eqnarray*}
with probability at least $1-c_{1}\exp(-c_{2}n)$, where $\underbar{\ensuremath{\alpha}}$,
$\bar{\alpha}$, and $\alpha^{'}$ only depend on $\Sigma_{X}$ and
$\sigma_{X}$.\\
\textbf{}\\
\textbf{Remark. }This lemma is Lemma 13 in Loh and Wainwright (2012).
The choice of $\underbar{\ensuremath{\alpha}}$, $\bar{\alpha}$,
and $\alpha^{'}$ depends on the problems.$\square$\\
\\
\textbf{Lemma B.3}: For all $u\geq r_{nj}$, we have 
\[
\mathbb{P}\left[\mathcal{A}_{j}(u)\,\vert\,\left\{ w_{i}^{T}\theta^{*}\right\} _{i=1}^{n}\right]\leq\exp\left(-nC_{j}^{*}u^{2}\right)
\]
 where 
\[
\mathcal{A}_{j}(u):=\left\{ \exists f\in\mathcal{\bar{F}}_{j}\cap\left\{ \left|f_{\theta^{*}}\right|_{n}\geq u\right\} \,:\,\frac{1}{n}\left|\sum_{i=1}^{n}v_{ij}f(w_{i}^{T}\theta^{*})\right|\geq2\left|f_{\theta^{*}}\right|_{n}u\right\} ,
\]
and 
\begin{equation}
C_{j}^{*}=c\frac{\gamma_{v_{j}}\wedge(B_{v_{j}}^{2}\vee B_{\eta}^{2})^{-1}}{\sigma_{v_{j}}^{2}\vee\sigma_{\eta}^{2}}\label{eq:36}
\end{equation}
when $v_{j}$ is either (i) sub-Gaussian with \textit{strongly log-concave}
distribution for some $\gamma_{v_{j}}>0$; or, (ii) a bounded vector
such that for every $i=1,...,n$, $v_{ij}$ is supported on the interval
$(a_{v_{j}}^{'},\, a_{v_{j}}^{''})$ with $B_{v_{j}}:=a_{v_{j}}^{''}-a_{v_{j}}^{'}$;
or, (iii) a mixture of (i) and (ii) in terms of its probability measure.
\\
\textbf{}\\
\textbf{Proof.} This proof is a modification of the proofs for Corollary
8.3 in van de Geer (2000) and Lemma 13.2 in Wainwright (2015). Suppose
that there exists some $f\in\mathcal{\bar{F}}_{j}$ with $\left|f_{\theta^{*}}\right|_{n}\geq u$
such that 
\begin{equation}
\frac{1}{n}\left|\sum_{i=1}^{n}v_{ij}f(w_{i}^{T}\theta^{*})\right|\geq2\left|f_{\theta^{*}}\right|_{n}u.\label{eq:37}
\end{equation}
Defining the function $\tilde{f}_{\theta^{*}}:=\frac{u}{\left|f_{\theta^{*}}\right|_{n}}f_{\theta^{*}}$,
observe that $\left|\tilde{f}_{\theta^{*}}\right|_{n}=u$. Since $f\in\mathcal{\bar{F}}_{j}$
and $u\leq\left|f_{\theta^{*}}\right|_{n}$ by construction, the star-shaped
assumption implies that $\tilde{f}\in\mathcal{\bar{F}}_{j}$. Consequently,
we have shown that if there exists a function $f$ satisfying inequality
(\ref{eq:37}), which occurs whenever the event $\mathcal{A}_{j}(u)$
is true, then there exists a function $\tilde{f}\in\mathcal{\bar{F}}_{j}$
with $\left|\tilde{f}_{\theta^{*}}\right|_{n}=u$ such that 
\[
\frac{1}{n}\left|\sum_{i=1}^{n}v_{ij}\tilde{f}(w_{i}^{T}\theta^{*})\right|=\frac{u}{\left|f_{\theta^{*}}\right|_{n}}\frac{1}{n}\left|\sum_{i=1}^{n}v_{ij}f(w_{i}^{T}\theta^{*})\right|\geq2u^{2}.
\]
Summarizing then, we have established the inequality 
\[
\mathbb{P}\left[\mathcal{A}_{j}(u)\,\vert\,\left\{ w_{i}^{T}\theta^{*}\right\} _{i=1}^{n}\right]\leq\mathbb{P}\left[A_{nj}(u)\geq2u^{2}\,\vert\,\left\{ w_{i}^{T}\theta^{*}\right\} _{i=1}^{n}\right],
\]
where 

\[
A_{nj}(u):=\sup_{\tilde{f}(\theta^{*})\in\Omega(u;\,\mathcal{F}_{j})}\frac{1}{n}\left|\sum_{i=1}^{n}v_{ij}\tilde{f}(w_{i}^{T}\theta^{*})\right|
\]
where 
\[
\Omega(u;\,\mathcal{F}_{j})=\left\{ f\,:\, f\in\mathcal{\bar{F}}_{j},\,\left|f_{\theta^{*}}\right|_{n}\leq u\right\} .
\]
By construction, $\mathbb{E}(v_{ij}\vert w_{i}^{T}\theta^{*})=0$
for all $i=1,...,n$. Consequently, conditioning on $\left\{ w_{i}^{T}\theta^{*}\right\} _{i=1}^{n}$,
for each fixed $\tilde{f}$, the variable $\frac{1}{n}\left|\sum_{i=1}^{n}v_{ij}\tilde{f}(w_{i}^{T}\theta^{*})\right|$
is zero-mean sub-Gaussian, so that the variable $A_{nj}(u)$ is the
supremum of a sub-Gaussian process. If we view this supremum as a
function of a standardized sub-Gaussian vector, then by standard maximal
inequality (e.g., Corollary 8.3 in van de Geer, 2000), we have 
\[
\mathbb{P}\left[A_{nj}(u)\geq\mathbb{E}_{v_{j}}\left[A_{nj}(u)\,\vert\,\left\{ w_{i}^{T}\theta^{*}\right\} _{i=1}^{n}\right]+b\wedge2\sigma_{v_{j}}^{2}\,\vert\,\left\{ w_{i}^{T}\theta^{*}\right\} _{i=1}^{n}\right]\leq\exp\left(-\frac{nC_{j}^{*}b^{2}}{u^{2}}\right).
\]
Consequently, for any $b=u^{2}$, 
\begin{equation}
\mathbb{P}\left[A_{nj}(u)\geq\mathbb{E}_{v_{j}}\left[A_{nj}(u)\,\vert\,\left\{ w_{i}^{T}\theta^{*}\right\} _{i=1}^{n}\right]+u^{2}\,\vert\,\left\{ w_{i}^{T}\theta^{*}\right\} _{i=1}^{n}\right]\leq\exp\left(-nC_{j}^{*}u^{2}\right).\label{eq:38}
\end{equation}
Now let us look at special cases that yield sharp constant $C_{j}^{*}$.\\
\textbf{\textit{}}\\
\textbf{\textit{Special case }}\textbf{(}\textbf{\textit{i}}\textbf{)}:
If $v_{j}$ is sub-Gaussian with \textit{strongly log-concave} distribution
for some $\gamma_{v_{j}}>0$, then by Cauchy-Schwarz inequality, it
can be verified that the associated Lipschitz constant is at most
$\frac{\sigma_{v_{j}}u}{\sqrt{n}}$. Consequently, by Lemma B.6, for
any $b>0$, we have the sub-Gaussian tail bound 
\[
\mathbb{P}\left[A_{nj}(u)\geq\mathbb{E}_{v_{j}}\left[A_{nj}(u)\,\vert\,\left\{ w_{i}^{T}\theta^{*}\right\} _{i=1}^{n}\right]+b\,\vert\,\left\{ w_{i}^{T}\theta^{*}\right\} _{i=1}^{n}\right]\leq\exp\left(-\frac{n\gamma_{v_{j}}b^{2}}{2u^{2}\sigma_{v_{j}}^{2}}\right).
\]
Setting $b=u^{2}$ yields 
\[
\mathbb{P}\left[A_{nj}(u)\geq\mathbb{E}_{v_{j}}\left[A_{nj}(u)\,\vert\,\left\{ w_{i}^{T}\theta^{*}\right\} _{i=1}^{n}\right]+u^{2}\,\vert\,\left\{ w_{i}^{T}\theta^{*}\right\} _{i=1}^{n}\right]\leq\exp\left(-\frac{n\gamma_{v_{j}}u^{2}}{2\sigma_{v_{j}}^{2}}\right).
\]
\textbf{\textit{Special case }}\textbf{(}\textbf{\textit{ii}}\textbf{)}:
If $v_{j}$ is a bounded vector such that for every $i=1,...,n$,
$v_{ij}$ is supported on the interval $(a_{v_{j}}^{'},\, a_{v_{j}}^{''})$
with $B_{v_{j}}:=a_{v_{j}}^{''}-a_{v_{j}}^{'}$, then by Lemma B.7,
following similar argument as above, we have 
\[
\mathbb{P}\left[A_{nj}(u)\geq\mathbb{E}_{v_{j}}\left[A_{nj}(u)\,\vert\,\left\{ w_{i}^{T}\theta^{*}\right\} _{i=1}^{n}\right]+u^{2}\,\vert\,\left\{ w_{i}^{T}\theta^{*}\right\} _{i=1}^{n}\right]\leq\exp\left(-\frac{nu^{2}}{2\sigma_{v_{j}}^{2}B_{v_{j}}^{2}}\right).
\]
Finally, by definition of $A_{nj}(u)$ and $\mathcal{G}_{n}(u;\,\mathcal{F}_{j})$,
we have $\mathbb{E}_{v_{j}}\left[A_{nj}(u)\,\vert\,\left\{ w_{i}^{T}\theta^{*}\right\} _{i=1}^{n}\right]=\sigma_{v_{j}}\mathcal{G}_{n}(u;\,\mathcal{F}_{j})$.
By Lemma B.8, the function $t\mapsto\frac{\mathcal{G}_{n}(t;\,\mathcal{F}_{j})}{t}$
is non-decreasing, and since $u\geq r_{nj}$ by assumption, we have
\[
\sigma_{v_{j}}\frac{\mathcal{G}_{n}(u;\,\mathcal{F}_{j})}{u}\leq\sigma_{v_{j}}\frac{\mathcal{G}_{n}(r_{nj};\,\mathcal{F}_{j})}{r_{nj}}\leq r_{nj},
\]
where the last inequality uses the definition of $r_{nj}$. Putting
everything together, we have shown that $\mathbb{E}_{v_{j}}\left[A_{nj}(u)\,\vert\,\left\{ w_{i}^{T}\theta^{*}\right\} _{i=1}^{n}\right]\leq ur_{nj}$.
Combined with the tail bound (\ref{eq:38}), respectively, we obtain
\begin{eqnarray*}
\mathbb{P}\left[A_{nj}(u)\geq2u^{2}\,\vert\,\left\{ w_{i}^{T}\theta^{*}\right\} _{i=1}^{n}\right] & \leq & \mathbb{P}\left[A_{nj}(u)\geq ur_{nj}+u^{2}\,\vert\,\left\{ w_{i}^{T}\theta^{*}\right\} _{i=1}^{n}\right]\\
 & \leq & \exp\left(-nC_{j}^{*}u^{2}\right)
\end{eqnarray*}
where the inequality uses the fact that $u^{2}\geq ur_{nj}$. $\square$\\
\\
\textbf{Lemma B.4}: There are universal positive constants $(c_{1},\, c_{2})$
such that for all $\Delta_{j}\in\left\{ f\in\mathcal{\bar{F}}_{j}\,:\,\left|f\right|_{\mathcal{F}_{j}}\leq3,\,\left|f_{\theta^{*}}\right|_{n}\geq u\right\} $,
\begin{equation}
\frac{1}{n}\left|\sum_{i=1}^{n}v_{ij}\Delta_{j}(w_{i}^{T}\theta^{*})\right|\leq6r_{nj}\left|\Delta_{j}(\theta^{*})\right|_{n}+\frac{1}{32}\left|\Delta_{j}(\theta^{*})\right|_{n}^{2}\label{eq:39}
\end{equation}
with probability at least $1-c_{1}\exp\left(-nC_{j}^{*}r_{nj}^{2}\right)$
and $C_{j}^{*}$ follows (\ref{eq:36}) under the same special cases
in Lemma B.3.\\
\\
\textbf{Proof.} To establish bound (\ref{eq:39}), we first consider
over the ball 
\[
\left\{ \Delta_{j}\in\mathcal{\bar{F}}_{j}\,:\,\left|\Delta_{j}\right|_{\mathcal{F}_{j}}\leq3\right\} \cap\left\{ \Delta_{j}\in\mathcal{\bar{F}}_{j}\,:\,\left|\Delta_{j}(\theta^{*})\right|_{n}\leq u\right\} 
\]
for some fixed radius $u\geq r_{nj}$. Later in the proof we extend
the bound to one that is uniform in $\left|\Delta_{j}(\theta^{*})\right|_{n}$
via a peeling argument. Define the random variable 
\[
A_{nj}^{'}(u):=\sup_{\Delta_{j}\in\Omega^{'}(u;\,\mathcal{F}_{j})}\frac{1}{n}\left|\sum_{i=1}^{n}v_{ij}\Delta_{j}(w_{i}^{T}\theta^{*})\right|.
\]
where $\Omega^{'}(u;\,\mathcal{F}_{j})=\left\{ \Delta_{j}\in\mathcal{\bar{F}}_{j}\,:\,\left|\Delta_{j}\right|_{\mathcal{F}_{j}}\leq3,\,\left|\Delta_{j}(\theta^{*})\right|_{n}\leq u\right\} $.
Following the same argument in Lemma B.3, we obtain 
\begin{equation}
\mathbb{P}\left[A_{nj}^{'}(u)\geq\mathbb{E}_{v_{j}}\left[A_{nj}^{'}(u)\,\vert\,\left\{ w_{i}^{T}\theta^{*}\right\} _{i=1}^{n}\right]+b\,\vert\,\left\{ w_{i}^{T}\theta^{*}\right\} _{i=1}^{n}\right]\leq\exp\left(-\frac{nC_{j}^{*}b^{2}}{u^{2}}\right).\label{eq:40}
\end{equation}
We first derive a bound for $u=r_{nj}$. By definition of $A_{nj}^{'}(u)$
and $\mathcal{G}_{n}\left(r_{nj};\,\mathcal{F}_{j}\right)$ and the
critical radius $r_{nj}$, we have 
\[
\mathbb{E}_{v_{j}}\left[A_{nj}^{'}(r_{nj})\,\vert\,\left\{ w_{i}^{T}\theta^{*}\right\} _{i=1}^{n}\right]=\tilde{\sigma}_{v_{j}}\mathcal{G}_{n}\left(r_{nj};\,\mathcal{F}_{j}\right)\leq3r_{nj}^{2}.
\]
Setting $b=3r_{nj}^{2}$ in the tail bound (\ref{eq:40}) yields 
\[
\mathbb{P}\left[A_{nj}^{'}(r_{nj})\geq6r_{nj}^{2}\,\vert\,\left\{ w_{i}^{T}\theta^{*}\right\} _{i=1}^{n}\right]\leq\exp\left(-nC_{j}^{*}r_{nj}^{2}\right).
\]
On the other hand, for any $u>r_{nj}$, we have 
\begin{eqnarray*}
\mathbb{E}_{v_{j}}\left[A_{nj}^{'}(u)\,\vert\,\left\{ w_{i}^{T}\theta^{*}\right\} _{i=1}^{n}\right] & = & \tilde{\sigma}_{v_{j}}\mathcal{G}_{n}\left(u;\,\mathbb{B}_{\mathcal{F}_{j}}(3)\right)=u\frac{\tilde{\sigma}_{v_{j}}\mathcal{G}_{n}\left(u;\,\mathbb{B}_{\mathcal{F}_{j}}(3)\right)}{u}\\
 & \leq & u\frac{\tilde{\sigma}_{v_{j}}\mathcal{G}_{n}\left(r_{nj};\,\mathbb{B}_{\mathcal{F}_{j}}(3)\right)}{r_{nj}}\leq3ur_{nj},
\end{eqnarray*}
where the first inequality follows from Lemma B.8 and the second inequality
follows from the critical radius $r_{nj}$. Setting $b=\frac{u^{2}}{128}$
in the tail bound (\ref{eq:40}) yields 
\begin{equation}
\mathbb{P}\left[A_{nj}^{'}(u)\geq3ur_{nj}+\frac{u^{2}}{128}\,\vert\,\left\{ w_{i}^{T}\theta^{*}\right\} _{i=1}^{n}\right]\leq\exp\left(-nC_{j}^{*}u^{2}\right).\label{eq:41}
\end{equation}
It remains to prove the bound (\ref{eq:39}) via a peeling argument.
Let $\mathcal{E}$ denote the event that the bound (\ref{eq:39})
is violated for some function $\Delta_{j}$ with $\left|\Delta_{j}\right|_{\mathcal{F}_{j}}\leq3$.
For real numbers $0\leq a_{1}<a_{2}$, let $\mathcal{E}(a_{1},\, a_{2})$
denote the event that it is violated for some function such that $\left|\Delta_{j}(\theta^{*})\right|_{n}\in[a_{1},\, a_{2}]$,
and $\left|\Delta_{j}\right|_{\mathcal{F}_{j}}\leq3$. For $m=0,\,1,\,2,\,...$,
define $t_{m}=2^{m}r_{nj}$. We then have the decomposition $\mathcal{E}=\mathcal{E}(0,\, t_{0})\cup\left(\cup_{m=0}^{\infty}\mathcal{E}(t_{m},\, t_{m+1})\right)$
and hence by union bound, 
\begin{equation}
\mathbb{P}\left[\mathcal{E}\right]\leq\mathbb{P}\left[\mathcal{E}(0,\, t_{0})\right]+\sum_{m=0}^{\infty}\mathbb{P}\left[\mathcal{E}(t_{m},\, t_{m+1})\right].\label{eq:42}
\end{equation}
The final step is to bound each of the terms in this summation. Since
$t_{0}=r_{nj}$, we have 
\begin{equation}
\mathbb{P}\left[\mathcal{E}(0,\, t_{0})\right]\leq\mathbb{P}\left[A_{nj}^{'}(r_{nj})\geq6r_{nj}^{2}\,\vert\,\left\{ w_{i}^{T}\theta^{*}\right\} _{i=1}^{n}\right]\leq\exp\left(-nC_{j}^{*}r_{nj}^{2}\right).\label{eq:43}
\end{equation}
On the other hand, suppose that $\mathcal{E}(t_{m},\, t_{m+1})$ holds,
meaning that there exists some function $\Delta_{j}$ with $\left|\Delta_{j}\right|_{\mathcal{F}_{j}}\leq3$
and $\left|\Delta_{j}(\theta^{*})\right|_{n}\leq t_{m+1}$ such that
\begin{eqnarray*}
\frac{1}{n}\left|\sum_{i=1}^{n}v_{ij}\Delta_{j}(w_{i}^{T}\theta^{*})\right| & > & 6r_{nj}\left|\Delta_{j}(\theta^{*})\right|_{n}+\frac{1}{32}\left|\Delta_{j}(\theta^{*})\right|_{n}^{2}\\
 & > & 6r_{nj}t_{m}+\frac{1}{32}t_{m}^{2}\\
 & = & 3r_{nj}t_{m+1}+\frac{1}{128}t_{m+1}^{2},
\end{eqnarray*}
where the second inequality follows since $\left|\Delta_{j}(\theta^{*})\right|_{n}\geq t_{m}$;
and the third inequality follows since $t_{m+1}=2t_{m}$. This lower
bound implies that 
\[
A_{nj}^{'}(t_{m+1})\geq3r_{nj}t_{m+1}+\frac{1}{128}t_{m+1}^{2},
\]
hence the bound (\ref{eq:41}) implies that 
\[
\mathbb{P}\left[\mathcal{E}(t_{m},\, t_{m+1})\right]\leq\exp\left(-\frac{nt_{m+1}^{2}}{c_{2}\tilde{\sigma}_{v_{j}}}\right)=\exp\left(-\frac{n2^{2(m+1)}r_{nj}^{2}}{c_{2}\tilde{\sigma}_{v_{j}}}\right).
\]
Combining this tail bound with our earlier bound (\ref{eq:43}), and
substituting into the union bound (\ref{eq:42}) yields 
\[
\mathbb{P}\left[\mathcal{E}\right]\leq\exp\left(-\frac{n\gamma_{v_{j}}r_{nj}^{2}}{c_{2}\tilde{\sigma}_{v_{j}}}\right)+\sum_{m=0}^{\infty}\exp\left(-\frac{n\gamma_{v_{j}}2^{2(m+1)}r_{nj}^{2}}{c_{2}\tilde{\sigma}_{v_{j}}}\right)\leq c_{1}\exp\left(-\frac{n\gamma_{v_{j}}r_{nj}^{2}}{c_{2}\tilde{\sigma}_{v_{j}}^{2}}\right).
\]
\textbf{Lemma B.5}: There are universal positive constants $(c_{1},\, c_{2})$
such that for all $\Delta_{j}\in\left\{ f\in\mathcal{\bar{F}}_{j}\,:\,\left|f\right|_{\mathcal{F}_{j}}\geq1\right\} $,
\[
\frac{1}{n}\left|\sum_{i=1}^{n}v_{ij}\Delta_{ij}(\theta^{*})\right|\leq2r_{nj}\left|\Delta_{j}(\theta^{*})\right|_{n}+2r_{nj}^{2}\left|\Delta_{j}\right|_{\mathcal{F}_{j}}+\frac{1}{16}\left|\Delta_{j}(\theta^{*})\right|_{n}^{2}
\]
with probability at least $1-c_{1}\exp\left(-nC_{j}^{*}r_{nj}^{2}\right)$
and $C_{j}^{*}$ follows (\ref{eq:36}) under the same special cases
in Lemma B.3.\\
\\
\textbf{Proof. }This proof is a modification of the proofs for Corollary
8.3 in van de Geer (2000) and Lemma 13.2 in Wainwright (2015). It
suffices to show that 
\begin{equation}
\frac{1}{n}\left|\sum_{i=1}^{n}v_{ij}f(w_{i}^{T}\theta^{*})\right|\leq2r_{nj}\left|f_{\theta^{*}}\right|_{n}+2r_{nj}^{2}+\frac{1}{32}\left|f_{\theta^{*}}\right|_{n}^{2},\quad\textrm{for all}\,\left|f\right|_{\mathcal{F}_{j}}=1.\label{eq:44}
\end{equation}
To see this, noting that for any function $\Delta$ with $\left|\Delta\right|_{\mathcal{F}_{j}}\geq\frac{1}{2}$,
we can define $f=\frac{\Delta(\theta^{*})}{\left|\Delta(\theta^{*})\right|_{\mathcal{F}_{j}}}$.
Substituting this definition and then multiplying both sides of inequality
(\ref{eq:44}) by $\left|\Delta(\theta^{*})\right|_{\mathcal{F}_{j}}$,
we obtain 
\begin{eqnarray*}
\frac{1}{n}\left|\sum_{i=1}^{n}v_{ij}\Delta_{ij}(\theta^{*})\right| & \leq & 2r_{nj}\left|\Delta_{j}(\theta^{*})\right|_{n}+2r_{nj}^{2}\left|\Delta_{j}(\theta^{*})\right|_{\mathcal{F}_{j}}+\frac{1}{32}\frac{\left|\Delta_{j}(\theta^{*})\right|_{n}^{2}}{\left|\Delta_{j}(\theta^{*})\right|_{\mathcal{F}_{j}}}\\
 & \leq & 2r_{nj}\left|\Delta_{j}(\theta^{*})\right|_{n}+2r_{nj}^{2}\left|\Delta_{j}(\theta^{*})\right|_{\mathcal{F}_{j}}+\frac{1}{16}\left|\Delta_{j}(\theta^{*})\right|_{n}^{2},
\end{eqnarray*}
where the second inequality follows by the fact that $\left|\Delta\right|_{\mathcal{F}_{j}}\geq\frac{1}{2}$.

To establish the bound (\ref{eq:44}), we first consider it over the
ball $\left\{ \left|f_{\theta^{*}}\right|_{n}\leq u\right\} $ for
some fixed radius $u\geq r_{nj}$. Later in the proof we extend bound
(\ref{eq:44}) to one that is uniform in $\left|f_{\theta^{*}}\right|_{n}$
via a peeling argument. Define the random variable 
\[
A_{nj}^{''}(u):=\sup_{f(\theta)\in\Omega^{''}(u;\,\mathcal{F}_{j})}\frac{1}{n}\left|\sum_{i=1}^{n}v_{ij}f(w_{i}^{T}\theta^{*})\right|.
\]
where $\Omega^{''}(u;\,\mathcal{F}_{j})=\left\{ f\in\mathcal{\bar{F}}_{j}\,:\,\left|f\right|_{\mathcal{F}_{j}}\leq1,\,\left|f_{\theta^{*}}\right|_{n}\leq u\right\} $. 

The next part of the proof is almost identical to that in the proof
for Lemma B.4 except that we replace $A_{nj}^{'}(u)$ with $A_{nj}^{''}(u)$
and $\Omega^{'}(u;\,\mathcal{F}_{j})$ with $\Omega^{''}(u;\,\mathcal{F}_{j})$.
We first derive a bound for $u=r_{nj}$. By definition of $A_{nj}^{''}(u)$
and $\mathcal{G}_{n}\left(r_{nj};\,\mathcal{F}_{j}\right)$ and the
critical radius $r_{nj}$, we have 
\[
\mathbb{E}_{v_{j}}\left[A_{nj}^{''}(r_{nj})\,\vert\,\left\{ w_{i}^{T}\theta^{*}\right\} _{i=1}^{n}\right]=\tilde{\sigma}_{v_{j}}\mathcal{G}_{n}\left(r_{nj};\,\mathbb{B}_{\mathcal{F}_{j}}(1)\right)\leq r_{nj}^{2}.
\]
Setting $b=r_{nj}^{2}$ in the tail bound (\ref{eq:40}) with $A_{nj}^{'}(u)$
replaced by $A_{nj}^{''}(u)$ and $\Omega^{'}(u;\,\mathcal{F}_{j})$
replaced by $\Omega^{''}(u;\,\mathcal{F}_{j})$ yields 
\begin{equation}
\mathbb{P}\left[A_{nj}^{''}(r_{nj})\geq2r_{nj}^{2}\,\vert\,\left\{ w_{i}^{T}\theta^{*}\right\} _{i=1}^{n}\right]\leq\exp\left(-nC_{j}^{*}r_{nj}^{2}\right).\label{eq:45}
\end{equation}
On the other hand, for any $u>r_{nj}$, we have 
\begin{eqnarray*}
\mathbb{E}_{v_{j}}\left[A_{nj}^{''}(u)\,\vert\,\left\{ w_{i}^{T}\theta^{*}\right\} _{i=1}^{n}\right] & = & \tilde{\sigma}_{v_{j}}\mathcal{G}_{n}\left(u;\,\mathbb{B}_{\mathcal{F}_{j}}(1)\right)=u\frac{\tilde{\sigma}_{v_{j}}\mathcal{G}_{n}\left(u;\,\mathbb{B}_{\mathcal{F}_{j}}(1)\right)}{u}\\
 & \leq & u\frac{\tilde{\sigma}_{v_{j}}\mathcal{G}_{n}\left(r_{nj};\,\mathbb{B}_{\mathcal{F}_{j}}(1)\right)}{r_{nj}}\leq ur_{nj},
\end{eqnarray*}
where the first inequality follows from Lemma B.8 and the second inequality
follows from the critical radius $r_{nj}$. Setting $b=\frac{u^{2}}{128}$
in the tail bound (\ref{eq:40}) with $A_{nj}^{'}(u)$ replaced by
$A_{nj}^{''}(u)$ and $\Omega^{'}(u;\,\mathcal{F}_{j})$ replaced
by $\Omega^{''}(u;\,\mathcal{F}_{j})$ yields 
\begin{equation}
\mathbb{P}\left[A_{nj}^{''}(u)\geq ur_{nj}+\frac{u^{2}}{128}\,\vert\,\left\{ w_{i}^{T}\theta^{*}\right\} _{i=1}^{n}\right]\leq\exp\left(-nC_{j}^{*}u^{2}\right).\label{eq:46}
\end{equation}
It remains to prove the bound (\ref{eq:44}) via a peeling argument.
Let $\mathcal{E}^{'}$ denote the event that the bound (\ref{eq:44})
is violated for some function $f$ with $\left|f\right|_{\mathcal{F}_{j}}=1$.
For real numbers $0\leq a_{1}<a_{2}$, let $\mathcal{E}^{'}(a_{1},\, a_{2})$
denote the event that it is violated for some function such that $\left|f_{\theta^{*}}\right|_{n}\in[a_{1},\, a_{2}]$,
and $\left|f\right|_{\mathcal{F}_{j}}=1$. For $m=0,\,1,\,2,\,...$,
define $t_{m}=2^{m}r_{nj}$. We then have the decomposition $\mathcal{E}^{'}=\mathcal{E}^{'}(0,\, t_{0})\cup\left(\cup_{m=0}^{\infty}\mathcal{E}^{'}(t_{m},\, t_{m+1})\right)$
and hence by union bound, 
\begin{equation}
\mathbb{P}\left[\mathcal{E}^{'}\right]\leq\mathbb{P}\left[\mathcal{E}^{'}(0,\, t_{0})\right]+\sum_{m=0}^{\infty}\mathbb{P}\left[\mathcal{E}^{'}(t_{m},\, t_{m+1})\right].\label{eq:47}
\end{equation}
The final step is to bound each of the terms in this summation. Since
$t_{0}=r_{nj}$, we have 
\begin{equation}
\mathbb{P}\left[\mathcal{E}^{'}(0,\, t_{0})\right]\leq\mathbb{P}\left[A_{nj}^{''}(r_{nj})\geq2r_{nj}^{2}\,\vert\,\left\{ w_{i}^{T}\theta^{*}\right\} _{i=1}^{n}\right]\leq\exp\left(-nC_{j}^{*}r_{nj}^{2}\right).\label{eq:48}
\end{equation}
On the other hand, suppose that $\mathcal{E}^{'}(t_{m},\, t_{m+1})$
holds, meaning that there exists some function $f$ with $\left|f\right|_{\mathcal{F}_{j}}=1$
and $\left|f_{\theta^{*}}\right|_{n}\leq t_{m+1}$ such that 
\begin{eqnarray*}
\frac{1}{n}\left|\sum_{i=1}^{n}v_{ij}f(w_{i}^{T}\theta^{*})\right| & > & 2r_{nj}\left|f_{\theta^{*}}\right|_{n}+2r_{nj}^{2}+\frac{1}{32}\left|f_{\theta^{*}}\right|_{n}^{2}\\
 & > & 2r_{nj}t_{m}+2r_{nj}^{2}+\frac{1}{32}t_{m}^{2}\\
 & = & r_{nj}t_{m+1}+2r_{nj}^{2}+\frac{1}{128}t_{m+1}^{2},
\end{eqnarray*}
where the second inequality follows since $\left|f_{\theta^{*}}\right|_{n}\geq t_{m}$;
and the third inequality follows since $t_{m+1}=2t_{m}$. This lower
bound implies that 
\[
A_{nj}^{''}(t_{m+1})\geq r_{nj}t_{m+1}+\frac{1}{128}t_{m+1}^{2},
\]
hence the bound (\ref{eq:46}) implies that 
\[
\mathbb{P}\left[\mathcal{E}^{'}(t_{m},\, t_{m+1})\right]\leq\exp\left(-\frac{nt_{m+1}^{2}}{c_{2}\tilde{\sigma}_{v_{j}}}\right)=\exp\left(-\frac{n2^{2(m+1)}r_{nj}^{2}}{c_{2}\tilde{\sigma}_{v_{j}}}\right).
\]
Combining this tail bound with our earlier bound (\ref{eq:48}), and
substituting into the union bound (\ref{eq:47}) yields 
\[
\mathbb{P}\left[\mathcal{E}^{'}\right]\leq\exp\left(-\frac{n\gamma_{v_{j}}r_{nj}^{2}}{2\tilde{\sigma}_{v_{j}}}\right)+\sum_{m=0}^{\infty}\exp\left(-\frac{n\gamma_{v_{j}}2^{2(m+1)}r_{nj}^{2}}{c_{2}\tilde{\sigma}_{v_{j}}}\right)\leq c_{1}\exp\left(-\frac{n\gamma_{v_{j}}r_{nj}^{2}}{c_{2}\tilde{\sigma}_{v_{j}}^{2}}\right).
\]
$\square$\\
\textbf{}\\
\textbf{Lemma B.6}: Suppose the variable $X$ has \textit{strongly
log-concave} distribution with parameter $\gamma>0$. Then for any
function $f:\,\mathbb{R}^{n}\rightarrow\mathbb{R}$ that is \textit{L}-Lipschtiz
with respect to Euclidean norm, we have 
\[
\mathbb{P}\left[\left|f(X)-\mathbb{E}\left[f(X)\right]\right|\geq t\right]\leq2\exp\left(-\frac{\gamma t^{2}}{2L^{2}}\right).
\]
\textbf{Remark}: This result was initiated by Maurey (1991) and further
developed by Bobkov and Ledoux (2000).\\
\textbf{}\\
\textbf{Definition B.1} (Separate convex functions). A function $f:\,\mathbb{R}^{n}\rightarrow\mathbb{R}$
is \textit{separately convex }if, for each $j\in\left\{ 1,\,2,\,...,\, n\right\} $,
the co-ordinate function $f_{j}:\,\mathbb{R}\rightarrow\mathbb{R}$
defined by varying only the $j^{th}$ co-ordinate, is a convex function
of $x_{j}$.\\
\textbf{}\\
\textbf{Lemma B.7}: Let $\left\{ X_{i}\right\} _{i=1}^{n}$ be independent
random variables, each supported on the interval $[a,\, b]$, and
let $f:\,\mathbb{R}^{n}\rightarrow\mathbb{R}$ be separately convex,
and $L-$Lipschitz with respect to the Euclidean norm. Then for all
$t>0$, we have\textbf{
\[
\mathbb{P}\left[f(X)-\mathbb{E}\left[f(X)\right]\geq t\right]\leq2\exp\left(-\frac{t^{2}}{2L^{2}(b-a)^{2}}\right).
\]
Remark}: The proof can be found in Wainwright (2015), which is based
on Ledoux (1996).\textbf{}\\
\\
\textbf{Lemma B.8}: For any star-shaped function class $\mathcal{F}$,
the function $t\mapsto\frac{\mathcal{G}_{n}(t;\,\mathcal{F})}{t}$
is non-decreasing on the interval $(0,\,\infty)$.\\
\\
\textbf{Proof.} This is Lemma 13.1 from Wainwright (2015).\\
\\
\textbf{Definition B.2} (Covering and packing numbers). Consider a
metric space consisting of a set $\mathcal{X}$ and a metric $\rho:\,\mathcal{X}\times\mathcal{X}\rightarrow\mathbb{R}_{+}$. 
\begin{enumerate}
\item An $t-$covering of $\mathcal{X}$ in the metric $\rho$ is a collection
$\left\{ \beta^{1},\,...,\,\beta^{N}\right\} \subset\mathcal{X}$
such that for all $\beta\in\mathcal{X}$, there exists some $i\in\left\{ 1,\,...,\, N\right\} $
with $\rho(\beta,\,\beta^{i})\leq t$. The $t-$covering number $N(t;\,\mathcal{X},\,\rho)$
is the cardinality of the smallest $t-$covering.
\item An $t^{'}-$packing of $\mathcal{X}$ in the metric $\rho$ is a collection
$\left\{ \beta^{1},\,...,\,\beta^{N}\right\} \subset\mathcal{X}$
such that $\rho(\beta^{i},\,\beta^{j})\geq t$ for all $i\neq j$.
The $t^{'}-$packing number $M(t^{'};\,\mathcal{X},\,\rho)$ is the
cardinality of the largest $t^{'}-$covering.
\end{enumerate}
\textbf{Lemma B.9}: Let $N_{n}(t;\,\mathbb{B}_{n}(r_{n};\,\mathcal{F}))$
denote the $t-$covering number of the set 
\[
\mathbb{B}_{n}(r_{n};\,\mathcal{F})=\left\{ f\in2\mathcal{F}\,:\,\left|f\right|_{n}\leq r_{n}\right\} 
\]
in the empirical $L^{2}(\mathbb{P}_{n})$ norm. Then the critical
radius condition (\ref{eq:18}) holds for any $r_{n}\in(0,\,\sigma]$
such that 
\[
\frac{32}{\sqrt{n}}\int_{\frac{r_{n}^{2}}{2\sigma}}^{r_{n}}\sqrt{\log N_{n}(t;\,\mathbb{B}_{n}(r_{n};\,\mathcal{F}))}dt\leq\frac{r_{n}^{2}}{\sigma}.
\]
\textbf{Remark.} This result is established by van der Vaart and Wellner
(1996), van de Geer (2000), Barlett and Mendelson (2002), Koltchinski
(2006), Wainwright (2015), etc.\\
\\
\textbf{Lemma B.10}: For $q\in(0,\,1]$, let 
\[
\theta\in\mathbb{\mathcal{B}}^{d}(R_{q}):=\left\{ \theta^{'}\in\mathbb{R}^{d}\,:\,\left|\theta^{'}\right|_{1}^{q}=\sum_{j=1}^{d}\left|\theta_{j}^{'}\right|^{q}\leq R_{q}\right\} 
\]
and $N_{2}(t;\,\mathbb{\mathcal{B}}^{d}(R_{q}))$ be the $t-$covering
number of the set $\mathbb{\mathcal{B}}^{d}(R_{q})$ in the $l_{2}-$norm.
Then there is a universal constant $c$ such that 
\[
\log N_{2}(t;\,\mathbb{\mathcal{B}}^{d}(R_{q}))\leq cR_{q}^{\frac{2}{2-q}}\left(\frac{1}{t}\right)^{\frac{2q}{2-q}}\log d\quad\textrm{for all}\, t\in(0,\, R_{q}^{\frac{1}{q}})
\]
Conversely, suppose in addition that $t<1$ and $t^{2}\geq c^{''}\left(\frac{\log d}{d^{\kappa}}\right)^{\frac{2-q}{q}}$
for some fixed $\kappa\in(0,\,1)$. Then there is a constant $c^{'}\leq c$
such that 
\[
\log N_{2}(t;\,\mathbb{\mathcal{B}}^{d}(R_{q}))\geq c^{'}R_{q}^{\frac{2}{2-q}}\left(\frac{1}{t}\right)^{\frac{2q}{2-q}}\log d.
\]
\textbf{Proof.} The result is obtained by inverting known results
on (dyadic) entropy numbers of $l_{q}-$balls as in Raskutii, Wainwright,
and Yu (2011). For a $d-$dimensional $l_{q}-$ball with $q\in(0,\, p)$,
it is known (e.g., Schütt, 1984; Guedon and E. Litvak, 2000; Kühn,
2001) that for all integers $k\in[\log d,\, d]$, the dyadic entropy
numbers $t_{k}$ of the $d-$dimensional $l_{q}-$ball of unit radius
$\mathbb{\mathcal{B}}_{q}^{d}(1)$ with respect to the $l_{p}-$norm
scale as 
\[
t_{k}\left(l_{q},\,\left|\cdot\right|_{p}\right)=C_{q,\, p}\left[\frac{\log(1+\frac{d}{k})}{k}\right]^{\frac{1}{q}-\frac{1}{p}}.
\]
Moreover, for $k\in[1,\,\log d]$, we have $t_{k}(l_{q})\leq C_{q,\, p}$.
We first establish the upper bound on the metric entropy for $q\in(0,\,1]$
and $p=2$. Since $d\geq2$, we have 
\[
t_{k}\left(l_{1},\,\left|\cdot\right|_{2}\right)\leq C_{q,\,2}\left(\frac{\log(1+\frac{d}{2})}{k}\right)^{\frac{1}{q}-\frac{1}{2}}\leq C_{q,\,2}\left(\frac{\log d}{k}\right)^{\frac{1}{q}-\frac{1}{2}}.
\]
Inverting this inequality for $k=\log N_{2}(t;\,\mathbb{\mathcal{B}}^{d}(R_{q}))$
and allowing for a ball radius $R_{q}$ yields 
\begin{equation}
\log N_{2}(t;\,\mathbb{\mathcal{B}}^{d}(R_{q}))\leq cR_{q}^{\frac{2}{2-q}}\left(\frac{1}{t}\right)^{\frac{2q}{2-q}}\log d.\label{eq:49}
\end{equation}
For the lower bound on the metric entropy, which requires the existence
of some fixed $\kappa\in(0,\,1)$ such that $k\leq d^{1-\kappa}$.
Under this assumption, we have $1+\frac{d}{k}\geq\frac{d}{k}\geq d^{\kappa}$,
and hence 
\[
C_{q,\,2}\left(\frac{\log(1+\frac{d}{k})}{k}\right)^{\frac{1}{q}-\frac{1}{2}}\geq C_{q,\,2}\left(\frac{\kappa\log d}{k}\right)^{\frac{1}{q}-\frac{1}{2}}.
\]
Accounting for the radius $R_{q}$ yields 
\[
\log N_{2}(t;\,\mathbb{\mathcal{B}}^{d}(R_{q}))\geq\kappa C_{q,\,2}R_{q}^{\frac{2}{2-q}}\left(\frac{1}{t}\right)^{\frac{2q}{2-q}}\log d.
\]
Finally, let us check that our assumptions on $k$ needed to perform
the inversion are ensured by the conditions that we have imposed on
$t$. The condition $k\geq\log d$ is ensured by setting $t<1$. Turning
to the condition $k\leq d^{1-\kappa}$, from bound (\ref{eq:49}),
it suffices to choose $t$ such that $\left(\frac{C_{q,\,2}}{t}\right)^{\frac{2q}{2-q}}\log d\leq d^{1-\kappa}$.
This condition is ensured by enforcing the bound $t^{2}\geq c^{''}\left(\frac{\log d}{d^{\kappa}}\right)^{\frac{2-q}{q}}$.
$\square$\\
\\
\textbf{Lemma B.11}: Consider the class of smooth functions 
\[
\mathcal{F}_{L}:=\left\{ f\,:\,[0,\,1]\rightarrow\mathbb{R}\,:\, f(0)=0,\,\left|f(x)-f(x^{'})\right|\leq L\left|x-x^{'}\right|\quad\forall x,\, x^{'}\in[0,\,1]\right\} .
\]
We have 
\begin{eqnarray*}
\log N_{\infty}(t;\,\mathcal{F}_{L}) & \asymp & \frac{L}{t},\\
\textrm{and }\log N_{2}(t;\,\mathcal{F}_{L}) & \asymp & \frac{L}{t}
\end{eqnarray*}
where $N_{\infty}(t;\,\mathcal{F}_{L})$ is the $t-$covering number
of the set $\mathcal{F}_{L}$ in the sup-norm and $N_{2}(t;\,\mathcal{F}_{L})$
is the $t-$covering number of the set $\mathcal{F}_{L}$ in the $L^{2}\left(\mathbb{P}\right)-$norm.\\
\\
\textbf{Remark.} The proof for $\log N_{\infty}(t;\,\mathcal{F}_{L})$
is a well-known result due to Kolmogorov and Tikhomirov (1961). To
proof for $\log N_{2}(t;\,\mathcal{F}_{L})$ is an application of
Varshamov-Gilbert lemma. $\square$\\
\\
\textbf{Lemma B.12}: For all $\epsilon>0$, the packing and covering
number of a set $\mathbb{T}$ with respect to a metric $\rho$ are
related as follows:
\[
M(2\epsilon;\,\mathbb{T},\,\rho)\leq N(\epsilon;\,\mathbb{T},\,\rho)\leq M(\epsilon;\,\mathbb{T},\,\rho).
\]
\textbf{Proof.} To prove the first inequality, let $\left\{ \theta^{1},\,...,\,\theta^{2M}\right\} $
be a maximal $2\epsilon-$packing, so that $\rho(\theta^{i},\,\theta^{j})>2\epsilon$
for all $i\neq j$. Given a $\epsilon-$cover, let $\gamma$ be an
arbitrary element. For any distinct pair $i\neq j$, we have 
\[
2\epsilon<\rho(\theta^{i},\,\theta^{j})\leq\rho(\theta^{i},\,\gamma)+\rho(\gamma,\,\theta^{j}),
\]
so that it is not possible to have $\rho(\gamma,\,\theta^{i})\leq\epsilon$
and $\rho(\gamma,\,\theta^{j})\leq\epsilon$ simultaneously. Conversely,
we must have at least one element in the cover for each element of
the packing set, showing that $M(2\epsilon;\,\mathbb{T},\,\rho)\leq N(\epsilon;\,\mathbb{T},\,\rho)$.
Turning to the second inequality, let $\left\{ \theta^{1},\,...,\,\theta^{2M}\right\} $
be a maximal $\epsilon-$packing. Maximality implies that for any
other $\theta$ not already in the packing set, we must have $\rho(\theta,\,\theta^{i})\leq\epsilon$
for some $\theta^{i}$. Thus, the given set forms a $\epsilon-$cover,
showing that $N(\epsilon;\,\mathbb{T},\,\rho)\leq M(\epsilon;\,\mathbb{T},\,\rho)$.
$\square$\\
\\
\textbf{Lemma B.13}: Let $k\geq1$. There exists a subset $W$ of
$\{0,\,1\}^{k}$ with $\left|W\right|>e^{k/8}$ such that the Hamming
distance, 
\[
d(\tau,\,\tau^{'}):=\sum_{j=1}^{k}\mathbb{I}\left\{ \tau_{i}\neq\tau_{i}^{'}\right\} >\frac{k}{4}
\]
for all $\tau,\,\tau^{'}\in W$ with $\tau\neq\tau^{'}$. \\
\\
\textbf{Remark.} This is Varshamov-Gilbert Lemma.\\
\textbf{}\\
\textbf{Lemma B.14}: For a given radius $r>0$, define the set 
\[
\mathcal{S}(s,\, r):=\left\{ \delta\in\mathbb{R}^{d}\,\vert\,\left|\delta\right|_{0}\leq2s,\:\left|\delta\right|_{2}\leq r\right\} ,
\]
there is a covering set $\left\{ \delta^{1},\,...,\,\delta^{N}\right\} \in\mathcal{S}(s,\, r)$
with cardinality $N=N(s,\, r,\,\epsilon)$ such that 
\[
\log N(s,\, r,\,\epsilon)\leq\log\left(\begin{array}{c}
d\\
2s
\end{array}\right)+2s\log\left(\frac{1}{\epsilon}\right).
\]
\textbf{Remark.} This result is established in Matousek (2002).\\
\\
\textbf{Lemma B.15}: There exists a subset $\tilde{\mathcal{M}}\subset\mathcal{M}$
with cardinality $\left|\tilde{\mathcal{M}}\right|\geq\exp\left(\frac{s}{2}\log\frac{d-s}{s/2}\right)$
such that $\rho_{\mathcal{M}}(z,\, z^{'})\geq\frac{s}{2}$ for all
$z,\, z^{'}\in\tilde{\mathcal{M}}$.\\
\textbf{}\\
\textbf{Remark.} This is Lemma 5 from Raskutii, Wainwright, and Yu
(2011). \\
\textbf{}\\
\textbf{Lemma B.16}: Let $X_{1},\,...,\, X_{n}$ be independent random
variables such that, for every $i$, $\mathbb{E}(X_{i})=0$ and $0<\mathbb{E}\left(\left|X_{i}\right|^{2+\delta}\right)<\infty$
for some $0<\delta<1$. Set $S_{n}=\sum_{i=1}^{n}X_{i}$, $B_{n}^{2}=\sum_{i=1}^{n}\mathbb{E}(X_{i}^{2})$,
$V_{n}^{2}=\sum_{i=1}^{n}X_{i}^{2}$, $L_{n,\delta}=\sum_{i=1}^{n}\mathbb{E}\left(\left|X_{i}\right|^{2+\delta}\right)$,
$b_{n,\delta}=\frac{B_{n}}{L_{n,\delta}^{1/(2+\delta)}}$. Then for
all $0\leq t\leq b_{n,\delta}$ and an absolute constant $a_{0}>0$
\[
\left|\mathbb{P}\left[\frac{S_{n}}{V_{n}}\geq t\right]-\left[1-\Phi(t)\right]\right|\leq a_{0}(1+t)^{1+\delta}\frac{e^{-\frac{t^{2}}{2}}}{b_{n,\delta}^{2+\delta}}.
\]
\textbf{Remark.} This is formula 2.11 from Jing, Shao, and Wang (2003)
on moderate deviations for self-normalized sums.\\
\\
\textbf{Lemma B.17}: Let $X_{k}$, $k\geq1$ be \textit{i.i.d.} sub-Gaussian
variables with parameters at most $\sigma$ and $Z:=\max_{k=1,...,N}\left|X_{k}\right|$.
For some integer $N\geq10$, we have $\mathbb{E}\left[Z\right]\leq3\sigma\sqrt{N}$.\\
\\
\textbf{Remark.} This is a well-known result on sub-Gaussian maxima
(e.g., van der Vaart and Wellner, 1996).